\documentclass[12pt]{article}
\usepackage{amssymb}
\usepackage{amsmath}
\usepackage{makeidx}
\usepackage[dvips]{graphicx}
\usepackage{sw20elba}

\newtheorem{theorem}{Theorem}[section]

\newtheorem{conclusion}[theorem]{Conclusion}
\newtheorem{condition}[theorem]{Condition}

\newtheorem{corollary}[theorem]{Corollary}

\newtheorem{definition}[theorem]{Definition}
\newtheorem{example}[theorem]{Example}

\newtheorem{lemma}[theorem]{Lemma}

\newtheorem{proposition}[theorem]{Proposition}
\newtheorem{remark}[theorem]{Remark}

\newenvironment{proof}[1][Proof]{\textbf{#1.} }{\ \rule{0.5em}{0.5em}}
\input{tcilatex}
\numberwithin{equation}{section}

\begin{document}

\title{Linear superposition in nonlinear wave dynamics }
\author{A. Babin and A. Figotin \\
Department of Mathematics, University of California at Irvine, CA 92697}
\maketitle

\begin{abstract}
We study nonlinear dispersive wave systems described by hyperbolic PDE's in $%
\mathbb{R}^{d}$ and difference equations on the lattice $\mathbb{Z}^{d}$.
The systems involve two small parameters: one is the ratio of the slow and
the fast time scales, and another one is the ratio of the small and the
large space scales. We show that a wide class of\ such systems, including
nonlinear Schrodinger and Maxwell equations, Fermi-Pasta-Ulam\ model and
many other not completely integrable systems, satisfy a superposition
principle. The principle essentially states that if a nonlinear evolution of
a wave starts initially as a sum of generic wavepackets (defined as almost
monochromatic waves), then this wave with a high accuracy remains a sum of
separate wavepacket\ waves undergoing independent nonlinear evolution. The
time intervals for which the evolution is considered are long enough to
observe fully developed nonlinear phenomena \ for involved wavepackets. In
particular, our approach provides a simple justification for numerically
observed effect of almost non-interaction\ of solitons passing through each
other without any recourse to the complete integrability. Our analysis does
not rely on any ansatz or common asymptotic expansions with respect to the
two small parameters but it uses rather explicit and constructive
representation for solutions as functions of the initial data in the form of
functional analytic series.
\end{abstract}

\section{Introduction}

The principal object of our studies here is a general nonlinear evolutionary
system which describes wave propagation in homogeneous media governed either
by a hyperbolic PDE's in $\mathbb{R}^{d}$ or by a difference equation on the
lattice $\mathbb{Z}^{d}$, $d=1,2,3,\ldots $ is the space dimension. We
assume the evolution to be governed by the following equation with constant
coefficients 
\begin{equation}
\partial _{\tau }\mathbf{U}=-\frac{\mathrm{i}}{\varrho }\mathbf{L}\left( -%
\mathrm{i}\nabla \right) \mathbf{U}+\mathbf{F}\left( \mathbf{U}\right) ,\
\left. \mathbf{U}\left( \mathbf{r},\tau \right) \right\vert _{\tau =0}=%
\mathbf{h}\left( \mathbf{r}\right) ,\ \mathbf{r}\in \mathbb{R}^{d},
\label{difeqintr}
\end{equation}%
where (i) $\mathbf{U}=\mathbf{U}\left( \mathbf{r},\tau \right) $, $\mathbf{r}%
\in \mathbb{R}^{d}$, $\mathbf{U}\in \mathbb{C}^{2J}$ is a $2J$ dimensional
vector; (ii) $\mathbf{L}\left( -\mathrm{i}\nabla \right) $ is a linear
self-adjoint differential (pseudodifferential) operator with constant
coefficients with the symbol $\mathbf{L}\left( \mathbf{k}\right) $, which is
a Hermitian $2J\times 2J$ matrix; (iii) $\ \mathbf{F}$ is a general
polynomial nonlinearity; (iv) $\varrho >0$ is a \emph{small parameter}. \
The form of the equation suggests that the processes described by it involve
two time scales. Since the nonlinearity $\mathbf{F}\left( \mathbf{U}\right) $
is of order one, nonlinear effects occur at times $\tau $ of order one,
whereas the natural time scale of linear effects, governed by the operator $%
\mathbf{L}$ with the coefficient $1/\varrho $, is of order $\varrho $.
Consequently, the small parameter $\varrho $ measures the ratio of the slow
(nonlinear effects) time scale and the fast (linear effects) time scale. A
typical example an equation of the form (\ref{difeqintr}) is nonlinear
Schrodinger equation (NLS) or a system of NLS. Another one is the Maxwell
equation in a periodic medium when truncated to a finite number of bands,
and more examples are discussed below.

We assume further that the initial data $\mathbf{h}$ for the evolution
equation (\ref{difeqintr}) to be the sum of a finite number of \emph{%
wavepackets} $\mathbf{h}_{l}$, $l=1,\ldots ,N$, i.e.%
\begin{equation}
\mathbf{h}=\mathbf{h}_{1}+\ldots +\mathbf{h}_{N}\   \label{hhh1}
\end{equation}%
where the monochromaticity of every wavepacket $\mathbf{h}_{l}$ is
characterized by another small parameter $\beta $.

The well known \emph{superposition principle} is a fundamental property of
every linear evolutionary system, stating that the solution $\mathbf{U}$
corresponding to the initial data $\mathbf{h}$ as in (\ref{hhh1}) equals%
\begin{equation}
\mathbf{U}=\mathbf{U}_{1}+\ldots +\mathbf{U}_{N},\text{ for }\mathbf{h}=%
\mathbf{h}_{1}+\ldots +\mathbf{h}_{N},  \label{hhh2}
\end{equation}
where $\mathbf{U}_{l}$ is the solution to the same linear problem with the
initial data $\mathbf{h}_{l}$.

Evidently the standard superposition principle can not hold exactly as a
general principle in the presence of a nonlinearity, and, at the first
glance, there is no expectation for it to hold even approximately. We have
discovered though that the \emph{superposition principle does hold with a
high accuracy for general dispersive nonlinear wave systems provided that
the initial data are a sum of generic wavepackets, and this constitutes the
subject of this paper}. Namely, the superposition principle for nonlinear
wave systems states that the solution $\mathbf{U}$ corresponding to the
multi-wavepacket initial data $\mathbf{h}$ as in (\ref{hhh1}) equals \emph{\ 
}%
\begin{equation*}
\mathbf{U}=\mathbf{U}_{1}+\ldots +\mathbf{U}_{N}+\mathbf{D},\text{ for }%
\mathbf{h}=\mathbf{h}_{1}+\ldots +\mathbf{h}_{N},\text{ where }\mathbf{D}%
\text{ is small.}
\end{equation*}

As to the particular form (\ref{difeqintr}) we chose to be our primary one,
we would like to point out that many important classes of problems involving
small parameters \ can be readily reduced to the framework of (\ref%
{difeqintr}) by a simple rescaling. It can be seen from the following
examples. First example is a system with a small factor before the
nonlinearity 
\begin{equation}
\partial _{t}\mathbf{v}=-\mathrm{i}\mathbf{Lv}+\alpha \mathbf{f}\left( 
\mathbf{v}\right) ,\ \left. \mathbf{v}\right\vert _{t=0}=\mathbf{h},\
0<\alpha \ll 1,  \label{smallnon}
\end{equation}%
where initial data are bounded uniformly in $\alpha $. \ Such problems are
reduced to (\ref{difeqintr}) by the time rescaling $\tau =t\alpha .\ $Note
that now $\varrho =\alpha $ \ and the finite time interval $0\leq \tau \leq
\tau _{\ast }$ corresponds to the \emph{long time interval} $0\leq t\leq
\tau _{\ast }/\alpha $.

Second example is a system with small initial data on a long time interval.
The system here is given and has no small parameters but the initial data
are small, namely 
\begin{gather}
\partial _{t}\mathbf{v}=-\mathrm{i}\mathbf{Lv}+\mathbf{f}_{0}\left( \mathbf{v%
}\right) ,\ \left. \mathbf{v}\right\vert _{t=0}=\alpha _{0}\mathbf{h},\
0<\alpha _{0}\ll 1,\text{ where}  \label{smallinit} \\
\mathbf{f}_{0}\left( \mathbf{v}\right) =\mathbf{f}_{0}^{\left( m\right)
}\left( \mathbf{v}\right) +\mathbf{f}_{0}^{\left( m+1\right) }\left( \mathbf{%
v}\right) +\ldots ,  \notag
\end{gather}%
where $\alpha _{0}$ is a small parameter and $\mathbf{f}^{\left( m\right)
}\left( \mathbf{v}\right) $ is a \ homogeneous polynomial of degree $m\geq 2$%
. After the rescaling $\mathbf{v}=\alpha _{0}\mathbf{V}$ we obtain the
following equation \ with a small nonlinearity 
\begin{equation}
\partial _{t}\mathbf{V}=-\mathrm{i}\mathbf{LV}+\alpha _{0}^{m-1}\left[ 
\mathbf{f}_{0}^{\left( m\right) }\left( \mathbf{V}\right) +\alpha _{0}%
\mathbf{f}^{0\left( m+1\right) }\left( \mathbf{V}\right) +\ldots \right] ,\
\left. \mathbf{V}\right\vert _{t=0}=\mathbf{h},  \label{smin1}
\end{equation}%
which is of the form of (\ref{smallnon}) with $\alpha =\alpha _{0}^{m-1}$.
Introducing the slow time variable $\tau =t\alpha _{0}^{m-1}$ we get from
the above an equation of the form (\ref{difeqintr}), namely%
\begin{equation}
\partial _{\tau }\mathbf{V}=-\frac{\mathrm{i}}{\alpha _{0}^{m-1}}\mathbf{LV}+%
\left[ \mathbf{f}^{\left( m\right) }\left( \mathbf{V}\right) +\alpha _{0}%
\mathbf{f}^{\left( m+1\right) }\left( \mathbf{V}\right) +\ldots \right] ,\
\left. \mathbf{V}\right\vert _{t=0}=\mathbf{h},  \label{smin2}
\end{equation}%
where the nonlinearity does not vanish as $\alpha _{0}\rightarrow 0$. In
this case $\varrho =\alpha _{0}^{m-1}$ and the finite time interval $0\leq
\tau \leq \tau _{\ast }$ corresponds to the long time interval $0\leq t\leq 
\frac{\tau _{\ast }}{\alpha _{0}^{m-1}}$ \ with small $\alpha _{0}\ll 1$.

Very often in theoretical studies of equations of the form (\ref{difeqintr})
or ones reducible to it a functional dependence between $\varrho $ and $%
\beta $ is imposed, resulting in a single small parameter. The most common
scaling is $\varrho =\beta ^{2}$. The nonlinear evolution of wavepackets for
a variety of equations which can be reduced to the form (\ref{difeqintr})
was studied in numerous physical and mathematical papers, mostly by
asymptotic expansions of solutions with respect to a single small parameter
similar to $\beta $, see \cite{BenYoussefLannes02}, \cite{BonaCL05}, \cite%
{ColinLannes}, \cite{CraigSulemS92}, \cite{GiaMielke}, \cite{JolyMR98}, \cite%
{KalyakinUMN}, \cite{Maslov83}, \cite{PW}, \cite{Schneider98a}, \cite%
{Schneider05} and references therein. Often the asymptotic expansions are
based on a specific ansatz prescribing a certain form\ to the solution. In
our studies here we do not use asymptotic expansions with respect to a small
parameter and do not prescribe a specific form to the solution, but we
impose conditions on the initial data requiring it to be a wavepacket or a
linear combination of wavepackets. Since we want to establish a general
property of a wide class of systems, we apply a general enough dynamical
approach. There is a number of general approaches developed for the studies
of high-dimensional and infinite-dimensional nonlinear evolutionary systems
\ of hyperbolic type, \cite{Bambusi03}, \cite{BM}, \cite{Craig}, \cite{GW}, 
\cite{Iooss}, \cite{Kuksin}, \cite{MSZ}, \cite{Schneider98a}, \cite{SU}, 
\cite{Weinstein}, \cite{W}) and references therein. \ We develop here an
approach which allows to exploit specific properties of a certain class of
initial data, namely wavepackets and their linear combinaions, which comply
with the symmetries of equations. Such a class of the initial data is
obviously lesser than all possible initial data. One of the key mathematical
tools developed here for the nonlinear studies is a \emph{refined implicit
function theorem} (Theorem \ref{Theorem Monomial Convergence}). This theorem
provides a constructive and rather explicit representation of the solution
to an abstract nonlinear equation in a Banach space as a certain functional
series. The representation is explicit enough to prove the superposition
principle and is general enough to carry out the studies of the problem
without imposing restrictions on dimension of the problem, structural
restrictions on nonlinearities or a functional dependence between the two
small parameters $\varrho ,\beta $.

As we have already stated the superposition principle holds with high
accuracy for linear combinations of wavepackets. A wavepacket $\mathbf{h}%
\left( \beta ,\mathbf{r}\right) $ can be most easily described in terms of
its Fourier transform $\mathbf{\tilde{h}}\left( \beta ,\mathbf{k}\right) $.
Simply speaking, wavepacket $\mathbf{\tilde{h}}\left( \beta ,\mathbf{k}%
\right) $ is a function which is localized in $\beta $-neighborhood of a
given wavevector $\mathbf{k}_{\ast }$ (the \emph{wavepacket center) }and as
a vector is an eigenfunction of the matrix $\mathbf{L}\left( \mathbf{k}%
\right) $, details of the definition of the wavepacket can be found in the
following Section 2. The simplest example of a wavepacket is a function of
the form%
\begin{equation}
\mathbf{\tilde{h}}\left( \beta ,\mathbf{k}\right) =\beta ^{-d}\hat{h}\left( 
\frac{\mathbf{k}-\mathbf{k}_{\ast }}{\beta }\right) \mathbf{g}_{n}\left( 
\mathbf{k}_{\ast }\right) ,\ \mathbf{k}\in \mathbb{R}^{d},  \label{hh}
\end{equation}%
where $\mathbf{g}_{n}\left( \mathbf{k}_{\ast }\right) $ is an eigenvector of
the matrix $\mathbf{L}\left( \mathbf{k}_{\ast }\right) $ and $\hat{h}\left( 
\mathbf{k}\right) $ is a Schwartz function (i.e. it is infinitely smooth and
rapidly decaying one). Note that the inverse Fourier transform $\mathbf{h}%
\left( \beta ,\mathbf{r}\right) $ of $\mathbf{\tilde{h}}\left( \beta ,%
\mathbf{k}\right) $ has the form 
\begin{equation}
\mathbf{h}\left( \beta ,\mathbf{r}\right) =h\left( \beta \mathbf{r}\right) 
\mathrm{e}^{\mathrm{i}\mathbf{k}_{\ast }\mathbf{r}}\mathbf{g}_{n}\left( 
\mathbf{k}_{\ast }\right) ,\ \mathbf{r}\in \mathbb{R}^{d},  \label{hr}
\end{equation}%
where $h\left( \mathbf{r}\right) $ is a Schwartz function, and obviously has
a large spatial extension of order $\beta ^{-1}$.

We study the nonlinear evolution equation (\ref{difeqintr}) on a finite time
interval 
\begin{equation}
0\leq \tau \leq \tau _{\ast },\text{ where }\tau _{\ast }>0\text{ is a fixed
number}  \label{taustar}
\end{equation}%
which may depend on the $L^{\infty }$ norm of the initial data $\mathbf{h}$
but, importantly, $\tau _{\ast }$\emph{\ does not depend on }$\varrho $. 
\emph{We consider classes of initial data such that wave evolution governed
by (\ref{difeqintr}) is significantly nonlinear on time interval }$\left[
0,\tau _{\ast }\right] $\emph{\ and the effect of the nonlinearity }$F\left( 
\mathbf{U}\right) $\emph{\ does not vanish as }$\varrho \rightarrow 0$\emph{%
. \ } We assume that $\beta ,\varrho $ satisfy%
\begin{equation}
0<\beta \leq 1,\ 0<\varrho \leq 1,\ \frac{\beta ^{2}}{\varrho }\leq C_{1}%
\text{ with some }C_{1}>0.  \label{rbb1}
\end{equation}%
The above condition on the dispersion parameter $\frac{\beta ^{2}}{\varrho }$
ensures that the dispersive effects are not dominant and do not suppress
nonlinear effects, see \cite{BF7} for a discussion.

To formulate the superposition principle more precisely we introduce first
the solution operator $\mathcal{S}\left( \mathbf{h}\right) \left( \tau
\right) :\mathbf{h}\rightarrow \mathbf{U}\left( \tau \right) $ which relates
to the initial data $\mathbf{h}$ of the nonlinear evolution equation (\ref%
{difeqintr}) the solution $\mathbf{U}\left( t\right) $ of this equation.
Suppose that the initial state is a \emph{multi-wavepacket}, namely $\mathbf{%
h}=\dsum \mathbf{h}_{l}$, with $\mathbf{h}_{l}$, $l=1,\ldots ,N$ being
"generic" wavepackets. Then for all times\emph{\ }$0\leq \tau \leq \tau
_{\ast }$ the following\emph{\ superposition principle} holds 
\begin{gather}
\mathcal{S}\left( \dsum\nolimits_{l=1}^{N}\mathbf{h}_{l}\right) \left( \tau
\right) =\dsum\nolimits_{l=1}^{N}\mathcal{S}\left( \mathbf{h}_{l}\right)
\left( \tau \right) +\mathbf{D}\left( \tau \right) ,  \label{Gsum1} \\
\left\Vert \mathbf{D}\left( \tau \right) \right\Vert _{E}=\sup\limits_{0\leq
\tau \leq \tau _{\ast }}\left\Vert \mathbf{D}\left( \tau \right) \right\Vert
_{L^{\infty }}\leq C_{\delta }\frac{\varrho }{\beta ^{1+\delta }}\text{ for
any small }\delta >0.  \label{Dbetkap}
\end{gather}%
Obviously, the right-hand side of (\ref{Dbetkap}) may be small \ only if $%
\varrho \leq C_{1}\beta $. There are examples (see \cite{BF7}) in which $%
\mathbf{D}\left( \tau \right) $ is not small for $\varrho =C_{1}\beta $. \
In what follows we refer to a linear combination of wavepackets as a\emph{\
multi-wavepacket}, and to wavepackets which constitutes the multi-wavepacket
as \emph{component wavepackets}.

The superposition principle implies, in particular, that in the process of
nonlinear evolution every single wavepacket propagates almost \
independently of other wavepackets even though they may "collide" in
physical space for a certain period of time and the exact solution equals
the sum of particular single wavepacket solutions with a high precision. In
particular, the dynamics of a solution with multi-wavepacket initial data is
reduced to dynamics of separate solutions with single wavepacket data. Note
that the nonlinear evolution of a single wavepacket solution for many
problems is studied in detail, namely it is well approximated by its own
nonlinear Schrodinger equation (NLS), see \cite{ColinLannes}, \cite%
{GiaMielke}, \cite{KalyakinUMN}, \cite{Kalyakin2}, \cite{Schneider98a}, \cite%
{Schneider05}, \cite{SU}, \cite{BF7} and references therein.

The superposition principle (\ref{Gsum1}), (\ref{Dbetkap}) can also be
looked at as a form of separation of variables. Such a form of separation of
variables is different from usual complete integrability, and its important
factor is the continuity of spectrum of the linear component of the system.
The approximate superposition principle imposes certain restrictions on
dynamics which differ from usual constraints imposed by the conserved
quantities as in completely integrable systems as well as from topological
constraints related to invariant tori as in KAM theory.

Now we present an elementary physical argument justifying the superposition
principle. If nonlinearity is absent, the superposition principle holds
exactly and any deviation from it is due to the nonlinear interactions
between wavepackets, so we need to estimate their impact. Suppose that
initially at time $\tau =0$ the spatial extension $\ s$ of every composite
wavepacket is characterized by the parameter $\beta ^{-1}$ as in (\ref{hr}).
Assume also (and it is quite an assumption) that the component wavepackets
during the nonlinear evolution maintain somehow their wavepacket identity,
group velocities and spatial extension. Then, consequently, the spatial
extension of every component wavepacket is propositional to $\beta ^{-1}$
and its group velocity $v_{j}$ is proportional to $\varrho ^{-1}$. The
difference $\Delta v$ between any two different component group velocities
is also proportional to $\varrho ^{-1}$. The time when two different
component wavepackets overlap in space is proportional to $s/\left\vert
\Delta v\right\vert $ and, hence, to $\varrho /\beta $. Since the nonlinear
term is of order one, the \ magnitude of the impact of the nonlinearity
during this time interval should be proportional to $\varrho /\beta $, which
results in the same order of magnitude of $\mathbf{D}$. This conclusion is
in agreement with the estimate of magnitude of $\mathbf{D}$ in (\ref{Dbetkap}%
) (if we set $\delta =0$).

The rigorous proof of the superposition principle we present in this paper
is not based on the above argument since it implicitly relies on a
superposition principle in the form of an assumption that component
wavepackets can somehow maintain their identity, group velocities and
spatial extension during nonlinear evolution which by no means is obvious.
In fact, the question if a wavepacket or a multi-wavepacket structure can be
preserved during nonlinear evolution is important and interesting question
on its own right. The answer to it under natural conditions is affirmative
as we have shown in \cite{BF7}. Namely, if initially solution was a
multi-wavepacket at $\tau =0$, it remains a multi-wavepacket at $\tau >0$,
and every component wavepacket maintains its identity. Therefore a
wavepacket can be interpreted as a quasi-particle which maintains its
identity and can interact with other quasi-particles. This property holds
also in the situation when there are stronger nonlinear interactions between
wavepacket components which do not allow the superposition principle to
hold, see \cite{BF7} for details.

The proof we present here is based on general algebraic-functional
considerations. The strategy of our proof is as follows. First, we prove
that the operator $\mathcal{S}\left( \mathbf{h}\right) $ in (\ref{Gsum1}) is
analytical, i.e. it can be written in the form of a convergent series%
\begin{equation*}
\mathcal{S}\left( \mathbf{h}\right) =\dsum\nolimits_{j=1}^{\infty }\mathcal{S%
}^{\left( j\right) }\left( \mathbf{h}^{j}\right) ,\ \mathbf{h}^{j}=\mathbf{h}%
,\ldots ,\mathbf{h}\text{ (}j\text{ copies of }\mathbf{h}\text{)},
\end{equation*}%
where $\mathcal{S}^{\left( j\right) }\left( \mathbf{h}^{j}\right) $ is a $j$%
-linear operator applied to $\mathbf{h}$. Now we substitute $\mathbf{h}$ in $%
\mathcal{S}^{\left( j\right) }$ with the sum of $\mathbf{h}_{l}$ as in (\ref%
{hhh1}). Considering for simplicity the case $N=2$ and using the
polylinearity of $\mathcal{S}^{\left( j\right) }$ we get%
\begin{equation*}
\mathcal{S}^{\left( 2\right) }\left( \left( \mathbf{h}_{1}+\mathbf{h}%
_{2}\right) ^{2}\right) =\mathcal{S}^{\left( 2\right) }\left( \left( \mathbf{%
h}_{1}\right) ^{2}\right) +2\mathcal{S}^{\left( 2\right) }\left( \mathbf{h}%
_{1}\mathbf{h}_{2}\right) +\mathcal{S}^{\left( 2\right) }\left( \left( 
\mathbf{h}_{2}\right) ^{2}\right) ,\ldots ,
\end{equation*}%
implying after the summation 
\begin{eqnarray*}
\mathcal{S}\left( \mathbf{h}\right) &=&\mathcal{S}^{\left( 2\right) }\left(
\left( \mathbf{h}_{1}\right) ^{2}\right) +\mathcal{S}^{\left( 3\right)
}\left( \left( \mathbf{h}_{1}\right) ^{3}\right) +\ldots +\mathcal{S}%
^{\left( 2\right) }\left( \left( \mathbf{h}_{2}\right) ^{2}\right) +\mathcal{%
S}^{\left( 3\right) }\left( \left( \mathbf{h}_{2}\right) ^{3}\right) +\ldots 
\mathcal{S}_{cr} \\
&=&\mathcal{S}\left( \mathbf{h}_{1}\right) +\mathcal{S}\left( \mathbf{h}%
_{2}\right) +\mathcal{S}_{cr},
\end{eqnarray*}%
where $\mathcal{S}_{cr}$ is a sum of all cross-terms such as $\mathcal{S}%
^{\left( 2\right) }\left( \mathbf{h}_{1}\mathbf{h}_{2}\right) $ etc. The
main part of the proof is to show that every term in $\mathcal{S}_{cr}$ is
small. An important step \ for that is based on the refined implicit
function theorem (Theorem \ref{Theorem Monomial Convergence}) which allows
to represent the operators $\mathcal{S}^{\left( j\right) }$ in the form of a
sum of certain composition monomials, which, in turn, have a relatively
simple oscillatory integral representation. Importantly, the relevant
oscillatory integrals involve the known initial data $\mathbf{h}_{l}$ rather
than unknown solution $\mathbf{U}$. The analysis of the oscillatory
integrals shows that there are two mechanisms responsible for the smallness
of the integrals. The first one is \emph{time averaging}, and the second one
is based on large group velocities (in the slow time scale) of wavepackets. 
\emph{Remarkably, if wavepackets satisfy proper genericity conditions, every
cross term is small due one of the above mentioned two mechanisms.
Importantly, the both mechanism are instrumental for the smallness of terms
in }$S_{cr}$\emph{, and the time averaging alone is not sufficient}. We
obtain estimates on terms in $\mathcal{S}_{cr}$ which ultimately yield the
estimate (\ref{Dbetkap}). \emph{Since the smallness of interactions between
waves under nonlinear evolution stems from high frequency oscillations in
time and space of functions involved in the interaction integrals, we can
interpret it as a result of the destructive wave interference}. The above
sketch shows that the mathematical tools we use in our studies are (i) the
theory of analytic functions and corresponding series of
infinite-dimensional (Banach) variable, and (ii) the theory of oscillatory
integrals.\ 

We would like to point out that the estimate (\ref{Dbetkap}) for the
remainder in the superposition principle is quite accurate. For example,
when the estimate is applied to the sine-Gordon equation with bimodal
initial data, it yields essentially optimal estimates for the magnitude of
the interaction of counterpropagating waves. These estimates are more
accurate than ones obtained by the well known \emph{ansatz method} as in 
\cite{PW}, and the comparative analysis is provided below in Example 1,
Section 2.2.

To summarize the above analysis we list important ingredients of our
approach.

\begin{itemize}
\item The spectrum of the underlying linear problem is \emph{continuous.}

\item The wave nonlinear evolution is analyzed based on the modal
decomposition with respect to the linear component of the system because
there is no exchange of energy between modes by linear mechanisms.
Wavepacket definition is based on the modal expansion determining, in
particular, its the spatial extension and the group velocity..

\item The problem involves two small parameters $\beta $ and $\varrho $
respectively in the initial data and coefficients of the equations. These
parameters scale respectively (i) the range of wavevectors involved in its
modal composition, with $\beta ^{-1}$ scaling its spatial extension, and
(ii) $\varrho $ scaling the ratio of the slow and the fast time scales. We
make no assumption on the functional dependence between $\beta $ and $%
\varrho $, which are essentially independent and are subject only to
inequalities.

\item The nonlinear \emph{evolution is studied for a finite time} $\tau
_{\ast }$ which may depend on, say, the amplitude of the initial excitation,
and, importantly, $\tau _{\ast }$ is long enough to observe appreciable
nonlinear phenomena which are not vanishingly small. The superposition
principle can be extended to longer time intervals up to blow-up time or
even infinity if relevant uniform in $\beta $ and $\varrho $ estimates of
solutions in appropriate norms are available.

\item Two fast wave processes (in the chosen slow time scale) attributed to
the linear operator $\mathbf{L}$ and having typical time scale of order $%
\varrho $ can be identified as responsible for the essential independence of
wavepackets: (i) fast time oscillations which lead to time averaging; (ii)
fast wavepacket propagation with large group velocities produce effective
weakening of interactions which are\ not subjected to time averaging.
\end{itemize}

The rest of the paper is organized as follows. In the following Section 2 we
formulate exact conditions and theorems for lattice equations and partial
differential equations and give examples. In Section 3 we recast the
original evolution equation in a convenient reduced form allowing, in
particular, to construct a representation of the solution in a form of
convergent functional operator series explicitly involving the equation
nonlinear term. In Section 4 we provide the detailed analysis of
function-analytic series used to get a constructive representation of the
solution. Section 5 is devoted to the analysis of certain oscillatory
integrals which are terms of the series representing the solution. Note that
when making estimations we use the same letter $C$ for different constants
in different statements. Finally, the proofs of Theorems \ref{Theorem
Superposition} and \ref{Theorem Superposition1} are provided in Section 6.
more examples and generalizations are given in Section 7. For reader's
convenience we provide a list of notations in the end of the paper.

\section{Statement of results}

In this section we consider two classes of problems: lattice equations and
partial differential equations. After Fourier transform they can be written
in the modal form which is essentially the same in both cases. We formulate
the exact conditions on the modal equations and present the main theorems on
the superposition principle. We also give examples of equations to which the
general theorems apply, in particular Fermi-Pasta-Ulam system and Nonlinear
Schrodinger equation.

\subsection{Main definitions, statements and examples for the lattice
equation}

The first class of evolutionary systems we consider involves systems of
equations describing coupled nonlinear oscillators on a lattice $\mathbb{Z}%
^{d}$, namely the following lattice system of ordinary differential
equations (ODE's) with respect to time%
\begin{equation}
\partial _{\tau }\mathbf{U}\left( \mathbf{\mathbf{m}},\tau \right) =-\frac{%
\mathrm{i}}{\varrho }\mathbf{LU}\left( \mathbf{m},\tau \right) +F\left( 
\mathbf{U}\right) \left( \mathbf{\mathbf{m}},\tau \right) ,\ \mathbf{U}%
\left( \mathbf{\mathbf{m}},0\right) =\mathbf{h}\left( \mathbf{m}\right) ,\ 
\mathbf{\mathbf{m}}\in \mathbb{Z}^{d},  \label{lat0}
\end{equation}%
where $\mathbf{L}$ is a linear operator, $F$ is a nonlinear operator and $%
\varrho >0$ is a small parameter (see \cite{BF6}). To analyze the evolution
equation (\ref{lat0}) it is instrumental to recast it in the modal form (the
wavevector domain), in other words, to apply to it the \emph{lattice Fourier
transform} as defined by the formula 
\begin{equation}
\mathbf{\tilde{U}}\left( \mathbf{\mathbf{k}}\right) =\sum_{\mathbf{\mathbf{m}%
}\in \mathbb{Z}^{d}}\mathbf{U}\left( \mathbf{\mathbf{m}}\right) \mathrm{e}^{-%
\mathrm{i}\mathbf{\mathbf{m}}\cdot \mathbf{\mathbf{\mathbf{\mathbf{k}}}}},\ 
\text{where }\mathbf{k}\in \left[ -\pi ,\pi \right] ^{d}\text{,}
\label{Fourintr}
\end{equation}%
$\mathbf{k}$ is called a wave vector. \ We assume that the Fourier
transformation of the original lattice evolutionary equation (\ref{lat0}) is
of the form 
\begin{equation}
\partial _{\tau }\mathbf{\tilde{U}}\left( \mathbf{k},\tau \right) =-\frac{%
\mathrm{i}}{\varrho }\mathbf{L}\left( \mathbf{k}\right) \mathbf{\tilde{U}}%
\left( \mathbf{k},\tau \right) +\tilde{F}\left( \mathbf{\tilde{U}}\right)
\left( \mathbf{k},\tau \right) ;\ \mathbf{\tilde{U}}\left( \mathbf{k}%
,0\right) =\mathbf{\tilde{h}}\left( \mathbf{k}\right) \ \text{for }\tau =0.
\label{eqFur}
\end{equation}%
Here, $\mathbf{\tilde{U}}\left( \mathbf{k},\tau \right) $ is $2J$- component
vector, $\mathbf{L}\left( \mathbf{k}\right) \mathbf{\ }$is a $\mathbf{k}$%
-dependent $2J\times 2J$ matrix that corresponds to the linear operator $%
\mathbf{L}$ and $\tilde{F}\left( \mathbf{\tilde{U}}\right) $ is a nonlinear
operator, which we describe later. The matrix $\mathbf{L}\left( \mathbf{k}%
\right) \mathbf{\ }$and the coefficients of the nonlinear operator $\tilde{F}%
\left( \mathbf{\tilde{U}}\right) $ in (\ref{eqFur}) are $2\pi $-periodic
functions of $\mathbf{k}$ and for that reason we assume that $\mathbf{k}$\
belongs to the torus $\mathbb{R}^{d}/\left( 2\pi \mathbb{Z}\right) ^{d}$
which we denote by $\left[ -\pi ,\pi \right] ^{d}$. The $\mathbf{k}$%
-dependent matrix $\mathbf{L}\left( \mathbf{k}\right) $\textbf{\ }determines
the linear operator $\mathbf{L}$ and plays an important role in the
analysis. We refer to $\mathbf{L}\left( \mathbf{k}\right) $ as to the \emph{%
linear symbol}. \ Since (\ref{eqFur}) \ describes evolution of the Fourier
modes \ of the solution, we call (\ref{eqFur}) \emph{modal evolution equation%
}.

We study the modal evolution equation (\ref{eqFur})\ on a finite time
interval 
\begin{equation}
0\leq \tau \leq \tau _{\ast },  \label{tstarrho}
\end{equation}%
where $\tau _{\ast }>0$ is a fixed number which, as we will see, may depend
on the magnitude of the initial data. The time $\tau _{\ast }$ does not
depend on small parameters, it is of order one and is determined by norms of
operators and initial data; it is almost optimal for general $F$ since there
are examples when $\ \tau _{\ast }$ is of the same order as the blow up time
of solutions. To make formulas and estimates simpler we assume without loss
of generality that 
\begin{equation}
\tau _{\ast }\leq 1.  \label{tau1}
\end{equation}%
For a number of reasons the modal form (\ref{eqFur}) of the evolution
equation is much more suitable for nonlinear analysis than the original
evolution equation (\ref{lat0}). This is why from now on \emph{we consider
the modal form of evolution equation (\ref{eqFur}) for the modal components }%
$\mathbf{\tilde{U}}\left( \mathbf{k},\tau \right) $\emph{\ as our primary
evolution equation}.

First, as an illustration, let us look at the simplest nontrivial example of
(\ref{eqFur}) with $J=1$ corresponding to two-component vector fields on the
lattice $\mathbb{Z}^{d}$. A two-component vector function $\mathbf{U}\left( 
\mathbf{\mathbf{m}}\right) $ of a discrete argument $\mathbf{\mathbf{m}}\in 
\mathbb{Z}^{d}$ has the form \ 
\begin{equation}
\mathbf{U\left( \mathbf{\mathbf{m}}\right) =}\left[ 
\begin{array}{c}
U_{+}\mathbf{\left( \mathbf{\mathbf{m}}\right) } \\ 
U_{-}\mathbf{\left( \mathbf{\mathbf{m}}\right) }%
\end{array}%
\right] ,\ \mathbf{\mathbf{m}}\in \mathbb{Z}^{d}.  \label{U2comp}
\end{equation}%
In this example $\mathbf{L}\left( \mathbf{k}\right) $ in (\ref{eqFur}) is a $%
2\times 2$ matrix, and we assume that for almost all $\mathbf{k}$ it has two
different real eigenvalues $\omega _{-}\left( \mathbf{k}\right) $ and $%
\omega _{+}\left( \mathbf{k}\right) $ (the dependence of $\omega _{\pm
}\left( \mathbf{k}\right) $ on $\mathbf{k}$ is called the \emph{dispersion
relation}) satisfying the relation $\omega _{-}\left( \mathbf{k}\right)
=-\omega _{+}\left( \mathbf{k}\right) $, namely, 
\begin{equation}
\mathbf{L}\left( \mathbf{k}\right) \mathbf{g}_{\zeta }\left( \mathbf{k}%
\right) =\omega _{\zeta }\left( \mathbf{k}\right) \mathbf{g}_{\zeta }\left( 
\mathbf{k}\right) ,\ \omega _{\zeta }\left( \mathbf{k}\right) =\zeta \omega
\left( \mathbf{k}\right) ,\ \zeta =\pm ,  \label{Omom}
\end{equation}%
where, evidently, $\mathbf{g}_{\zeta }\left( \mathbf{k}\right) $ are the
eigenvectors of $\mathbf{L}\left( \mathbf{k}\right) $. These eigenvalues $%
\omega _{\zeta }\left( \mathbf{k}\right) $,$\ \zeta =\pm $, are $2\pi $%
-periodic real valued functions 
\begin{equation}
\omega _{\zeta }\left( k_{1}+2\pi ,k_{2},\ldots ,k_{d}\right) =\ldots
=\omega _{\zeta }\left( k_{1},k_{2},\ldots ,k_{d}+2\pi \right) =\omega
_{\zeta }\left( k_{1},k_{2},\ldots ,k_{d}\right) .  \label{omper}
\end{equation}%
The simplest nonlinearity in (\ref{eqFur}) is a quadratic nonlinear operator 
$\tilde{F}\left( \mathbf{\tilde{U}}\right) =\tilde{F}^{\left( 2\right)
}\left( \mathbf{\tilde{U}}^{2}\right) $ which is given by the following
convolution integral%
\begin{equation}
\tilde{F}^{\left( 2\right) }\left( \mathbf{\tilde{U}}_{1}\mathbf{\tilde{U}}%
_{2}\right) \left( \mathbf{\mathbf{k}}\right) =\frac{1}{\left( 2\pi \right)
^{d}}\int\limits_{\mathbf{k}^{\prime }\in \left[ -\pi ,\pi \right] ^{d};\ 
\mathbf{k}^{\prime }+\mathbf{k}^{\prime \prime }=\mathbf{k}}\chi ^{\left(
2\right) }\left( \mathbf{\mathbf{k}},\vec{k}\right) \left( \mathbf{\tilde{U}}%
_{1}\left( \mathbf{k}^{\prime }\right) \mathbf{\tilde{U}}_{2}\left( \mathbf{k%
}^{\prime \prime }\right) \right) \,\mathrm{d}\mathbf{k}^{\prime },
\label{Fintr}
\end{equation}%
where $\vec{k}=\left( \mathbf{\mathbf{k}}^{\prime },\mathbf{\mathbf{k}}%
^{\prime \prime }\right) $, $\chi ^{\left( 2\right) }\left( \mathbf{\mathbf{k%
}},\vec{k}\right) $ is a quadratic tensor (susceptibility) which acts on
vectors $\mathbf{\tilde{U}}_{1},\mathbf{\tilde{U}}_{2}$. We refer to the \
case $J=1$ as the \emph{one-band case} since the corresponding linear
operator is described by a single function $\omega \left( \mathbf{k}\right) $%
.

A particular example of (\ref{eqFur}) is obtained as a Fourier transform of
the following \emph{Fermi-Pasta-Ulam equation} (FPU) (see \cite{BermanI}, 
\cite{Pankov}, \cite{Weissert97}) describing a nonlinear system of coupled\
oscillators:%
\begin{gather}
\partial _{\tau }x_{n}=\frac{1}{\varrho }\left( y_{n}-y_{n-1}\right) ,
\label{FPMint} \\
\partial _{\tau }y_{n}=\frac{1}{\varrho }\left( x_{n+1}-x_{n}\right) +\alpha
_{2}\left( x_{n+1}-x_{n}\right) ^{2}+\alpha _{3}\left( x_{n+1}-x_{n}\right)
^{3},\ n\in \mathbb{Z}.  \notag
\end{gather}%
Note that an equivalent form of (\ref{FPMint}) (with $\alpha _{2}=0$) is the
second order equation 
\begin{equation}
\partial _{\tau }^{2}x_{n}=\frac{1}{\varrho ^{2}}\left(
x_{n-1}-2x_{n}+x_{n+1}\right) +\frac{\alpha _{3}}{\varrho }\left( \left(
x_{n+1}-x_{n}\right) ^{3}-\left( x_{n}-x_{n-1}\right) ^{3}\right) .
\label{FPMin2}
\end{equation}%
In this example $d=1$, $\mathbf{k}=k$ and elementary computations show that
the Fourier transform of the FPU equation (\ref{FPMint}) has the form of the
modal evolution equation (\ref{eqFur}), (\ref{Fintr}) where 
\begin{gather}
\mathbf{\tilde{U}}=\left[ 
\begin{array}{c}
\tilde{x} \\ 
\tilde{y}%
\end{array}%
\right] ,\ \mathrm{i}\mathbf{L}\left( k\right) =\left[ 
\begin{array}{cc}
0 & -\left( 1-\mathrm{e}^{-\mathrm{i}k}\right) ^{\ast } \\ 
\left( 1-\mathrm{e}^{-\mathrm{i}k}\right) & 0%
\end{array}%
\right] ,\ \omega _{\zeta }\left( k\right) =2\zeta \left\vert \sin \frac{k}{2%
}\right\vert ,\   \label{FPMom} \\
\chi ^{\left( 2\right) }\left( k,k^{\prime },k^{\prime \prime }\right) 
\mathbf{\tilde{U}}_{1}\left( k^{\prime }\right) \mathbf{\tilde{U}}_{2}\left(
k^{\prime \prime }\right) =\alpha _{2}\left( 1-\mathrm{e}^{-\mathrm{i}%
k^{\prime }}\right) \left( 1-\mathrm{e}^{-\mathrm{i}k^{\prime \prime
}}\right) \left[ 
\begin{array}{c}
0 \\ 
\tilde{x}_{1}\left( k^{\prime }\right) \tilde{x}_{2}\left( k^{\prime \prime
}\right)%
\end{array}%
\right] ,  \notag
\end{gather}%
and a similar formula for $\chi ^{\left( 3\right) }$ (see (\ref{chi3})).

Now let us consider the general multi-component vector case with $J>1$ which
we refer to as $J$\emph{-band case} for which the system (\ref{eqFur}) has $%
2J$ components, and instead of (\ref{Omom}) we assume that $\mathbf{L}\left( 
\mathbf{\mathbf{k}}\right) $ has eigenvalues and eigenvectors as follows:%
\begin{equation}
\mathbf{L}\left( \mathbf{\mathbf{k}}\right) \mathbf{g}_{n,\zeta }\left( 
\mathbf{k}\right) =\omega _{n,\zeta }\left( \mathbf{k}\right) \mathbf{g}%
_{n,\zeta }\left( \mathbf{k}\right) ,\ \omega _{n,\zeta }\left( \mathbf{k}%
\right) =\zeta \omega _{n}\left( \mathbf{k}\right) ,\ \zeta =\pm ,\
n=1,\ldots ,J,  \label{OmomL}
\end{equation}%
where $\omega _{n}\left( \mathbf{k}\right) $ are real-valued, \ continuous
for all $\mathbf{k}$ \ functions, and eigenvectors $\mathbf{g}_{n,\zeta
}\left( \mathbf{k}\right) \in \mathbb{C}^{2J}$ have unit length in the
standard Euclidean norm. We also suppose that the eigenvalues are numbered
so that 
\begin{equation}
\omega _{n+1}\left( \mathbf{k}\right) \geq \omega _{n}\left( \mathbf{k}%
\right) \geq 0,\ n=1,\ldots ,J-1,  \label{omgr0}
\end{equation}%
and we call $n$ the \emph{band index}. Note that the presence of $\zeta =\pm 
$ reflects a symmetry of the system allowing it, in particular, to have
real-valued solutions. Such a symmetry of dispersion relation $\omega
_{n}\left( \mathbf{k}\right) $ occurs in photonic crystals and many other
physical problems.

Note that (\ref{OmomL}) implies that the following symmetry relation hold:%
\begin{equation}
\omega _{n,-\zeta }\left( \mathbf{k}\right) =-\omega _{n,\zeta }\left( 
\mathbf{k}\right) ,n=1,\ldots ,J.  \label{invsym}
\end{equation}%
We also always assume that the following inversion symmetry holds:%
\begin{equation}
\omega _{n,\zeta }\left( -\mathbf{k}\right) =\omega _{n,\zeta }\left( 
\mathbf{k}\right) .  \label{omeven}
\end{equation}

\begin{remark}
Assuming (\ref{invsym})and \ (\ref{omeven}) we suppose that the dispersion
relations $\omega _{\zeta }\left( \mathbf{k}\right) $ have the same symmetry
properties as the dispersion relations of Maxwell equations in periodic
media, see \cite{BF1}-\cite{BF3}, \cite{BF5}. We would like to stress that 
\emph{these} \emph{symmetry conditions are not imposed for technical reasons}
but because they are consequences of fundamental symmetries of physical
media. Such symmetries arise in many problems including, for instance, the
Fermi-Pasta-Ulam equation, or when $\mathbf{L}\left( \mathbf{k}\right) $
originates from a Hamiltonian $H\left( p,q\right) =\frac{1}{2}\left(
H_{1}\left( p^{2}\right) \right) +\frac{1}{2}H_{2}\left( q^{2}\right) $. In
the opposite case if it is assumed that (\ref{invsym})and (\ref{omeven})
never hold, the results of this paper hold and the proofs, in fact, are
simpler. The case with the symmetry is more difficult and delicate because
of a possibility of resonant nonlinear interactions.
\end{remark}

There are values of $\mathbf{\mathbf{k}}$ \ for which inequalities (\ref%
{omgr0}) turn into equalities, these points require special treatment.

\begin{definition}[band-crossing points]
\label{Definition band-crossing point copy(1)} We call $\mathbf{k}_{0}$ a 
\emph{band-crossing point} if $\omega _{n+1}\left( \mathbf{k}_{0}\right)
=\omega _{n}\left( \mathbf{k}_{0}\right) $ \ for some $n$ or $\omega
_{1}\left( \mathbf{k}_{0}\right) =0$ and denote the set of \ band-crossing
points by $\sigma $.
\end{definition}

Everywhere in this paper we assume that the following condition is satisfied.

\begin{condition}
\label{Definition band-crossing point} The set $\sigma $ of \ band-crossing
points\ is a closed nowhere dense set in $\mathbb{R}^{d}$ with zero Lebesgue
measure, the entries of the matrix $\mathbf{L}\left( \mathbf{\mathbf{k}}%
\right) $ are infinitely smooth functions of $\mathbf{k}\notin \sigma $ and $%
\omega _{n}\left( \mathbf{k}\right) $ are continuous functions of $\mathbf{k}
$\ for all $\mathbf{k}$ and are infinitely smooth when $\mathbf{k}\notin
\sigma $. \ 
\end{condition}

Observe that for $\mathbf{k}\notin \sigma $ all the eigenvalues of the
matrix $\mathbf{L}\left( \mathbf{\mathbf{k}}\right) $ are different and the
corresponding eigenvectors $\mathbf{g}_{n,\zeta }\left( \mathbf{k}\right) $
of $\mathbf{L}\left( \mathbf{\mathbf{k}}\right) $\ can be locally defined as
smooth functions of $\mathbf{k}\notin \sigma $ as long as $\mathbf{L}\left( 
\mathbf{\mathbf{k}}\right) $ is smooth.

\begin{remark}
The band-crossing points are discussed in more details in \cite{BF1}, \cite%
{BF2}. Here we only note that generically the singular set $\sigma $ is a
manifold of the dimension $d-2$, see \cite{BF1}, \cite{BF2}). A simple
example of a band-crossing point is $k=0$ in (\ref{FPMom}).
\end{remark}

Since we do not assume the matrix $\mathbf{L}\left( \mathbf{\mathbf{k}}%
\right) $ to be Hermitian, we impose the following condition on its
eigenfunctions which guarantees its uniform diagonalization.

\begin{condition}
\label{Diagonalization} We assume that the $2J\times 2J$ matrix formed by
the eigenvectors $\mathbf{g}_{n,\zeta }\left( \mathbf{k}\right) $ of $%
\mathbf{L}\left( \mathbf{\mathbf{k}}\right) $, namely, 
\begin{equation*}
\Xi \left( \mathbf{k}\right) =\left[ \mathbf{g}_{1,+}\left( \mathbf{k}%
\right) ,\mathbf{g}_{1,-}\left( \mathbf{k}\right) ,\ldots ,\mathbf{g}%
_{J,+}\left( \mathbf{k}\right) ,\mathbf{g}_{J,+}\left( \mathbf{k}\right) %
\right]
\end{equation*}%
is uniformly bounded together with its inverse%
\begin{equation}
\sup_{\mathbf{k}\notin \sigma }\left\Vert \Xi \left( \mathbf{k}\right)
\right\Vert ,\ \sup_{\mathbf{k}\notin \sigma }\left\Vert \Xi ^{-1}\left( 
\mathbf{k}\right) \right\Vert \leq C_{\Xi }\text{ for some constant }C_{\Xi
}.  \label{supxi}
\end{equation}%
Here and everywhere we use the standard Euclidean norm in $\mathbb{C}^{2J}$.
\end{condition}

Note that if the matrix $\mathbf{L}\left( \mathbf{\mathbf{k}}\right) $ is
Hermitian for every $\mathbf{\mathbf{k}}$, the eigenvectors form an
orthonormal system. \ Then the matrix $\Xi $, which diagonalizes $\mathbf{L}$%
, is unitary and (\ref{supxi}) is satisfied with $C_{\Xi }=1$. \emph{%
Everywhere throughout the paper we assume that Condition \ref%
{Diagonalization} is satisfied.}

We introduce for vectors $\mathbf{\tilde{u}}\in \mathbb{C}^{2J}$ their
expansion with respect to the basis $\mathbf{g}_{n,\zeta }$: 
\begin{equation}
\mathbf{\tilde{u}}\left( \mathbf{k}\right) =\sum_{n=1}^{J}\sum_{\zeta =\pm }%
\tilde{u}_{n,\zeta }\left( \mathbf{k}\right) \mathbf{g}_{n,\zeta }\left( 
\mathbf{k}\right) =\sum_{n=1}^{J}\sum_{\zeta =\pm }\mathbf{\tilde{u}}%
_{n,\zeta }\left( \mathbf{k}\right) ,  \label{Uboldj}
\end{equation}%
and we refer to it as the \emph{modal decomposition} of $\mathbf{\tilde{u}}%
\left( \mathbf{k}\right) $, and call the coefficients $\tilde{u}_{n,\zeta
}\left( \mathbf{k}\right) $ the \emph{modal coefficients} of $\mathbf{\tilde{%
u}}\left( \mathbf{k}\right) $. In this expansion we assign to every $n,\zeta 
$ a linear projection $\Pi _{n,\zeta }\left( \mathbf{\mathbf{k}}\right) $ in 
$\mathbb{C}^{2J}$ corresponding to $\mathbf{g}_{n,\zeta }\left( \mathbf{k}%
\right) $, namely 
\begin{equation}
\Pi _{n,\zeta }\left( \mathbf{\mathbf{k}}\right) \mathbf{\tilde{u}}\left( 
\mathbf{k}\right) =\tilde{u}_{n,\zeta }\left( \mathbf{k}\right) \mathbf{g}%
_{n,\zeta }\left( \mathbf{k}\right) =\mathbf{\tilde{u}}_{n,\zeta }\left( 
\mathbf{k}\right) ,\ n=1,\ldots ,J,\ \zeta =\pm .  \label{Pin}
\end{equation}%
Note that these \emph{projections may be not orthogonal} if $\mathbf{L}%
\left( \mathbf{\mathbf{k}}\right) $ is not Hermitian. Evidently the
projections $\Pi _{n,\zeta }\left( \mathbf{\mathbf{k}}\right) $ are
determined by the matrix $\mathbf{L}\left( \mathbf{\mathbf{k}}\right) $ and
therefore do not depend on the choice of the basis $\mathbf{g}_{n,\zeta
}\left( \mathbf{\mathbf{k}}\right) $. Projections $\Pi _{n,\zeta }\left( 
\mathbf{\mathbf{k}}\right) $ depend smoothly on $\mathbf{\mathbf{k}}\notin
\sigma $ \ (note that we do not assume that the basis elements $\mathbf{g}%
_{n,\zeta }\left( \mathbf{k}\right) $ are defined globally as smooth
functions for all $\mathbf{\mathbf{k}}\notin \sigma $, in fact band-crossing
points may be branching points for eigenfunctions, see for example \cite{BF1}%
.) They are also uniformly bounded thanks to Condition \ref{Diagonalization}:%
\begin{equation}
C_{\Xi }^{-1}\left\vert \mathbf{V}\right\vert \leq \left( \sum_{n,\zeta
}\left\vert \Pi _{n,\zeta }\left( \mathbf{\mathbf{k}}\right) \mathbf{V}%
\right\vert ^{2}\right) ^{1/2}\leq C_{\Xi }\left\vert \mathbf{V}\right\vert
,\ \mathbf{V}\in \mathbb{C}^{2J},\ \mathbf{k}\notin \sigma .  \label{normP}
\end{equation}%
We would like to point out that most of the quantities are defined outside
of the singular set $\sigma $ of band-crossing points. It is sufficient
since we consider $\mathbf{\tilde{U}}\left( \mathbf{k}\right) $ as an
element of the space $L_{1}$ of Lebesgue integrable functions and the set $%
\sigma $ has zero Lebesgue measure.

The class of nonlinearities $\tilde{F}$ in (\ref{eqFur}) which we consider
can be described as follows. $\tilde{F}$ is a general polynomial
nonlinearity of the form 
\begin{equation}
\tilde{F}\left( \mathbf{\tilde{U}}\right) =\sum_{m=2}^{m_{F}}\tilde{F}%
^{\left( m\right) }\left( \mathbf{\tilde{U}}^{m}\right) ,\text{ with }%
m_{F}\geq 2\text{,}  \label{Fseries}
\end{equation}%
where $m$-linear operators $\tilde{F}^{\left( m\right) }$ are represented by
integral convolution formulas similar to (\ref{Fintr}), namely%
\begin{equation}
\tilde{F}^{\left( m\right) }\left( \mathbf{\tilde{U}}_{1},\ldots ,\mathbf{%
\tilde{U}}_{m}\right) \left( \mathbf{k},\tau \right) =\int_{\mathbb{D}%
_{m}}\chi ^{\left( m\right) }\left( \mathbf{\mathbf{k}},\vec{k}\right) 
\mathbf{\tilde{U}}_{1}\left( \mathbf{k}^{\prime }\right) \ldots \mathbf{%
\tilde{U}}_{m}\left( \mathbf{k}^{\left( m\right) }\left( \mathbf{k},\vec{k}%
\right) \right) \,\mathrm{\tilde{d}}^{\left( m-1\right) d}\vec{k},
\label{Fmintr}
\end{equation}%
where the domain%
\begin{equation}
\mathbb{D}_{m}=\left[ -\pi ,\pi \right] ^{\left( m-1\right) d},  \label{Dm}
\end{equation}%
and we use notation%
\begin{equation}
\mathrm{\tilde{d}}^{\left( m-1\right) d}\vec{k}=\frac{1}{\left( 2\pi \right)
^{\left( m-1\right) d}}\,\mathrm{d}\mathbf{k}^{\prime }\ldots \,\mathrm{d}%
\mathbf{k}^{\left( m-1\right) }  \label{dtild}
\end{equation}%
and%
\begin{equation}
\mathbf{k}^{\left( m\right) }\left( \mathbf{k},\vec{k}\right) =\mathbf{k}-%
\mathbf{k}^{\prime }-\ldots -\mathbf{k}^{\left( m-1\right) },\ \vec{k}%
=\left( \mathbf{k}^{\prime },\ldots ,\mathbf{k}^{\left( m\right) }\right) .
\label{kkar}
\end{equation}

\begin{condition}[nonlinearity regularity]
\label{cnonreg}The nonlinear operator $\tilde{F}\left( \mathbf{\tilde{U}}%
\right) $ defined by (\ref{Fseries}) satisfy 
\begin{equation}
\left\Vert \chi ^{\left( m\right) }\right\Vert =\frac{1}{\left( 2\pi \right)
^{\left( m-1\right) d}}\sup_{\mathbf{\mathbf{k}},\mathbf{k}^{\prime },\ldots
,\mathbf{k}^{\left( m\right) }}\left\Vert \chi ^{\left( m\right) }\left( 
\mathbf{\mathbf{k}},\mathbf{k}^{\prime },\ldots ,\mathbf{k}^{\left( m\right)
}\right) \right\Vert \leq C_{\chi },\ m=2,3,\ldots ,  \label{chiCR}
\end{equation}%
where, without loss of generality, we can assume that $C_{\chi }\geq 1$. The
norm $\left\vert \chi ^{\left( m\right) }\left( \mathbf{k},\vec{k}\right)
\right\vert $ of the tensor $\chi ^{\left( m\right) }$ with a fixed $\vec{k}$
as a $m$-linear operator from $\left( \mathbb{C}^{2J}\right) ^{m}$ into $%
\left( \mathbb{C}^{2J}\right) $ is defined by 
\begin{equation}
\left\vert \chi ^{\left( m\right) }\left( \mathbf{k},\vec{k}\right)
\right\vert =\sup_{\left\vert \mathbf{x}_{j}\right\vert \leq 1}\left\vert
\chi _{\ }^{\left( m\right) }\left( \mathbf{k},\vec{k}\right) \left( \mathbf{%
x}_{1},\ldots ,\mathbf{x}_{m}\right) \right\vert ,  \label{normchi0}
\end{equation}%
where as always, $\left\vert \cdot \right\vert $ stands for the standard
Euclidean norm. The tensors $\chi ^{\left( m\right) }\left( \mathbf{\mathbf{k%
}},\vec{k}\right) $ are assumed to be smooth functions of $\mathbf{k},%
\mathbf{k}^{\prime },\ldots ,\mathbf{k}^{\left( m\right) }\notin \sigma $,
namely for every compact $K\subset \mathbb{R}^{d}\setminus \sigma $ and for
all $m=2,3,\ldots .$ 
\begin{equation}
\left\vert \nabla ^{l}\chi ^{\left( m\right) }\left( \mathbf{\mathbf{k}},%
\mathbf{k}^{\prime },\ldots ,\mathbf{k}^{\left( m\right) }\right)
\right\vert \leq C_{K,l}\text{ if }\mathbf{\mathbf{k}},\mathbf{k}^{\prime
},\ldots ,\mathbf{k}^{\left( m\right) }\in K,l=1,2,\ldots ,  \label{gradchi}
\end{equation}%
where $\nabla ^{l}\chi ^{\left( m\right) }$ is the vector composed of \ all
partial derivatives of order $l$ of all components of \ the tensor $\chi
^{\left( m\right) }$ with respect to the variables $\mathbf{\mathbf{k}},%
\mathbf{k}^{\prime },\ldots ,\mathbf{k}^{\left( m\right) }$.
\end{condition}

>From now on all the nonlinear operators we consider are assumed to satisfy
the nonlinearity regularity Condition \ref{cnonreg}.

\begin{remark}
At first sight, since $\varrho $ is a small parameter, one might think that
the linear term in (\ref{lat0}) with the factor $\frac{1}{\varrho }$ is
dominant. But it is not that simple. Indeed, since all eigenvalues of $%
\mathbf{L}\left( \mathbf{k}\right) $ are purely imaginary the magnitude of $%
\mathrm{e}^{-\frac{\mathrm{i}}{\varrho }\mathbf{L}\left( \mathbf{k}\right) }%
\mathbf{\tilde{h}}\left( \mathbf{k}\right) \ $ which represents the solution
of a linear equation (with $\tilde{F}=0$) is bounded uniformly in $\varrho $%
. A nonlinearity $\tilde{F}$ alters the solution for a bounded time $\tau
_{\ast }$ which is not small for small $\varrho $. Therefore the influence
of the nonlinearity can be significant. This phenomenon can be illustrated
by the following toy model. Let us consider the partial differential
equation for a scalar function $y\left( x,\tau \right) $:%
\begin{equation*}
\partial _{\tau }y=-\frac{1}{\varrho }\partial _{x}y+y^{2},\;y\left(
x,0\right) =h\left( x\right) .
\end{equation*}%
Its solution is of the form 
\begin{equation}
y\left( x,\tau \right) =\frac{h\left( x-\frac{\tau }{\varrho }\right) }{%
1-\tau h\left( x-\frac{\tau }{\varrho }\right) },  \label{yxt}
\end{equation}%
and regularly it exists only for a finite time. The solution (\ref{yxt})
shows that the \ large coefficient $\frac{1}{\varrho }$ enters it so that
the corresponding wave moves faster with the velocity $\frac{1}{\varrho }$
along the $x$-axis but the wave's shape does not depend on $\varrho $ at
all. For the NLS with the initial data $\mathbf{\tilde{h}}\left( \mathbf{k}%
\right) =\mathbf{\tilde{h}}\left( \mathbf{k},\beta \right) $, $\varrho
=\beta ^{2}$, and the coefficient $\frac{1}{\varrho }$ at the linear part,
the nonlinearity balances the effect of dispersion leading to emergence of
solitons, see \cite{BF6} \ for a discussion.
\end{remark}

To formulate our results we introduce a Banach space $E=C\left( \left[
0,\tau _{\ast }\right] ,L_{1}\right) $ of functions $\mathbf{\tilde{v}}%
\left( \mathbf{k},\tau \right) $, $0\leq \tau \leq \tau _{\ast }$, with the
norm 
\begin{equation}
\left\Vert \mathbf{\tilde{v}}\left( \mathbf{k},\tau \right) \right\Vert
_{E}=\left\Vert \mathbf{\tilde{v}}\left( \mathbf{k},\tau \right) \right\Vert
_{C\left( \left[ 0,\tau _{\ast }\right] ,L_{1}\right) }=\sup_{0\leq \tau
\leq \tau _{\ast }}\int_{\left[ -\pi ,\pi \right] ^{d}}\left\vert \mathbf{%
\tilde{v}}\left( \mathbf{k},\tau \right) \right\vert \,\mathrm{d}\mathbf{k}.
\label{Elat}
\end{equation}%
Here $L_{1}$ is the Lebesgue function space with the standard norm defined
by the formula 
\begin{equation}
\left\Vert \mathbf{\tilde{v}}\left( \mathbf{\cdot }\right) \right\Vert
_{L_{1}}=\int_{\left[ -\pi ,\pi \right] ^{d}}\left\vert \mathbf{\tilde{v}}%
\left( \mathbf{k}\right) \right\vert \mathrm{d}\mathbf{k}.  \label{L1}
\end{equation}%
The following theorem guarantees the existence and the uniqueness of a
solution to the modal evolution equation (\ref{eqFur}) on a time interval
which does not depend on $\varrho $ (see Theorem \ref{Imfth1} for details).

\begin{theorem}[existence and uniqueness]
\label{Theorem Existence}\ Let the modal evolution equation (\ref{eqFur})
satisfy the Condition \ref{Diagonalization}, and let $\mathbf{\tilde{h}}\in
L_{1}$, $\left\Vert \mathbf{\tilde{h}}\right\Vert _{L_{1}}\leq R$. Then
there exists a unique solution $\mathbf{\tilde{U}}=\mathcal{G}\left( \mathbf{%
\tilde{h}}\right) $ of (\ref{eqFur}) which belongs to $C^{1}\left( \left[
0,\tau _{\ast }\right] ,L_{1}\right) $. The number $\tau _{\ast }>0$ depends
on $R$, $C_{\chi }$ and $C_{\Xi }$ and it does not depend on $\varrho $.
\end{theorem}

Now we would like to formulate the main result of this paper, a theorem on
the superposition principle, showing that the generic wavepackets evolve
almost independently for the case of lattice equations. To do that, first,
we define an important concept of \emph{wavepacket. }

\begin{definition}[wavepacket]
\label{dwavepack} A function $\mathbf{\tilde{h}}\left( \beta ,\mathbf{k}%
\right) $ \ which depends on a parameter $0<\beta <1$, is called a \emph{\
wavepacket with a center }$\mathbf{k}_{\ast }$ if it satisfies the following
conditions:

\begin{enumerate}
\item[(i)] It is bounded in $L_{1}$ uniformly in $\beta $, i.e.%
\begin{equation}
\left\Vert \mathbf{\tilde{h}}\left( \beta ,\mathbf{\cdot }\right)
\right\Vert _{L_{1}}\leq C_{h}.  \label{L1b}
\end{equation}

\item[(ii)] It is composed of modes from essentially a single band $n$,
namely for any $0<\epsilon <1$ there is a constant $C_{\epsilon }>0$ such
that 
\begin{equation}
\left\Vert \mathbf{\tilde{h}}\left( \mathbf{k}\right) -\mathbf{\tilde{h}}%
_{-}\left( \mathbf{k}\right) -\mathbf{\tilde{h}}_{+}\left( \mathbf{k}\right)
\right\Vert _{L_{1}}\leq C_{\epsilon }\beta ,\ \mathbf{\tilde{h}}_{\zeta
}\left( \mathbf{k}\right) =\Pi _{n,\zeta }\mathbf{\tilde{h}}\left( \mathbf{k}%
\right) ,\ \zeta =\pm ,  \label{hbold}
\end{equation}%
and $\mathbf{\tilde{h}}_{\zeta }\left( \beta ,\mathbf{k}\right) $ is
essentially supported in a small vicinity of $\zeta \mathbf{k}_{\ast }$,
where $\mathbf{k}_{\ast }$ is the \emph{wavepacket center}, namely 
\begin{equation}
\int_{\left\vert \mathbf{k-}\zeta \mathbf{k}_{\ast }\right\vert \geq \beta
^{1-\epsilon }}\left\vert \mathbf{\tilde{h}}_{\zeta }\left( \beta ,\mathbf{k}%
\right) \right\vert \mathrm{d}\mathbf{k}\leq C_{\epsilon }\beta .
\label{sourloc}
\end{equation}

\item[(iii)] The wavepacket center $\mathbf{k}_{\ast }$ is not a
band-crossing point, that is $\mathbf{k}_{\ast l}\notin \sigma $, and the
following regularity condition holds: 
\begin{equation}
\int_{\left\vert \mathbf{k-}\zeta \mathbf{k}_{\ast }\right\vert \leq \beta
^{1-\epsilon }}\left\vert \nabla _{\mathbf{k}}\mathbf{\tilde{h}}_{\zeta
}\left( \beta ,\mathbf{k}\right) \right\vert \,\mathrm{d}\mathbf{k}\leq
C_{\epsilon }\beta ^{-1-\epsilon }.  \label{hbder}
\end{equation}

In the above conditions (ii) and (iii) $\ C_{\epsilon }$ does not depend on $%
\beta ,$ $0<\beta <1$.
\end{enumerate}
\end{definition}

The simplest example of a wavepacket in the sense of Definition \ref%
{dwavepack} is a function of the form%
\begin{equation}
\mathbf{\tilde{h}}_{\zeta }\left( \beta ,\mathbf{k}\right) =\beta ^{-d}\hat{h%
}_{\zeta }\left( \frac{\mathbf{k}-\zeta \mathbf{k}_{\ast }}{\beta }\right) 
\mathbf{g}_{n,\zeta }\left( \mathbf{k}\right) ,\ \zeta =\pm ,  \label{h0}
\end{equation}%
where $\hat{h}_{\zeta }\left( \mathbf{k}\right) $ is a Schwartz function,
that is an infinitely smooth, rapidly decaying function. Another typical and
natural example of a wavepacket $\mathbf{\tilde{h}}$ centered at $\mathbf{k}%
_{\ast }$ is readily provided by 
\begin{equation}
\mathbf{\tilde{h}}\left( \beta ,\mathbf{k}\right) =\Pi _{n,+}\left( \mathbf{%
\mathbf{k}}\right) \mathbf{\tilde{h}}_{0,+}\left( \beta ,\mathbf{k}\right)
+\Pi _{n,-}\left( \mathbf{\mathbf{k}}\right) \mathbf{\tilde{h}}_{0,-}\left(
\beta ,\mathbf{k}\right)  \label{ex0}
\end{equation}%
where $\ \mathbf{\tilde{h}}_{0,\zeta }\left( \beta ,\mathbf{k}\right) $ is
the lattice Fourier transform of the following function 
\begin{equation}
\mathbf{h}_{0,\zeta }\left( \mathbf{m},\beta \right) \mathbf{\ }=\mathrm{e}^{%
\mathrm{i}\zeta \mathbf{k}_{\ast }\cdot \mathbf{m}}\Phi _{\zeta }\left(
\beta \mathbf{m}-\mathbf{r}_{0}\right) \mathbf{g},\ \zeta =\pm ,  \label{ex}
\end{equation}%
where $\mathbf{g}$ is a vector in $\mathbb{C}^{2J}$, projection $\Pi
_{n,\zeta }$ is as in (\ref{Pin}) with some $n$, vectors $\mathbf{m},\ 
\mathbf{r}_{0}\in \mathbb{R}^{d}$ and $\Phi _{\zeta }\left( \mathbf{r}%
\right) $ being an \emph{arbitrary} Schwartz function (see Lemma \ref{Lemma
201}).

Our special interest is in the waves that are finite sums of wavepackets and
we refer to them as \emph{multi-wavepackets}.

\begin{definition}[multi-wavepacket]
\label{dmwavepack} A function $\mathbf{\tilde{h}}\left( \beta ,\mathbf{k}%
\right) $, $0<\beta <1$, is called a \emph{multi-wavepacket} if it is a
finite sum of wavepackets $\mathbf{\tilde{h}}_{l}$ as defined in Definition %
\ref{dwavepack}, namely 
\begin{equation}
\mathbf{\tilde{h}}\left( \beta ,\mathbf{k}\right) =\sum_{l=1}^{N_{h}}\mathbf{%
\tilde{h}}_{l}\left( \beta ,\mathbf{k}\right) ,  \label{JJ1}
\end{equation}%
and we call the set $\left\{ \mathbf{k}_{\ast l}\right\} $ of all the
centers $\mathbf{k}_{\ast l}$ of involved wavepackets center set of $\mathbf{%
\tilde{h}}$.
\end{definition}

In what follows we will be interested in \emph{generic multi-wavepackets}
such that their centers are \emph{generic}. The exact meaning of this is
provided below in the following conditions.

\begin{condition}[non-zero frequency]
\label{Condition non-zero} We assume that every center $\ \mathbf{\mathbf{k}}%
_{\ast l}$ of a wavepacket satisfies the following condition 
\begin{equation}
\omega _{n_{l}}\left( \mathbf{\mathbf{k}}_{\ast l}\right) \neq 0,\
l=1,\ldots ,N_{h}.  \label{omne00}
\end{equation}
\end{condition}

\begin{condition}[group velocity]
\label{Condition GV} We assume that all centers $\ \mathbf{\mathbf{k}}_{\ast
l}$, $l=1,\ldots ,N_{h}$, of the multi-wavepacket $\mathbf{\tilde{h}}$ as
defined in Definition \ref{dmwavepack} are not band-crossing points, and the
gradients $\nabla _{\mathbf{k}}\omega _{n_{l_{j}}}\left( \mathbf{k}_{\ast
l_{j}}\right) $ (called group velocities) at these points satisfy the
following condition 
\begin{equation}
\left\vert \nabla _{\mathbf{k}}\omega _{n_{l_{1}}}\left( \mathbf{k}_{\ast
l_{1}}\right) -\nabla _{\mathbf{k}}\omega _{n_{l_{2}}}\left( \mathbf{k}%
_{\ast l_{2}}\right) \right\vert \neq 0\text{ when \ }l_{1}\neq l_{2},
\label{NGVM}
\end{equation}%
indicating that the group velocities are different.
\end{condition}

We also want the functions (dispersion relations) $\omega _{n_{l}}\left( 
\mathbf{k}\right) $ to be non-degenerate in the sense that they are not
exactly linear, below we give exact conditions.

Consider the following equation for $n$ and $\theta $%
\begin{equation}
\theta \omega _{n_{l}}\left( \mathbf{k}_{\ast }\right) -\zeta \omega
_{n}\left( \theta \mathbf{k}_{\ast }\right) =0,\ \ \zeta =\pm 1,
\label{zomomz0}
\end{equation}%
where the admissible $\theta $ have the form 
\begin{equation}
\theta =\sum_{j=1}^{m}\zeta ^{\left( j\right) },\ \zeta ^{\left( j\right)
}=\pm 1,\;m\leq m_{F},  \label{suml}
\end{equation}%
$m_{F}$ is the same as in (\ref{Fseries}). In the case when in the series (%
\ref{Fseries}) some terms $\tilde{F}^{\left( m\right) }$ vanish, we take in (%
\ref{suml}) only $m$ corresponding to non-zero $\tilde{F}^{\left( m\right) }$%
.

\begin{condition}[non-degeneracy]
\label{Condition generic1} Given a point \ $\mathbf{k}_{\ast }=\mathbf{k}%
_{\ast l}$ and band $n_{l}$ we assume that dispersion relations $\omega
_{n}\left( \mathbf{k}\right) $ are such that all solutions $n$,$\theta $ of (%
\ref{zomomz0}) are necessarily of the form 
\begin{equation}
n=n_{l},\ \theta =\zeta .  \label{FMeq}
\end{equation}
\end{condition}

\begin{definition}[generic multi-wavepackets]
\label{Condition generic} A multi-wavepacket $\mathbf{\tilde{h}}$ as defined
in Definition \ref{dmwavepack} is called generic if the centers $\ \mathbf{%
\mathbf{k}}_{\ast l}$, $l=1,\ldots ,N_{h}$, of all wavepackets satisfy
Conditions \ref{Condition non-zero} and \ref{Condition GV}; and the
dispersion relations $\omega _{n}\left( \mathbf{k}\right) $ at every $%
\mathbf{\mathbf{k}}_{\ast l}$ and band $n_{l}$ satisfy Condition \ref%
{Condition generic1}.
\end{definition}

We introduce now the \emph{solution operator} $\mathcal{G}$ mapping the
initial data $\mathbf{\tilde{h}}$ into the solution $\mathbf{\tilde{U}}=%
\mathcal{G}\left( \mathbf{\tilde{h}}\right) $ of the modal evolution
equation (\ref{eqFur}); this operator is defined for $\left\Vert \mathbf{%
\tilde{h}}\right\Vert \leq R$ according to Theorem \ref{Theorem Existence}.
The main result of this paper for the lattice case is the following
statement.

\begin{theorem}[superposition principle for lattice equations]
\label{Theorem Superposition}\ Suppose that the initial data $\mathbf{\tilde{%
h}}$ of (\ref{eqFur}) is a multi-wavepacket of the form%
\begin{equation}
\mathbf{\tilde{h}}=\sum_{l=1}^{N_{h}}\mathbf{\tilde{h}}_{l},\;N_{h}\max_{l}%
\left\Vert \mathbf{\tilde{h}}_{l}\right\Vert _{L_{1}}\leq R,  \label{hsumR}
\end{equation}%
satisfying Definition \ref{dmwavepack}, where $\mathbf{\tilde{h}}$\ is
generic in the sense of Definition \ref{Condition generic}. Let us assume
that 
\begin{equation}
\frac{\beta ^{2}}{\varrho }\leq C,\ \text{with \ some \ }C,\text{\ \ }%
0<\beta \leq \frac{1}{2},\ 0<\varrho \leq \frac{1}{2}.  \label{scale1}
\end{equation}
Then the solution $\mathbf{\tilde{U}}=\mathcal{G}\left( \mathbf{\tilde{h}}%
\right) $ to the evolution equation (\ref{eqFur}) satisfies the following
approximate superposition principle 
\begin{equation}
\mathcal{G}\left( \sum_{l=1}^{N_{h}}\mathbf{\tilde{h}}_{l}\right)
=\sum_{l=1}^{N_{h}}\mathcal{G}\left( \mathbf{\tilde{h}}_{l}\right) +\mathbf{%
\tilde{D}},  \label{Gsum}
\end{equation}%
with a small remainder $\mathbf{\tilde{D}}\left( \tau \right) $ satisfying
the following estimate 
\begin{equation}
\sup_{0\leq \tau \leq \tau _{\ast }}\left\Vert \mathbf{\tilde{D}}\left( \tau
\right) \right\Vert _{L_{1}}\leq C_{\epsilon }\frac{\varrho }{\beta
^{1+\epsilon }}\left\vert \ln \beta \right\vert ,  \label{rem}
\end{equation}%
where $\epsilon $ is the same as in Definition \ref{dwavepack} and can be
arbitrary small, $\tau _{\ast }$ does not depend on $\beta ,\varrho $ \ and $%
\epsilon $.
\end{theorem}

The most common case when (\ref{scale1}) \ holds is $\varrho =\beta ^{2}$, a
discussion of different scalings is provided in \cite{BF6} and \cite{BF7}.

Observe that solutions to the original evolution equation (\ref{lat0}) with
the initial data (\ref{JJ1}), (\ref{ex}) satisfy the superposition principle
if the wave vectors $\mathbf{k}_{\ast l}$ in (\ref{ex}) satisfy (\ref{NGVM}%
), (\ref{zomomz0}) and $\Phi _{l}$ are Schwartz functions. It turns out,
that the evolution of every coefficient $\tilde{u}_{n,\zeta }\left( \mathbf{k%
}\right) $ of the solution as defined by (\ref{Uboldj}) can be accurately
approximated by a solution a relevant Nonlinear Schrodinger equation (NLS),
see \cite{GiaMielke}. Therefore Theorem \ref{Theorem Superposition} provides
a reduction of multi-wavepacket problem to several single-wavepacket
problems.

We also would like to stress that though $\beta $ is small the nonlinear
effects are not small. Namely, there can be a significant difference between
solutions of a nonlinear and the corresponding linear (with $F\left( \mathbf{%
U}\right) $ being set zero) equations with the same initial data for times $%
\tau =\tau _{\ast }$.

Recall that up to now we analyzed the nonlinear evolution in the modal form (%
\ref{eqFur}) for $\mathbf{\tilde{U}}\left( \mathbf{\mathbf{k}},\tau \right) $%
. To make a statement on the nonlinear evolution for the original evolution
equation (\ref{lat0}), i.e. in terms of the quantities $\mathbf{U}\left( 
\mathbf{\mathbf{m}},\tau \right) $, we introduce $\mathbf{U}\left( \mathbf{h}%
\right) \left( \mathbf{\mathbf{m}}\right) $ as the inverse Fourier transform
of the solution $\mathcal{G}\left( \mathbf{\tilde{h}}\right) \left( \mathbf{%
\mathbf{k}}\right) $ of the modal evolution equation (\ref{eqFur}). Recall
that the inverse Fourier transform corresponding to (\ref{Fourintr}) is
given by the formula%
\begin{equation}
\mathbf{U}\left( \mathbf{\mathbf{m}}\right) =\left( 2\pi \right) ^{-d}\int_{%
\left[ -\pi ,\pi \right] ^{d}}\mathrm{e}^{\mathrm{i}\mathbf{\mathbf{m}}\cdot 
\mathbf{\mathbf{\mathbf{k}}}}\mathbf{\tilde{U}}\left( \mathbf{\mathbf{k}}%
\right) \,\mathrm{d}\mathbf{k},  \label{LF}
\end{equation}%
and when applying the inverse Fourier transform we get back the original
lattice system (\ref{lat0}) from its modal form (\ref{eqFur}). The
convolution form of the nonlinearity makes the lattice system invariant with
respect to translations on the lattice $\mathbb{Z}^{d}$. Using Theorem \ref%
{Theorem Superposition} and applying the inverse Fourier transform together
with the inequality 
\begin{equation}
\left\Vert \mathbf{U}\right\Vert _{L_{\infty }}\leq \left( 2\pi \right)
^{-d}\left\Vert \mathbf{\tilde{U}}\right\Vert _{L_{1}}  \label{Linf}
\end{equation}%
we obtain the following statement.

\begin{corollary}
\label{Corollary Uh} Let the evolution equation (\ref{lat0}) be obtained as
the lattice Fourier transform of (\ref{eqFur}). If $\mathbf{h}$ is given by (%
\ref{ex}) where every $\Phi _{l,\zeta }\left( \mathbf{r}\right) $ is a\emph{%
\ } Schwartz function (that is an infinitely smooth, rapidly decaying
function) then $\mathbf{U}\left( \mathbf{h}\right) $ is a solution to the
evolution equation (\ref{lat0}). If $\mathbf{h}=\mathbf{h}_{1}+\ldots +%
\mathbf{h}_{N_{h}}$ and every $\mathbf{h}_{l}$ is given by (\ref{ex}) then
the \emph{approximate superposition principle holds:}%
\begin{equation}
\mathbf{U}\left( \mathbf{h}\right) =\mathbf{U}\left( \mathbf{h}_{1}\right)
+\ldots +\mathbf{U}\left( \mathbf{h}_{N_{h}}\right) +\mathbf{D},
\label{sumU}
\end{equation}%
with a small coupling remainder $\mathbf{D}\left( \tau \right) $ satisfying 
\begin{equation}
\sup_{0\leq \tau \leq \tau _{\ast }}\left\Vert \mathbf{D}\left( \tau \right)
\right\Vert _{L_{\infty }}\leq C_{\delta }^{\prime }\frac{\varrho }{\beta
^{1+\delta }},  \label{Linf1}
\end{equation}%
where $\delta >0$ can be taken arbitrary small.
\end{corollary}

As an application of Theorem \ref{Theorem Superposition} let us consider the 
\emph{Fermi-Pasta-Ulam equation} (\ref{FPMint}). We impose the initial
condition for (\ref{FPMint}) 
\begin{equation}
x_{n}\left( 0\right) =\sum_{l=1}^{n_{h}}\Psi _{0l}\left( \beta
n-r_{l}\right) \mathrm{e}^{\mathrm{i}\mathbf{k}_{\ast l}n}+\limfunc{cc},\
y_{n}\left( 0\right) =\sum_{l=1}^{n_{h}}\Psi _{1l}\left( \beta
n-r_{l}\right) \mathrm{e}^{\mathrm{i}\mathbf{k}_{\ast l}n}+\limfunc{cc},\
n\in \mathbb{Z},  \label{inFPMN}
\end{equation}%
where $\Psi _{0l}\left( r\right) ,\Psi _{1l}\left( r\right) $ are arbitrary
Schwartz functions, and $r_{l}$ are arbitrary real numbers, $\limfunc{cc}$
means complex conjugate to the preceding terms and assume that $\varrho
,\beta $ satisfy (\ref{scale1}). For any given $k_{\ast l}$ there are two
eigenvectors$\ \mathbf{g}_{\pm }\left( k_{\ast l}\right) $ of the matrix $%
\mathbf{L}\left( k_{\ast l}\right) $ \ in (\ref{FPMom}) given by (\ref{Gzet}%
) and corresponding terms in (\ref{inFPMN}) can be written as 
\begin{equation*}
\left[ 
\begin{array}{c}
\Psi _{0l} \\ 
\Psi _{1l}%
\end{array}%
\right] \mathrm{e}^{\mathrm{i}\mathbf{k}_{\ast l}n}=\left[ \Phi _{-,l}%
\mathbf{g}_{-}\left( k_{\ast l}\right) +\Phi _{+,l}\mathbf{g}_{+}\left(
k_{\ast l}\right) \right] \mathrm{e}^{\mathrm{i}\mathbf{k}_{\ast l}n}.
\end{equation*}%
In this case all requirements of Definition \ref{dmwavepack} are fulfilled,
and (\ref{inFPMN}) defines a multiwavepacket. Note that the multiwavepacket (%
\ref{inFPMN}) involves $N_{h}=2n_{h}$ wavepackets \ with $2n_{h}$ wavepacket
centers $\vartheta k_{\ast l},\vartheta =\pm $. To satisfy Condition \ref%
{Condition GV} the wavepacket centers $k_{\ast l}$ must satisfy 
\begin{equation}
\frac{\cos \frac{k_{\ast l}}{2}}{\left\vert \sin \frac{k_{\ast l}}{2}%
\right\vert }\neq \frac{\cos \frac{k_{\ast j}}{2}}{\left\vert \sin \frac{%
k_{\ast j}}{2}\right\vert }\text{ if }l\neq j.  \label{kFPU}
\end{equation}%
To check if the centers $k_{\ast l}$ satisfy Condition \ref{Condition
generic1} we consider the equation \ 
\begin{equation}
z\left\vert \sin \frac{k_{\ast l}}{2}\right\vert -\zeta \left\vert \sin
\left( z\frac{k_{\ast l}}{2}\right) \right\vert =0,\ z=\sum_{j=1}^{3}\zeta
^{\left( j\right) },\ \zeta ^{\left( j\right) }=\pm 1.  \label{gen0}
\end{equation}%
Evidently the possible values of $z$ are $-3,-1,1,3$. Since the equation $%
3\left\vert \sin \phi \right\vert =\left\vert \sin \left( 3\phi \right)
\right\vert $ has the only solution $\ \phi =0$ on $\left[ 0,\pi /2\right] $
the equation (\ref{gen0}) has the only solution $z=\zeta $. Consequently,
all points $k_{\ast l}\neq 0$ satisfy Condition \ref{Condition generic1},
and Theorem \ref{Theorem Superposition} applies. The initial data for a
single wavepacket solution \ have the form 
\begin{equation}
\left[ 
\begin{array}{c}
x_{\vartheta ,n,l}\left( 0\right) \\ 
y_{\vartheta ,n,l}\left( 0\right)%
\end{array}%
\right] =\Phi _{\vartheta ,l}\left( \beta n-r_{l}\right) \mathbf{g}%
_{\vartheta }\left( k_{\ast l}\right) +\limfunc{cc},\;n\in \mathbb{Z}%
,\;\vartheta =\pm .  \label{ins}
\end{equation}%
According to this theorem and Corollary \ref{Corollary Uh} the solution to (%
\ref{FPMint}), (\ref{inFPMN}) equals the sum of solutions of (\ref{FPMint})
with single wavepacket initial data, that is 
\begin{equation}
x_{n}\left( \tau \right) =\sum_{\vartheta =\pm
}\sum_{l=1}^{n_{h}}x_{\vartheta ,n,l}\left( \tau \right) +D_{1,n}\left( \tau
\right) ,\;y_{n}\left( \tau \right) =\sum_{\vartheta =\pm
}\sum_{l=1}^{n_{h}}y_{\vartheta ,n,l}\left( \tau \right) +D_{2,n}\left( \tau
\right) .  \label{sum}
\end{equation}%
where $D_{n}$ is a small remainder satisfying%
\begin{equation}
\sup_{0\leq \tau \leq \tau _{\ast }}\sup_{n}\left[ \left\vert D_{1,n}\left(
\tau \right) \right\vert +\left\vert D_{2,n}\left( \tau \right) \right\vert %
\right] \leq C_{\delta }\frac{\varrho }{\beta ^{1+\delta }}  \label{RFPU}
\end{equation}%
with arbitrarily small positive $\delta $. Hence, the following statement
holds.

\begin{theorem}[superposition for Fermi-Pasta-Ulam equation]
\label{FPU} If $\mathbf{\ }$every $\Phi _{l,\zeta }\left( \mathbf{r}\right) $
is a\emph{\ } Schwartz function, and the wavevectors $k_{\ast l}\neq 0$
satisfy (\ref{kFPU}), then the solution $x_{n}\left( \tau \right)
,y_{n}\left( \tau \right) $ of the initial value problem for the
Fermi-Pasta-Ulam equation (\ref{FPMint}) with multi-wavepacket initial
condition (\ref{inFPMN}) is a linear superposition of solutions $%
x_{n,l}\left( \tau \right) ,y_{n,l}\left( \tau \right) $ of the same
equation with single-wavepacket initial condition (\ref{ins}) up to a small
coupling term $D_{1,n}\left( \tau \right) ,D_{2,n}\left( \tau \right) $
satisfying (\ref{sum}), (\ref{RFPU}) with arbitrary small $\delta >0$ and $%
\tau _{\ast }$ which do not depend on $\beta ,\varrho ,\delta $.
\end{theorem}

\begin{figure}[tbph]
\scalebox{0.3}{\includegraphics[viewport=-360 60 750
600,clip]{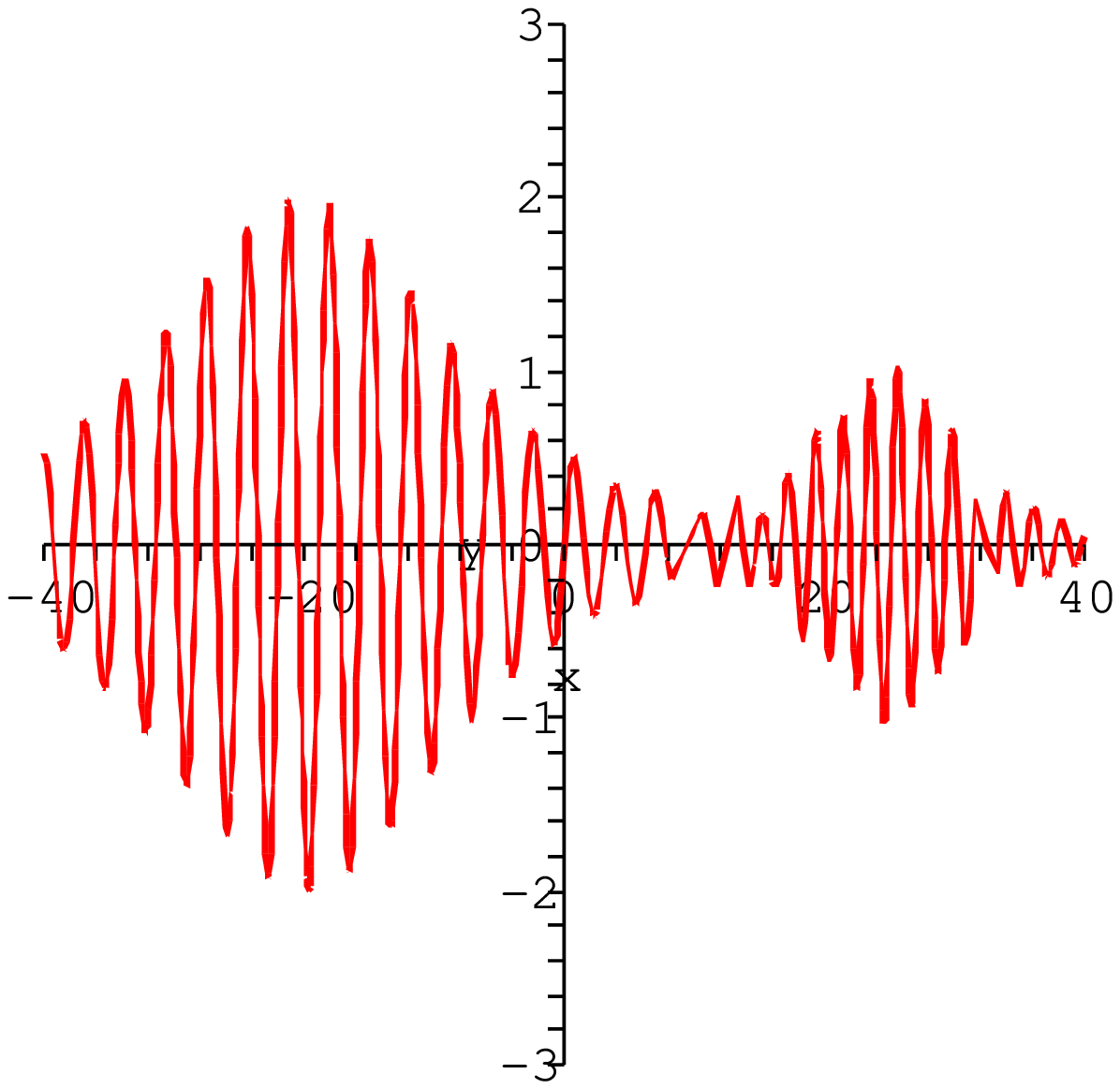}}
\caption{ In this picture two wavepackets are shown with different "centers" 
$k_{\ast 1}$ \ and $k_{\ast 2}$. The values of $k_{\ast 1}$ \ and $k_{\ast
2} $ are proportional to the frequences of spatial oscillations. Though the
wavepackets overlap in physical space, they pass one through another in the
process of nonlinear evolution almost without interaction if their group
velocities are different. }
\label{figtwopack}
\end{figure}

Note that solutions $x_{\vartheta ,n,l}\left( \tau \right) $ with different $%
\vartheta ,l$ \ resemble $2n_{h}$ solitons which originate at different
points $r_{l}$ and propagate with different group velocities. According to (%
\ref{sum}), (\ref{RFPU}) all these \emph{soliton-like wavepackets pass
through one another with very little interaction}, see Fig. \ref{figtwopack}%
. 
Note that Theorem \ref{Theorem Superposition} shows that this phenomenon is
robust in the class of general difference equations on the lattice $\mathbb{Z%
}$, and that it persists under polynomial perturbations of the nonlinearity
as well as perturbations of the linear part of the equation (\ref{FPMin2})
as long as they leave the linear difference operator non-positive and
self-adjoint. Observe also that the evolution of every single wavepacket is
nonlinear, and it is well-approximated by a properly constructed NLS (we
intend to write a proof of this statement for general lattice systems in
another article; see \cite{GiaMielke} for a particular case). For example,\
for a special choice of $\Psi _{jl}$ the solution $x_{n,l}\left( \tau
\right) $ can be well approximated by a soliton solution of a corresponding
NLS.

\subsection{Main statements and examples for semilinear systems of
hyperbolic PDE}

In this subsection we consider nonlinear evolution equation involving
partial differential (and pseudodifferential) operators with respect to
spatial variables with constant coefficients in the entire space $\mathbb{R}%
^{d}$. There is a great deal of similarity between such nonlinear evolution
PDE and the lattice nonlinear evolution equations considered in the previous
section. In particular, we study first not the original PDE but its Fourier
transform, modal evolution equation, and the results concerning the original
PDE are obtained by applying the inverse Fourier transform.

Recall that for functions $\mathbf{U}\left( \mathbf{r}\right) $ from $%
L_{1}\left( \mathbb{R}^{d}\right) $ the \emph{Fourier transform} and its
inverse are defined by the formulas%
\begin{equation}
\mathbf{\hat{U}}\left( \mathbf{k}\right) =\int_{\mathbb{R}^{d}}\mathbf{U}%
\left( \mathbf{r}\right) \mathrm{e}^{-\mathrm{i}\mathbf{r}\cdot \mathbf{k}}%
\mathrm{d}\mathbf{r},\text{ \ where }\mathbf{k}\in \mathbb{R}^{d},
\label{FoR}
\end{equation}%
\begin{equation}
\mathbf{U}\left( \mathbf{r}\right) =\frac{1}{\left( 2\pi \right) ^{d}}\int_{%
\mathbb{R}^{d}}\mathbf{\hat{U}}\left( \mathbf{k}\right) \mathrm{e}^{\mathrm{i%
}\mathbf{r}\cdot \mathbf{k}}\mathrm{d}\mathbf{r},\ \text{where }\mathbf{r}%
\in \mathbb{R}^{d}.  \label{Finv}
\end{equation}%
Similarly to (\ref{eqFur}) we introduce the following modal evolution
equation 
\begin{equation}
\partial _{\tau }\mathbf{\hat{U}}\left( \mathbf{\mathbf{k}},\tau \right) =-%
\frac{\mathrm{i}}{\varrho }\mathbf{L}\left( \mathbf{\mathbf{k}}\right) 
\mathbf{\hat{U}}\left( \mathbf{\mathbf{k}},\tau \right) +\hat{F}\left( 
\mathbf{\hat{U}}\right) \left( \mathbf{\mathbf{k}},\tau \right) ,\ \mathbf{%
\hat{U}}\left( \mathbf{\mathbf{k}},0\right) =\mathbf{\hat{h}}\left( \mathbf{%
\mathbf{k}}\right) ,\mathbf{k}\in \mathbb{R}^{d},  \label{difF}
\end{equation}%
where (i) $\mathbf{\hat{U}}\left( \mathbf{\mathbf{k}},\tau \right) $ is a $%
2J $-component vector-function of $\mathbf{\mathbf{k}}$, $\tau $, (ii) $%
\mathbf{L}\left( \mathbf{\mathbf{k}}\right) $ is a $2J\times 2J$ matrix
function of $\mathbf{\mathbf{k}}$, and (iii) $\hat{F}\left( \mathbf{\hat{U}}%
\right) $ is the nonlinearity. We assume that the $2J\times 2J$ matrix $%
\mathbf{L}\left( \mathbf{\mathbf{k}}\right) $, $\mathbf{k}\in \mathbb{R}^{d}$%
, has exactly $2J $ eigenvectors $\mathbf{g}_{n,\zeta }\left( \mathbf{k}%
\right) $ with corresponding $2J$ real eigenvalues $\omega _{n,\zeta }\left( 
\mathbf{k}\right) $ satisfying the relation (\ref{OmomL}), (\ref{omgr0}), (%
\ref{invsym}), (\ref{omeven}), (\ref{supxi}). We also assume the matrix $%
\mathbf{L}\left( \mathbf{\mathbf{k}}\right) $, $\mathbf{k}\in \mathbb{R}^{d}$%
, to satisfy the polynomial bound 
\begin{equation}
\left\vert \mathbf{L}\left( \mathbf{\mathbf{k}}\right) \right\vert \leq
C\left( 1+\left\vert \mathbf{\mathbf{k}}\right\vert ^{p}\right) .
\label{Lpol}
\end{equation}%
The singular set $\sigma $ for $\mathbf{L}\left( \mathbf{\mathbf{k}}\right) $
is as in Definition \ref{Definition band-crossing point} with the only
difference that functions $\omega _{n,\zeta }\left( \mathbf{k}\right) $ are
defined over $\mathbb{R}^{d}$ rather than the torus $\left[ -\pi ,\pi \right]
^{d}$, and, consequently they are not periodic. The nonlinearity $\hat{F}%
\left( \mathbf{\hat{U}}\right) $ has a form entirely similar to (\ref%
{Fseries}):%
\begin{equation}
\hat{F}\left( \mathbf{\hat{U}}\right) =\sum_{m=2}^{m_{F}}\hat{F}^{\left(
m\right) }\left( \mathbf{\hat{U}}^{m}\right) ,  \label{Fhat}
\end{equation}%
with $\hat{F}^{\left( m\right) }$ being $m$-linear operators with the
following representation similar to (\ref{Fmintr}):%
\begin{equation}
\hat{F}^{\left( m\right) }\left( \mathbf{\hat{U}}_{1},\ldots ,\mathbf{\hat{U}%
}_{m}\right) \left( \mathbf{k}\right) =\int_{\mathbb{D}_{m}}\chi ^{\left(
m\right) }\left( \mathbf{k},\vec{k}\right) \mathbf{\hat{U}}_{1}\left( 
\mathbf{k}^{\prime }\right) \ldots \mathbf{\hat{U}}_{m}\left( \mathbf{k}%
^{\left( m\right) }\left( \mathbf{k},\vec{k}\right) \right) \mathrm{\tilde{d}%
}^{\left( m-1\right) d}\vec{k},  \label{FY}
\end{equation}%
where $\mathbf{k}^{\left( m\right) }\left( \mathbf{k},\vec{k}\right) $ is
defined by the convolution equation (\ref{kkar}), $\mathrm{\tilde{d}}$ is
defined by (\ref{dtild}) and $\mathbb{D}_{m}$ in (\ref{FY}) is now defined
not by (\ref{Dm}) but by 
\begin{equation}
\mathbb{D}_{m}=\mathbb{R}^{\left( m-1\right) d}.  \label{DmR}
\end{equation}%
The difference with (\ref{eqFur}) now is that the involved functions of $%
\mathbf{k}$, $\mathbf{k}^{\prime }$ etc., are not $2\pi $-periodic, $\mathbb{%
D}_{m}$ in (\ref{FY}) is defined by (\ref{DmR}) instead of (\ref{Dm}), and
the tensors $\chi ^{\left( m\right) }\left( \mathbf{k},\vec{k}\right) $
satisfy the nonlinear regularity Condition \ref{cnonreg} without the
periodicity assumption. The functions $\mathbf{\hat{U}}_{l}\left( \mathbf{k}%
^{\left( l\right) }\right) $ in (\ref{FY}) are assumed to be from the space $%
L_{1}=L_{1}\left( \mathbb{R}^{d}\right) $ with the norm%
\begin{equation}
\left\Vert \mathbf{\hat{U}}\left( \mathbf{\cdot }\right) \right\Vert
_{L_{1}}=\int_{\mathbb{R}^{d}}\left\vert \mathbf{\tilde{v}}\left( \mathbf{k}%
\right) \right\vert \mathrm{d}\mathbf{k}.  \label{L1R}
\end{equation}%
We seek solutions to (\ref{difF}) in the space $C^{1}\left( \left[ 0,\tau
_{\ast }\right] ,L_{1}\right) $ with $0<\tau _{\ast }\leq 1$.

Applying the inverse Fourier transform to the modal evolution equation (\ref%
{difF}) we obtain a hyperbolic $2J$-component systems in $\mathbb{R}^{d}$ of
the form%
\begin{equation}
\partial _{\tau }\mathbf{U}\left( \mathbf{\mathbf{r}},\tau \right) =-\frac{%
\mathrm{i}}{\varrho }\mathbf{L}\left( -\mathrm{i}\nabla _{\mathbf{r}}\right) 
\mathbf{U}\left( \mathbf{r},\tau \right) +F\left( \mathbf{U}\right) \left( 
\mathbf{\mathbf{r}},\tau \right) ,\ \mathbf{U}\left( \mathbf{\mathbf{r}}%
,0\right) =\mathbf{h}\left( \mathbf{\mathbf{r}}\right) .  \label{dif}
\end{equation}%
Note that since $\mathbf{L}\left( \mathbf{\mathbf{k}}\right) $ satisfies the
polynomial bound \ref{Lpol}) we can define the action of the operator $%
\mathbf{L}\left( -\mathrm{i}\nabla _{\mathbf{r}}\right) $ on any Schwartz
function $\mathbf{Y}\left( \mathbf{r}\right) $ by the formula%
\begin{equation}
\widehat{\mathbf{L}\left( -\mathrm{i}\nabla _{\mathbf{r}}\right) \mathbf{Y}}%
\left( \mathbf{\mathbf{k}}\right) =\mathbf{L}\left( \mathbf{\mathbf{k}}%
\right) \mathbf{\hat{Y}}\left( \mathbf{\mathbf{k}}\right) ,  \label{Ldiff}
\end{equation}%
where, in view of (\ref{Lpol}), the order of $\mathbf{L}$ does not exceed $p$%
. If all the entries of $\mathbf{L}\left( \mathbf{\mathbf{k}}\right) $ are
polynomials, such a definition coincides with the common definition of the
action of a differential operator $\mathbf{L}\left( -\mathrm{i}\nabla _{%
\mathbf{r}}\right) $. In this case $\mathbf{L}\left( -\mathrm{i}\nabla _{%
\mathbf{r}}\right) $ defined by (\ref{Ldiff}) is a differential operator
with constant coefficients of order not greater than $p$.

The properties of the modal evolution equation (\ref{difF}) are completely
similar to its lattice counterpart and are as follows. The existence and
uniqueness theorem is similar to Theorem \ref{Theorem Existence}.

\begin{theorem}[existence and uniqueness]
\label{Theorem Existence1}\ Let equation (\ref{difF}) satisfy conditions (%
\ref{supxi}) and (\ref{chiCR}) and $\mathbf{h}\in L_{1}=L_{1}\left( \mathbb{R%
}^{d}\right) $,$\left\Vert \mathbf{\tilde{h}}\right\Vert _{L_{1}}\leq R$.
Then there exists a unique solution to the modal evolution equation (\ref%
{difF}) in the functional space $C^{1}\left( \left[ 0,\tau _{\ast }\right]
,L_{1}\right) $. The number $\tau _{\ast }$ depends on $R$, $C_{\chi }$ and $%
C_{\Xi }$.
\end{theorem}

Here is the main result for the semilinear hyperbolic systems of PDE which
is completely similar to Theorem \ref{Theorem Superposition}.

\begin{theorem}[principle of superposition for PDE systems]
\label{Theorem Superposition1} Let the initial data of the modal evolution
equation (\ref{difF}) be a multi-wavepacket, i.e. the sum of $N_{h}$
wavepackets $\mathbf{\hat{h}}_{l}$ as in (\ref{hsumR}) satisfying
Definitions \ref{dwavepack}, \ref{dmwavepack}. Suppose that $\varrho ,\beta $
satisfy condition (\ref{scale1}). \ Assume also that $\mathbf{\hat{h}}$ is
generic in the sense of Definition \ref{Condition generic}. Then the
solution $\mathbf{\hat{U}}=\mathcal{G}\left( \mathbf{\hat{h}}\right) $ to
the modal evolution equation (\ref{difF}) satisfies the \emph{approximate
linear superposition principle}, namely%
\begin{equation}
\mathcal{G}\left( \sum_{l=1}^{N_{h}}\mathbf{\hat{h}}_{l}\right)
=\sum_{l=1}^{N_{h}}\mathcal{G}\left( \mathbf{\hat{h}}_{l}\right) +\mathbf{%
\hat{D}},  \label{GsumR}
\end{equation}%
with a small remainder $\mathbf{\hat{D}}\left( \tau \right) $%
\begin{equation}
\sup_{0\leq \tau \leq \tau _{\ast }}\left\Vert \mathbf{\hat{D}}\left( \tau
\right) \right\Vert _{L_{1}}\leq C_{\epsilon }\frac{\varrho }{\beta
^{1+\epsilon }}\left\vert \ln \beta \right\vert ,  \label{Dr}
\end{equation}%
where $\epsilon $ is the same as in Definition \ref{dwavepack}, $\tau _{\ast
}$ does not depend on $\beta ,\varrho $ and $\epsilon $. The solutions $%
\mathbf{U}\left( \mathbf{h}\right) \left( \mathbf{r},\tau \right) $ of the
space evolution equation (\ref{dif}) are obtained as the inverse Fourier
transform of$\ \mathcal{G}\left( \mathbf{\hat{h}}\right) $ and they satisfy
the \emph{approximate linear superposition principle}, namely 
\begin{equation}
\mathbf{U}\left( \mathbf{h}\right) =\mathbf{U}\left( \mathbf{h}_{1}\right)
+\ldots +\mathbf{U}\left( \mathbf{h}_{N_{h}}\right) +\mathbf{D},
\label{UsumR}
\end{equation}%
with a small coupling remainder $\mathbf{D}\left( \tau \right) $ satisfying 
\begin{equation}
\sup_{0\leq \tau \leq \tau _{\ast }}\left\Vert \mathbf{D}\left( \tau \right)
\right\Vert _{L_{\infty }}\leq C_{\epsilon }\frac{\varrho }{\beta
^{1+\epsilon }}\left\vert \ln \beta \right\vert ,  \label{estDPDE}
\end{equation}%
where $\epsilon >0$ is the same as in Definition \ref{dwavepack} and can be
arbitrary small.
\end{theorem}

\paragraph{Example 1: Sine-Gordon and Klein-Gordon equations with small
initial data}

Let us consider the \emph{sine-Gordon} equation (see \cite{InfeldR}) 
\begin{equation}
\partial _{t}^{2}u=\partial _{r}^{2}u-\sin u  \label{sgu1}
\end{equation}%
with small initial data 
\begin{equation}
u\left( r,0\right) =\beta b_{0},\ \partial _{t}u\left( r,0\right) =\beta
b_{1},\;\beta \ll 1.  \label{insineG}
\end{equation}%
First, we recast this the equation into our framework by rescaling the
variables 
\begin{equation}
u=\beta U_{1},\;\beta ^{2}t=\tau .  \label{resc}
\end{equation}%
Since $\sin \beta U_{1}=\beta U_{1}-\frac{1}{6}\beta ^{3}U_{1}^{3}+\beta
^{5}f\left( U_{1}\right) $, where evidently $f\left( U_{1}\right) $ is an
enitire function, we can recast the equation (\ref{sgu1}) into the following
form 
\begin{equation}
\partial _{\tau }^{2}U_{1}=\frac{1}{\beta ^{4}}\left[ \partial
_{x}^{2}U_{1}-U_{1}\right] +\frac{1}{\beta ^{2}}\left[ qU_{1}^{3}+\beta
^{2}f\left( U_{1}\right) \right] .  \label{sineGres}
\end{equation}%
We introduce then a linear pseudodifferential operator $\ A=\left(
I-\partial _{x}^{2}\right) ^{1/2}$ with the symbol $\left( 1+k^{2}\right)
^{1/2}$ and rewrite the equation (\ref{sineGres}) as the following system%
\begin{equation}
\partial _{\tau }U_{1}=\frac{1}{\beta ^{2}}AU_{2},\ \partial _{\tau }U_{2}=-%
\frac{1}{\beta ^{2}}AU_{1}+A^{-1}\left[ qU_{1}^{3}+\beta ^{2}f\left(
U_{1}\right) \right] ,  \label{sineGsys}
\end{equation}%
with the\ initial data 
\begin{equation}
U_{1}\left( 0\right) =h_{0},\ U_{2}\left( 0\right) =h_{1},  \label{inSGsys}
\end{equation}%
where $h_{0}$ and $h_{1}$ are assumed to be of the form 
\begin{equation}
z\left( \mathbf{r},0\right) =h_{0},\ p\left( \mathbf{r},0\right)
=h_{1},h_{j}=\sum_{l=1}^{n_{h}}\Psi _{jl}\left( \beta \mathbf{r}-\mathbf{r}%
_{l}\right) \mathrm{e}^{\mathrm{i}\mathbf{k}_{\ast l}\cdot \mathbf{r}}+%
\limfunc{cc},\ j=0,1,  \label{inhyp0}
\end{equation}%
in one-dimentional case with$\ \mathbf{r}=r,$ $\mathbf{k}=k$. \ Evidently,
the relations with the initial data of (\ref{sgu1}) are 
\begin{equation*}
b_{0}=h_{0},\ b_{1}=Ah_{1}.
\end{equation*}%
Notice that the system (\ref{sineGsys}) is of the form (\ref{dif}) \ with 
\begin{gather}
\varrho =\beta ^{2},\;\mathbf{LU}=\left[ 
\begin{array}{c}
AU_{2} \\ 
-AU_{1}%
\end{array}%
\right] ,\ F\left( \mathbf{U}\right) =F_{0}\left( \mathbf{U}\right) +\beta
^{2}F_{1}\left( \mathbf{U}\right) ,  \label{FsinG} \\
F_{0}\left( \mathbf{U}\right) =A^{-1}\left[ 
\begin{array}{c}
0 \\ 
qU_{1}^{3}%
\end{array}%
\right] ,\ F_{1}\left( \mathbf{U}\right) =A^{-1}\left[ 
\begin{array}{c}
0 \\ 
f\left( U_{1}\right)%
\end{array}%
\right] .  \notag
\end{gather}%
Observe now that $\mathbf{L}$ has only one spectral band with the dispersion
relation and eigenvectors given by 
\begin{equation*}
\omega \left( k\right) =\left( I+k^{2}\right) ^{1/2},\ \mathbf{g}_{\vartheta
}\left( k\right) =\mathbf{g}_{\vartheta }=2^{-1/2}\left[ 
\begin{array}{c}
-\mathrm{i}\vartheta \\ 
1%
\end{array}%
\right] ,\;\ \vartheta =\pm 1,
\end{equation*}%
and there is no band-crossing points. We use expansion in the basis $\mathbf{%
g}_{\pm }$ 
\begin{equation}
\left[ 
\begin{array}{c}
\Psi _{0l} \\ 
\Psi _{1l}%
\end{array}%
\right] \mathrm{e}^{\mathrm{i}\mathbf{k}_{\ast l}\cdot \mathbf{r}}=\left[
\Phi _{+,l}\mathbf{g}_{+}+\Phi _{-,l}\mathbf{g}_{-}\right] \mathrm{e}^{%
\mathrm{i}\mathbf{k}_{\ast l}\cdot \mathbf{r}}.  \label{Psibas0}
\end{equation}%
to represent initial data (\ref{inSGsys}) and (\ref{inhyp0}). The equation (%
\ref{zomomz0}) takes here the form 
\begin{equation*}
\left( 1+k_{\ast l}^{2}\right) ^{1/2}\lambda =\zeta \left( 1+\lambda
^{2}k_{\ast l}^{2}\right) ^{1/2},\;\zeta =\pm 1.
\end{equation*}%
Obviously, this equation has only solutions $\ \lambda =\zeta $ and
Condition \ref{Condition generic1} is fulfilled. Condition \ref{Condition GV}
holds if%
\begin{equation}
\frac{\vartheta k_{\ast l}}{\left( 1+k_{\ast l}^{2}\right) ^{1/2}}\neq \frac{%
\vartheta ^{\prime }k_{\ast l^{\prime }}}{\left( 1+k_{\ast l^{\prime
}}^{2}\right) ^{1/2}}\text{ for \ }l\neq l^{\prime }\text{ or }\vartheta
\neq \vartheta ^{\prime }
\end{equation}%
which is equivalent to 
\begin{equation}
k_{\ast l^{\prime }}\neq k_{\ast l}\text{ for }l^{\prime }\neq l,\text{ and
\ }k_{\ast l}\neq 0\ \text{for all }l.  \label{SGgen}
\end{equation}%
\ Equation (\ref{sineGsys}) can be written in the integral form (\ref{varcu}%
) with $m_{F}=\infty $ and by Theorem \ref{Imfth1} \ it has unique solution $%
\mathbf{U}$ for $\tau \leq \tau _{\ast }$. If we replace $F\left( \mathbf{U}%
\right) $ in (\ref{FsinG}) by $F_{0}\left( \mathbf{U}\right) $ \ we obtain 
\begin{equation}
\partial _{\tau }V_{1}=\frac{1}{\beta ^{2}}AV_{2},\;\partial _{\tau }V_{2}=-%
\frac{1}{\beta ^{2}}AV_{1}+A^{-1}qV_{1}^{3},  \label{sinGcube}
\end{equation}%
where we take the initial data to be\ as in (\ref{inSGsys}), namely \ 
\begin{equation}
V_{1}\left( 0\right) =h_{0},\ V_{2}\left( 0\right) =h_{1}.  \label{inV}
\end{equation}%
Equations (\ref{sinGcube}) can be obtained by replacing $\sin u$ in (\ref%
{sgu1}) by the cubic polynomial $u-u^{3}/6$ \ producing the quasilinear 
\emph{Klein-Gordon} equation (see \cite{Nayfeh}). Observe that the solutions
to the sine-Gordon and the Klein-Gordon equations with small initial data
are very close. To see that, note that the operator $\widehat{f\left(
U\right) }$ $\left( k\right) $ is bounded \ in $L_{1}$ for $\widehat{U}%
\left( k\right) $ which are bounded in $L_{1}$. Therefore the norm of the
neglected term is small, namely $\left\Vert \beta ^{2}\widehat{f\left(
U\right) }\right\Vert _{L_{1}}\leq C\beta ^{2}$. \ Thus, by Remark \ref%
{Remark Lip}, the solutions of (\ref{sineGsys}) and (\ref{sinGcube})\ are
close, namely \ 
\begin{equation}
\ \left\Vert U_{1}-V_{1}\right\Vert _{L_{\infty }}+\left\Vert
U_{2}-V_{2}\right\Vert _{L_{\infty }}\leq C\beta ^{2},\;0\leq \tau \leq \tau
_{\ast .}  \label{UVbet2}
\end{equation}%
\ According to Theorem \ref{Theorem Superposition1} the superposition
principle is applicable to \ the equation (\ref{sinGcube}) \ with initial
data as in (\ref{inV}), and the following statements hold.

\begin{theorem}[Superposition for Klein-Gordon]
\label{TheoremKG}Assume that the initial data $h_{0},h_{1}$ in (\ref{inV})
are as in (\ref{inhyp0}). Then the solution $\left\{ V_{1},V_{2}\right\} $
to the system (\ref{sinGcube}) satisfies the linear superposition principle,
namely 
\begin{equation}
V_{1}\left( \mathbf{r},\tau \right) =\sum_{\vartheta =\pm
}\sum_{l=1}^{n_{h}}V_{1,\vartheta ,l}\left( \mathbf{r},\tau \right) +\mathbf{%
D}_{1}\left( \mathbf{r},\tau \right) ,\ V_{2}\left( \mathbf{r},\tau \right)
=\sum_{\vartheta =\pm }\sum_{l=1}^{n_{h}}V_{2,\vartheta ,l}\left( \mathbf{r}%
,\tau \right) +\mathbf{D}_{2}\left( \mathbf{r},\tau \right) ,  \label{KGsum}
\end{equation}%
where $\left\{ V_{1,\vartheta ,l}\left( \mathbf{r},\tau \right)
,V_{2,\vartheta ,l}\left( \mathbf{r},\tau \right) \right\} $ is a solution
to (\ref{sinGcube}) with the one-wavepacket initial condition 
\begin{equation}
\left[ 
\begin{array}{c}
V_{1,\vartheta ,l}\left( \mathbf{r},0\right) \\ 
V_{2,\vartheta ,l}\left( \mathbf{r},0\right)%
\end{array}%
\right] =\Phi _{\vartheta ,l}\left( \beta \mathbf{r}-\mathbf{r}_{l}\right) 
\mathbf{g}_{\vartheta }\mathrm{e}^{\mathrm{i}\mathbf{k}_{\ast l}\cdot 
\mathbf{r}}+\limfunc{cc},  \label{inKG1}
\end{equation}%
where $\Phi _{\vartheta ,l}\left( \mathbf{r}\right) $ are arbitrary Schwartz
functions. If\textbf{\ }(\ref{SGgen}) holds, the coupling terms $\mathbf{D}%
_{1},\mathbf{D}_{2}$ satisfy the bound 
\begin{equation}
\sup_{0\leq \tau \leq \tau _{\ast }}\left[ \left\Vert \mathbf{D}_{1}\left(
\tau \right) \right\Vert _{L_{\infty }}+\left\Vert \mathbf{D}_{2}\left( \tau
\right) \right\Vert _{L_{\infty }}\right] \leq C_{\delta }^{\prime }\frac{%
\varrho }{\beta ^{1+\delta }}=C_{\delta }^{\prime }\beta ^{1-\delta },
\label{D1D2}
\end{equation}%
where $\tau _{\ast }$ and $C_{\delta }^{\prime }$ do not depend on $\beta $,
and $\delta $ can be taken arbitrary small.
\end{theorem}

Using (\ref{UVbet2}) we obtain a similar superposition theorem for the
sine-Gordon equation.

\begin{theorem}[Superposition for sine-Gordon]
\label{TheoremSG}Assume that the initial data $h_{0},h_{1}$ in (\ref{inSGsys}%
) are as in (\ref{inhyp0}). Then the solution $\left\{ U_{1},U_{2}\right\} $
to (\ref{sineGsys}), (\ref{inSGsys}) satisfies the linear superposition
principle, namely 
\begin{equation*}
U_{1}\left( \mathbf{r},\tau \right) =\sum_{\vartheta =\pm
}\sum_{l=1}^{n_{h}}U_{1,\vartheta ,l}\left( \mathbf{r},\tau \right) +\mathbf{%
D}_{1}\left( \mathbf{r},\tau \right) ,\ U_{2}\left( \mathbf{r},\tau \right)
=\sum_{\vartheta =\pm }\sum_{l=1}^{n_{h}}U_{2,\vartheta ,l}\left( \mathbf{r}%
,\tau \right) +\mathbf{D}_{2}\left( \mathbf{r},\tau \right) ,
\end{equation*}%
where $U_{1,\vartheta ,l}\left( \mathbf{r},\tau \right) ,U_{2,\vartheta
,l}\left( \mathbf{r},\tau \right) $ is a solution of (\ref{sineGsys}) with
the one-wavepacket initial condition 
\begin{equation*}
\left[ 
\begin{array}{c}
U_{1,\vartheta ,l}\left( \mathbf{r},0\right) \\ 
U_{2,\vartheta ,l}\left( \mathbf{r},0\right)%
\end{array}%
\right] =\Phi _{\vartheta ,l}\left( \beta \mathbf{r}-\mathbf{r}_{l}\right) 
\mathbf{g}_{\vartheta }\mathrm{e}^{\mathrm{i}\mathbf{k}_{\ast l}\cdot 
\mathbf{r}}+\limfunc{cc},\vartheta =\pm ,
\end{equation*}%
where $\Phi _{\vartheta ,,l}\left( \mathbf{r}\right) $ are arbitrary
Schwartz functions. If\textbf{\ }(\ref{SGgen}) holds, \ the coupling terms $%
\mathbf{D}_{1},\mathbf{D}_{2}$ satisfy the bound (\ref{D1D2}).
\end{theorem}

Note that a theorem completely similar to Theorem \ref{TheoremKG} holds also
for a generalized Klein-Gordon equation where $qV_{1}^{3}$ is replaced by an
arbitrary polynomial $P\left( V_{1}\right) $. Hence, the superposition
principle holds for the sine-Gordon equation (\ref{sgu1}) with a small
initial data and a strongly perturbed nonlinearity as, for example, when $%
\sin u$ is replaced by $\sin u+\beta ^{-1}u^{4}+\beta ^{-2}u^{5}$.

We would like to compare now our results and methods with that of \cite{PW}
where the interaction of counterpropagating waves is studied by the ansatz
method. Pierce and Wayne considered in \cite{PW} the sine-Gordon equation in
the case of small initial data which have the form of a bimodal wavepacket.
In our notation it corresponds to the case when $\varrho =\beta ^{2},$ $%
n_{h}=1$ in (\ref{inhyp0}), when two wavepackets, corresponding to $%
\vartheta =+$ \ and $\vartheta =-$, have exactly opposite group velocities.
They proved that the bimodal wavepacket data generate two waves which are
described by two uncoupled nonlinear Schrodinger equations \ with a small
error. The magnitude of the error given in \cite{PW} \ (which we formulate
here for the solution $U_{1}$ of the rescaled equation \ (\ref{sineGres}) )
is estimated by $C\beta ^{1/2}$ on the time interval $0\leq \tau \leq \tau
_{0}$ ( or $0\leq t\leq \tau _{0}\beta ^{-2}$). Note that our general
Theorem \ref{Theorem Superposition1} when applied to the special case of the
sine-Gordon equation (\ref{sineGres}) provides a better estimate of the
coupling error, namely $C\varrho /\beta ^{1+\delta }=C\beta ^{1-\delta }$ in
(\ref{D1D2}) with arbitrary small $\delta $, for the same time interval.
Notice that the estimate (\ref{estDPDE}) given in Theorem \ref{Theorem
Superposition1} \ is almost optimal, since it is possible to construct
examples when the coupling error is greater than $c\beta ^{1+\delta }$ with
arbitrary small $\delta $.

We would like to point out\ that the general mechanism responsible for the
wavepacket decoupling is the destructive wave interference, this mechanism
is subtle though general. We treat the destructive wave interference by \
taking into account explicitly all nonlinear interactions of high-frequency
waves. In our approach we use the exact representation of a general solution
in the form of a functional-analytic operator monomial series, every term of
the series is explicitly given as a multilinear oscillatory integral
operator applied to the initial data. A key advantage of such an approach is
that it allows to estimate wavepacket coupling as a sum of contributions of
highly oscillatory terms and to get a precise estimate of magnitude of every
term. In contrast, the well known "ansatz" approach as, for instance, in 
\cite{PW} and \cite{KSM}, \ requires to find a clever ansatz with consequent
estimations of the \ "residuum" in an appropriate norm. Our approach can
naturally treat general tensorial polynomial nonlinearities $F$ of arbitrary
large degree $N_{F}$ and any number\ of wavepackets, whereas finding a good
ansatz which allows to estimate the residuum in such a general situation
would be difficult. For readers interested in detailed features of
one-wavepacket solutions to the sine-Gordon equations, we refer to \cite{KSM}%
, \cite{PW} \ and \cite{Schneider98a}. \ \ 

\paragraph{Example 2: Nonlinear Schrodinger equation.}

The Nonlinear Schrodinger equation (NLS) with $d$ spatial variables (\cite%
{Sulem}, \cite{Caz}, \cite{Bourgain}) has the form 
\begin{equation}
\partial _{\tau }z\left( \mathbf{r},\tau \right) =\mathrm{i}\frac{1}{\varrho 
}\gamma \left( -\mathrm{i}\nabla \right) z\left( \mathbf{r},\tau \right) \
+\alpha \left\vert z\right\vert ^{2}z\left( \mathbf{r},\tau \right) ,\
z\left( \mathbf{r},0\right) =h\left( \mathbf{r}\right) ,\;\mathbf{r}\in 
\mathbb{R}^{d},  \label{NLS}
\end{equation}%
where $\alpha $ is a complex constant, $\gamma \left( -\mathrm{i}\nabla
\right) \ $is a second-order differential operator, its symbol $\gamma
\left( \mathbf{k}\right) $ is a real, symmetric quadratic form%
\begin{equation*}
\gamma \left( \mathbf{k}\right) =\gamma \left( \mathbf{k},\mathbf{k}\right)
=\sum \gamma _{ij}k_{i}k_{j},\ \gamma \left( -\mathrm{i}\nabla \right)
z=-\sum \gamma _{ij}\partial _{r_{i}}\partial _{r_{j}}z.
\end{equation*}%
To put the NLS into the framework of this paper we introduce the following
two-component system%
\begin{gather}
\partial _{\tau }z_{+}\left( \mathbf{r},\tau \right) =\mathrm{i}\frac{1}{%
\varrho }\gamma \left( -\mathrm{i}\nabla \right) z_{+}\left( \mathbf{r},\tau
\right) +\alpha z_{-}z_{+}^{2}\left( \mathbf{r},\tau \right) ,  \label{NLS2}
\\
\partial _{\tau }z_{-}\left( \mathbf{r},\tau \right) =-\mathrm{i}\frac{1}{%
\varrho }\gamma \left( \mathrm{i}\nabla \right) z_{-}\left( \mathbf{r},\tau
\right) +\alpha ^{\ast }z_{+}z_{-}^{2}\left( \mathbf{r},\tau \right) , 
\notag \\
z_{+}\left( \mathbf{r},0\right) =h\left( \mathbf{r}\right) ,\ z_{-}\left( 
\mathbf{r},0\right) =h^{\ast }\left( \mathbf{r}\right) ,\ \mathbf{r}\in 
\mathbb{R}^{d},  \notag
\end{gather}%
where $\alpha ^{\ast }$ denotes complex conjugate to $\alpha $. Obviously if 
$z\left( \mathbf{r},\tau \right) $ is a solution of (\ref{NLS}) then $%
z_{+}\left( \mathbf{r},\tau \right) =z\left( \mathbf{r},\tau \right) $, $%
z_{-}\left( \mathbf{r},\tau \right) =z^{\ast }\left( \mathbf{r},\tau \right) 
$ gives a solution of (\ref{NLS2}). Using the Fourier transform we get from (%
\ref{NLS}) 
\begin{equation}
\partial _{\tau }\hat{z}\left( \mathbf{k},\tau \right) =\mathrm{i}\frac{1}{%
\varrho }\gamma \left( \mathbf{k}\right) \hat{z}\left( \mathbf{k},\tau
\right) +\alpha \widehat{\left( z^{\ast }z^{2}\right) }\left( \mathbf{k}%
,\tau \right) \ ,\mathbf{k}\in \mathbb{R}^{d}.  \label{NLSk}
\end{equation}%
Now the band-crossing set $\sigma =\left\{ \mathbf{k}\in \mathbb{R}%
^{d}:\gamma \left( \mathbf{k}\right) =0\right\} $. \ We assume that the
quadratic form $\gamma \ $is not identically zero. $\ $ The Fourier
transform of (\ref{NLS2}) takes the form of (\ref{dif}) \ with 
\begin{gather*}
\mathbf{\hat{U}=}\left[ 
\begin{array}{c}
\hat{U}_{+} \\ 
\hat{U}_{-}%
\end{array}%
\right] ,\ \mathbf{L}\left( \mathbf{k}\right) \mathbf{\hat{U}}=\left[ 
\begin{array}{cc}
\gamma \left( \mathbf{k}\right) & 0 \\ 
0 & -\gamma \left( -\mathbf{k}\right)%
\end{array}%
\right] \left[ 
\begin{array}{c}
\hat{U}_{+} \\ 
\hat{U}_{-}%
\end{array}%
\right] , \\
\omega \left( \mathbf{k}\right) =\left\vert \gamma \left( \mathbf{k}\right)
\right\vert ,\hat{F}^{\left( 3\right) }\left( \mathbf{\hat{U}}^{3}\right) =%
\left[ 
\begin{array}{c}
\alpha \widehat{\left( \hat{z}_{+}\left( \mathbf{\hat{U}}\right) \hat{z}%
_{+}\left( \mathbf{\hat{U}}\right) \hat{z}_{-}\left( \mathbf{\hat{U}}\right)
\right) } \\ 
\alpha ^{\ast }\widehat{\left( \hat{z}_{-}\left( \mathbf{\hat{U}}\right) 
\hat{z}_{-}\left( \mathbf{\hat{U}}\right) \hat{z}_{+}\left( \mathbf{\hat{U}}%
\right) \right) }%
\end{array}%
\right] ,
\end{gather*}%
To satisfy the requirements of Condition \ref{Condition generic} we have to
take the wave vectors $\mathbf{k}_{\ast l}\notin \sigma $ \ so that%
\begin{equation}
\nabla \left\vert \gamma \left( \mathbf{k}_{\ast l}\right) \right\vert =%
\frac{2\gamma \left( \mathbf{k}_{\ast l}\right) }{\left\vert \gamma \left( 
\mathbf{k}_{\ast l}\right) \right\vert }\gamma \left( \mathbf{k}_{\ast
l},\cdot \right) \neq \frac{2\gamma \left( \mathbf{k}_{\ast l^{\prime
}}\right) }{\left\vert \gamma \left( \mathbf{k}_{\ast l^{\prime }}\right)
\right\vert }\gamma \left( \mathbf{k}_{\ast l^{\prime }},\cdot \right) \text{
if \ }l\neq l^{\prime },  \label{NLSNGV}
\end{equation}%
which provides (\ref{NGVM}). Since 
\begin{equation*}
\left\vert \gamma \left( \mathbf{k}_{\ast l}\right) \right\vert \lambda
-\zeta \left\vert \gamma \left( \lambda \mathbf{k}_{\ast l}\right)
\right\vert =\left\vert \gamma \left( \mathbf{k}_{\ast l}\right) \right\vert %
\left[ \lambda -\zeta \left\vert \lambda \right\vert ^{2}\right] ,
\end{equation*}%
and $\lambda $ is odd, every point $\mathbf{k}_{\ast l}\notin \sigma $
satisfies Condition \ref{Condition generic1}. If the quadratic form $\gamma $
is not singular, that is $\det \gamma \neq 0$, then condition (\ref{NLSNGV}%
), which ensures that group velocities of wavepackets are different, holds
when%
\begin{equation*}
\frac{\gamma \left( \mathbf{k}_{\ast l}\right) }{\left\vert \gamma \left( 
\mathbf{k}_{\ast l}\right) \right\vert }\mathbf{k}_{\ast l}\neq \frac{\gamma
\left( \mathbf{k}_{\ast l^{\prime }}\right) }{\left\vert \gamma \left( 
\mathbf{k}_{\ast l^{\prime }}\right) \right\vert }\mathbf{k}_{\ast l^{\prime
}}\text{ if \ }l\neq l^{\prime }.
\end{equation*}%
In this case Theorem \ref{Theorem Superposition1} is applicable, and generic
wavepacket solutions of the NLS are linearly superposed and propagate almost
independently with coupling $O\left( \beta \right) $. More precisely, as a
corollary of Theorem \ref{Theorem Superposition1} we obtain the following
statement.

\begin{theorem}[Superposition for NLS]
\label{Theorem Superposition NLS} Assume that initial data of the NLS (\ref%
{NLS}) have the form $h=\ h_{1}+\ldots +h_{N_{h}}$ 
\begin{equation*}
h_{l}\left( \mathbf{r}\right) =\mathrm{e}^{\mathrm{i}\mathbf{k}_{\ast
l}\cdot \mathbf{m}}\Phi _{l,+}\left( \beta \mathbf{r}-\mathbf{r}_{0}\right) +%
\mathrm{e}^{-\mathrm{i}\mathbf{k}_{\ast l}\cdot \mathbf{m}}\Phi _{l,-}\left(
\beta \mathbf{r}-\mathbf{r}_{0}\right) ,l=1,\ldots ,N_{h}
\end{equation*}%
where $\Phi _{l,\zeta }\left( \mathbf{r}\right) $ are arbitrary Schwartz
functions. Assume also that $\det \gamma \neq 0$ and the vectors $\mathbf{k}%
_{\ast l}$ satisfy conditions 
\begin{equation*}
\gamma \left( \mathbf{k}_{\ast l}\right) \neq 0,\;l=1,\ldots ,N_{h};\ 
\mathbf{k}_{\ast l}\neq \mathbf{k}_{\ast l^{\prime }}\ \text{ if \ }l\neq
l^{\prime }.
\end{equation*}%
Then \ solution $z=z\left( h\right) $ is a linear superposition 
\begin{equation*}
z\left( h\right) =z\left( h_{1}\right) +\ldots +z\left( h_{N_{h}}\right) +D
\end{equation*}%
with a small coupling term $D$ 
\begin{equation*}
\sup_{0\leq \tau \leq \tau _{\ast }}\left\Vert D\left( \tau \right)
\right\Vert _{L_{\infty }\left( \mathbb{R}^{d}\right) }\leq C_{\delta }\frac{%
\varrho }{\beta ^{1+\delta }},
\end{equation*}%
where $\delta >0$ can be taken arbitrary small.
\end{theorem}

We note in conclusion, that the superposition principle reduces dynamics of
multi-wavepacket solutions to dynamics of single-wavepacket solutions; we do
not study dynamics of single-wavepacket solutions in this paper. Note that
the theory of NLS-type approximations of one-wavepacket solutions of
hyperbolic PDE is well-developed, see \cite{KalyakinUMN}, \cite{Kalyakin2}, 
\cite{ColinLannes}, \cite{Schneider05}, \cite{SU}, \cite{BF5} and references
therein. Relevance of different group velocities of wavepackets for
smallness of their interaction was noted in \cite{KalyakinUMN}.

\subsection{Generalizations}

Note that in a degenerate case when the function $\omega _{n_{l}}\left( 
\mathbf{k}\right) $ is linear in the direction of $\mathbf{\mathbf{k}}_{\ast
}$ the equation (\ref{zomomz0}) for $\zeta =1$ has many solutions for which $%
\theta \neq \pm 1$\ and Condition \ref{Condition generic1} does not hold. It
turns out, that if \ Condition \ref{Condition generic1} for dispersion
relations $\omega _{n}\left( \mathbf{k}\right) $ at $\ \mathbf{\mathbf{k}}%
_{\ast }$ is not satisfied, still we can prove our results under the
following alternative condition. We consider here the case of PDE\ in the
entire space $\mathbb{R}^{d}$ \ and $\mathbf{k}\in \mathbb{R}^{d}$.

\begin{condition}[complete degeneracy]
\label{Condition generic2} The series (\ref{Fseries}) has only $\tilde{F}%
^{\left( m\right) }$ with odd $m$. \ The wavevectors $\mathbf{k}_{\ast l}$
and functions $\omega _{n_{l}}\left( \mathbf{k}\right) $, $l=1,\ldots ,N_{h}$%
, \ have the following three properties:

(i) There exists $\delta >0$ such that for every $l_{1}\neq l_{2}$, the
following inequality holds: 
\begin{equation}
\left\vert \nabla _{\mathbf{k}}\omega _{n_{l_{1}}}\left( \nu _{1}\mathbf{k}%
_{\ast l_{1}}\right) -\nabla _{\mathbf{k}}\omega _{n_{l_{2}}}\left( \nu _{2}%
\mathbf{k}_{\ast l_{2}}\right) \right\vert \geq \delta \text{,}  \label{csg1}
\end{equation}%
for any odd integers $\nu _{1},\nu _{2}=1,3,\ldots $.

(ii) There exists $\delta >0$ such that $\nu \mathbf{k}_{\ast l}$ does not
get in a $\delta $-neighborhood of $\sigma $ for any \ odd integer $\nu $ \
and any $l=1,\ldots ,N_{h}$.

(iii) For any positive integer odd number $\theta $ and any $\mathbf{k}%
_{\ast l}$, for any $n$\ the following identities hold: 
\begin{equation}
\nabla _{\mathbf{k}}\omega _{n}\left( \theta \mathbf{k}_{\ast l}\right)
=\nabla _{\mathbf{k}}\omega _{n}\left( \mathbf{k}_{\ast l}\right) ,
\label{nd22}
\end{equation}%
\begin{equation}
\;\omega _{n}\left( \theta \mathbf{k}_{\ast l}\right) =\theta \omega
_{n}\left( \mathbf{k}_{\ast l}\right) .  \label{nd23}
\end{equation}
\end{condition}

A nontrivial examples, where the above Condition \ref{Condition generic2} is
satisfied, is given below.

We give here a generalization of Definition \ref{Condition generic}.

\begin{definition}[generic multi-wavepackets]
\label{Condition generic22} A multi-wavepacket $\mathbf{\hat{h}}$ as defined
in Definition \ref{dmwavepack} is called generic if (i) the centers $\ 
\mathbf{\mathbf{k}}_{\ast l}$, $l=1,\ldots ,N_{h}$, of all wavepackets
satisfy Conditions \ref{Condition non-zero} and \ref{Condition GV}; (ii)
either the dispersion relations $\omega _{n}\left( \mathbf{k}\right) $ at
every $\mathbf{\mathbf{k}}_{\ast l}$ and band $n_{l}$ satisfy Condition \ref%
{Condition generic1} or \ \ they satisfy Condition \ref{Condition generic2}.
\end{definition}

The statement of Theorem \ref{Theorem Superposition1} remains true if
Condition \ref{Condition generic} is replaced by less restrictive Condition %
\ref{Condition generic22}, namely the following theorem holds.

\begin{theorem}
\label{Theorem Superposition11} Let the initial data of the modal evolution
equation (\ref{difF}) be a multi-wavepacket, i.e. the sum of $N_{h}$
wavepackets $\mathbf{\hat{h}}_{l}$ as in (\ref{hsumR}) satisfying
Definitions \ref{dwavepack}, \ref{dmwavepack}. Suppose that (\ref{scale1})
holds. \ Assume also that $\mathbf{\hat{h}}$ is generic in the sense of
Definition \ref{Condition generic22}. Then the solution $\mathbf{\hat{U}}=%
\mathcal{G}\left( \mathbf{\hat{h}}\right) $ to the modal evolution equation (%
\ref{difF}) satisfies the approximate linear superposition principle, namely
(\ref{GsumR}), (\ref{Dr}), (\ref{UsumR}) and (\ref{estDPDE}) hold.
\end{theorem}

The proofs we give in this paper directly apply to more general Theorem \ref%
{Theorem Superposition11}.

Another generalization concerns the possibility to shift independently
initial wavepackets. If initial data involve parameters $\mathbf{r}_{l}$ as
in (\ref{inhyp0}) it is possible to prove that $C_{\epsilon }$ in (\ref{rem}%
), (\ref{Dr}) and (\ref{estDPDE}) does not depend on $\mathbf{r}_{l}\in 
\mathbb{R}^{d}$ if the functions $\Psi _{jl}$ are Schwartz functions. Most
of the proofs remain the same, but several statements have to be modified,
and we present proofs in a subsequent paper.

One more generalization concerns the smoothness of initial data. It is
possible to take initial data $\mathbf{h}_{l}\left( \mathbf{r}\right) $ with
a finite smoothness \ rather than from Schwartz class. Namely, consider
weighted spaces $L_{1,a}$ with the norm 
\begin{equation}
\left\Vert \mathbf{\hat{v}}\right\Vert _{L_{1,a}}=\int_{\mathbb{R}%
^{d}}\left( 1+\left\vert \mathbf{k}\right\vert \right) ^{a}\left\vert 
\mathbf{\hat{v}}\left( \mathbf{k}\right) \right\vert \,\mathrm{d}\mathbf{k}%
,\ a\geq 0.  \label{L1a}
\end{equation}%
Obviously, large $a$ corresponds to high smoothness of the inverse Fourier
transform $\mathbf{v}\left( \mathbf{r}\right) $. Then if functions $\mathbf{%
\hat{h}}_{l,\zeta }\left( \mathbf{\mathbf{k}}\right) $ have the form (\ref%
{h0}) with $\ \hat{h}_{\zeta }\left( \mathbf{k}\right) =\hat{h}_{l,\zeta
}\left( \mathbf{k}\right) $ \ from the class $L_{1,a}$ \ the inequality (\ref%
{Dr}) can be replaced by 
\begin{equation}
\sup_{0\leq \tau \leq \tau _{\ast }}\left\Vert \mathbf{\hat{D}}\left( \tau
\right) \right\Vert _{L_{1}}\leq C_{\epsilon }\frac{\varrho }{\beta
^{1+\epsilon }}\left\vert \ln \beta \right\vert +C_{\epsilon }\beta ^{s},
\label{ests}
\end{equation}%
where $s>0$ and $\epsilon >0$ have to satisfy restriction $\frac{s}{\epsilon 
}<a$. This generalization requires minor modifications in the proofs and in
conditions (\ref{hbold}) and (\ref{sourloc}) $C_{\epsilon }\beta $ has to be
replaced by $C_{\epsilon }\beta ^{s}$. In particular, if $a=1$, $\varrho
=\beta ^{2}$ and $s=1/2$\ the right-hand side of (\ref{ests}) can be
estimated by $C_{\epsilon _{1}}\beta ^{1/2-\epsilon _{1}}$ with arbitrary
small $\epsilon _{1}$.

More generalizations which involve the structure of equations are discussed
in Sections 7.3 and 7.4. Now we give an example where Condition \ref%
{Condition generic2} is applicable.

\paragraph{Example 3: Semilinear wave equation.}

Let us consider a semilinear wave equation with $d$ spatial variables 
\begin{equation}
\partial _{\tau }^{2}z\left( \mathbf{r},\tau \right) =\frac{1}{\varrho ^{2}}%
\Delta z\left( \mathbf{r},\tau \right) +\frac{\alpha }{\varrho }\partial
_{x_{1}}z^{3}\left( \mathbf{r},\tau \right) ,\ \mathbf{r}\in \mathbb{R}^{d},
\label{hyp2}
\end{equation}%
where $\Delta $ is the Laplace operator, $\alpha \ $is an arbitrary complex
constant, $\varrho =\beta ^{2}$. We introduce the operator $A=\sqrt{-\Delta }
$ which is defined in terms of the Fourier transform, it has symbol $%
\left\vert \mathbf{k}\right\vert $. We rewrite (\ref{hyp2}) in the form of a
first-order system 
\begin{gather}
\partial _{\tau }z\left( \mathbf{r},\tau \right) =\frac{1}{\varrho }Ap\left( 
\mathbf{r},\tau \right) ,\ \mathbf{r}\in \mathbb{R}^{d};  \label{hyp21} \\
\partial _{\tau }p\left( \mathbf{r},\tau \right) =-\frac{1}{\varrho }%
Az\left( \mathbf{r},\tau \right) +\alpha A^{-1}\partial _{x_{1}}z^{3}\left( 
\mathbf{r},\tau \right) .  \notag
\end{gather}%
The linear operator $A^{-1}\partial _{x_{1}}$ has the symbol $\frac{-\mathrm{%
i}k_{1}}{\left\vert \mathbf{k}\right\vert }$, it is a zero order operator.
We rewrite (\ref{hyp21}) in the form of (\ref{dif}) where \ 
\begin{equation*}
\mathbf{U=}\left[ 
\begin{array}{c}
z \\ 
p%
\end{array}%
\right] ,\ -\mathrm{i}\mathbf{L}\left( -\mathrm{i}\nabla _{\mathbf{r}%
}\right) \mathbf{U=}\left[ 
\begin{array}{cc}
0 & A \\ 
-A & 0%
\end{array}%
\right] \left[ 
\begin{array}{c}
z \\ 
p%
\end{array}%
\right] ,\ F\left( \left[ 
\begin{array}{c}
z \\ 
p%
\end{array}%
\right] \right) =\alpha \left[ 
\begin{array}{c}
0 \\ 
-A^{-1}\partial _{x_{1}}z^{3}%
\end{array}%
\right] .
\end{equation*}%
Using the Fourier transform we get (\ref{difF}) with 
\begin{equation}
\mathbf{\hat{U}=}\left[ 
\begin{array}{c}
\hat{z} \\ 
\hat{p}%
\end{array}%
\right] ,\ -\mathrm{i}\mathbf{L}\left( \mathbf{k}\right) \mathbf{\hat{U}}=%
\left[ 
\begin{array}{cc}
0 & \left\vert \mathbf{k}\right\vert  \\ 
-\left\vert \mathbf{k}\right\vert  & 0%
\end{array}%
\right] \left[ 
\begin{array}{c}
\hat{z} \\ 
\hat{p}%
\end{array}%
\right] ,\ \hat{F}^{\left( 3\right) }\left( \mathbf{\hat{U}}^{3}\right) =%
\frac{-\mathrm{i}\alpha k_{1}}{\left\vert \mathbf{k}\right\vert }\widehat{%
\left( z^{3}\right) }\left[ 
\begin{array}{c}
0 \\ 
1%
\end{array}%
\right] ,  \notag
\end{equation}%
\begin{equation}
\widehat{\left( z^{3}\right) }\left( \mathbf{k}\right) =\frac{1}{\left( 2\pi
\right) ^{2d}}\int\limits_{\mathbf{k}^{\prime },\mathbf{k}^{\prime \prime
}\in \mathbb{R}^{2d};\mathbf{k}^{\prime }+\mathbf{k}^{\prime \prime }+%
\mathbf{k}^{\prime \prime \prime }=\mathbf{k}}\hat{z}\left( \mathbf{k}%
^{\prime }\right) \hat{z}\left( \mathbf{k}^{\prime \prime }\right) \hat{z}%
\left( \mathbf{k}^{\prime \prime \prime }\right) \,\mathrm{d}\mathbf{k}%
^{\prime }\,\mathrm{d}\mathbf{k}^{\prime \prime }.  \notag
\end{equation}%
Since the factor $\frac{k_{1}}{\left\vert \mathbf{k}\right\vert }$ is
uniformly bounded \ and smooth for $\left\vert \mathbf{k}\right\vert \neq 0$
conditions (\ref{chiCR}) and (\ref{gradchi}) are satisfied. The eigenvalues
\ and corresponding eigenvectors of $\ \mathbf{L}$ are given explicitly: 
\begin{equation}
\omega _{+}\left( \mathbf{k}\right) =\left\vert \mathbf{k}\right\vert
,\omega _{-}\left( \mathbf{k}\right) =-\left\vert \mathbf{k}\right\vert ,\ 
\mathbf{g}_{+}\left( \mathbf{k}\right) =2^{-1/2}\left[ 
\begin{array}{c}
-\mathrm{i} \\ 
1%
\end{array}%
\right] ,\mathbf{g}_{-}\left( \mathbf{k}\right) =2^{-1/2}\left[ 
\begin{array}{c}
\mathrm{i} \\ 
1%
\end{array}%
\right] .  \label{hypg}
\end{equation}%
Since the matrix $\mathbf{L}\left( \mathbf{k}\right) $ is Hermitian,
Condition \ref{Diagonalization} is satisfied. The singular set $\sigma $
consists of the single point $\mathbf{k}=\mathbf{0}$. Note that conclusions
of Theorem \ref{Theorem Superposition1} are applicable to equation (\ref%
{hyp21}) and consequently to (\ref{hyp2}). For instance, we take the initial
data for (\ref{hyp21}) \ in the form (\ref{inhyp0}) 
\begin{equation}
z\left( \mathbf{r},0\right) =h_{0},\ p\left( \mathbf{r},0\right)
=h_{1},h_{j}=\sum_{l=1}^{n_{h}}\Psi _{jl}\left( \beta \mathbf{r}-\mathbf{r}%
_{l}\right) \mathrm{e}^{\mathrm{i}\mathbf{k}_{\ast l}\cdot \mathbf{r}}+%
\limfunc{cc},\ j=0,1,  \label{inhyp}
\end{equation}%
where $\Psi _{0l}\left( \mathbf{r}\right) $, $\Psi _{1l}\left( \mathbf{r}%
\right) $ are arbitrary Schwartz functions, $\limfunc{cc}$ means complex
conjugate to the preceding terms. The points $\mathbf{r}_{l}$ are arbitrary.
Note that terms corresponding to $\mathbf{k}_{\ast l}$ can be written using
the basis (\ref{hypg}) as 
\begin{equation}
\left[ 
\begin{array}{c}
\Psi _{0l} \\ 
\Psi _{1l}%
\end{array}%
\right] \mathrm{e}^{\mathrm{i}\mathbf{k}_{\ast l}\cdot \mathbf{r}}=\left[
\Phi _{+,l}\mathbf{g}_{+}+\Phi _{-,l}\mathbf{g}_{-}\right] \mathrm{e}^{%
\mathrm{i}\mathbf{k}_{\ast l}\cdot \mathbf{r}}.  \label{Psibas}
\end{equation}%
In this case all requirements of Definition \ref{dwavepack} are fulfilled.
The number of initial wavepackets for the first-order system (\ref{hyp21})
corresponding to initial data (\ref{inhyp}) ) equals $N_{h}=2n_{h}$ and
there are $2N_{h}$ wavepacket centers $\vartheta \mathbf{k}_{\ast
l},\vartheta =\pm $. \ To satisfy the requirements of Condition \ref%
{Condition generic} we have to take the wave vectors $\mathbf{k}_{\ast
l}\neq 0$ \ so that%
\begin{equation*}
\frac{\vartheta \mathbf{k}_{\ast l}}{\left\vert \mathbf{k}_{\ast
l}\right\vert }\neq \frac{\vartheta ^{\prime }\mathbf{k}_{\ast l^{\prime }}}{%
\left\vert \mathbf{k}_{\ast l^{\prime }}\right\vert }\text{ if \ }l\neq
l^{\prime }\text{ \ or }\vartheta \neq \vartheta ^{\prime },
\end{equation*}%
which provides (\ref{NGVM}). Since 
\begin{equation*}
\left\vert \mathbf{k}_{\ast l}\right\vert \lambda -\zeta \left\vert \lambda 
\mathbf{k}_{\ast l}\right\vert =\left\vert \mathbf{k}_{\ast l}\right\vert
\left( \lambda -\zeta \left\vert \lambda \right\vert \right) ,
\end{equation*}%
equation (\ref{zomomz0}) has solutions $\lambda \neq \zeta $ and every point 
$\mathbf{k}_{\ast l}$ does not satisfy Condition \ref{Condition generic1}.
This is the property of the very special, purely homogeneous $\omega \left( 
\mathbf{k}\right) =\left\vert \mathbf{k}\right\vert $. Checking the second
alternative, namely Condition \ref{Condition generic2} we observe that 
\begin{equation*}
\nabla _{\mathbf{k}}\left\vert \nu \mathbf{k}_{\ast l}\right\vert =\frac{\nu 
\mathbf{k}_{\ast l}}{\left\vert \nu \mathbf{k}_{\ast l}\right\vert }=\frac{%
\nu }{\left\vert \nu \right\vert }\frac{\mathbf{k}_{\ast l}}{\left\vert 
\mathbf{k}_{\ast l}\right\vert }.
\end{equation*}%
Hence, if 
\begin{equation}
\frac{\vartheta \mathbf{k}_{\ast l}}{\left\vert \mathbf{k}_{\ast
l}\right\vert }\neq \frac{\vartheta ^{\prime }\mathbf{k}_{\ast l^{\prime }}}{%
\left\vert \mathbf{k}_{\ast l^{\prime }}\right\vert }\text{ for \ }l\neq
l^{\prime }\text{ or }\vartheta \neq \vartheta ^{\prime }\;\text{ and if }\;%
\mathbf{k}_{\ast l}\neq 0  \label{kjl}
\end{equation}%
then Condition \ref{Condition generic2} is satisfied and Superposition
Theorem \ref{Theorem Superposition1} is applicable. As a corollary of
Theorem \ \ref{Theorem Superposition1} \ applied to (\ref{hyp2}) we obtain
that if the initial data for (\ref{hyp2}) \ equal the sum of wavepackets,
then the solution equals the sum of separate solutions plus a small
remainder, more precisely we have the following theorem.

\begin{theorem}[superposition principle for wave equation]
\label{Superposition wave} Assume that the initial data for (\ref{hyp21}) to
be a multi-wavepacket of the form (\ref{inhyp}) and (\ref{scale1}) holds.
Then the solution $z\left( \mathbf{r},\tau \right) $ to (\ref{hyp21}), (\ref%
{inhyp}) satisfy the superposition principle, namely 
\begin{equation*}
z\left( \mathbf{r},\tau \right) =\sum_{\vartheta =\pm
}\sum_{l=1}^{n_{h}}z_{\vartheta ,l}\left( \mathbf{r},\tau \right) +\mathbf{D}%
_{1}\left( \mathbf{r},\tau \right) ,p\left( \mathbf{r},\tau \right)
=\sum_{\vartheta =\pm }\sum_{l=1}^{n_{h}}p_{\vartheta ,l}\left( \mathbf{r}%
,\tau \right) +\mathbf{D}_{2}\left( \mathbf{r},\tau \right)
\end{equation*}%
where $z_{\vartheta ,l}\left( \mathbf{r},\tau \right) ,p_{\vartheta
,l}\left( \mathbf{r},\tau \right) $ is a solution of (\ref{hyp21}) with the
initial condition 
\begin{equation}
\left[ 
\begin{array}{c}
z_{\vartheta ,l}\left( \mathbf{r},0\right) \\ 
p_{\vartheta ,l}\left( \mathbf{r},0\right)%
\end{array}%
\right] =\Phi _{\vartheta ,l}\left( \beta \mathbf{r}-\mathbf{r}_{l}\right) 
\mathbf{g}_{\vartheta }\mathrm{e}^{\mathrm{i}\mathbf{k}_{\ast l}\cdot 
\mathbf{r}}+\limfunc{cc},  \label{zpw}
\end{equation}%
with $\Phi _{\vartheta ,l}\left( \mathbf{r}\right) $ being arbitrary
Schwartz functions. If\textbf{\ }(\ref{kjl}) holds, \ the coupling terms $%
\mathbf{D}_{1}$ and $\mathbf{D}_{2}$ satisfy the bound 
\begin{equation}
\sup_{0\leq \tau \leq \tau _{\ast }}\left[ \left\Vert \mathbf{D}_{1}\left(
\tau \right) \right\Vert _{L_{\infty }}+\left\Vert \mathbf{D}_{2}\left( \tau
\right) \right\Vert _{L_{\infty }}\right] \leq C_{\delta }^{\prime }\frac{%
\varrho }{\beta ^{1+\delta }},  \label{D1D21}
\end{equation}%
where $\tau _{\ast }$ and $C_{\delta }^{\prime }$ do not depend on $\beta $,$%
\varrho $ and $\delta $ can be taken arbitrary small..
\end{theorem}

In the following sections we introduce concepts and develop analytic tools
allowing to prove the approximate linear superposition principle as stated
in Theorems \ref{Theorem Superposition}, \ref{Theorem Superposition1} and %
\ref{Theorem Superposition11}.

\section{Reduced evolution equation}

Since the properties of the evolution equations (\ref{eqFur}) and (\ref{difF}%
) are very similar, we consider here in detail the lattice evolution
equation (\ref{eqFur}) with understanding that all the statements apply to
the PDE (\ref{difF}) if we replace $\mathbf{\tilde{U}}$ with $\mathbf{\hat{U}%
}$, $\left[ -\pi ,\pi \right] ^{d}$ with $\mathbb{R}^{d}$, the function
space $L_{1}=L_{1}\left( \left[ -\pi ,\pi \right] ^{d}\right) $ with $%
L_{1}=L_{1}\left( \mathbb{R}^{d}\right) $ and so on.

First, using the variation of constants formula we recast the modal
evolution equation (\ref{eqFur}) into the following equivalent integral form 
\begin{equation}
\mathbf{\tilde{U}}\left( \mathbf{k},\tau \right) =\int_{0}^{\tau }\mathrm{e}%
^{\frac{-\mathrm{i}\left( \tau -\tau ^{\prime }\right) }{\varrho }\mathbf{L}%
\left( \mathbf{k}\right) }\tilde{F}\left( \mathbf{\tilde{U}}\right) \left( 
\mathbf{k},\tau \right) d\tau ^{\prime }+\mathrm{e}^{\frac{-\mathrm{i}\zeta
\tau }{\varrho }\mathbf{L}\left( \mathbf{k}\right) }\mathbf{\tilde{h}}\left( 
\mathbf{k}\right) ,\ \tau \geq 0.  \label{varc}
\end{equation}%
Then we introduce for $\mathbf{\tilde{U}}\left( \mathbf{k},\tau \right) $
its \emph{two-time-scale representation} (with respectively slow and fast
times $\tau $ and $t=\frac{\tau }{\varrho }$) 
\begin{equation}
\mathbf{\tilde{U}}\left( \mathbf{k},\tau \right) =\mathrm{e}^{-\frac{\mathrm{%
i}\tau }{\varrho }\mathbf{L}\left( \mathbf{k}\right) }\mathbf{\tilde{u}}%
\left( \mathbf{k},\tau \right) ,\ \mathbf{\tilde{U}}_{n,\zeta }\left( 
\mathbf{k},\tau \right) =\mathbf{\tilde{u}}_{n,\zeta }\left( \mathbf{k},\tau
\right) \mathrm{e}^{-\frac{\mathrm{i}\tau }{\varrho }\zeta \omega _{n}\left( 
\mathbf{k}\right) },  \label{Uu}
\end{equation}%
where $\mathbf{\tilde{u}}_{n,\zeta }\left( \mathbf{k},\tau \right) $ are the
modal coefficients of $\mathbf{\tilde{u}}\left( \mathbf{k},\tau \right) $
(see (\ref{Uboldj})); note that $\mathbf{\tilde{u}}_{n,\zeta }\left( \mathbf{%
k},\tau \right) $ may depend on $\varrho $, therefore (\ref{Uu}) is just a
change of variables. Consequently we obtain the following \emph{reduced
evolution equation} for $\mathbf{\tilde{u}}=\mathbf{\tilde{u}}\left( \mathbf{%
k},\tau \right) $, $\tau \geq 0$, 
\begin{gather}
\mathbf{\tilde{u}}\left( \mathbf{k},\tau \right) =\mathcal{F}\left( \mathbf{%
\tilde{u}}\right) \left( \mathbf{k},\tau \right) +\mathbf{\tilde{h}}\left( 
\mathbf{k}\right) ,\ \mathcal{F}\left( \mathbf{\tilde{u}}\right)
=\sum_{m=2}^{m_{F}}\mathcal{F}^{\left( m\right) }\left( \mathbf{\tilde{u}}%
^{m}\left( \mathbf{k},\tau \right) \right) ,  \label{varcu} \\
\mathcal{F}^{\left( m\right) }\left( \mathbf{\tilde{u}}^{m}\right) \left( 
\mathbf{k},\tau \right) =\int_{0}^{\tau }\mathrm{e}^{\frac{\mathrm{i}\tau
^{\prime }}{\varrho }\mathbf{L}\left( \mathbf{k}\right) }\tilde{F}^{\left(
m\right) }\left( \left( \mathrm{e}^{\frac{-\mathrm{i}\tau ^{\prime }}{%
\varrho }\mathbf{L}\left( \cdot \right) }\mathbf{\tilde{u}}\right)
^{m}\right) \left( \mathbf{k},\tau ^{\prime }\right) \mathrm{d}\tau ^{\prime
},  \label{Fu}
\end{gather}%
where the quantities $\tilde{F}^{\left( m\right) }$ are defined by (\ref%
{Fseries}) and (\ref{Fmintr}) in terms of the susceptibilities $\chi
^{\left( m\right) }$.

The norm of the oscillatory integral $\mathcal{F}^{\left( m\right) }$ \ in (%
\ref{Fu}) is estimated in terms of the norm of the tensor $\chi ^{\left(
m\right) }\left( \mathbf{\mathbf{k}},\vec{k}\right) $ defined in (\ref{chiCR}%
), (\ref{normchi0}). The operator $\mathcal{F}^{\left( m\right) }$ is shown
to be a bounded one from $\left( E\right) ^{m}$ into $E$, see Lemma \ref%
{Lemma 8} for details The proof of this property is based on the following
Young inequality for the convolution 
\begin{equation}
\left\Vert \mathbf{\tilde{u}}\ast \mathbf{\tilde{v}}\right\Vert _{L_{1}}\leq
\left\Vert \mathbf{\tilde{u}}\right\Vert _{L_{1}}\left\Vert \mathbf{\tilde{v}%
}\right\Vert _{L_{1}}.  \label{Yconv}
\end{equation}

For a detailed analysis of \ solutions of (\ref{varcu}) we recast the
equation (\ref{varcu}) for $\mathbf{\tilde{u}}\left( \mathbf{k},\tau \right) 
$ using projections (\ref{Pin}) as the following \emph{expanded reduced
evolution equation} 
\begin{equation}
\mathbf{\tilde{u}}_{n,\zeta }\left( \mathbf{k},\tau \right)
=\sum_{m=2}^{\infty }\sum_{\vec{n},\vec{\zeta}}\mathcal{F}_{n,\zeta ,\vec{n},%
\vec{\zeta}}^{\left( m\right) }\left( \mathbf{\tilde{u}}^{m}\right) \left( 
\mathbf{k},\tau \right) +\mathbf{h}_{n,\zeta }\left( \mathbf{k}\right) ,\
\tau \geq 0,  \label{equfa}
\end{equation}%
for the modal coefficient $\mathbf{\tilde{u}}_{n,\zeta }\left( \mathbf{k}%
,\tau \right) $. In the above formula and elsewhere we use notations 
\begin{equation}
\vec{n}=\left( n^{\prime },\ldots ,n^{\left( m\right) }\right) ,\ \vec{\zeta}%
=\left( \zeta ^{\prime },\ldots ,\zeta ^{\left( m\right) }\right) ,\ \vec{k}%
=\left( \mathbf{\mathbf{k}}^{\prime },\ldots ,\mathbf{\mathbf{k}}^{\left(
m\right) }\right) .  \label{zetaar}
\end{equation}%
The operators $\mathcal{F}_{n,\zeta ,\vec{n},\vec{\zeta}}^{\left( m\right) }$
are $m$-linear \emph{oscillatory integral operators} defined by the formulas 
\begin{gather}
\mathcal{F}_{n,\zeta ,\vec{n},\vec{\zeta}}^{\left( m\right) }\left( \mathbf{%
\tilde{u}}_{1}\ldots \mathbf{\tilde{u}}_{m}\right) \left( \mathbf{k},\tau
\right) =\int_{0}^{\tau }\int_{\mathbb{D}_{m}}\exp \left\{ \mathrm{i}\phi
_{n,\zeta ,\vec{n},\vec{\zeta}}\left( \mathbf{\mathbf{k}},\vec{k}\right) 
\frac{\tau _{1}}{\varrho }\right\}  \label{Fm} \\
\chi _{n,\zeta ,\vec{n},\vec{\zeta}}^{\left( m\right) }\left( \mathbf{%
\mathbf{k}},\vec{k}\right) \left[ \mathbf{\tilde{u}}_{1}\left( \mathbf{k}%
^{\prime },\tau _{1}\right) ,\ldots ,\mathbf{\tilde{u}}_{m}\left( \mathbf{k}%
^{\left( m\right) }\left( \mathbf{\mathbf{k}},\vec{k}\right) ,\tau
_{1}\right) \right] \mathrm{\tilde{d}}^{\left( m-1\right) d}\vec{k}\mathrm{d}%
\tau _{1},  \notag
\end{gather}%
where we use notations (\ref{Dm}), (\ref{dtild}), (\ref{kkar}). In (\ref{Fm}%
) the \emph{interaction phase function }$\phi $ is defined by 
\begin{equation}
\phi _{n,\zeta ,\vec{n},\vec{\zeta}}\left( \mathbf{\mathbf{k}},\vec{k}%
\right) =\zeta \omega _{n}\left( \mathbf{k}\right) -\zeta ^{\prime }\omega
_{n^{\prime }}\left( \mathbf{k}^{\prime }\right) -\ldots -\zeta ^{\left(
m\right) }\omega _{n^{\left( m\right) }}\left( \mathbf{k}^{\left( m\right)
}\right) ,\ \mathbf{k}^{\left( m\right) }=\mathbf{k}^{\left( m\right)
}\left( \mathbf{\mathbf{k}},\vec{k}\right)  \label{phim}
\end{equation}%
and the \emph{susceptibilities} $\chi _{n,\zeta ,\vec{n},\vec{\zeta}%
}^{\left( m\right) }\left( \mathbf{\mathbf{k}},\vec{k}\right) $ are$\ m$%
-linear symmetric\ tensors (i.e. mappings from $\left( \mathbb{C}%
^{2J}\right) ^{m}$ into $\mathbb{C}^{2J}$) defined for almost all $\mathbf{%
\mathbf{k}},\vec{k}$ by the following formula 
\begin{gather}
\chi _{n,\zeta ,\vec{n},\vec{\zeta}}^{\left( m\right) }\left( \mathbf{%
\mathbf{k}},\vec{k}\right) \left[ \mathbf{\tilde{u}}_{1}\left( \mathbf{k}%
^{\prime }\right) ,\ldots ,\mathbf{\tilde{u}}_{m}\left( \mathbf{k}^{\left(
m\right) }\right) \right] =  \label{chim} \\
\Pi _{n,\zeta }\left( \mathbf{\mathbf{k}}\right) \chi ^{\left( m\right)
}\left( \mathbf{\mathbf{k}},\vec{k}\right) \left[ \Pi _{n^{\prime },\zeta
^{\prime }}\left( \mathbf{k}^{\prime }\right) \mathbf{\tilde{u}}_{1}\left( 
\mathbf{k}^{\prime }\right) ,\ldots ,\Pi _{n^{\left( m\right) },\zeta
^{\left( m\right) }}\left( \mathbf{k}^{\left( m\right) }\left( \mathbf{k},%
\vec{k}\right) \right) \mathbf{\tilde{u}}_{m}\left( \mathbf{k}^{\left(
m\right) }\left( \mathbf{k},\vec{k}\right) \right) \right] .  \notag
\end{gather}%
For the lattice equation $\chi _{\ n,\zeta ,\vec{n},\vec{\zeta}}^{\left(
m\right) }\left( \mathbf{\mathbf{k}},\vec{k}\right) $ is $2\pi $-periodic
with respect to every variable $\mathbf{\mathbf{k}},\mathbf{\mathbf{k}}%
^{\prime },\ldots ,\mathbf{\mathbf{k}}^{\left( m\right) }$. Note that
operators $\mathcal{F}^{\left( m\right) }\left( \mathbf{u}^{m}\right) $ in (%
\ref{varcu}) can be rewritten using (\ref{Fm}) as 
\begin{equation}
\mathcal{F}^{\left( m\right) }\left( \mathbf{u}^{m}\right) =\sum_{\vec{n},%
\vec{\zeta}}\mathcal{F}_{n,\zeta ,\vec{n},\vec{\zeta}}^{\left( m\right)
}\left( \mathbf{\tilde{u}}^{m}\right) .  \label{Tplusmin}
\end{equation}%
We also call operators $\mathcal{F}_{n,\zeta ,\vec{n},\vec{\zeta}}^{\left(
m\right) }$ \emph{decorated operators. }

\begin{remark}
The expanded reduced evolution equation (\ref{equfa}) is instrumental to the
nonlinear analysis. Its very form, a convergent series of multilinear forms
which are oscillatory integrals (\ref{Fm}), is already a significant step in
the analysis of the solution accomplishing several tasks: (i) it suggests a
constructive representation for the solution; (ii) every term $\mathcal{F}%
_{n,\zeta ,\vec{n},\vec{\zeta}}^{\left( m\right) }$ can be naturally
interpreted as nonlinear interaction of the underlying linear modes; (iii)
the representation of $\mathcal{F}_{n,\zeta ,\vec{n},\vec{\zeta}}^{\left(
m\right) }$ as the oscillatory integral (\ref{Fm}) involving the interaction
phase $\phi _{n,\zeta ,\vec{n},\vec{\zeta}}$ and the susceptibilities $\chi
_{n,\zeta ,\vec{n},\vec{\zeta}}^{\left( m\right) }\left( \mathbf{\mathbf{k}},%
\vec{k}\right) $ directly relates $\mathcal{F}_{n,\zeta ,\vec{n},\vec{\zeta}%
}^{\left( m\right) }$ to the terms of the original evolution equation as
well as to physically significant quantities. We can also add that since we
consider $\varrho \rightarrow 0$ the interaction phase function $\phi
_{n,\zeta ,\vec{n},\vec{\zeta}}\left( \mathbf{\mathbf{k}},\vec{k}\right) $
plays the decisive role in the analysis of nonlinear interactions of
different modes.
\end{remark}

The analysis of fundamental properties of the reduced evolution equation (%
\ref{equfa}), including, in particular, the linear modal superposition
principle, involves and combines the following three components: (i) \emph{%
the linear spectral theory component} in the form of \emph{the modal
decomposition} of the solution and introduction of \emph{wavepackets} as
elementary waves; (ii) \emph{function-analytic component} which deals with
the structure of series similar to the one in (\ref{equfa}) and its
dependence on the nonlinearity of the original evolution equation; (iii) 
\emph{asymptotic analysis of oscillatory integrals} (\ref{Fm}) which allows
to estimate the magnitude of nonlinear interactions between different modes
and, in particular, to show that generically different modes almost do not
interact leading to the superposition principle.

Sometimes it is convenient to rewrite (\ref{Fm}) in a slightly different
form. The convolution integral (\ref{Fm}) according to (\ref{kkar}) involves
the following \emph{phase matching condition} 
\begin{equation}
\mathbf{k}^{\prime }+\ldots +\mathbf{k}^{\left( m\right) }=\mathbf{k}.
\label{Phmc}
\end{equation}%
Using the following notation for the integral over the plane (\ref{Phmc}) 
\begin{gather}
\int\limits_{\mathbf{k}^{\prime },\mathbf{\ldots ,k}^{\left( m-1\right) }\in %
\left[ -\pi ,\pi \right] ^{\left( m-1\right) d};\mathbf{k}^{\prime }+\mathbf{%
\ldots +k}^{\left( m\right) }=\mathbf{k}}f\left( \mathbf{\mathbf{k}},\vec{k}%
\right) \,\mathrm{d}\mathbf{k}^{\prime }\ldots \,\mathrm{d}\mathbf{k}%
^{\left( m-1\right) }=  \label{delk} \\
\int_{\left[ -\pi ,\pi \right] ^{md}}f\left( \mathbf{\mathbf{k}},\vec{k}%
\right) \delta \left( \mathbf{k}-\mathbf{k}^{\prime }-\ldots -\mathbf{k}%
^{\left( m\right) }\right) \mathrm{d}\mathbf{k}^{\prime }\ldots \,\mathrm{d}%
\mathbf{k}^{\left( m\right) }  \notag
\end{gather}%
in terms of a delta-function we can rewrite (\ref{Fm}) in the form 
\begin{gather}
\mathcal{F}_{n,\zeta ,\vec{n},\vec{\zeta}}^{\left( m\right) }\left( \mathbf{%
\tilde{u}}_{1}\ldots \mathbf{\tilde{u}}_{m}\right) \left( \mathbf{k},\tau
\right) =\frac{1}{\left( 2\pi \right) ^{m\left( d-1\right) }}\int_{0}^{\tau
}\int_{\left[ -\pi ,\pi \right] ^{md}}\exp \left\{ \mathrm{i}\phi _{n,\zeta ,%
\vec{n},\vec{\zeta}}\left( \mathbf{\mathbf{k}},\vec{k}\right) \frac{\tau _{1}%
}{\varrho }\right\}   \label{Fmdel} \\
\delta \left( \mathbf{k}-\mathbf{k}^{\prime }-\ldots -\mathbf{k}^{\left(
m\right) }\right) \chi _{n,\zeta ,\vec{n},\vec{\zeta}}^{\left( m\right)
}\left( \mathbf{\mathbf{k}},\vec{k}\right) \mathbf{\tilde{u}}_{1,\zeta
^{\prime }}\left( \mathbf{k}^{\prime }\right) \ldots \mathbf{\tilde{u}}%
_{m,\zeta ^{\left( m\right) }}\left( \mathbf{k}^{\left( m\right) }\right) \,%
\mathrm{d}\mathbf{k}^{\prime }\ldots \,\mathrm{d}\mathbf{k}^{\left( m\right)
}d\tau _{1}.  \notag
\end{gather}

\section{Function-analytic operator series}

In this section necessary algebraic concepts required for the analysis are
introduced. We study the reduced evolution equation (\ref{varcu}) as a
particular case of the following \emph{abstract nonlinear equation} in a
Banach space 
\begin{equation}
\mathbf{\mathbf{\mathbf{u}}}=\mathcal{F}\left( \mathbf{\mathbf{\mathbf{u}}}%
\right) +\mathbf{x},\ \mathcal{F}\left( \mathbf{\mathbf{\mathbf{u}}}\right)
=\sum_{s=2}^{\infty }\mathcal{F}^{\left( s\right) }\left( \mathbf{x}%
^{s}\right)  \label{eqF}
\end{equation}%
where the \emph{nonlinearity} $\mathcal{F}\left( \mathbf{\mathbf{\mathbf{u}}}%
\right) $ is an analytic operator represented by a convergent operator
series. It is well known (see \cite{HPh}) that the solution $\mathbf{\mathbf{%
\mathbf{u}}}=\mathcal{G}\left( \mathbf{x}\right) $ of such equation can be
represented as a convergent series in terms of$\ m$-linear operators $%
\mathcal{G}_{m}$ which are constructed based on $\mathcal{F}$: $\ $%
\begin{gather*}
\mathcal{G}\left( \mathbf{x}\right) =\mathcal{G}\left( \mathcal{F},\mathbf{x}%
\right) =\sum_{m=1}^{\infty }\mathcal{G}^{\left( m\right) }\left( \mathbf{x}%
^{m}\right) ,\ \mathcal{G}^{\left( m\right) }\left( \mathbf{x}^{m}\right) =%
\mathcal{G}^{\left( m\right) }\left( \mathcal{F},\mathbf{x}^{m}\right) ,%
\text{ where} \\
\mathbf{x}^{m}\underset{m\text{ times}}{=\;\underbrace{\mathbf{x}\ldots 
\mathbf{x}}}\;=\;\underset{m\text{ times}}{\underbrace{\mathbf{x}\ldots 
\mathbf{x}}}.
\end{gather*}%
Using the multilinearity of $\mathcal{G}^{\left( m\right) }$ we readily
obtain the formula 
\begin{gather}
\mathcal{G}\left( \mathbf{x}_{1}+\ldots +\mathbf{x}_{N}\right)
=\sum_{m=1}^{\infty }\mathcal{G}^{\left( m\right) }\left( \left( \mathbf{x}%
_{1}+\ldots +\mathbf{x}_{N}\right) ^{m}\right)  \label{GCI} \\
=\sum_{m=1}^{\infty }\mathcal{G}\left( \left( \mathbf{x}_{1}\right)
^{m}\right) +\ldots +\sum_{m=1}^{\infty }\mathcal{G}\left( \left( \mathbf{x}%
_{N}\right) ^{m}\right) +\mathcal{G}_{\text{CI}}\left( \mathbf{x}_{1},\ldots
,\mathbf{x}_{N}\right) ,  \notag
\end{gather}%
where $\mathbf{x}=\mathbf{x}_{1}+\ldots +\mathbf{x}_{N}$ represents a
multi-wavepacket and $\mathcal{G}_{\text{CI}}\left( \mathbf{x}_{1},\ldots ,%
\mathbf{x}_{N}\right) $ collects all "cross terms" and describes the "cross
interaction" (CI) of involved wavepackets $\mathbf{x}_{1},\ldots ,\mathbf{x}%
_{N}$. We will find in sufficient detail the dependence of the solution
operators $\mathcal{G}_{m}$ on the nonlinearity $\mathcal{F}$ and prepare a
basis for the consequent estimation of nonlinear interactions between
different modes and wavepackets. \ Then combining the facts about the
structure of the solution operators $\mathcal{G}^{\left( m\right) }$ with
asymptotic estimates of relevant oscillatory integrals we show that for a
multi-wavepacket $\mathbf{x}=\mathbf{x}_{1}+\ldots +\mathbf{x}_{N}$ the
cross-interaction term satisfies the following estimate 
\begin{equation*}
\left\Vert \mathcal{G}_{\text{CI}}\left( \mathbf{x}_{1},\ldots ,\mathbf{x}%
_{N}\right) \right\Vert =O\left( \beta \right) +O\left( \varrho \left\vert
\ln \beta \right\vert /\beta ^{1+\epsilon }\right) ,\ \beta ,\varrho
\rightarrow 0,
\end{equation*}%
implying the modal superposition principle.

\subsection{Multilinear forms and polynomial operators}

The analysis of nonlinear equations of the form (\ref{varcu}) requires the
use of appropriate Banach spaces of time dependent fields, as well as
multilinear and analytic functions in those spaces. It also uses an
appropriate version of the implicit function theorem. For the reader's
convenience we collect in this section the known concepts and statements on
the above-mentioned subjects needed for our analysis. We in this section
consider functional-analytic operators which are defined in a ball in a
Banach space $X$ with the norm $\left\Vert \mathbf{x}\right\Vert _{X}$. In
our treatment of the analytic functions in infinitely-dimensional Banach
spaces we follow to \cite[Section 3]{HPh}, \cite{Din}.

\begin{definition}[polylinear operator]
Suppose that $\mathbf{x}_{1},\mathbf{x}_{2},\ldots ,\mathbf{x}_{n}$ are
vectors in a Banach space $X$. Let a function $\ \mathcal{F}^{\left(
n\right) }\left( \vec{x}\right) $, $\vec{x}=\left( \mathbf{x}_{1},\ldots ,%
\mathbf{x}_{n}\right) $, take values in $X$ and be defined for all $\vec{x}%
\in X^{n}$. Such a function $\mathcal{F}^{\left( n\right) }$ is called a $n$%
\emph{-linear operator} if it is linear in each variable, and it is said to
be bounded if its following norm is finite%
\begin{equation}
\left\Vert \mathcal{F}^{\left( n\right) }\right\Vert =\sup_{\left\Vert 
\mathbf{x}_{1}\right\Vert _{X}=\ldots =\left\Vert \mathbf{x}_{n}\right\Vert
_{X}=1}\left\Vert \mathcal{F}^{\left( n\right) }\left( \mathbf{x}_{1}\mathbf{%
x}_{2}\ldots \mathbf{x}_{n}\right) \right\Vert _{X}<\infty .  \label{normXY}
\end{equation}
\end{definition}

\begin{definition}[polynomial]
A function $P\left( x\right) $ from $X$ to $X$ defined for all $x\in X$ is
called a \emph{polynomial} in $\mathbf{x}$ of \emph{degree} $n$ if for all $%
\mathbf{a},\mathbf{h}\in X$ and all complex $\alpha $%
\begin{equation*}
P\left( \mathbf{a}+\alpha \mathbf{h}\right) =\sum_{\nu =0}^{n}P_{\nu }\left( 
\mathbf{a},\mathbf{h}\right) \alpha ^{\nu },
\end{equation*}%
where $P_{\nu }\left( a,h\right) \in X$ are independent of $\alpha \ $. The
degree of $P_{n}$ is exactly $n$ if $P_{n}\left( a,h\right) $ is not
identically zero. A polynomial $\mathcal{F}\left( \mathbf{x}\right) $ is a 
\emph{homogeneous polynomial} of a degree $n$ if for all $c\in \mathbb{C}$ 
\begin{equation*}
\mathcal{F}\left( c\mathbf{x}\right) =c^{n}\mathcal{F}\left( \mathbf{x}%
\right) .
\end{equation*}%
Then $n$ is called also the \emph{homogeneity index} of $\mathcal{F}\left( 
\mathbf{x}\right) $. A homogeneous polynomial $\mathcal{F}$ is called
bounded if its norm 
\begin{equation}
\left\Vert \mathcal{F}\right\Vert _{\ast }=\sup_{\left\Vert x\right\Vert
_{X}=1}\left\{ \left\Vert \mathcal{F}\left( \mathbf{x}\right) \right\Vert
_{X}\right\}  \label{normF*}
\end{equation}%
is finite. For a given $n$-linear operator $\mathcal{F}^{\left( n\right)
}\left( \vec{x}\right) =\ \mathcal{F}^{\left( n\right) }\left( \mathbf{x}_{1}%
\mathbf{x}_{2}\ldots \mathbf{x}_{n}\right) $ we denote by $\mathcal{F}%
^{\left( n\right) }\left( \mathbf{x}^{n}\right) $ a homogeneous of degree $n$
polynomial \ from $X$ to $X$: 
\begin{equation}
\mathcal{F}^{\left( n\right) }\left( \mathbf{x}^{n}\right) =\mathcal{F}%
^{\left( n\right) }\left( \mathbf{x}\ldots \mathbf{x}\right) .
\label{normF*1}
\end{equation}
\end{definition}

Note the norm definitions (\ref{normXY}), (\ref{normF*}) and (\ref{normF*1})
readily imply%
\begin{equation}
\left\Vert \mathcal{F}^{\left( n\right) }\right\Vert _{\ast }\leq \left\Vert 
\mathcal{F}^{\left( n\right) }\right\Vert .  \label{normn}
\end{equation}

\begin{definition}[analyticity class 1]
\label{defA*}Let a function $\mathcal{F}$ \ be defined by the following
convergent series 
\begin{equation}
\mathcal{F}\left( x\right) =\sum_{m=2}^{\infty }\mathcal{F}^{\left( m\right)
}\left( \mathbf{x}^{m}\right) \text{ for }\left\Vert \mathbf{x}\right\Vert
_{X}<R_{\ast \mathcal{F}},  \label{deff}
\end{equation}%
where $\mathcal{F}^{\left( m\right) }\left( \mathbf{x}^{m}\right) $, $%
m=2,3,\ldots $ is a sequence of bounded $m$-homogenious polynomials
satisfying 
\begin{equation}
\left\Vert \mathcal{F}^{\left( m\right) }\right\Vert _{\ast }\leq C_{\ast 
\mathcal{F}}R_{\ast \mathcal{F}}^{-m},\ m=2,3,\ldots ,  \label{deff1}
\end{equation}%
Then we say that $\mathcal{F}\left( x\right) $ belongs to the \emph{%
analyticity class} $A_{\ast }\left( C_{\ast \mathcal{F}},R_{\ast \mathcal{F}%
}\right) $ and write $\mathcal{F}\in A_{\ast }\left( C_{\ast \mathcal{F}%
},R_{\ast \mathcal{F}}\right) .$
\end{definition}

Notice that for $\left\Vert \mathbf{x}\right\Vert _{X}<R_{\ast \mathcal{F}}$
we have 
\begin{equation}
\left\Vert \mathcal{F}\left( \mathbf{x}\right) \right\Vert _{X}\leq C_{\ast 
\mathcal{F}}\sum_{n=2}^{\infty }\left\Vert \mathbf{x}\right\Vert
_{X}^{n}R_{\ast \mathcal{F}}^{-n}\leq C_{\ast \mathcal{F}}\frac{\left\Vert 
\mathbf{x}\right\Vert _{X}^{n_{0}}R_{\ast \mathcal{F}}^{-n_{0}}}{%
1-\left\Vert \mathbf{x}\right\Vert _{X}R_{\ast \mathcal{F}}^{-1}},
\label{sumineq}
\end{equation}
implying, in particular, the convergence of the series (\ref{deff}).

\begin{definition}[analyticity class 2]
\label{defan}If $\mathcal{F}^{\left( m\right) }\left( \vec{x}\right) $, $%
m=2,3,\ldots $, is a sequence of bounded $m$-linear operators from $X^{m}$
to $X$ and 
\begin{equation}
\left\Vert \mathcal{F}^{\left( m\right) }\right\Vert \leq C_{\mathcal{F}}R_{%
\mathcal{F}}^{-m},\ m=2,3\ldots .,  \label{normFj}
\end{equation}%
we say that a function $\mathcal{F}$ defined by the series (\ref{deff}) for $%
\left\Vert x\right\Vert _{X}<R_{\mathcal{F}}$ belongs to the analyticity
class $A\left( C_{\mathcal{F}},R_{\mathcal{F}}\right) $ and write $\mathcal{F%
}\in A\left( C_{\mathcal{F}},R_{\mathcal{F}}\right) .$
\end{definition}

In this paper we will use operators from the classes $A\left( C_{\mathcal{F}%
},R_{\mathcal{F}}\right) $ based on multilinear operators.

Note that evidently $A\left( C_{\mathcal{F}},R_{\mathcal{F}}\right) \subset
A_{\ast }\left( C_{\mathcal{F}},R_{\mathcal{F}}\right) $. One can construct
a polynomial based on a multilinear operator according to the formula (\ref%
{normF*1}). Conversely, the construction of a multilinear operator, called 
\emph{polar form}, based on a given homogeneous polynomial is described by
the following statement, \cite[Section 1.1, 1.3]{Din}, \cite[Section 26.2]%
{HPh}.

\begin{proposition}[polar form]
\label{polar}For any homogeneous polynomial $P^{\left( n\right) }\left(
x\right) $ of degree $n$ there is a unique symmetric $n$-linear operator $%
\tilde{P}^{\left( n\right) }\left( \mathbf{x}_{1}\mathbf{x}_{2}\ldots 
\mathbf{x}_{n}\right) $, called the \emph{polar form} of $P_{n}\left( 
\mathbf{x}\right) $, such that $P^{\left( n\right) }\left( \mathbf{x}\right)
=\tilde{P}^{\left( n\right) }\left( \mathbf{x}\ldots \mathbf{x}\right) $. It
is defined by the following \emph{polarization formula}:%
\begin{equation}
\tilde{P}^{\left( n\right) }\left( \mathbf{x}_{1}\mathbf{x}_{2}\ldots 
\mathbf{x}_{n}\right) =\frac{1}{2^{n}n!}\sum_{\xi _{j}=\pm 1}P^{\left(
n\right) }\left( \sum_{j=1}^{n}\xi _{j}\mathbf{x}_{j}\right) .  \label{Ppol}
\end{equation}%
In addition to that, the following estimate holds: 
\begin{equation}
\left\Vert P_{n}\right\Vert _{\ast }\leq \left\Vert \tilde{P}_{n}\right\Vert
\leq \frac{n^{n}}{n!}\left\Vert P_{n}\right\Vert _{\ast }\leq \mathrm{e}%
^{n}\left\Vert P^{\left( n\right) }\right\Vert _{\ast }.  \label{Ppol1}
\end{equation}
\end{proposition}

Since by Definition \ref{defan} functions from $A\left( C,R\right) $ have
zero of the second order at zero, their Lipschitz constant is small in a
vicinity of zero. More exactly, the following statement holds.

\begin{lemma}[Lipschitz estimate]
\label{Lemma 101Continuity} \ If $\mathcal{F}\in A\left( C_{\mathcal{F}},R_{%
\mathcal{F}}\right) $ then 
\begin{equation}
\left\Vert \mathcal{F}\left( \mathbf{x}\right) -\mathcal{F}\left( \mathbf{y}%
\right) \right\Vert \leq C_{\mathcal{F}}C\left\Vert \mathbf{x}-\mathbf{y}%
\right\Vert \left( \left\Vert \mathbf{x}\right\Vert +\left\Vert \mathbf{y}%
\right\Vert \right) \text{ for }\left\Vert \mathbf{x}\right\Vert ,\left\Vert 
\mathbf{y}\right\Vert \leq R_{\mathcal{F}}^{\prime }<R_{\mathcal{F}},
\label{Lipf}
\end{equation}%
where $C>0$\ depends on $R_{\mathcal{F}}^{\prime }$ and $R_{\mathcal{F}}$.
\end{lemma}

\subsection{Implicit Function Theorem and expansion of operators into
composition monomials}

Here we provide a version of the Implicit Function Theorem, first we
formulate classical implicit function theorem for equations $\mathbf{\mathbf{%
\mathbf{u}}}=\mathcal{F}\left( \mathbf{\mathbf{\mathbf{u}}}\right) +\mathbf{x%
}$ with analytic function $\mathcal{F}$ and then we present a refined
implicit function theorem. The refined implicit function theorem we prove
here produces expansion of the solution $\mathbf{\mathbf{\mathbf{u}}}$ into
a sum of terms which are multi-linear not only with respect to $\mathbf{x}$
but also with respect to $\mathcal{F}$. The formulation of the theorem and
the proof involve convenient labeling of the terms of the expansion (called
composition monomials), and we use properly introduced trees to this end.
The explicit expansion produced by the refined implicit function theorem is
required to be able to take into account rather subtle mechanisms which lead
to the superposition principle.

Let us consider the abstract nonlinear equation (\ref{eqF}) and its solution 
$\mathbf{u}=\mathbf{u}\left( \mathbf{x}\right) $ for small $\left\Vert 
\mathbf{x}\right\Vert $ when the nonlinear operator $\mathcal{F}\ $belongs
to the class $A\left( C_{\mathcal{F}},R_{\mathcal{F}}\right) $. We seek the
solution $\mathbf{u}$ in the following form%
\begin{equation}
\mathbf{u}=\mathcal{G}\left( \mathcal{F},\mathbf{x}\right)
=\sum_{m=1}^{\infty }\mathcal{G}^{\left( m\right) }\left( \mathbf{x}%
^{m}\right) \text{ for sufficiently small }\left\Vert \mathbf{x}\right\Vert ,
\label{Gser}
\end{equation}%
and we call $\mathcal{G}$ the \emph{solution operator} for (\ref{eqF}). It
readily follows from (\ref{eqF}) that 
\begin{equation}
\mathcal{G}\left( \mathcal{F},\mathbf{x}\right) =\mathbf{x}+\mathcal{F}%
\left( \mathcal{G}\left( \mathcal{F},\mathbf{x}\right) \right)  \label{FG}
\end{equation}%
and 
\begin{equation}
\sum_{m=1}^{\infty }\mathcal{G}^{\left( m\right) }\left( \mathbf{\mathbf{x}}%
^{m}\right) =\mathbf{\mathbf{x}}+\sum_{s=2}^{\infty }\mathcal{F}^{\left(
s\right) }\left( \left( \sum_{m=1}^{\infty }\mathcal{G}^{\left( m\right)
}\left( \mathbf{x}^{m}\right) \right) ^{s}\right) .  \label{FGser}
\end{equation}%
>From the above equation we can deduce recurrent formulas for multilinear
operators $\mathcal{G}^{\left( m\right) }$. Indeed for $m=1$\ the linear
term is the identity operator 
\begin{equation}
\mathcal{G}^{\left( 1\right) }\left( \mathbf{x}\right) =\mathcal{F}^{\left(
1\right) }\left( \mathbf{x}\right) \equiv \mathbf{x}.  \label{G1}
\end{equation}%
\ For $m\geq 2$ we write the following recurrent formula 
\begin{equation}
\mathcal{G}^{\left( m\right) }\left( \mathbf{x}_{1}\ldots \mathbf{x}%
_{m}\right) =\sum_{s=2}^{m}\sum_{i_{1}+\ldots +i_{s}=m}\mathcal{F}^{\left(
s\right) }\left( \mathcal{G}^{\left( i_{1}\right) }\left( \mathbf{x}%
_{1}\ldots \mathbf{x}_{i_{1}}\right) \ldots \mathcal{G}^{\left( i_{s}\right)
}\left( \mathbf{x}_{m-i_{s}+1}\ldots \mathbf{x}_{m}\right) \right) .
\label{recG0}
\end{equation}%
By the construction, if multilinear operators $\mathcal{G}^{\left( i\right)
} $ are defined by (\ref{recG0}), then (\ref{FGser}) is satisfied. Namely,
expanding right-hand side of (\ref{FGser}) using multi-linearity of $%
\mathcal{F}^{\left( s\right) }$ we obtain a sum of expressions \ as in
right-hand side of (\ref{recG0}), and since (\ref{recG0}) holds, terms in
the left-hand side of (\ref{FGser}) with given homogeneity index $p$ cancel
with the terms in the right-hand side with the same homogeneity. \ Note that
in (\ref{recG0}) \ we do not assume that the operators $\mathcal{F}^{\left(
s\right) }$ and $\mathcal{G}^{\left( i\right) }$ are symmetrized and the
order of variables is important; \ we prefer to treat $\mathcal{F}^{\left(
s\right) }$ and $\ \mathcal{G}^{\left( m\right) }$ as multilinear operators
of $s$ and $m$ variables respectively. Though, when we apply constructed $%
\mathcal{G}^{\left( i\right) }$ to solve (\ref{eqF}), \ we set $\mathbf{x}%
_{1}=\ldots =\mathbf{x}_{m}$.

The following implicit function theorem holds (see \cite{BF4} and Theorem %
\ref{Theorem Monomial Convergence} below with a similar proof).

\begin{theorem}[Implicit Function Theorem]
\label{Imfth} Let $\ \mathcal{F}\in A\left( C_{\mathcal{F}},R_{\mathcal{F}%
}\right) $. Then there exists a solution $\mathbf{u}=\mathbf{x}+\mathcal{G}%
\left( \mathcal{F},\mathbf{x}\right) $ of the equation (\ref{eqF}) $\mathbf{%
\mathbf{\mathbf{u}}}=\mathbf{x}+\mathcal{F}\left( \mathbf{\mathbf{\mathbf{u}}%
}\right) $, given by the solution operator $\mathcal{G}\in A\left( C_{%
\mathcal{G}},R_{\mathcal{G}}\right) $, where we can take 
\begin{equation}
C_{\mathcal{G}}=\frac{R_{\mathcal{F}}^{2}}{2\left( C_{\mathcal{F}}+R_{%
\mathcal{F}}\right) },\ R_{\mathcal{G}}=\frac{R_{\mathcal{F}}^{2}}{4\left(
C_{\mathcal{F}}+R_{\mathcal{F}}\right) },  \label{CRC}
\end{equation}%
the series (\ref{Gser}) converges for $\left\Vert \mathbf{x}\right\Vert
_{X}<R_{\mathcal{G}}$. The multilinear operators $\ \mathcal{G}^{\left(
m\right) }\left( \vec{x}\right) $ satisfy the recursive relations (\ref{G1}%
), (\ref{recG0}).
\end{theorem}

Note that uniqueness of the solution and continuous dependence on parameters
follows from Lemma \ref{Lemma 101Continuity} and from a standard observation
which we formulate in the following remark.

\begin{remark}
\label{Remark Lip}If $\mathbf{u}_{1},\mathbf{u}_{2}$ are two solutions of
the equation (\ref{eqF}) with $\mathbf{x}=\mathbf{x}_{1},\mathbf{x}_{2}\ $%
respectively\ and $\left\Vert \mathbf{u}_{1}\right\Vert ,\left\Vert \mathbf{u%
}_{2}\right\Vert \leq R$, and $\mathcal{F}\left( \mathbf{u}\right) $ is
Lipschitz continuous for $\left\Vert \mathbf{u}\right\Vert \leq R$ with a
Lipschitz constant $q<1$ then $\left\Vert \mathbf{u}_{1}-\mathbf{u}%
_{2}\right\Vert \leq \left( 1-q\right) ^{-1}\left\Vert \mathbf{h}_{1}-%
\mathbf{h}_{2}\right\Vert $ . If $\ \mathbf{u}_{1},\mathbf{u}_{2}$ are two
solutions of the equation (\ref{eqF}) with $\mathcal{F}=\mathcal{F}_{0}%
\mathbf{\ \ }$and\ $\mathcal{F}=\mathcal{F}_{0}+\mathcal{F}_{1}$ \
respectively, $\left\Vert \mathbf{u}_{1}\right\Vert ,\left\Vert \mathbf{u}%
_{2}\right\Vert \leq R$, and $\mathcal{F}\left( \mathbf{u}\right) $ is
Lipschitz continuous for $\left\Vert \mathbf{u}\right\Vert \leq R$ with a
Lipschitz constant $q<1$ and $\mathcal{F}_{1}\left( \mathbf{u}\right) \leq
\epsilon $ when $\left\Vert \mathbf{u}\right\Vert \leq R$ then $\left\Vert 
\mathbf{u}_{1}-\mathbf{u}_{2}\right\Vert \leq \epsilon \left( 1-q\right)
^{-1}$ .
\end{remark}

Observe that every term $\mathcal{G}^{\left( i_{l}\right) }$ in (\ref{recG0}%
), in turn, can be recast as a sum (\ref{recG0}) with $m$ replaced by $%
i_{l}<m$. Evidently applying the recurrent representation (\ref{recG0}) and
multilinearity of $\mathcal{F}^{\left( s\right) }$ we can get a formula for $%
\mathcal{G}^{\left( m\right) }$ as a sum of terms involving exclusively (i)
the symbols $\mathcal{F}^{\left( m\right) }$, (ii) variables $\mathbf{x}_{j}$
\ and (iii) parentheses. We will refer to the terms of such a formula as 
\emph{composition monomials}. To be precise we give below a formal recursive
definition of composition monomials. The monomials are expressions which
involve variables $\mathbf{\mathbf{u}}_{j}$,$\ j=1,2,\ldots $, and $m$%
-linear operators $\mathcal{F}^{\left( m\right) }$, $m=2,3\ldots $, and are
constructed by induction as follows.

\begin{definition}[composition monomials]
\label{Definition Rank} Let $\left\{ \mathcal{F}^{\left( s\right) }\right\}
_{s=2}^{\infty }$ be a sequence of $s$-linear operators which act on
variables $\mathbf{\mathbf{u}}_{j}$, $j=1,2,\ldots $. A \emph{composition
monomial} $M\ $ of \emph{rank} $0$ is the identity operator, namely $M\left( 
\mathbf{\mathbf{u}}_{j}\right) =\mathbf{\mathbf{u}}_{j}$, and its
homogeneity index is $1$. A composition monomial $M$ of a non-zero rank $%
r\geq 1$ has the form 
\begin{equation}
M\left( \mathbf{u}_{i_{0}}\ldots \mathbf{u}_{i_{s}}\right) =\mathcal{F}%
^{\left( s\right) }\left( M_{1}\left( \mathbf{u}_{i_{0}}\ldots \mathbf{u}%
_{i_{1}}\right) \ldots M_{s}\left( \mathbf{u}_{i_{s-1}+1}\ldots \mathbf{u}%
_{i_{s}}\right) \right) ,  \label{submon}
\end{equation}%
where $M_{1}\left( \mathbf{u}_{i_{0}}\ldots \mathbf{u}_{i_{1}}\right) $, $%
M_{2}\left( \mathbf{u}_{i_{1}+1}\ldots \mathbf{u}_{i_{2}}\right) $,..., $%
M_{s}\left( \mathbf{u}_{i_{s-1}+1}\ldots \mathbf{u}_{i_{s}}\right) $, with $%
1\leq i_{0}<i_{1}<\ldots <i_{s}$, are composition monomials of ranks not
exceeding $r-1$ (submonomials) and at least one of the rank $r-1$, the
homogeneity index of $M_{j}$ equals $i_{j}-i_{j-1}$. \ For a composition
monomial $M$ the operator $\mathcal{F}^{\left( s\right) }$ in its
representation (\ref{submon}) is called its \emph{root operator}. The index
of homogeneity of $M$ defined by (\ref{submon}) equals $i_{m}-i_{0}+1$. We
call the labeling of the arguments of a composition monomial $M$ defined by (%
\ref{submon}) \ by consecutive integers \emph{standard labeling} if $i_{0}=1$%
.
\end{definition}

If the monomials $M_{1},.,M_{s}$ have the respective \ homogeneity indexes $%
\nu \left( M_{i}\right) $ then we readily get that the homogeneity index of
the monomial $M$ satisfies the identity%
\begin{equation}
\nu \left( M\right) =\nu \left( M_{1}\right) +\ldots +\nu \left(
M_{s}\right) .  \label{sumhom}
\end{equation}%
Using the formula (\ref{submon}) inductively we find that any composition
monomial $M$ \ is given by a formula which involves symbols from the set $%
\left\{ \mathcal{F}^{\left( s\right) }\right\} _{s=2}^{\infty }$, arguments $%
\mathbf{u}_{i}$ \ and parentheses, and if $s$-linear operators are
substituted as $\mathcal{F}^{\left( s\right) }$ we obtain the terms
contained in the expansion of $\mathcal{G}^{\left( m\right) }$.

\begin{definition}[incidence number]
\label{defsubmon}The total number of symbols $\mathcal{F}^{\left( q\right) }$
involved in $M$ is called the \emph{incidence number} for $M$.
\end{definition}

For instance, the expression of the form 
\begin{equation}
M=\mathcal{F}^{\left( 4\right) }\left( \mathbf{\mathbf{u}}_{1}\mathbf{%
\mathbf{u}}_{2}\mathbf{\mathbf{u}}_{3}\mathcal{F}^{\left( 3\right) }\left( 
\mathbf{\mathbf{u}}_{4}\mathcal{F}^{\left( 2\right) }\left( \mathbf{\mathbf{u%
}}_{5}\mathbf{\mathbf{u}}_{6}\right) \mathcal{F}^{\left( 3\right) }\left( 
\mathbf{\mathbf{u}}_{7}\mathbf{\mathbf{u}}_{8}\mathbf{\mathbf{u}}_{9}\right)
\right) \right)  \label{supmon}
\end{equation}%
is an example of a composition monomial $M$ of rank 3, incidence number 4
and homogeneity index 9. It has three submonomials. Namely, the first one is 
$\mathcal{F}^{\left( 3\right) }\left( \mathbf{\mathbf{u}}_{4}\mathcal{F}%
^{\left( 2\right) }\left( \mathbf{\mathbf{u}}_{5}\mathbf{\mathbf{u}}%
_{6}\right) \mathcal{F}^{\left( 3\right) }\left( \mathbf{\mathbf{u}}_{7}%
\mathbf{\mathbf{u}}_{8}\mathbf{\mathbf{u}}_{9}\right) \right) $ of rank 2
and incidence number 3. The second submonomial $\mathcal{F}^{\left( 2\right)
}\left( \mathbf{\mathbf{u}}_{5}\mathbf{\mathbf{u}}_{6}\right) $ has rank 1
and incidence number 1, and the third one is$\ \mathcal{F}^{\left( 3\right)
}\left( \mathbf{\mathbf{u}}_{7}\mathbf{\mathbf{u}}_{8}\mathbf{\mathbf{u}}%
_{9}\right) $ of rank 1 and incidence number 1.

When analyzing the structure of composition monomials we use basic concepts
and notation from the graph theory, namely, nodes, trees and subtrees.

\begin{definition}[nodes, tree, subtree]
\label{Tree}\textbf{\ }A (finite) directed graph $T$ consists of \emph{nodes}
$N_{i}\in N_{T}$ where $N_{T}$ is the set (finite) of nodes of $T$ and a set
of edges $N_{i}N_{j}\in N_{T}\times N_{T}$. An edge $N_{i}N_{j}$ connects $%
N_{i}$ with $N_{j}$, it is an outcoming edge of $N_{i}$ and an incoming edge
of $N_{j}$. A \emph{tree} (more precisely a \emph{rooted tree}, we only
consider rooted trees) is a directed connected graph which is cycle-free and
has a selected \emph{root node}, that is a node $N_{\ast }$ which has no
incoming edges. If a \emph{node} $N$ has an outcoming edge $NN_{j}$ the node 
$N_{j}$ is called a \emph{child node} of $N$; if a node $N$ has an incoming
edge $N_{j}N$ the node $N_{j}$ is called the \emph{parent node} of $N$. We
denote the parent node of $N$ by $p\left( N\right) $. If a node does not
have children it is called an \emph{end node} (or a leaf). For every node $N$
we denote by $\mu \left( N\right) $ the \emph{number of child nodes} of the
node $N.$ If a \emph{path } connects two nodes we call the number of edges
in the path its \emph{length}. We denote by $l\left( N\right) $ the length
of a path which connects $N_{\ast }$ with $N$. Every node $N$ of the tree $T$
can be taken as a root node of a \emph{subtree} which involves all
descendent nodes of $N$ and connecting edges; we denote this maximal subtree 
$T^{\prime }\left( N\right) $. Since we consider only maximal subtrees we
simply call them subtrees. We call by the \emph{rank of a tree} the maximal
length of a path from its root node to an end node and denote it by $r\left(
T\right) $. We call by the \emph{rank of a node} $N$ of the tree $T$ the
rank of the subtree $T^{\prime }\left( N\right) $.
\end{definition}

\begin{definition}[tree incidence number and homogeneity index]
\label{Definition incidence} For a tree $T$ we call the number of non-end
nodes \emph{incidence number} $i=i\left( T\right) $. We denote the \emph{%
number of end nodes} of the tree by $\nu \left( T\right) $ and call it \emph{%
homogeneity index}.
\end{definition}

\paragraph{Elementary properties of trees.}

Since a tree does not have cycles, the path connecting two nodes on a tree
is unique. The root node $N_{\ast }$ does not have a parent node, and since
it is connected with every other node, every non-root node has a parent
node. The end nodes have zero rank. The only node with rank $r\left(
T\right) $ is the root node. The total number of nodes of a tree $T$ equals $%
m\left( T\right) +i\left( T\right) $.

\begin{definition}[ordered tree]
\label{Definition ordered}\textbf{\ }A tree is called an \emph{ordered tree}
if for every node $N$ all child nodes of $N$ are labeled by consecutive
positive integers (which may start not from $1$). Hence, for any node $%
N^{\prime }\neq N_{\ast }$ there is the parent node $N=p\left( N^{\prime
}\right) $ and the order number (label) $o\left( N^{\prime }\right) $, $%
i_{1}\leq o\left( N^{\prime }\right) \leq i_{1}+\mu \left( N\right) -1$. Two
trees are equal if there is one-to-one mapping $\Theta $ between the nodes
which preserves edges, maps the root node into the root node and preserves
the order of children of every node up to a shift: if $\Theta \left(
N\right) =\tilde{N}$ and $p\left( N_{1}\right) =p\left( N_{2}\right) =N$
then $o\left( N_{1}\right) -o\left( N_{2}\right) =o\left( \Theta \left(
N_{1}\right) \right) -o\left( \Theta \left( N_{2}\right) \right) $.
\end{definition}

Since we use in this paper only ordered trees we simply call them \emph{trees%
}.

\paragraph{Standard node labeling and ordering.}

We use the following way of labeling and ordering of end nodes of a given
ordered tree $T$. Let $\hat{r}$ be the rank of $T$. For any end node $N$ we
take the unique path $N_{\ast }N_{1}\ldots N_{l\left( N\right) -1}N$ of
length $l\left( N\right) \leq \hat{r}$ connecting it to the root. Since the
tree is ordered, every node $N_{j}$ in the path has an order number $o\left(
N_{j}\right) $. These order numbers form a word $w\left( N\right) $ of
length $l\left( N\right) $. If $l\left( N\right) <\hat{r}$ we complete $%
w\left( N\right) $ to the length $\hat{r}$ adding several symbols $\infty $
and assuming that $\infty >n$ for $n=1,2,\ldots $. After that we order words 
$w\left( N\right) $ in the lexicographic order. We obtain the ordered list $%
w_{1}\left( N_{1}\right) ,\ldots ,w_{\nu \left( T\right) }\left( N_{\nu
\left( T\right) }\right) $. We take this ordering and labeling of the end
nodes $N_{1},\ldots ,N_{\nu \left( T\right) }$ as a standard ordering and
denote by $o_{0}\left( N\right) $ the consecutive number with respect to
this labeling: $j=o_{0}\left( N_{j}\right) $. To label the nodes with rank $%
r $ we delete all the nodes of rank less than $r$ together with the incoming
edges and nodes of rank $\ r$ become end nodes. We apply to them the
described labeling and denote the indexes obtained by $o_{r}\left( N\right) $%
. Hence, every node $N$ of the tree $T$ \ has two integer numbers assigned: $%
r\left( N\right) $ and $o_{r\left( N\right) }\left( N\right) $. We introduce
the standard labeling of all nodes of $T$ by applying the lexicographic
ordering to pairs $\left( r\left( N\right) ,o_{r\left( N\right) }\left(
N\right) \right) $, and denote the corresponding number $o\left( N\right) $, 
$1\leq o\left( N\right) \leq m\left( T\right) +i\left( T\right) $.

The following statement follows straightforwardly from the definition of the
standard ordering.

\begin{proposition}
\label{Proposition orinterval} If a tree $T$ has a subtree $T^{\prime }$ and
the standard labeling of end nodes is used, then all the end nodes of the
subtree $T^{\prime }$ fill an interval $j_{1}\leq o_{0}\left( N\right) \leq
j_{2}$ for some $j_{1}$ and $j_{2}$.
\end{proposition}

\begin{theorem}
\label{Theorem montree} Let $\mathcal{T}_{2}$ be the set of ordered trees
such that each node of a tree which is not an end node has at least two
children nodes. The set of\emph{\ }composition monomials based on $\left\{ 
\mathcal{F}^{\left( s\right) },s=2,3,\ldots \right\} $ is in one-to-one
correspondence with the set $\mathcal{T}_{2}$. The correspondence has the
following properties. The monomials of rank $r$ correspond to trees of rank $%
r$. The root node of the tree $T$ corresponds to the root operator of the
composition monomial. The end nodes correspond to variables $\mathbf{\mathbf{%
u}}_{j}$, $j=1,\ldots ,\nu \left( T\right) $. The standard labeling of end
nodes coincides with the consecutive labeling of the variables $\mathbf{%
\mathbf{u}}_{j}$ of monomial from left to right.\ The homogeneity index of a
monomial equals the homogeneity index of the corresponding tree. The
incidence number of a monomial equals the incidence number of a tree, and
the rank of a monomial equals the rank of a tree.
\end{theorem}

\begin{proof}
For a given $\left\{ \mathcal{F}^{\left( s\right) }\right\} $ the set of
monomials with rank $r$ is finite, the set of trees with rank $r$ is finite
too. Therefore, to prove one-to-one correspondence of the two sets it is
sufficient to construct two one-to-one mappings from the first set into the
second and from the second into the first. First of all, using the induction
with respect to $r$ we construct for every monomial the corresponding tree.
Let $r=0$. A monomial of rank $0$ has the form $\mathbf{\mathbf{u}}_{1}$,
and it corresponds to a tree involving one node. The tree has no edges and
the node is the both the root and the end node; its incidence number is zero
and homogeneity power is one. Assume now that we have defined a tree for any
monomial of rank not greater than $r-1$. A monomial of rank $r$ has the form 
$\mathcal{F}^{\left( m\right) }\left( M_{1}\ldots M_{m}\right) $ where
monomials $M_{1}\ldots M_{m}$ \ have rank not greater than $r-1$. Every
monomial $M_{1}\ldots M_{m}$ corresponds to an ordered tree $T_{1},\ldots
,T_{m}$ with the root nodes $N_{\ast 1},\ldots N_{\ast m}$. We form the tree 
$T$ as a union of the nodes of $T_{1},\ldots ,T_{m}$ \ and add one more node 
$N_{\ast }$ which corresponds to the root operator $\mathcal{F}^{\left(
m\right) }$ and it becomes the root node of $T$. We take the union of edges
from $T_{1},\ldots ,T_{m}$ \ and add $m$ more edges connecting $N_{\ast }$
with the nodes $N_{\ast 1},\ldots N_{\ast m}$, the order of the nodes
corresponds to ordering of $M_{1}\ldots M_{m}$ from left to right. The first
mapping is constructed.

Now let us define for every ordered tree $T$ the corresponding monomial $%
M\left( \mathcal{F},T\right) $. If we have a tree $T$ of rank zero we set $%
M\left( \mathcal{F},T\right) =\mathbf{\mathbf{u}}_{j}$ and \ $j=1$ if we use
the standard labeling. Now we do induction step from $r-1$ to $r$. If we
have a tree of rank $r$ we take the root node $N_{\ast }$ and its children $%
N_{\ast 1},\ldots ,N_{\ast s}$, $s=\mu \left( N_{\ast }\right) $. The
subtrees $T^{\prime }\left( N_{\ast 1}\right) $,..., $T^{\prime }\left(
N_{\ast s}\right) $ have rank not greater than $r-1$ and the monomials $%
M\left( \mathcal{F},T^{\prime }\left( N_{\ast 1}\right) \right) $ ,..., $%
M\left( \mathcal{F},T^{\prime }\left( N_{\ast s}\right) \right) $ are
defined according to induction assumption, let $m\left( T^{\prime }\left(
N_{\ast 1}\right) \right) $,...,$m\left( T^{\prime }\left( N_{\ast s}\right)
\right) $ be their homogeneity indices. We set $m\left( T\right) =m\left(
T^{\prime }\left( N_{\ast 1}\right) \right) +$...$+m\left( T^{\prime }\left(
N_{\ast s}\right) \right) $. We denote the variables of every monomial $%
M\left( \mathcal{F},T^{\prime }\left( N_{\ast j}\right) \right) $ by $%
\mathbf{\mathbf{u}}_{j,1},\ldots ,\mathbf{\mathbf{u}}_{j,m\left( T^{\prime
}\left( N_{\ast j}\right) \right) \ }$counting from left to right, and then
labeling all the variables $\mathbf{\mathbf{u}}_{j,l}$ using the
lexicographic ordering of pairs $j,l$ we obtain variables $\mathbf{\mathbf{u}%
}_{1},\ldots ,\mathbf{\mathbf{u}}_{m\left( T\right) }$ and monomials 
\begin{equation*}
M\left( \mathcal{F},T^{\prime }\left( N_{\ast 1}\right) \right) \left( 
\mathbf{\mathbf{u}}_{1},\ldots ,\mathbf{\mathbf{u}}_{m\left( T^{\prime
}\left( N_{\ast 1}\right) \right) \ }\right) ,M\left( \mathcal{F},T^{\prime
}\left( N_{\ast 2}\right) \right) \left( \mathbf{\mathbf{u}}%
_{m_{1}+1},\ldots ,\mathbf{\mathbf{u}}_{m_{1}+m_{2}\ \ }\right) ,\ 
\end{equation*}%
etc., where $m_{j}=m\left( T^{\prime }\left( N_{\ast j}\right) \right) $.
After that we set 
\begin{gather*}
M\left( \mathcal{F},T\right) \left( \mathbf{\mathbf{u}}_{1},\ldots ,\mathbf{%
\mathbf{u}}_{m\left( T\right) }\right) = \\
\mathcal{F}^{\left( s\right) }\left( M\left( \mathcal{F},T^{\prime }\left(
N_{\ast 1}\right) \right) \left( \mathbf{\mathbf{u}}_{1},\ldots ,\mathbf{%
\mathbf{u}}_{m\left( T^{\prime }\left( N_{\ast 1}\right) \right) }\right)
,\ldots ,M\left( \mathcal{F},T^{\prime }\left( N_{\ast s}\right) \right)
\left( \mathbf{\mathbf{u}}_{m\left( T\right) -m_{s-1}+1},\ldots ,\mathbf{%
\mathbf{u}}_{m\left( T\right) }\right) \right) .
\end{gather*}%
Note that the homogeneity index for the monomial $M$ equals the sum of the
indices for submonomials $M_{1}\ldots M_{m}$, the homogeneity index for the
tree $T$ equals the sum of the indices for subtrees $T_{1},\ldots ,T_{m}$,
this implies their equality by induction. The incidence number for the
monomial $M$ equals the sum of the numbers for submonomials $M_{1}\ldots
M_{m}$ plus one; the incidence number for the tree $T$ equals the sum of the
numbers for submonomials $T_{1},\ldots ,T_{m}$ plus one. Therefore, these
quantities for monomials and trees are $\ $ equal by induction. Induction is
completed. Therefore we constructed the two mappings, one can easily check
that they are one-to-one and have all required properties.
\end{proof}

\begin{definition}[monomial to a tree]
\label{Definition monoftree} For a tree $T\in \mathcal{T}_{2}$ we denote by $%
M\left( \mathcal{F},T\right) $ the monomial which is constructed in Theorem %
\ref{Theorem montree}.
\end{definition}

\begin{conclusion}
\label{Tree monomial}\textbf{\ }The above construction shows that the
structure of every composition monomial is completely described by an
(ordered) tree $T$ with nodes $N_{i}$ corresponding to the operators$\ 
\mathcal{F}^{\left( m_{i}\right) }$. At such \ a node $N_{i}$ (i) the number 
$m_{i}$ of outcoming edges equals the homogeneity index of $\mathcal{F}%
^{\left( m_{i}\right) }$; (ii) the outcoming edges are in one-to-one
correspondence with the arguments of \ $\mathcal{F}^{\left( m_{i}\right) }$,
and the ordering of the child nodes coincides with the ordering of arguments
of $\mathcal{F}^{\left( m_{i}\right) }$ from left to right. The value of $%
m_{i}$ may be different for different nodes. A node corresponding to $%
\mathcal{F}^{\left( m\right) }$ is connected by edges with $m$ child nodes
corresponding to the arguments of $\mathcal{F}^{\left( m\right) }$. \ Every
node $N$ of the tree $T$ can be taken as a root node of a subtree $T^{\prime
}\left( N\right) $ which correspond to a submonomial $M\left( \mathcal{F}%
,T^{\prime }\left( N\right) \right) $.$\ $Conversely, every submonomial of $%
M\left( \mathcal{F},T\right) $ equals $M\left( \mathcal{F},T^{\prime }\left(
N\right) \right) $ for some mode $N$. If $m>1$ the submonomial has a nonzero
rank. The number of non-end nodes equals to the number of symbols $\mathcal{F%
}^{\left( m\right) }$ used in $\mathcal{F}$-represenation of the monomial
which is the incidence number of the monomial. The total number of end nodes
of an $m$-homogenious operator equals to $m=\nu \left( T\right) $. The rank
of a node $N$ equals the rank of the corresponding submonomial $M\left( 
\mathcal{F},T^{\prime }\left( N\right) \right) $. The arguments $\mathbf{%
\mathbf{u}}_{1},\ldots \mathbf{\mathbf{u}}_{s}$ of a monomial correspond to
the end nodes of the tree. \ The standard labeling of nodes of $T$ \ agrees
with the standard labeling (from left to right) of the$\mathcal{\ }$%
arguments of the composition monomial $M\left( \mathcal{F},T\right) $. The
number of end nodes of the tree $T$ equals the homogeneity index of
corresponding monomial.\ If the root mode of the tree $T$ of a monomial $M$
has $\mu \left( N_{\ast }\right) =m$ edges which are connected to child
nodes $N_{1},\ldots N_{m}$ \ then there is a node $\mathcal{F}^{\left(
m_{j}\right) }$, $j=1,\ldots ,n$ at the end of every edge such that $M$ has
the form%
\begin{equation}
\mathcal{F}^{\left( m\right) }\left( \mathcal{F}^{\left( \mu \left(
N_{1}\right) \right) }\left( \ldots \right) ,\ldots ,\mathcal{F}^{\left( \mu
\left( N_{m}\right) \right) }\left( \ldots \right) \right) .  \label{MG0}
\end{equation}
\end{conclusion}

\begin{example}
The tree corresponding to $\mathcal{F}^{\left( 3\right) }\left( \mathbf{%
\mathbf{u}}_{1}\mathbf{\mathbf{u}}_{2}\mathcal{F}\left( \mathbf{\mathbf{u}}%
_{1}\mathbf{\mathbf{u}}_{2}\mathbf{\mathbf{u}}_{3}\right) \right) $ has two
nodes of non-zero rank, the root node of rank 2, one non-end node of rank 1
and five end nodes of rank 0. \ Another example, the monomial (\ref{supmon})
has the root node corresponding to $\mathcal{F}^{\left( 4\right) }$, four
edges lead respectively to nodes corresponding to the end nodes with $%
\mathbf{\mathbf{u}}_{1}$, $\mathbf{\mathbf{u}}_{2}$, $\mathbf{\mathbf{u}}%
_{3}\ $and to the non-end node with $\mathcal{F}^{\left( 3\right) }$, see
Fig. 2.
\end{example}

\begin{figure}[tbph]
\scalebox{0.5}{\includegraphics[viewport=-160 60 750 600,clip]{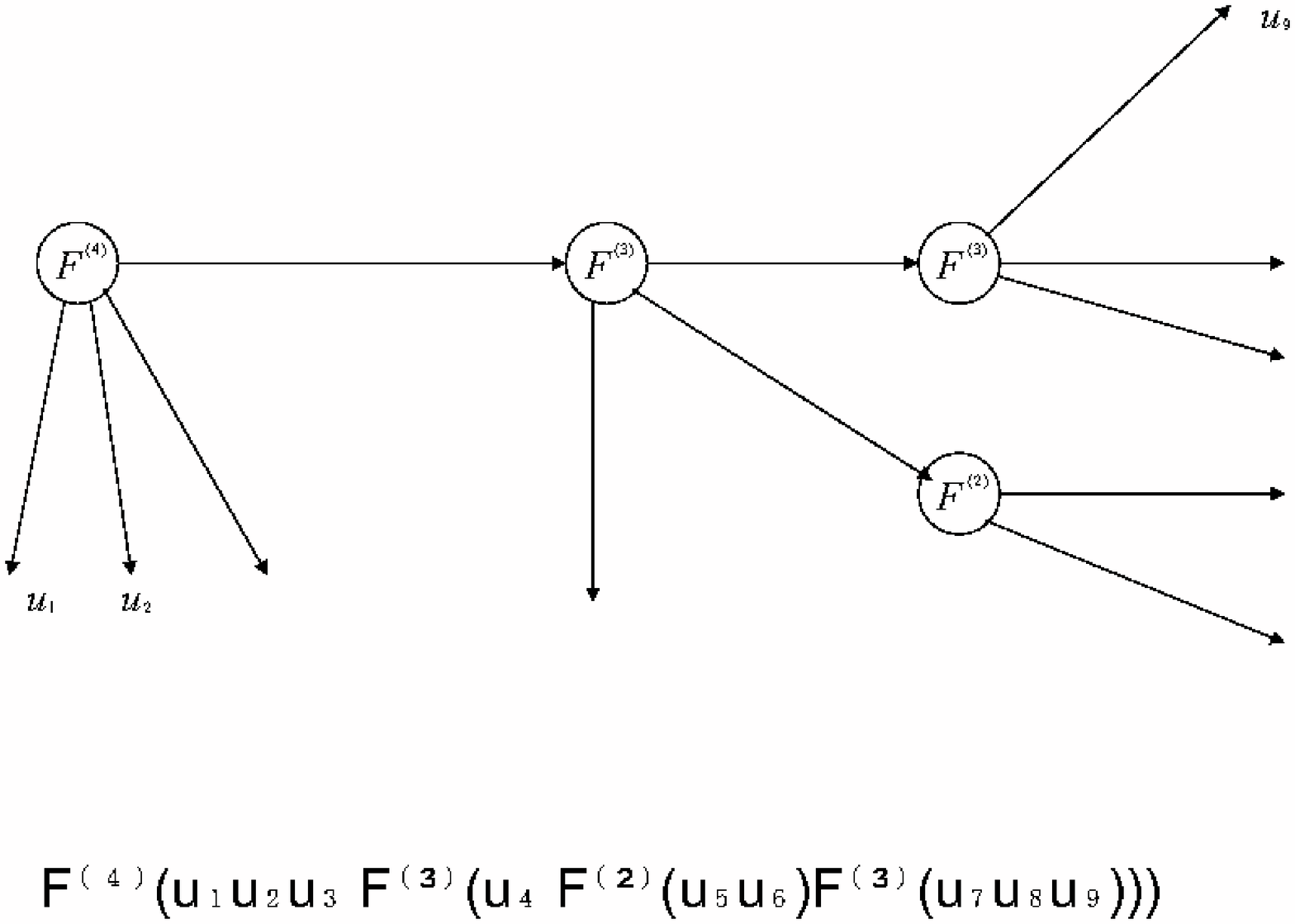}}
\caption{In this picture a tree corresponding to a monomial is drawn.}
\end{figure}

\begin{remark}
Since all operators in the set $\left\{ \mathcal{F}^{\left( s\right)
}\right\} _{s=2}^{\infty }$ in (\ref{recG0}) have the homogeneity index at
least two, the trees of monomials generated by recurrent relations (\ref%
{recG0}) have a special property: every non-end mode has at least two
children.
\end{remark}

Sometimes it is convenient to use monomials involving several types of
operators. To describe such a situation we introduce for a given tree a
decorated monomial.

\begin{definition}[decorated monomial of a tree]
\label{Generalized monomial of a tree}\textbf{\ }Assume that we have several
formal series \ $\left\{ \mathcal{F}_{1},\ldots ,\mathcal{F}_{l}\right\} $
where $\mathcal{F}_{i}$ is represented by a formal series$\ \mathcal{F}%
_{l}=\sum_{m}\mathcal{F}_{i}^{\left( m\right) }$, $i=1,\ldots ,l$. \ We call
the set $\left\{ \mathcal{F}\right\} =\left\{ \mathcal{F}_{j},\ j=1,\ldots
,S\right\} $\ the \emph{operator alphabet}, and $j$ is called the \emph{%
decoration index}. We consider a function $\Gamma \left( N\right) $, $N\in T$%
, defined on the nodes of the tree $T$ and taking values in the set $\left\{
1,\ldots ,l\right\} $ of the decoration indices, and call such a function a 
\emph{decoration function} on the tree $T$. Then for a decoration function $%
\Gamma \left( N\right) $ we define the \emph{decorated monomial }$M\left(
\left\{ \mathcal{F}\right\} ,\Gamma ,T\right) $ of the tree\ $T$ by picking
operators $\mathcal{F}_{j}^{\left( m\right) }$ with $j$ defined by $\Gamma $%
. For every node $N$ the homogeneity index $m=\mu \left( N\right) $ of the
operator $\mathcal{F}_{j}^{\left( m\right) }$ equals to the number of
children of $N$ and $j$ is defined by $\Gamma $, namely $\mathcal{F}_{j}$, $%
j=\Gamma \left( N\right) $.
\end{definition}

Hence, a decorated monomial $M\left( \left\{ \mathcal{F}\right\} ,\Gamma
,T\right) $ has instead of (\ref{MG0}) the following form 
\begin{equation}
\mathcal{F}_{\Gamma \left( N\right) }^{\left( m\right) }\left( \mathcal{F}%
_{\Gamma \left( N_{1}\right) }^{\left( \mu \left( N_{1}\right) \right)
}\left( \ldots \right) ,\ldots ,\mathcal{F}_{\Gamma \left( N_{m}\right)
}^{\left( \mu \left( N_{m}\right) \right) }\left( \ldots \right) \right) .
\label{MG}
\end{equation}%
When $\mathcal{F}_{i}^{\left( m\right) }$ are multilinear operators, a
monomial $M\left( \left\{ \mathcal{F}\right\} ,T,\Gamma \right) $ is also a
multilinear operator, its homogeneity index $m$ equals $\nu \left( T\right) $
and we denote its arguments by $\left( \mathbf{x}_{1}\ldots \mathbf{x}%
_{m}\right) $. Respectively, if \ $\mathbf{x}_{1}\ldots \mathbf{x}_{\nu }$
are arguments of a monomial $M\left( \left\{ \mathcal{F}\right\} ,T,\Gamma
\right) $ and we use the standard labeling of the nodes then according to
Proposition \ref{Proposition orinterval} a submonomial $M\left( \left\{ 
\mathcal{F}\right\} ,T,\Gamma \right) $ has arguments $\mathbf{x}_{\varkappa
\left( T^{\prime }\right) },\ldots ,\mathbf{x}_{\varkappa \left( T^{\prime
}\right) +\nu \left( T^{\prime }\right) -1}$ which are labeled
constructively.

Now we would like to describe elementary properties of composition monomials
and the related trees. Note that for every $N\in T$ a composition monomial
is a linear function of operator $\mathcal{F}_{\Gamma \left( N\right) }^{\mu
\left( N\right) }$. Consequently, the concept of the decorated composition
monomial can be naturally extended to monomials associated with the
following family of operators 
\begin{equation*}
\left\{ \mathcal{F}\right\} =\left\{ \mathcal{F}:\mathcal{F}=c_{1}\mathcal{F}%
_{1}+\ldots +c_{l}\mathcal{F}_{l},\ c_{i}\in \mathbb{C}\right\} .
\end{equation*}%
For a given tree $T$ the submonomial $M\left( \left\{ \mathcal{F}\right\}
,\Gamma ,T\right) $ is represented as a function on the tree $T$ with values
in $\left\{ \mathcal{F}\right\} $, this is an $i$-linear function of $%
\mathcal{F}$ where $i$ is the incidence number of $T$.

There are elementary relations between the incidence number $i\left(
T\right) $, the rank $r\left( T\right) $, the number of edges of a tree $T$
\ which do not end at an end node $e^{0}\left( T\right) $ and the
homogeneity index $m$ of a tree $T$, and corresponding monomial $M\left(
\left\{ \mathcal{F}\right\} ,\Gamma ,T\right) $. For example, $e^{0}\left(
T\right) =i\left( T\right) -1$. Some useful relations expressed by
inequalities are given in the following lemma.

\begin{lemma}
\label{Lemma nodes} Let us consider trees $T$ for which every non-end node
has at least two children, $\mu \left( N\right) \geq 2$ for all $N\in T$.
Let for any $i$ the number $m\left( i\right) $ be the minimum number of the
end nodes $\nu \left( T\right) $ for all trees $T$ with given incidence
number $i$. Then 
\begin{equation}
m\left( i\right) \geq i+1.  \label{mofs}
\end{equation}%
Similarly for any given $r$ let\ $m\left( r\right) $ be the minimum number
of end nodes with given rank $r$. Then%
\begin{equation}
m\left( r\right) \geq r+1.  \label{mofr}
\end{equation}%
Let $e^{0}\left( T\right) $ be the number of edges of a tree $T$ which do
not end at end nodes. For any given $e$ let $m\left( e\right) $ be the
minimum number of end nodes with $e^{0}\left( T\right) =e$. Then%
\begin{equation}
m\left( e^{0}\right) >e^{0}+1.  \label{mofe}
\end{equation}%
\ 
\end{lemma}

\begin{proof}
For $i=1$ (\ref{mofs}) is true. Let the statement be true for $i=i_{0}$. Let 
$T$ be a tree\ with the minimum number of end nodes $m\left( i_{0}\right) =m$%
. We delete one of the end nodes together with the edge leading to it from
its parent obtaining a tree with $m\left( i_{0}\right) -1$ end node. If the
tree remains in the same class, then $m\left( i_{0}\right) $ is reduced by
one contradicting the minimality. Hence, the deletion of the edge created a
node with only one child. Such a node can be replaced by an edge leading
from its parent to its child and reducing the incidence number by one. Using
the induction assumption we get 
\begin{equation}
m\left( i_{0}\right) -1\geq m\left( i_{0}-1\right) \geq \left(
i_{0}-1\right) +1  \label{ms}
\end{equation}%
that completes the induction and proves (\ref{mofs}) for all $i$. Similar
induction proves (\ref{mofr}). For $r=1$ (\ref{mofr}) is true. Let $T$ be a
tree\ with the minimum number of end nodes $m\left( r_{0}\right) =m$. As
above, by deleting an end node and using the minimality we reduce the tree $%
T $ to a tree $T^{\prime }$ with a smaller rank. Since only one non-end node
is eliminated, the rank of $T^{\prime }$ is $r_{0}-1$ and we get (\ref{mofr}%
). Inequality (\ref{mofe}) holds for $e=0$ since $m\left( 0\right) \geq 2$.
Let $T$ be a tree\ with the \ minimum number of end nodes $m\left(
e_{0}\right) =m$. We again delete one of the end nodes together with the
edge joining it to its parent and obtain a tree with $m\left( e_{0}\right)
-1 $ end nodes and the same number of edges which do not end at an end node.
The minimality implies that the parent node has only one another child and
removing it we get either $e_{0}$ or $e_{0}-1$ edges which do not go to end
nodes. We use the induction as in (\ref{ms}) obtaining (\ref{mofe}).
\end{proof}

\paragraph{Monomial expansion in the Implicit Function Theorem}

If operators $\mathcal{G}^{m}\left( \mathbf{x}_{1}\ldots \mathbf{x}%
_{m}\right) $ are determined by the recurrent formulas (\ref{recG0}) it is
obvious that every $\mathcal{G}^{m}$ can be represented in terms of $%
\mathcal{F}=\left\{ \mathcal{F}^{\left( s\right) }\right\} $ using the
recurrence and multilinearity of $\mathcal{F}^{\left( s\right) }$. More
precisely the following representation holds%
\begin{equation}
\mathcal{G}^{\left( m\right) }\left( \mathcal{F},\mathbf{x}_{1}\ldots 
\mathbf{x}_{m}\right) =\sum_{T\in T_{m}}c_{T}M\left( \mathcal{F},T\right)
\left( \mathbf{x}_{1}\ldots \mathbf{x}_{m}\right) ,  \label{treesum}
\end{equation}%
where (i) $M\left( \mathcal{F},T\right) \in \mathcal{T}_{2}$ is a
composition monomial corresponding to a tree $T$ and $T_{m}\subset \mathcal{T%
}_{2}$ stands for the set of trees with $m$ end nodes; (ii) the
integer-valued \emph{multiplicity coefficient} $c_{T}\geq 0$ counts the
multiplicity of the related monomial $M\left( \mathcal{F},T\right) $ in the
expansion of (\ref{recG0}); for some trees $T$ its multiplicity coefficient $%
c_{T}$ may be zero. \ The expansion (\ref{treesum}) is obtained by an
inductive process with respect to $m$ since (\ref{recG0}) expresses $%
\mathcal{G}^{m}$ in terms of $\mathcal{G}^{\left( i_{j}\right) }$ with $%
2\leq i_{j}<m$. Notice that for a given operator $\mathcal{F}=\left\{ 
\mathcal{F}^{\left( s\right) }\right\} $ the monomial $M\left( \mathcal{F}%
,T\right) $ considered as an operator can be the same for different $T$, the
monomials and the multiplicity coefficients are determined purely
algebraically.

\begin{remark}
The expression (\ref{treesum}) for $\mathcal{G}^{\left( m\right) }$ as a
linear combination of composition monomials $M\left( \mathcal{F},T\right) $,
in particular the multiplicity coefficients $c_{T}$, \ does not depend on a
specific form of the operator $\ \mathcal{F}$. It is the same for a solution 
$\mathbf{z}=\mathbf{x}+\mathcal{G}\left( \mathcal{F},\mathbf{x}\right) $ of
the general functional equation (\ref{eqF}) and for an elementary algebraic
equation $u=\mathcal{F}\left( u\right) +x$ with $u,x\in \mathbb{C}$ and with
a scalar analytic function $\mathcal{F}$ \ of one complex variable.
\end{remark}

If all $\mathcal{F}_{i}^{\left( m\right) }$ are bounded multilinear
operators then a decorated monomial $M\left( \mathcal{F},T,\Gamma \right) $
is also a bounded multilinear operator as it follows from the following
statement.

\begin{lemma}
\label{Lemma Normmon} Let $M\left( \left\{ \mathcal{F}\right\} ,T,\Gamma
\right) $ be a decorated monomial of the homogeneity index $\nu \left(
T\right) =m$ and all $\mathcal{F}_{i}^{\left( s\right) }$ be bounded
operators from $E^{s}$ into $E$ for a Banach space $E$. Then the following
estimate holds 
\begin{equation}
\left\Vert M\left( \left\{ \mathcal{F}\right\} ,T,\Gamma \right) \left( 
\mathbf{x}_{1}\ldots \mathbf{x}_{m}\right) \right\Vert _{E}\leq
\dprod\limits_{N\in T\ ,r\left( N\right) >0}\left\Vert \mathcal{F}_{\Gamma
\left( N\right) }^{\left( \mu \left( N\right) \right) }\right\Vert
\dprod\limits_{j=1\ }^{m}\left\Vert \mathbf{x}_{j}\right\Vert _{E}.
\label{normsimple}
\end{equation}
\end{lemma}

\begin{proof}
Notice that 
\begin{equation}
\left\Vert \mathcal{F}^{\left( m\right) }\left( M_{1}\ldots M_{m}\right)
\right\Vert _{E}\leq \left\Vert \mathcal{F}^{\left( m\right) }\right\Vert
\left\Vert M_{1}\right\Vert _{E}\ldots \left\Vert M_{m}\right\Vert _{E}
\label{Fmnorm}
\end{equation}%
where $M_{j}$ are submonomials. Applying the above inequality repeatedly we
obtain (\ref{normsimple}).
\end{proof}

The next statement provides a bound for the norm of a decorated monomial
which involves as a factor the norm of a submonomial.

\begin{lemma}
\label{Norm submonomial}. Let $M\left( \left\{ \mathcal{F}\right\} ,T,\Gamma
\right) $ be a decorated monomial evaluated at $\mathbf{x}_{1}\ldots \mathbf{%
x}_{m}$. Let all $\mathcal{F}^{\left( s\right) }$ be bounded operators from $%
E^{s}$ \ into Banach space $E$. \ Then \ for every evaluated submonomial $%
M\left( \left\{ \mathcal{F}\right\} ,T^{\prime }\left( N_{0}\right) ,\Gamma
\right) $ we have an estimate 
\begin{gather}
\left\Vert M\left( \left\{ \mathcal{F}\right\} ,T,\Gamma \right) \left( 
\mathbf{x}_{1}\ldots \mathbf{x}_{m}\right) \right\Vert _{E}\leq \left\Vert
M\left( \left\{ \mathcal{F}\right\} ,T^{\prime }\left( N_{0}\right) ,\Gamma
\right) \left( \mathbf{x}_{\varkappa },\ldots ,\mathbf{x}_{\varkappa +\nu
\left( T^{\prime }\left( N\right) \right) -1}\right) \right\Vert _{E}
\label{MFr} \\
\dprod\limits_{N\in T\setminus T^{\prime }\left( N_{0}\right) ,r\left(
N\right) >0}\left\Vert \mathcal{F}_{\Gamma \left( N\right) }^{\left( \mu
\left( N\right) \right) }\right\Vert \dprod\limits_{j<\varkappa }\left\Vert 
\mathbf{x}_{j}\right\Vert \dprod\limits_{j\geq \varkappa +\nu \left(
T^{\prime }\left( N_{0}\right) \right) }\left\Vert \mathbf{x}_{j}\right\Vert
.  \notag
\end{gather}%
where $\mathbf{x}_{\varkappa \ },\ldots ,\mathbf{x}_{\varkappa \ +\nu \left(
T^{\prime }\left( N\right) \right) -1}$ are the arguments of the submonomial 
$M\left( \left\{ \mathcal{F}\right\} ,T^{\prime }\left( N_{0}\right) ,\Gamma
\right) $.
\end{lemma}

\begin{proof}
The proof uses the induction with respect to the length $l\left(
N_{0}\right) $. For $l\left( N_{0}\right) =0$ the statement is obvious.
Assuming that the statement is true for $l\left( N\right) <l_{0}$ \ we
consider the case when $l\left( N_{0}\right) =l_{0}$. Notice that 
\begin{equation*}
\left\Vert \mathcal{F}_{\Gamma \left( N_{\ast }\right) }^{\left( \mu \left(
N_{\ast }\right) \right) }\left( M_{1}\ldots M_{\mu \left( N\right) }\right)
\right\Vert _{E}\leq \left\Vert \mathcal{F}_{\Gamma \left( N_{\ast }\right)
}^{\left( \mu \left( N_{\ast }\right) \right) }\right\Vert \left\Vert
M_{1}\right\Vert _{E}\ldots \left\Vert M_{\mu \left( N\right) }\right\Vert
_{E},
\end{equation*}%
where $M_{j}=M\left( \left\{ \mathcal{F}\right\} ,T^{\prime }\left( N_{\ast
j}\right) ,\Gamma \right) $, $N_{\ast j}$ are child nodes of $N_{\ast }$.
One of the submonomials $M_{1}\ldots M_{\mu \left( N\right) }$ contains $%
M\left( \left\{ \mathcal{F}\right\} ,T^{\prime }\left( N_{0}\right) ,\Gamma
\right) $ as a submonomial, and let it be $M\left( \left\{ \mathcal{F}%
\right\} ,T^{\prime }\left( N_{\ast j_{0}}\right) ,\Gamma \right) $. The
length of the path from $N_{0}$ to $N_{\ast j}$ is less than $l_{0}$ and we
can use the induction hypothesis to estimate the norm of $M\left( \left\{ 
\mathcal{F}\right\} ,T^{\prime }\left( N_{\ast j_{0}}\right) ,\Gamma \right) 
$. The norms of $M_{j}$ with $j\neq j_{0}$ are estimated using (\ref%
{normsimple}). The labels of the arguments of the submonomial fill an
interval according to Proposition \ref{Proposition orinterval}.
\end{proof}

The following theorem gives a needed refinement of the Implicit Function
Theorem \ref{Imfth}.

\begin{theorem}[refined Implicit Function Theorem]
\label{Theorem Monomial Convergence} Let $\mathcal{F}\in A\left( C_{\mathcal{%
F}},R_{\mathcal{F}}\right) $. Let $\ \mathcal{G}\in A\left( C_{\mathcal{G}%
},R_{\mathcal{G}}\right) $ be the analytic solution operator constructed in
Theorem \ref{Imfth} which solves (\ref{eqF}). Then the expansion of $%
\mathcal{G}\left( \mathcal{F},\mathbf{\mathbf{x}}\right) $ into composition
monomials%
\begin{equation}
\mathcal{G}\left( \mathcal{F},\mathbf{\mathbf{x}}\right) =\sum_{m=1}^{\infty
}\sum_{T\in T_{m}}c_{T}M\left( \mathcal{F},T\right) \left( \mathbf{\mathbf{x}%
}^{m}\right)  \label{Gtreeexp}
\end{equation}%
converges for $\left\Vert \mathbf{\mathbf{x}}\right\Vert <R_{\mathcal{G}}$,
and the following estimates hold%
\begin{equation}
\sum_{T\in T_{m}}c_{T}\left\Vert M\left( \mathcal{F},T\right) \left( \mathbf{%
\mathbf{x}}^{m}\right) \right\Vert \leq C_{\mathcal{G}}R_{\mathcal{G}%
}^{-m}\left\Vert \mathbf{\mathbf{x}}\right\Vert ^{m},\;m=2,\ldots ,
\label{cT}
\end{equation}%
\begin{equation*}
\sum_{m=2}^{\infty }\sum_{T\in T_{m}}c_{T}\left\Vert M\left( \mathcal{F}%
,T\right) \left( \mathbf{\mathbf{x}}^{m}\right) \right\Vert \leq C_{\mathcal{%
G}}\frac{\left\Vert x\right\Vert _{X}^{2}R_{\mathcal{G}}^{-2}}{1-\left\Vert
x\right\Vert _{X}R_{\mathcal{G}}^{-1}},
\end{equation*}%
where $C_{\mathcal{G}}$ and $R_{\mathcal{G}}$ depend only on $C_{\mathcal{F}%
} $ and $R_{\mathcal{F}}$ and satisfy 
\begin{equation*}
C_{\mathcal{G}}=\frac{R_{\mathcal{F}}^{2}}{2\left( C_{\mathcal{F}}+R_{%
\mathcal{F}}\right) },\ R_{\mathcal{G}}=\frac{R_{\mathcal{F}}^{2}}{4\left(
C_{\mathcal{F}}+R_{\mathcal{F}}\right) }.
\end{equation*}%
The multiplicity coefficients $c_{T}\geq 0$ satisfy the inequality%
\begin{equation}
\sum_{T\in T_{m}}c_{T}\leq \frac{1}{4}8^{m}.  \label{ctes}
\end{equation}
\end{theorem}

The proof of this statement is given in Appendix B.

\subsection{Decorated expansions}

In this section we develop a formalism for treating linear operators with
several invariant subspaces which span the entire space as, for example, in
the case of projections (\ref{Pin}). The decomposition into related
invariant subspaces is very important for the analysis. The general setting
is as follows. Suppose that a Banach space $E$ has several projection
operators $\Pi _{\lambda \ }$, $\lambda \in \Lambda $, where $\Lambda $ is a
finite set of indices, we call this set \emph{decoration set}. We assume
that the sum of the projections equals the identical operator, i.e. 
\begin{equation}
\sum_{\lambda \in \Lambda }\Pi _{\lambda }=\text{Id, where Id is the
identity operator,}  \label{Ps}
\end{equation}%
and 
\begin{equation}
\Pi _{\lambda }\Pi _{\lambda }=\Pi _{\lambda },\ \Pi _{\lambda ^{\prime
}}\Pi _{\lambda }=0\text{ if }\lambda ^{\prime }\neq \lambda ,\ \lambda
^{\prime },\ \lambda \in \Lambda .  \label{PP}
\end{equation}%
We call such projections \emph{decoration projections}. For example, let us
look at projections $\Pi _{n,\zeta }\left( \mathbf{\mathbf{k}}\right) $, $%
n=1,\ldots ,J$, $\zeta =\pm $ defined by (\ref{Pin}). These projections
define bounded operators $\Pi _{n,\zeta }$ acting on (i) functions of $%
\mathbf{\mathbf{k}}$ in the space $L_{1}$; (ii) functions of $\mathbf{%
\mathbf{k}},\tau $ in the space $E=C\left( \left[ 0,\tau _{\ast }\right]
,L_{1}\right) $. In another example based on (\ref{Pin}) we fix $n_{0}$ and
define 
\begin{equation}
\Pi _{\zeta }\left( \mathbf{\mathbf{k}}\right) =\Pi _{n_{0},\zeta }\left( 
\mathbf{\mathbf{k}}\right) ,\ \zeta =\pm ,\ \Pi _{\infty }\left( \mathbf{%
\mathbf{k}}\right) =\sum_{n\neq n_{0},\zeta =\pm }\Pi _{n,\zeta }\left( 
\mathbf{\mathbf{k}}\right) .  \label{Pinfk}
\end{equation}%
Using (\ref{Ps}) we expand vectors $\mathbf{x}\in E$ as follows 
\begin{equation}
\mathbf{x}=\sum_{\lambda \in \Lambda }\Pi _{\lambda }\mathbf{x}%
=\sum_{\lambda \in \Lambda }\mathbf{x}_{\lambda },\ \mathbf{x}_{\lambda
}=\Pi _{\lambda }\left( \mathbf{x}\right) .  \label{xP}
\end{equation}%
We also use notation 
\begin{equation}
\mathcal{F}_{\lambda }^{\left( n\right) }=\Pi _{\lambda }\mathcal{F}^{\left(
n\right) }  \label{Flam}
\end{equation}%
Often in applications the number of elements in $\Lambda $ is either $2$ or $%
3$. In the case when $\Lambda $ has three elements we set%
\begin{equation}
\Lambda =\left\{ +,-,\infty \right\} ,\ \Pi _{+}+\Pi _{-}+\Pi _{\infty }=%
\text{Id,}  \label{PPPI}
\end{equation}%
and 
\begin{equation}
\mathbf{x=x}_{+}+\mathbf{x}_{-}+\mathbf{x}_{\infty },\ \mathcal{F}\left( 
\mathbf{x}\right) =\mathcal{F}_{+}\left( \mathbf{x}\right) +\mathcal{F}%
_{-}\left( \mathbf{x}\right) +\mathcal{F}_{\infty }\left( \mathbf{x}\right) .
\label{xFinf}
\end{equation}%
Using the decomposition (\ref{Ps}) we introduce for $m$-linear operators $%
\mathcal{F}^{\left( n\right) }\left( \mathbf{x}_{1}\ldots \mathbf{x}%
_{n}\right) $ the corresponding \emph{decorated operators} $\mathcal{F}%
_{\lambda ,\vec{\zeta}}^{\left( n\right) }$ as follows: 
\begin{equation}
\mathcal{F}_{\lambda ,\vec{\zeta}}^{\left( n\right) }\left( \mathbf{x}%
_{1}\ldots \mathbf{x}_{n}\right) =\Pi _{\lambda }\mathcal{F}^{\left(
n\right) }\left( \Pi _{\zeta ^{\prime }}\mathbf{x}_{1}\ldots \Pi _{\zeta
^{\left( n\right) }}\mathbf{x}_{n}\right) =\mathcal{F}_{\lambda }^{\left(
n\right) }\left( \Pi _{\zeta ^{\prime }}\mathbf{x}_{1}\ldots \Pi _{\zeta
^{\left( n\right) }}\mathbf{x}_{n}\right) ,  \label{Flamz}
\end{equation}%
where $\vec{\zeta}$ is defined in (\ref{zetaar}). Obviously, we have 
\begin{equation}
\mathcal{F}^{\left( n\right) }\left( \mathbf{x}_{1}\ldots \mathbf{x}%
_{n}\right) =\sum_{\lambda \in \Lambda ,\ \vec{\zeta}\in \Lambda ^{n}}%
\mathcal{F}_{\lambda ,\vec{\zeta}}^{\left( n\right) }\left( \mathbf{x}%
_{1}\ldots \mathbf{x}_{n}\right) .  \label{Fsumlamz}
\end{equation}%
An example of expansion (\ref{Fsumlamz}) is given by (\ref{Tplusmin}).

\subsection{Decorated composition monomials}

We assume that operators $\mathcal{F}^{\left( n\right) }$ act in the space
allowing a decomposition into three components as in (\ref{PPPI}). Let $%
M\left( \mathcal{F},T\right) $ be a composition monomial of the homogeneity
index $m$, and assume that the corresponding tree $T$ has the incidence
number $i$, the rank $r$, and $e$ edges. Suppose also that every operator $%
\mathcal{F}^{\left( n\right) }$ is expanded into a sum of decorated
operators as in (\ref{Fsumlamz}) , (\ref{Flamz}). Using the linearity of $%
M\left( \mathcal{F},T\right) $ with respect to operators $\mathcal{F}%
^{\left( n\right) }$ we get 
\begin{gather}
M\left( \mathcal{F},T\right) =\mathcal{F}^{\left( n\right) }\left( \mathcal{F%
}^{\left( m_{1}\right) }\left( \ldots \right) \ldots \mathcal{F}^{\left(
m_{n}\right) }\left( \ldots \right) \right) =  \label{FFlam} \\
\sum_{\lambda \in \Lambda ,\ \vec{\lambda}\in \Lambda ^{i-1},\ \vec{\zeta}%
_{j},\ j=1,\ldots ,e}\mathcal{F}_{\lambda }^{\left( n\right) }\left( 
\mathcal{F}_{\lambda _{j_{1}},\vec{\zeta}_{j_{1}}}^{\left( m_{1}\right)
}\left( \ldots \right) \ldots \mathcal{F}_{\lambda _{j_{n}},\vec{\zeta}%
_{j_{n}}}^{\left( m_{n}\right) }\left( \ldots \right) \right) \text{,} 
\notag
\end{gather}%
where submonomials $\mathcal{F}_{\lambda _{1},\vec{\zeta}_{1}}^{\left(
m_{1}\right) }\left( \ldots \right) $,...,$\mathcal{F}_{\lambda _{n},\vec{%
\zeta}_{n}}^{\left( m_{n}\right) }\left( \ldots \right) $ have ranks not
exceeding $r-1$. We expanded repeatedly the expression in the left-hand side
of (\ref{FFlam}) as long as submonomials of non-zero rank were present
resulting in \emph{an expansion involving only decorated operators} $%
\mathcal{F}_{\lambda ,\vec{\zeta}}^{\left( n\right) }$.

\begin{remark}
\label{decproj}Note that 
\begin{equation}
\mathcal{F}_{\lambda }^{\left( n\right) }\left( \mathcal{F}_{\lambda _{1},%
\vec{\zeta}_{1}}^{\left( m_{1}\right) }\left( \ldots \right) \ldots \mathcal{%
F}_{\lambda _{n},\vec{\zeta}_{n}}^{\left( m_{n}\right) }\left( \ldots
\right) \right) =\mathcal{F}_{\lambda }^{\left( n\right) }\left( \Pi
_{\lambda _{1}}\mathcal{F}_{\lambda _{1},\vec{\zeta}_{1}}^{\left(
m_{1}\right) }\left( \ldots \right) \ldots \Pi _{\lambda _{n}}\mathcal{F}%
_{\lambda _{n},\vec{\zeta}_{n}}^{\left( m_{n}\right) }\left( \ldots \right)
\right)  \label{FPlam}
\end{equation}%
Since projections $\Pi _{\zeta }$ satisfy the identities (\ref{PP}) if a
vector $\vec{\zeta}=\left( \zeta ^{\prime },\ldots ,\zeta ^{\left( n\right)
}\right) $ and indices $\lambda _{1},\ldots ,\lambda _{n}$ are given, then
we have the identity 
\begin{equation}
\mathcal{F}_{\lambda ,\vec{\zeta}}^{\left( n\right) }\left( \mathcal{F}%
_{\lambda _{1},\vec{\zeta}_{1}}^{\left( m_{1}\right) }\ldots \mathcal{F}%
_{\lambda _{n},\vec{\zeta}_{n}}^{\left( m_{n}\right) }\right) =0\text{ \
when \ }\lambda _{i}\neq \zeta ^{\left( i\right) }\mathcal{\ }\text{\ for
some }i.  \label{FFlam0}
\end{equation}%
Hence, for \ non-zero terms in the expansion (\ref{FFlam}) if indices $%
\lambda _{1},\ldots ,\lambda _{n}$ for $\mathcal{F}_{\lambda _{1},\vec{\zeta}%
_{1}}^{\left( m_{1}\right) }$,...,$\mathcal{F}_{\lambda _{n},\vec{\zeta}%
_{n}}^{\left( m_{n}\right) }$ are given the vector $\vec{\zeta}$ in $%
\mathcal{F}_{\lambda ,\vec{\zeta}}^{\left( n\right) }$\ is determined by them%
\begin{equation}
\zeta ^{\left( i\right) }=\lambda _{i},\;i=1,\ldots ,n.  \label{zlam}
\end{equation}
\end{remark}

Note that according to (\ref{FFlam0}) \ and (\ref{zlam}) we have 
\begin{equation}
\mathcal{F}_{\lambda ,\vec{\lambda}}^{\left( n\right) }\left( \mathcal{F}%
^{\left( m_{1}\right) }\left( \ldots \right) \ldots \mathcal{F}^{\left(
m_{n}\right) }\left( \ldots \right) \right) =\mathcal{F}_{\lambda }^{\left(
n\right) }\left( \mathcal{F}_{\lambda _{1}}^{\left( m_{1}\right) }\left(
\ldots \right) \ldots \mathcal{F}_{\lambda _{n}}^{\left( m_{n}\right)
}\left( \ldots \right) \right) .  \label{FFlam1}
\end{equation}

According to (\ref{FFlam}), (\ref{FFlam1}) for every tree $T$ of the
homogeneity index $m$ and the incidence number $i$, we get an expansion into
a sum of monomials of the form%
\begin{gather}
M\left( \mathcal{F},T,\vec{\lambda},\vec{\zeta}\right) \left( \mathbf{x}_{1}%
\mathbf{x}_{2}\ldots \mathbf{x}_{m}\right) =M\left( \left\{ \mathcal{F}%
\right\} ,\Gamma ,T\right) \left( \mathbf{x}_{1}\mathbf{x}_{2}\ldots \mathbf{%
x}_{m}\right) ,  \label{MFsplit1} \\
\left\{ \mathcal{F}\right\} =\left\{ \mathcal{F}_{\lambda ,\vec{\zeta}%
}^{\left( n\right) }:\lambda \in \Lambda ,\ \vec{\zeta}\in \Lambda ^{n},\
n=2,3,\ldots \right\} .  \notag
\end{gather}%
Namely, if a monomial$\ M\left( \mathcal{F},T\right) $ has at a node $N$
operator $\mathcal{F}^{\left( m\left( N\right) \right) }$ then $M\left(
\left\{ \mathcal{F}\right\} ,\Gamma ,T\right) $ at this node has operator $%
\mathcal{F}_{\Gamma \left( N\right) }^{\left( m\left( N\right) \right) }$.
We call a composition monomial of the form (\ref{MFsplit1}), where (\ref%
{zlam}) is assumed, a \emph{decorated composition monomial}. Using the
standard labeling of nodes, for a given function $\Gamma $ on the tree $T$\
with values in $\Lambda $ we find the vectors $\vec{\lambda}\in \Lambda ^{i}$%
, $\vec{\zeta}\in \Lambda ^{m}$, with $\ i\ $ being the incidence number of
the tree $T$, and using (\ref{zlam}) we rewrite (\ref{FFlam}) in the form 
\begin{equation}
M\left( \mathcal{F},T\right) \left( \mathbf{x}_{1}\mathbf{x}_{2}\ldots 
\mathbf{x}_{m}\right) =\sum_{\vec{\lambda}\in \Lambda ^{i},\ \vec{\zeta}\in
\Lambda ^{m}}M\left( \mathcal{F},T,\vec{\lambda},\vec{\zeta}\right) \left( 
\mathbf{x}_{1}\mathbf{x}_{2}\ldots \mathbf{x}_{m}\right) .  \label{MFsplit}
\end{equation}%
where $\vec{\zeta}$ is determined by values of $\Gamma $ on the end nodes.
The sum (\ref{MFsplit}) contains at most $3^{i+m}$ non-zero terms, where $3$
is the number of elements in $\Lambda $. Combining (\ref{MFsplit}) with (\ref%
{Gtreeexp}) we obtain 
\begin{equation}
\mathcal{G}^{\left( m\right) }\left( \mathbf{x}^{m}\right) =\sum_{T\in
T_{m}}\sum_{\vec{\lambda}\in \Lambda ^{i\left( T\right) },\ \vec{\zeta}\in
\Lambda ^{m}}c_{T}M\left( \mathcal{F},T,\vec{\lambda},\vec{\zeta}\right)
\left( \mathbf{x}^{m}\right) .  \label{Gtreesum}
\end{equation}

\section{Expansions of solutions for oscillatory integral equation}

In this section we apply general concepts introduced in previous sections to
oscillatory integrals involving operators $\mathcal{F}$ as in (\ref{varcu}),
(\ref{Fu}). Based on projections $\Pi _{n,\zeta }\left( \mathbf{\mathbf{k}}%
\right) $ in (\ref{Pin}) for given $n=n_{0}$ we define as in (\ref{Pinfk})
decoration projections in $L_{1}$ which satisfy (\ref{PPPI}): 
\begin{equation}
\Pi _{\zeta }\mathbf{\tilde{u}}\left( \mathbf{\mathbf{k}}\right) =\Pi
_{n_{0},\zeta }\left( \mathbf{\mathbf{k}}\right) \mathbf{\tilde{u}}\left( 
\mathbf{\mathbf{k}}\right) ,\ \zeta =\pm ,\ \Pi _{\infty }=\sum_{n\neq
n_{0}}\sum_{\zeta =\pm }\Pi _{n,\zeta }.  \label{PPosc}
\end{equation}

\subsection{Boundedness of oscillatory integral operators}

In this subsection we estimate norms of multilinear operators $\mathcal{F}=%
\mathcal{F}^{\left( m\right) }$ defined by (\ref{Fu}) and the related
composition monomials. The operators $\mathcal{F}^{\left( m\right) }$ have
the form (\ref{Fu}) where $\mathbb{D}_{m}=\mathbb{R}^{d\left( m-1\right) }$
as in (\ref{DmR}) or $\mathbb{D}_{m}=\left[ -\pi ,\pi \right] ^{d\left(
m-1\right) }$ as in (\ref{Dm}). The both cases are completely similar since
we use the same properties of the spaces $L_{1}=$ $L_{1}\left( \left[ -\pi
,\pi \right] ^{d\ }\right) $ or $L_{1}=$ $L_{1}\left( \mathbb{R}^{d\
}\right) $, and we do not use in our proofs the boundedness and compactness
of the domain $\left[ -\pi ,\pi \right] ^{d}$. Hence, we will everywhere
consider the periodic case $\mathbb{\ }\left[ -\pi ,\pi \right] ^{d}$ which
corresponds to lattice equations and without further comment apply the
results to the case $\mathbb{R}^{d}$.

\begin{lemma}
\label{Lemma 8} The operator $\mathcal{F}^{\left( m\right) }$ defined by (%
\ref{Fu}), (\ref{Fmintr}) is bounded from $E=C\left( \left[ 0,\tau _{\ast }%
\right] ,L_{1}\right) $ into $C^{1}\left( \left[ 0,\tau _{\ast }\right]
,L_{1}\right) $ and its norm is estimated as follows 
\begin{equation}
\left\Vert \mathcal{F}^{\left( m\right) }\left( \mathbf{\tilde{u}}_{1}\ldots 
\mathbf{\tilde{u}}_{m}\right) \right\Vert _{E}\leq \tau _{\ast }C_{\Xi
}^{2m+1}\left\Vert \chi ^{\left( m\right) }\right\Vert
\dprod\limits_{j=1}^{m}\left\Vert \mathbf{\tilde{u}}_{j}\right\Vert _{E},
\label{dtf}
\end{equation}%
\begin{equation}
\left\Vert \partial _{\tau }\mathcal{F}^{\left( m\right) }\left( \mathbf{%
\tilde{u}}_{1}\ldots \mathbf{\tilde{u}}_{m}\right) \right\Vert _{E}\leq
C_{\Xi }^{2m+1}\left\Vert \chi ^{\left( m\right) }\right\Vert
\dprod\limits_{j}\left\Vert \mathbf{\tilde{u}}_{j}\right\Vert _{E}.
\label{dtf1}
\end{equation}
\end{lemma}

\begin{proof}
According to Condition \ref{Diagonalization} we can diagonalize the matrix $%
\exp \left\{ -\mathrm{i}\mathbf{L}\left( \mathbf{\mathbf{k}}\right) \frac{%
\tau _{1}}{\varrho }\right\} $ and its norm is bounded\ uniformly in $%
\mathbf{\mathbf{k}},\tau _{1}$ and $\varrho $:%
\begin{equation}
\left\Vert \exp \left\{ -\mathrm{i}\mathbf{L}\left( \mathbf{\mathbf{k}}%
\right) \frac{\tau _{1}}{\varrho }\right\} \right\Vert \leq C_{\Xi }^{2}\
\forall \mathbf{\mathbf{k}}\in \mathbb{R}^{d},\varrho >0,\tau _{1}\geq 0.
\label{expL}
\end{equation}%
By (\ref{Fu}), (\ref{Yconv}) and (\ref{Fmintr})%
\begin{gather*}
\left\Vert \mathcal{F}^{\left( m\right) }\left( \mathbf{\tilde{u}}_{1}\ldots 
\mathbf{\tilde{u}}_{m}\right) \left( \mathbf{\cdot },\tau \right)
\right\Vert _{L_{1}}\leq C_{\Xi }^{2m+1}\sup_{\ \mathbf{\mathbf{k}},\vec{k}%
}\left\vert \chi _{\ }^{\left( m\right) }\left( \mathbf{\mathbf{k}},\vec{k}%
\right) \right\vert \int \int_{0}^{\tau }\int_{\mathbb{D}_{m}} \\
\left\vert \mathbf{\tilde{u}}_{1}\left( \mathbf{k}^{\prime }\right)
\right\vert \ldots \left\vert \mathbf{\tilde{u}}_{m}\left( \mathbf{k}%
^{\left( m\right) }\left( \mathbf{k},\vec{k}\right) \right) \right\vert 
\mathrm{d}\mathbf{k}^{\prime }\ldots \mathrm{d}\mathbf{k}^{\left( m-1\right)
}\mathrm{d}\tau _{1}\mathrm{d}\mathbf{k}\leq \\
C_{\Xi }^{2m+1}\left\Vert \chi ^{\left( m\right) }\right\Vert \int_{0}^{\tau
}\left\Vert \mathbf{\tilde{u}}_{1}\left( \tau _{1}\right) \right\Vert
_{L_{1}}\ldots \left\Vert \mathbf{\tilde{u}}_{m}\left( \tau _{1}\right)
\right\Vert _{L_{1}}\mathrm{d}\tau _{1}\leq \tau _{\ast }C_{\Xi
}^{2m+1}\left\Vert \chi ^{\left( m\right) }\right\Vert \left\Vert \mathbf{%
\tilde{u}}_{1}\right\Vert _{E}\ldots \left\Vert \mathbf{\tilde{u}}%
_{m}\right\Vert _{E}.
\end{gather*}%
Similarly,%
\begin{gather*}
\left\Vert \partial _{\tau }\mathcal{F}^{\left( m\right) }\left( \mathbf{%
\tilde{u}}_{1}\ldots \mathbf{\tilde{u}}_{m}\right) \left( \mathbf{\cdot }%
,\tau \right) \right\Vert _{L_{1}}\leq \\
C_{\Xi }^{2m+1}\left\Vert \chi ^{\left( m\right) }\right\Vert \int \int_{%
\mathbb{D}_{m}}\left\vert \mathbf{\tilde{u}}_{1}\left( \mathbf{k}^{\prime
}\right) \right\vert \ldots \left\vert \mathbf{\tilde{u}}_{m}\left( \mathbf{k%
}^{\left( m\right) }\left( \mathbf{k},\vec{k}\right) \right) \right\vert 
\mathrm{d}\mathbf{k}^{\prime }\ldots \,\mathrm{d}\mathbf{k}^{\left(
m-1\right) }\mathrm{d}\mathbf{k}\,\leq \\
\left\Vert \chi ^{\left( m\right) }\right\Vert \left\Vert \mathbf{\tilde{u}}%
_{1}\right\Vert _{E}\ldots \left\Vert \mathbf{\tilde{u}}_{m}\right\Vert _{E}.
\end{gather*}
\end{proof}

\begin{corollary}
\label{Corollary normM}\textbf{\ }If $M\left( \mathcal{F},T,\vec{\lambda},%
\vec{\zeta}\right) \left( \mathbf{x}_{1}\ldots \mathbf{x}_{m}\right) $ is a
decorated composition monomial defined by (\ref{recG0}) and $\mathcal{F}$ is
defined by (\ref{varcu}), (\ref{Fu}) then%
\begin{equation}
\left\Vert M\left( \mathcal{F},T,\vec{\lambda},\vec{\zeta}\right) \left( 
\mathbf{x}_{1}\ldots \mathbf{x}_{m}\right) \right\Vert _{E}\leq C_{\Xi
}^{2e+i}\tau _{\ast }^{i}\dprod\limits_{N\in T}\left\Vert \chi _{\ }^{\left(
\mu \left( N\right) \right) }\right\Vert \dprod\limits_{l=1}^{m}\left\Vert 
\mathbf{x}_{l}\right\Vert _{E},  \label{gq}
\end{equation}%
\begin{equation}
\left\Vert \partial _{\tau }M\left( \mathcal{F},T,\vec{\lambda},\vec{\zeta}%
\right) \left( \mathbf{x}_{1}\ldots \mathbf{x}_{m}\right) \right\Vert
_{E}\leq C_{\Xi }^{2e+i}\tau _{\ast }^{i-1}\dprod\limits_{N\in T}\left\Vert
\chi _{\ }^{\left( \mu \left( N\right) \right) }\right\Vert
\dprod\limits_{l=1}^{m}\left\Vert \mathbf{x}_{l}\right\Vert _{E},
\label{gq1}
\end{equation}%
where $i$ is the incidence number of the tree $T$, and $e$ is the number of
edges of $T$.
\end{corollary}

\begin{proof}
We estimate the norm of the monomial $M=\mathcal{F}^{\left( m\right) }\left(
M_{1}\ldots M_{m}\right) $ and its time derivative applying Lemma \ref{Lemma
8}. Then we use (\ref{dtf}) to estimate $\left\Vert M_{j}\right\Vert
_{C\left( \left[ 0,\tau _{\ast }\right] ,L_{1}\right) }$. The formal proof
is straightforward and uses the induction with respect to the incidence
number of a monomial.
\end{proof}

Using boundedness of operators $\mathcal{F}^{\left( m\right) }$ we obtain in
a standard way uniqueness of solution of (\ref{varcu}).

\begin{lemma}
\label{Lemma uniqueness} \ If $\mathbf{\tilde{u}}_{1},\mathbf{\tilde{u}}%
_{2}\in C\left( \left[ 0,\tau _{0}\right] ,L_{1}\right) $ \ with $\tau
_{0}>0 $ \ are two solutions of (\ref{varcu}) \ with the same $\mathbf{%
\tilde{h}}$\ then $\mathbf{\tilde{u}}_{1}=\mathbf{\tilde{u}}_{2}$.
\end{lemma}

\begin{proof}
\ Applying Lemma \ref{Lemma 101Continuity},\ we conclude that \ 
\begin{equation*}
\left\Vert \mathcal{F}\left( \mathbf{\tilde{u}}_{1}\right) -\mathcal{F}%
\left( \mathbf{\tilde{u}}_{2}\right) \right\Vert _{C\left( \left[ 0,\tau _{1}%
\right] ,L_{1}\right) }\leq C\tau _{1}\left\Vert \mathcal{F}\left( \mathbf{%
\tilde{u}}_{1}\right) -\mathcal{F}\left( \mathbf{\tilde{u}}_{2}\right)
\right\Vert _{C\left( \left[ 0,\tau _{1}\right] ,L_{1}\right) },\;0<\tau
_{1}\leq \tau _{0}.
\end{equation*}%
Deriving the above inequality we use that since $N_{F}<\infty $ \ the radius 
$R_{\mathcal{F}}$ in Lemma \ref{Lemma 101Continuity} is arbitrary large and $%
C_{\mathcal{F}}$ in (\ref{Lipf}) \ according to (\ref{dtf}) is proportional
to $\tau _{1}$. When the Lipschitz constant $C\tau _{1}<1$, in a standard
way we obtain that $\mathbf{\tilde{u}}_{1}\left( \tau \right) =\mathbf{%
\tilde{u}}_{2}\left( \tau \right) $ for $0\leq \tau \leq \tau _{1}$. Since
this statement can be applied to $\mathbf{\tilde{u}}_{1}\left( \tau -\tau
_{1}\right) $ and $\mathbf{\tilde{u}}_{2}\left( \tau -\tau _{1}\right) $ we
obtain that solutions coincide for $0\leq \tau \leq \tau _{0}$.
\end{proof}

\subsection{Function-analytic expansion of solutions for modal integral
evolution equation}

The reduced evolution equation (\ref{varcu}) has the form 
\begin{equation}
\mathbf{\mathbf{\tilde{u}}}=\mathcal{F}\left( \mathbf{\mathbf{\tilde{u}}}%
\right) +\mathbf{\tilde{x}},  \label{UF}
\end{equation}%
where $\mathbf{\mathbf{\tilde{u}}}$, $\mathbf{\tilde{x}}$ are functions of $%
\left( \mathbf{k},\tau \right) $. The nonlinear operator $\mathcal{F}$ in
the right-hand side of (\ref{UF}) \ is determined by (\ref{Fu}), $\mathbf{%
\tilde{x}}\left( \mathbf{k},\tau \right) =\mathbf{\tilde{h}}\left( \mathbf{k}%
\right) $ as in (\ref{varcu}). We look for the solution operator $\mathcal{G}
$ in the form of operator series 
\begin{equation}
\mathbf{\mathbf{\tilde{u}}}=\mathcal{G}\left( \mathbf{\tilde{x}}\right)
=\sum_{m=1}^{\infty }\mathcal{G}^{\left( m\right) }\left( \mathbf{\tilde{x}}%
^{\left( m\right) }\right) .  \label{uTj}
\end{equation}%
The questions related to the existence and the convergence of such series
are addressed in Theorem \ref{Imfth}. As a direct corollary of Theorem \ref%
{Imfth} and Lemma \ref{Lemma uniqueness}\ if\ applied to the reduced
evolution equation (\ref{varcu}) we obtain the following theorem.

\begin{theorem}
\label{Imfth1}\textbf{\ }Let 
\begin{equation}
\left\Vert \mathbf{\tilde{x}}\right\Vert _{E\ }<R_{\mathcal{G}}=\left( \tau
_{\ast }C_{\chi }C_{\Xi }^{2m_{F}+1}\right) ^{-1/\left( m_{F}-1\right)
}/8,\;\tau _{\ast }\leq C_{\Xi }^{-3}C_{\chi }^{-1}.  \label{Jc}
\end{equation}%
with $C_{\chi }$ as in (\ref{chiCR}), $C_{\Xi }$ as in (\ref{supxi}). Then
the series (\ref{uTj}) converges in $E=C\left( \left[ 0,\tau _{\ast }\right]
,L_{1}\right) $. The solution operator $\mathcal{G}\left( \mathbf{\tilde{x}}%
\right) =\mathbf{\tilde{u}}$ determines the solution to (\ref{UF}) and the
operators $\mathcal{G}^{\left( m\right) }$ in series (\ref{uTj}) satisfy the
recursive relations \ (\ref{recG0}).
\end{theorem}

\begin{proof}
>From (\ref{chiCR}) and (\ref{dtf}) we infer that $\mathcal{F}$ defined by (%
\ref{Fseries}) belongs to the class $A\left( C_{\mathcal{F}},R_{\mathcal{F}%
}\right) $ if 
\begin{equation*}
\tau _{\ast }C_{\chi }C_{\Xi }^{2m+1}\leq C_{\mathcal{F}}R_{\mathcal{F}%
}^{-m},\ m=2,\ldots ,m_{F}.
\end{equation*}%
\ If $C_{\Xi }^{-2}R_{\mathcal{F}}^{-1}\leq 1$ it is sufficient to verify
the above condition at $m=m_{F}$ only. After this we apply Theorem \ref%
{Imfth} where according to (\ref{CRC}) we can take 
\begin{equation}
C_{\mathcal{G}}=\frac{R_{\mathcal{F}}^{2}}{2\left( C_{\mathcal{F}}+R_{%
\mathcal{F}}\right) }\ ,\ R_{\mathcal{G}}=\frac{R_{\mathcal{F}}^{2}}{4\left(
C_{\mathcal{F}}+R_{\mathcal{F}}\right) }.  \label{CRC1}
\end{equation}%
We take 
\begin{equation}
C_{\mathcal{F}}=R_{\mathcal{F}}=\left( \tau _{\ast }C_{\chi }C_{\Xi
}^{2m_{F}+1}\right) ^{-1/\left( m_{F}-1\right) },\ C_{\mathcal{G}}=2R_{%
\mathcal{G}}=R_{\mathcal{F}}/4  \label{RF}
\end{equation}%
and apply Theorem \ref{Imfth}. \ Note that $C_{\Xi }^{-2}R_{\mathcal{F}%
}^{-1}\leq 1$ if $\tau _{\ast }\leq C_{\Xi }^{-3}C_{\chi }^{-1}$. $\ $
\end{proof}

>From Theorem \ref{Imfth1} (observing that by (\ref{RF}) $R_{\mathcal{F}%
}\rightarrow \infty $ when $\tau _{\ast }\rightarrow 0$) we obtain Theorems %
\ref{Theorem Existence} and Theorem \ref{Theorem Existence1}.

To prove Theorem \ref{Theorem Superposition} on the superposition principle
we apply the solution operator $\mathcal{G}$ to a sum of wavepackets $%
\mathbf{\tilde{h}}_{l}\left( \mathbf{k},\beta \right) $ as in Definition \ref%
{dwavepack}. For technical reasons we have to modify the wavepackets using
cut-off functions described below.

\paragraph{Cutoff functions.}

We often use an infinitely smooth cutoff function $\Psi \left( \mathbf{\eta }%
\right) $, $\mathbf{\eta }\in \mathbb{R}^{d}$, satisfying the following
relations 
\begin{eqnarray}
0 &\leq &\Psi \left( \mathbf{\eta }\right) \leq 1,\ \Psi \left( -\mathbf{%
\eta }\right) =\Psi \left( \mathbf{\eta }\right) ,\   \label{j0} \\
\Psi \left( \mathbf{\eta }\right) &=&1\text{ for }\left\vert \mathbf{\eta }%
\right\vert \leq \pi _{0}/2,\ \Psi \left( \mathbf{\eta }\right) =0\text{ for 
}\left\vert \mathbf{\eta }\right\vert \geq \pi _{0},  \notag
\end{eqnarray}%
where $\pi _{0}\leq 1$ is a sufficiently small number which satisfies the
inequality 
\begin{equation}
0<\pi _{0}<\frac{1}{2}\min_{l}\limfunc{dist}\left\{ \mathbf{k}_{\ast
l},\sigma \right\} .  \label{kpi0}
\end{equation}%
Using $\Psi $ we introduce cutoff functions $\Psi _{l,\zeta }\left( \mathbf{k%
},\beta \right) $ with support near $\zeta \mathbf{k}_{\ast l}$ defined as
follows: 
\begin{equation}
\Psi _{l,\zeta }\left( \mathbf{k},\beta \right) =\Psi \left( \frac{\mathbf{k}%
-\zeta \mathbf{k}_{\ast l}}{\beta ^{1-\epsilon }}\right) ,\ l=1,\ldots
,N_{h}.  \label{Psilz}
\end{equation}%
Here $\epsilon $ is a small number, $1/2>\epsilon >0$; we take the same $%
\epsilon $ as in Definition \ref{dwavepack}.

Given a wavepacket $\mathbf{\tilde{h}}_{l}\left( \mathbf{k},\beta \right) $
we introduce a modified wavepacket 
\begin{equation}
\mathbf{\tilde{h}}_{l}^{\Psi }\left( \mathbf{k},\beta \right) =\mathbf{%
\tilde{h}}_{l,+}^{\Psi }\left( \mathbf{k},\beta \right) +\mathbf{\tilde{h}}%
_{l,-}^{\Psi }\left( \mathbf{k},\beta \right) ,\ \mathbf{\tilde{h}}_{l,\zeta
}^{\Psi }\left( \mathbf{k},\beta \right) =\Psi _{l,\zeta }\left( \mathbf{k}%
,\beta \right) \mathbf{\tilde{h}}_{l,\zeta }\left( \mathbf{k},\beta \right) ,
\label{hpsi}
\end{equation}%
where $\Psi _{l,\zeta }$ are defined by (\ref{Psilz}).

\begin{proposition}
\label{Prophpsi} If $\mathbf{\tilde{h}}_{l}^{\ }\left( \mathbf{k},\beta
\right) $ is a wavepacket in the sense of Definition \ref{dwavepack} then $%
\mathbf{\tilde{h}}_{l}^{\Psi }\left( \mathbf{k},\beta \right) $ defined by (%
\ref{hpsi}) and (\ref{Psilz}) is also a wavepacket in the sense of
Definition \ref{dwavepack} and, in addition to that, 
\begin{equation}
\mathbf{\tilde{h}}_{l,\zeta }^{\Psi }\left( \mathbf{k},\beta \right) =0\text{
if\ }\left\vert \mathbf{k}-\zeta \mathbf{k}_{\ast l}\right\vert \geq \pi
_{0}\beta ^{1-\epsilon },  \label{hpsi0}
\end{equation}%
\begin{equation}
\left\Vert \mathbf{\tilde{h}}_{l}-\mathbf{\tilde{h}}_{l}^{\Psi }\right\Vert
_{L_{1}}\leq C\beta .  \label{L1hpsi}
\end{equation}
\end{proposition}

\begin{proof}
To obtain (\ref{L1hpsi}) we note that (\ref{sourloc}) and (\ref{j0}) imply: 
\begin{equation}
\left\Vert \left( 1-\Psi _{l,\zeta }\right) \mathbf{\tilde{h}}_{l,\zeta
}\right\Vert _{L_{1}}=\int \left\vert \left( 1-\Psi _{l,\zeta }\left( 
\mathbf{k},\beta \right) \right) \tilde{h}_{l,\zeta }\left( \mathbf{k}%
\right) \right\vert \mathrm{d}\mathbf{k}\leq C\beta ,  \label{1minpsi}
\end{equation}%
and (\ref{L1hpsi}) follows. Remaining statements are obtained by a
straightforward verification.
\end{proof}

The following lemma shows that we can replace $\mathbf{\tilde{h}}_{l}$ by $%
\mathbf{\tilde{h}}_{l}^{\Psi }\mathbf{\ }$in the statement of Theorem \ref%
{Theorem Superposition}, in particular in (\ref{Gsum}), (\ref{rem}).

\begin{lemma}
\label{Lemma 1minpsi}\textbf{\ }Let $\mathbf{\tilde{h}}_{l,\zeta }$ satisfy (%
\ref{sourloc}) and $\mathbf{\tilde{h}}_{l}^{\Psi }\left( \mathbf{k},\beta
\right) $ be defined by (\ref{hpsi}). Let 
\begin{equation}
\left\Vert \mathbf{\tilde{h}}_{l}\right\Vert \leq R,\ l=1,\ldots ,N_{h}\text{
\ where \ }N_{h}R<R_{\mathcal{G}}.  \label{Nhh}
\end{equation}%
Then the difference 
\begin{equation}
\left[ \mathcal{G}\left( \sum_{l=1}^{N_{h}}\mathbf{\tilde{h}}_{l}\right)
-\sum_{l=1}^{N_{h}}\mathcal{G}\left( \mathbf{\tilde{h}}_{l}\right) \right] -%
\left[ \mathcal{G}\left( \sum_{l=1}^{N_{h}}\mathbf{\tilde{h}}_{l}^{\Psi
}\right) -\sum_{l=1}^{N_{h}}\mathcal{G}\left( \mathbf{\tilde{h}}_{l}^{\Psi
}\right) \right] =B_{\Psi },  \label{Bpsi}
\end{equation}%
is small, namely 
\begin{equation}
\left\Vert B_{\Psi }\right\Vert _{E}\leq C\left( R\right) \beta .
\label{Bpsin}
\end{equation}
\end{lemma}

\begin{proof}
Note that since $0\leq \Psi _{l}\leq 1$ we have 
\begin{equation}
\left\Vert \Psi _{l,\zeta }\mathbf{\tilde{h}}_{l,\zeta }\right\Vert
_{L_{1}}\leq \left\Vert \mathbf{\tilde{h}}_{l,\zeta }\right\Vert _{L_{1}},\
\left\Vert \left( 1-\Psi _{l,\zeta }\right) \mathbf{\tilde{h}}_{l,\zeta
}\right\Vert _{L_{1}}\leq \left\Vert \mathbf{\tilde{h}}_{l,\zeta
}\right\Vert _{L_{1}},  \label{psile1}
\end{equation}%
and (\ref{1minpsi}). Using the Lipschitz continuity of the solution operator 
$\mathcal{G}$ (see \ref{Lemma 101Continuity}) and (\ref{L1hpsi}) we obtain (%
\ref{Bpsin}).
\end{proof}

\paragraph{Truncation.}

We will truncate the infinite series (\ref{uTj}). To this end we define an
integer $m=m\left( \beta ^{q}\right) $ as a solution of the inequality \ 
\begin{equation}
\frac{2\left\vert \ln \beta ^{q}\right\vert }{\left\vert \ln R_{\mathcal{G}%
}\right\vert }<m\left( \beta ^{q}\right) \leq \frac{2\left\vert \ln \beta
^{q}\right\vert }{\left\vert \ln R_{\mathcal{G}}\right\vert }+1,
\label{mofrho}
\end{equation}%
where $R_{\mathcal{G}}$ is the same as in (\ref{Jc}).\ We consider then the
following partial sum of the expansion (\ref{uTj}) 
\begin{equation}
\mathcal{G}_{m\left( \beta ^{q}\right) }\left( \mathbf{\tilde{h}}\right)
=\sum_{m=1}^{m\left( \beta ^{q}\right) }\mathcal{G}^{\left( m\right) }\left( 
\mathbf{\tilde{h}}^{\left( m\right) }\right)  \label{uGjm}
\end{equation}%
and readily conclude that the following statement holds.

\begin{lemma}
\label{Lemma Truncation}\textbf{\ }Let\textbf{\ }$\mathcal{G}$ be defined by
(\ref{uTj}), then \textbf{\ } 
\begin{equation}
\left\Vert \mathcal{G}\left( \mathbf{\tilde{h}}\right) -\mathcal{G}_{m\left(
\beta \right) }\left( \mathbf{\tilde{h}}\right) \right\Vert _{E}\leq C\left(
R\right) \beta \text{ when \ }\left\Vert \mathbf{\tilde{h}}\right\Vert
_{E}\leq R<R_{\mathcal{G}}.  \label{uGmjd}
\end{equation}
\end{lemma}

\subsubsection{SI-CI splitting for evaluated monomials}

We consider a function $\mathbf{\tilde{h}}$ which is a sum of the form (\ref%
{JJ1}) and the solution $\mathcal{G}\left( \mathcal{F},\mathbf{\tilde{h}}%
\right) $. Expanding $\mathcal{G}^{\left( m\right) }\left( \mathbf{\tilde{h}}%
^{\left( m\right) }\right) $ into composition monomials as in (\ref{Gtreeexp}%
) we obtain a sum of composition monomials $M\left( \mathcal{F},T\right)
\left( \mathbf{\tilde{h}}^{m}\right) $. Then we look at the $m$-linear
monomial $M\left( \mathcal{F},T\right) \left( \mathbf{\tilde{h}}^{m}\right) $
where $\mathbf{\tilde{h}}$ equals a sum of $N_{h}$ one-band wavepacket $\ 
\mathbf{\tilde{h}}_{l}$ as in (\ref{JJ1}). Using the linearity with respect
to each argument we expand the monomial into a sum of $N_{h}^{m}$
expressions (\emph{evaluated monomials})%
\begin{equation}
M\left( \mathcal{F},T\right) \left( \sum_{l=1}^{N_{h}}\mathbf{\tilde{h}}%
_{l}\right) ^{m}=\sum_{l_{1},\ldots ,l_{m}}M\left( \mathcal{F},T\right)
\left( \mathbf{\tilde{h}}_{l_{1}}\ldots \mathbf{\tilde{h}}_{l_{m}}\right)
=\sum_{l_{1},\ldots ,l_{m}}M\left( \mathcal{F},T\right) \left(
\dprod\limits_{i}\mathbf{\tilde{h}}_{l_{i}}\right) .  \label{hlexp}
\end{equation}%
The sum contains evaluated monomials of two kinds: (i) ones which involve
the same wavepacket; and (ii) one corresponding to the cross-terms (terms
involving different wavepackets). To be precise, we introduce the following
definition.

\begin{definition}[SI and CI]
\label{Definition GVM} \ We say that an evaluated monomial $M\left( \mathcal{%
F},T\right) \left( \mathbf{\tilde{h}}_{l_{1}}\ldots \mathbf{\tilde{h}}%
_{l_{m}}\right) $ with the argument multiindex $l_{1},\ldots ,l_{m}\in
\left\{ 1,\ldots ,N\right\} ^{m}$ in the expansion (\ref{hlexp}) is \emph{%
self-interacting} (SI) if 
\begin{equation}
l_{1}=l_{2}=\ldots =l_{m}.  \label{lGVM}
\end{equation}%
Otherwise we say that $M\left( \mathcal{F},T\right) \left( \mathbf{\tilde{h}}%
_{l_{1}}\ldots \mathbf{\tilde{h}}_{l_{m}}\right) $\ is \emph{%
cross-interacting} (CI).
\end{definition}

Using this notation we rewrite (\ref{hlexp}): 
\begin{equation}
M\left( \mathcal{F},T\right) \left( \left( \sum_{l=1}^{N_{h}}\mathbf{\tilde{h%
}}_{l}\right) ^{m}\right) =\sum_{l=1}^{N_{h}}M\left( \mathcal{F},T\right)
\left( \left( \mathbf{\tilde{h}}_{l}\right) ^{m}\right) +\sum_{l_{1},\ldots
,l_{m}\text{ is CI}}M\left( \mathcal{F},T\right) \left( \mathbf{\tilde{h}}%
_{l_{1}}\ldots \mathbf{\tilde{h}}_{l_{m}}\right) .  \label{Msum}
\end{equation}%
Substituting this expression into (\ref{Gtreeexp}) we obtain the expansion 
\begin{gather}
\mathcal{G}\left( \mathbf{\tilde{h}}_{1}+\ldots +\mathbf{\tilde{h}}%
_{N_{h}}\right) =\sum_{m=1}^{\infty }\mathcal{G}_{m}\left( \left( \mathbf{%
\tilde{h}}_{1}+\ldots +\mathbf{\tilde{h}}_{N_{h}}\right) ^{m}\right)
\label{Gofsum} \\
=\sum_{m=1}^{\infty }\mathcal{G}\left( \left( \mathbf{\tilde{h}}_{1}\right)
^{m}\right) +\ldots +\sum_{m=1}^{\infty }\mathcal{G}\left( \left( \mathbf{%
\tilde{h}}_{N_{h}}\right) ^{m}\right) +\mathcal{G}_{\text{CI}}\left( \mathbf{%
\tilde{h}}_{1},\ldots ,\mathbf{\tilde{h}}_{N_{h}}\right) ,  \notag
\end{gather}%
where $\mathcal{G}_{\text{CI}}$ contains only CI monomials with cross-terms.

\begin{proposition}
\label{minsub}\ \ Every evaluated CI monomial $M\left( \mathcal{F},T\right)
\left( \mathbf{\tilde{h}}_{1},\ldots ,\mathbf{\tilde{h}}_{N_{h}}\right) $
has a submonomial of the form 
\begin{equation}
\mathcal{F}^{\left( s\right) }\left( M\left( \mathcal{F},T_{1}\right) \left( 
\mathbf{\tilde{h}}_{l_{1}}\ldots \mathbf{\tilde{h}}_{l_{1}}\right) \ldots
M\left( \mathcal{F},T_{s}\right) \left( \mathbf{\tilde{h}}_{l_{s}}\ldots 
\mathbf{\tilde{h}}_{l_{s}}\right) \right)  \label{FNGVM}
\end{equation}%
where all $M\left( \mathcal{F},T_{1}\right) \left( \mathbf{\tilde{h}}%
_{l_{1}}\ldots \mathbf{\tilde{h}}_{l_{1}}\right) $,..., $M\left( \mathcal{F}%
,T_{s}\right) \left( \mathbf{\tilde{h}}_{l_{s}}\ldots \mathbf{\tilde{h}}%
_{l_{s}}\right) $ are SI, and there are at least two indices $i$ and $j$
such that $\mathbf{\tilde{h}}_{l_{i}}\neq \mathbf{\tilde{h}}_{l_{j}}$. We
call such a monomial a minimal CI\ monomial.
\end{proposition}

\begin{proof}
The set of CI submonomials of $M\left( \mathcal{F},T\right) $ is finite and
it is non-empty since $M\left( \mathcal{F},T\right) $ itself is a CI
monomial. We take CI submonomial of $M\left( \mathcal{F},T\right) $ with a
minimal rank.\ Its rank is non-zero since every zero rank submonomial is SI.
Since the rank is minimal all submonomials are SI. Hence it has the form (%
\ref{FNGVM}).
\end{proof}

\subsection{Properties of SI monomials}

According to Definition \ref{Definition GVM} for a SI\ evaluated monomial we
have $\mathbf{\tilde{h}}_{l_{1}}=\ldots =\mathbf{\tilde{h}}_{l_{m}}$.
Observe also that in view of Definition \ref{dwavepack} every single-band
wavepacket $\mathbf{\tilde{h}}_{l}$ has its band number, and $n^{\prime
}=n^{\prime \prime }=\ldots =n^{\left( m\right) }$, that is the band $%
n_{l}=n_{0}$ is the same for all $\mathbf{\tilde{h}}_{l}$. Similarly, $%
\mathbf{k}_{\ast l_{1}}=\ldots =\mathbf{k}_{\ast l_{m}}$. Having these
properties we often omit in this section indices $n_{i}$, $l_{i}$ and skip $%
\vec{n}$ for notational brevity, writing, for example,%
\begin{equation*}
\omega _{n,\zeta }\left( \mathbf{k}\right) =\omega _{\zeta }\left( \mathbf{k}%
\right) ,\ \mathbf{\tilde{u}}_{n,\zeta }\left( \mathbf{k}\right) =\mathbf{%
\tilde{u}}_{\zeta }\left( \mathbf{k}\right) ,\ \chi _{\ n,\zeta ,\vec{n},%
\vec{\zeta}}^{\left( m\right) }=\chi _{\zeta ,\vec{\zeta}}^{\left( m\right)
}.
\end{equation*}

\subsubsection{Monomials applied to a single-band wavepacket.}

Here we consider monomials based on oscillatory integral operators and which
are applied to a single-band wavepacket. We recall that according to (\ref%
{hbold}) a single-band wavepacket $\mathbf{\tilde{h}}$ involves two
components $\mathbf{\tilde{h}}_{+}$ and $\mathbf{\tilde{h}}_{-}$ and a small
complement component $\mathbf{\tilde{h}}_{\infty }$.

\begin{definition}[frequency matching]
\label{Definition FM} We call a decorated composition monomial $M\left( 
\mathcal{F},T,\vec{\lambda},\vec{\zeta}\right) $ \emph{frequency matched}
(FM) if for every non-end node $N\in T$ the corresponding decorated
submonomial $M^{\prime }=\mathcal{F}_{\lambda }^{\left( m^{\prime }\right)
}\left( M_{1,\zeta ^{\prime }}\ldots M_{m^{\prime },\zeta ^{\left( m^{\prime
}\right) }}\right) $ satisfies the following conditions:%
\begin{equation}
\lambda \neq \infty ,\ \zeta ^{\left( j\right) }\neq \infty ,\ j=1,\ldots
,m^{\prime },  \label{neinf}
\end{equation}%
and%
\begin{equation}
\sum_{j=1}^{m^{\prime }}\zeta ^{\left( j\right) }=\lambda ,  \label{fmL}
\end{equation}%
where $\lambda ,\zeta ^{\left( j\right) }\in \Lambda $ defined by (\ref{PPPI}%
), we identify $\pm $ with $\pm 1$. A decorated composition monomial which
does not satisfy the above conditions is called \emph{not frequency matched}
(NFM) monomial.
\end{definition}

Collecting separately FM and NFM terms in the expression (\ref{MFsplit}) we
obtain 
\begin{gather}
M\left( \mathcal{F},T\right) \left( \mathbf{x}_{1}\mathbf{x}_{2}\ldots 
\mathbf{x}_{m}\right) =\sum_{\text{FM\ }\vec{\lambda},\vec{\zeta}}M\left( 
\mathcal{F},T,\vec{\lambda},\vec{\zeta}\right) \left( \mathbf{x}_{1}\mathbf{x%
}_{2}\ldots \mathbf{x}_{m}\right)  \label{GFMN} \\
+\sum_{\text{NFM\ }\vec{\lambda},\vec{\zeta}}M\left( \mathcal{F},T,\vec{%
\lambda},\vec{\zeta}\right) \left( \mathbf{x}_{1}\mathbf{x}_{2}\ldots 
\mathbf{x}_{m}\right) .  \notag
\end{gather}

\begin{remark}
Any SI evaluated monomial is either FM or NFM. We do not define for CI\
evaluated monomials \ if they are FM or NFM.
\end{remark}

Below we show that FM decorated monomials \ have the following properties
which can be briefly stated as follows.

\textbf{Property 1}. If $\mathbf{\tilde{h}}\left( \mathbf{k}\right) $ is a
wavepacket in the sense of Definition \ref{dwavepack} centered around $\pm 
\mathbf{k}_{\ast }$\ then FM\ monomial $M\left( \mathcal{F},T,\vec{\lambda},%
\vec{\zeta}\right) \left( \mathbf{\tilde{h}}^{m}\right) \left( \mathbf{k}%
\right) $ is also localized about $\pm \mathbf{k}_{\ast }$. This property is
proved below in Corollary \ref{Corollary 100}.

\textbf{Property 2}.\textbf{\ }The most important property concerning FM-NFM
splitting is that the result of a NFM\ monomial application to a wavepacket
has magnitude $O\left( \varrho \right) $, that is $O\left( \beta ^{2}\right) 
$ for the scaling (\ref{scale1}). Consequently, all NFM\ terms in (\ref{GFMN}%
) are small (see Lemma \ref{Lemma NFM} below) and they give contribution
only to the remainder $\mathbf{\tilde{D}}$ in (\ref{Gsum}).

Now we formulate exact statements clarifying the above properties. The
following two statements show, in particular, that an FM monomial transforms
a function supported in a vicinity of $\mathbf{k}_{\ast }$ into a similar
function.

\begin{lemma}[operator support ]
\label{Lemma Fsupport} If $\mathbf{\tilde{u}}_{1,\zeta ^{\prime }}\ldots 
\mathbf{\tilde{u}}_{m,\zeta ^{\left( m\right) }}$ \ are such that 
\begin{equation*}
\mathbf{\tilde{u}}_{\zeta ^{\left( l\right) }}\left( \mathbf{k}^{\left(
l\right) }\right) =0\text{ \ when }\left\vert \mathbf{k}^{\left( l\right)
}-\zeta ^{\left( l\right) }\mathbf{k}_{\ast }\right\vert >\delta _{l},\
l=1,\ldots ,m,
\end{equation*}%
and 
\begin{equation}
\mathbf{k}_{\vec{\zeta}}=\left( \zeta ^{\prime }+\ldots +\zeta ^{\left(
m\right) }\right) \mathbf{k}_{\ast }.  \label{kz}
\end{equation}%
then $\mathcal{F}^{\left( m\right) }\left( \mathbf{\tilde{u}}_{1,\zeta
^{\prime }}\ldots \mathbf{\tilde{u}}_{m,\zeta ^{\left( m\right) }}\right)
\left( \mathbf{k},\tau \right) $ given \ by (\ref{Fu}), \ satisfies 
\begin{equation}
\mathcal{F}_{\zeta }^{\left( m\right) }\left( \mathbf{\tilde{u}}_{1,\zeta
^{\prime }}\ldots \mathbf{\tilde{u}}_{m,\zeta ^{\left( m\right) }}\right)
\left( \mathbf{k},\tau \right) =0\text{ if }\left\vert \mathbf{k}-\mathbf{k}%
_{\vec{\zeta}}\right\vert >\delta _{1}+\ldots +\delta _{m}.  \label{FMz0}
\end{equation}%
In particular, if the binary indices $\ \zeta ,\vec{\zeta}_{\left( m\right)
} $ are frequency matched (FM), that is 
\begin{equation}
\zeta =\zeta ^{\prime }+\ldots +\zeta ^{\left( m\right) },\text{ where }%
\zeta ^{\left( j\right) },\ \zeta =\pm 1,  \label{FMz}
\end{equation}%
then (\ref{FMz0}) holds with $\mathbf{k}_{\vec{\zeta}}=\zeta \mathbf{k}%
_{\ast }$.
\end{lemma}

\begin{proof}
>From (\ref{Fm}) and (\ref{FMz}) \ we obtain the equality 
\begin{equation*}
\mathbf{k}-\zeta \mathbf{k}_{\ast }=\left( \mathbf{\mathbf{k}^{\prime }}%
-\zeta ^{\prime }\mathbf{k}_{\ast }\right) +\ldots +\left( \mathbf{k}%
^{\left( m\right) }-\zeta ^{\left( m\right) }\mathbf{k}_{\ast }\right)
\end{equation*}%
which implies lemma's statement.
\end{proof}

\begin{corollary}[support of a monomial]
\label{Corollary 100}\ If $M\left( \mathcal{F},T,\vec{\lambda},\vec{\zeta}%
\right) \left( \mathbf{\tilde{h}}_{1}\ldots \mathbf{\tilde{h}}_{m}\right) $
is a decorated composition monomial and 
\begin{equation}
\mathbf{\tilde{h}}_{l,\zeta ^{\left( l\right) }}=0\text{ \ when }\left\vert 
\mathbf{k}^{\left( l\right) }-\zeta ^{\left( l\right) }\mathbf{k}_{\ast
}\right\vert >\delta _{0},\ l=1,\ldots ,m,  \label{xeq0}
\end{equation}%
then 
\begin{equation}
\ M\left( \mathcal{F},T,\vec{\lambda},\vec{\zeta}\right) \left( \mathbf{%
\tilde{h}}_{1}\ldots \mathbf{\tilde{h}}_{m}\right) \left( \mathbf{k}\right)
=0\text{\ if }\left\vert \mathbf{k}-\mathbf{k}_{\vec{\zeta}}\right\vert
>m\delta _{0},  \label{GFMz}
\end{equation}%
where $\mathbf{k}_{\vec{\zeta}}$ is defined by (\ref{kz}). In particular, if 
$M\left( \mathcal{F},T,\vec{\lambda},\vec{\zeta}\right) \left( \mathbf{%
\tilde{h}}_{1}\ldots \mathbf{\tilde{h}}_{m}\right) $ is a FM decorated
composition monomial, then 
\begin{equation}
\ M\left( \mathcal{F},T,\vec{\lambda},\vec{\zeta}\right) \left( \mathbf{%
\tilde{h}}_{1}\ldots \mathbf{\tilde{h}}_{m}\right) \left( \mathbf{k}\right)
=0\text{\ if }\left\vert \mathbf{k}-\zeta \mathbf{k}_{\ast }\right\vert
>m\delta _{0}.  \label{GFM1}
\end{equation}%
where $\zeta $ satisfies (\ref{FMz}). In particular, if $\delta _{0}=\beta
^{1-\epsilon }$ and $m\leq C\ln \beta $\ then for any $\delta _{1}>0$ there
exists $\beta _{0}$ such that for $\beta <\beta _{0}$ we have $C\pi
_{0}\beta ^{1-\epsilon }\ln \beta <\delta _{1}$ and 
\begin{equation}
M\left( \mathcal{F},T,\vec{\lambda},\vec{\zeta}\right) \left( \mathbf{\tilde{%
h}}_{1}\ldots \mathbf{\tilde{h}}_{m}\right) \left( \mathbf{k}\right) =0\text{%
\ when }\left\vert \mathbf{k}-\zeta \mathbf{k}_{\ast }\right\vert >C\pi
_{0}\beta ^{1-\epsilon }\ln \beta .  \label{GFM0}
\end{equation}
\end{corollary}

\begin{proof}
To obtain (\ref{GFMz}) we apply Lemma \ref{Lemma Fsupport} \ and use\ the
induction with respect to the rank of a monomial.
\end{proof}

\begin{remark}
If $M\left( \mathcal{F},T,\vec{\lambda},\vec{\zeta}\right) $ is NFM and $%
\mathbf{\tilde{h}}\left( \mathbf{k}\right) $ \ is a wavepacket localized
near $\pm \mathbf{k}_{\ast }$, \ then $M\left( \mathcal{F},T,\vec{\lambda},%
\vec{\zeta}\right) \left( \mathbf{\tilde{h}}^{m}\right) \left( \mathbf{k}%
\right) $ is localized near the point $\mathbf{k}_{\vec{\zeta}}$. \ As $\vec{%
\zeta}$ vary over $\left\{ -1,1\right\} ^{m}$ such points $\mathbf{k}_{\vec{%
\zeta}}$ lie on a straight line parallel to $\mathbf{k}_{\ast }$. For $%
m\rightarrow \infty $ the closure of the set of such $\mathbf{k}_{\vec{\zeta}%
}$ with a generic $\mathbf{k}_{\ast }$ can be the entire torus $\left[ -\pi
,\pi \right] ^{d}$, whereas for the case of \ $\vec{\zeta}$ corresponding to
an FM monomial the closure is just two points $\pm \mathbf{k}_{\ast }$.\
Hence Property 1 is very useful and, in particular, allows to avoid small
denominators in coupling terms.
\end{remark}

The following lemma shows that the FM interaction phase function of a single
wavepacket has a critical point at its center, or, in other words, FM
monomials satisfy the group velocity matching condition (see \cite{BF3}, 
\cite{BF6}).

\begin{lemma}
\label{Lemma FMgr} If a decorated operator $\mathcal{F}_{\zeta ,\vec{\zeta}%
_{\left( m\right) }}^{\left( m\right) }$ is FM then the interaction phase
function $\phi $ in (\ref{Fm}) has a critical point: 
\begin{equation}
\nabla _{\mathbf{\mathbf{k}}}\phi _{n,\zeta ,\vec{n},\vec{\zeta}}\left(
\zeta \mathbf{\mathbf{k}}_{\ast },\vec{k}_{\ast }\right) =0\ \text{at}\ \vec{%
k}_{\ast }=\left( \zeta ^{\prime }\mathbf{k}_{\ast },\ldots ,\zeta ^{\left(
m\right) }\mathbf{k}_{\ast }\right) .  \label{gr0}
\end{equation}
\end{lemma}

\begin{proof}
For FM decorated operator all indices $\zeta ^{\left( j\right) }=\pm $ \ and 
\begin{equation}
n=n^{\prime }=\ldots =n^{\left( m\right) }\text{and }\zeta =\zeta ^{\prime
}+\ldots +\zeta ^{\left( m\right) }.  \label{sumz}
\end{equation}%
Hence we obtain from (\ref{phim}) that%
\begin{equation*}
\nabla _{\mathbf{\mathbf{k}}}\phi _{n,\zeta ,\vec{n},\vec{\zeta}}\left( 
\mathbf{\mathbf{k}},\vec{k}\right) =\zeta \nabla _{\mathbf{\mathbf{k}}%
}\omega \left( \mathbf{k}\right) -\zeta ^{\left( m\right) }\nabla _{\mathbf{%
\mathbf{k}}}\omega \left( \mathbf{k}-\mathbf{k}^{\prime }-\ldots -\mathbf{k}%
^{\left( m-1\right) }\right) .
\end{equation*}%
Since $\zeta \mathbf{k}_{\ast }-\zeta ^{\prime }\mathbf{k}_{\ast }^{\prime
}-\ldots -\zeta ^{\left( m-1\right) }\mathbf{k}_{\ast }^{\left( m-1\right)
}=\zeta ^{\left( m\right) }\mathbf{k}_{\ast }^{\left( m\right) }$ and (\ref%
{omeven}) implies 
\begin{equation}
\zeta \nabla _{\mathbf{\mathbf{k}}}\omega \left( \zeta \mathbf{\mathbf{k}}%
_{\ast }\right) =\zeta ^{\left( m\right) }\nabla _{\mathbf{\mathbf{k}}%
}\omega \left( \zeta ^{\left( m\right) }\zeta \mathbf{\mathbf{k}}_{\ast
}\right) \text{ for }\zeta =\pm ,\zeta ^{\left( m\right) }=\pm ,
\label{pmom}
\end{equation}%
we obtain the desired (\ref{gr0}).
\end{proof}

Now we consider NFM monomials and prove the Property 2. \ \ First we note
that (\ref{omne00}) implies%
\begin{equation}
\omega _{n_{l}}\left( \mathbf{\mathbf{k}}_{\ast l}\right) \geq \omega _{\ast
}>0,\ l=1,\ldots ,N_{h}.  \label{omne0}
\end{equation}%
If $\mathbf{k}_{\ast l}=\mathbf{k}_{\ast },$ $n_{l}=n_{0}$ satisfy Condition %
\ref{Condition generic1} then if (\ref{FMeq}) does not hold, (\ref{zomomz0})
does not hold too, hence for $m\leq m_{F}$%
\begin{equation}
\left\vert \sum_{j=1}^{m}\zeta ^{\left( j\right) }\omega _{n_{0}}\left( 
\mathbf{k}_{\ast }\right) -\zeta \omega _{n}\left( \mathbf{k}_{\vec{\zeta}%
}\right) \right\vert \geq \omega _{\ast }>0,\;\mathbf{k}_{\vec{\zeta}%
}=\sum_{j=1}^{m}\zeta ^{\left( j\right) }\mathbf{k}_{\ast },  \label{zomomz}
\end{equation}%
where $\omega _{\ast }>0$ is a positive number (we take for notation
simplicity the same small enough constsant in (\ref{omne0}) and (\ref{zomomz}%
).

The following Lemma, which is a version of the standard statement of the
stationary phase method, shows that the action of an NFM monomial on a
wavepacket produces a wave of a small amplitude.

\begin{lemma}
\label{Lemma NFM} Let the decoration projections be defined by (\ref{PPosc}%
). Assume that Condition \ref{Condition generic1} holds. Let $\ $indices $%
\zeta ,\zeta ^{\prime },\ldots ,\zeta ^{\left( m\right) }$ \ be NFM, that is
either one of them is $\infty \ $or 
\begin{equation}
\zeta \neq \zeta ^{\prime }+\ldots +\zeta ^{\left( m\right) },\ \zeta
^{\left( j\right) }=\pm 1,\ \zeta =\pm 1.  \label{NFMz}
\end{equation}%
Let $\delta _{NFM}>0$ be small enough to satisfy 
\begin{equation}
\delta _{NFM}\max_{\left\vert \mathbf{k}_{\ast l}-\mathbf{k}\right\vert \leq
\delta _{NFM}}\left\vert \nabla \omega _{l}\left( \mathbf{\mathbf{k}}\right)
\right\vert \leq \frac{1}{4}\omega _{\ast },\ l=1,\ldots ,N_{h},
\label{desmall}
\end{equation}%
where $\omega _{\ast }$ is given in (\ref{zomomz}). Let $\mathbf{k}$,$%
\mathbf{k}^{\left( j\right) }$ satisfy (\ref{Phmc}) and be such that \ 
\begin{equation}
\sum_{j=1}^{m}\left\vert \mathbf{k}^{\left( j\right) }-\zeta ^{\left(
j\right) }\mathbf{k}_{\ast }\right\vert \leq \delta _{NFM},\ \left\vert 
\mathbf{k}-\mathbf{k}_{\vec{\zeta}}\right\vert \leq \delta _{NFM},
\label{kmik}
\end{equation}%
where $\mathbf{k}_{\vec{\zeta}}$ is defined by (\ref{kz}) and $\mathbf{k}%
_{\ast }=\mathbf{k}_{\ast l}$ satisfy the conditions (\ref{omne0}) and (\ref%
{zomomz}). Let the functions $\tilde{u}_{j,\zeta ^{\left( j\right) }}\left( 
\mathbf{k},\tau \right) $ satisfy the condition%
\begin{equation}
\text{ }\tilde{u}_{j,\zeta ^{\left( j\right) }}\left( \mathbf{k},\tau
\right) =0\text{ \ when }\zeta ^{\left( j\right) }=\infty \text{ and }\tilde{%
u}_{j,\zeta ^{\left( j\right) }}\left( \zeta ^{\left( j\right) }\mathbf{k}%
_{\ast }+\mathbf{s},\tau \right) =0\ \text{when }\left\vert \mathbf{s}%
\right\vert \geq \delta _{NFM}\text{.}  \label{cond0}
\end{equation}%
Then 
\begin{gather}
\left\Vert \mathcal{F}_{\zeta ,\zeta ^{\prime },\ldots ,\zeta ^{\left(
m\right) }}^{\left( m\right) }\left( \mathbf{\tilde{u}}_{1,\zeta ^{\prime
}}\ldots \mathbf{\tilde{u}}_{m,\zeta ^{\left( m\right) }}\right) \right\Vert
_{E}\leq \frac{4\varrho }{\omega _{\ast }}\left\Vert \chi ^{\left( m\right)
}\right\Vert C_{\Xi }^{2m+1}\dprod\limits_{j}\left\Vert \mathbf{\tilde{u}}%
_{j}\right\Vert _{E}  \label{NFMF} \\
+\frac{2\varrho \tau _{\ast }}{\omega _{\ast }}C_{\Xi }^{2m+1}\left\Vert
\chi ^{\left( m\right) }\right\Vert \sum_{i}\left\Vert \partial _{\tau }%
\mathbf{\tilde{u}}_{i}\right\Vert _{E}\dprod\limits_{j\neq i}\left\Vert 
\mathbf{\tilde{u}}_{j}\right\Vert _{E}.  \notag
\end{gather}
\end{lemma}

\begin{proof}
If one of the indices $\zeta ^{\prime },\ldots ,\zeta ^{\left( m\right) }$
equals $\infty $ by (\ref{cond0}) $\mathcal{F}_{\zeta ,\zeta ^{\prime
},\ldots ,\zeta ^{\left( m\right) }}^{\left( m\right) }=0$ and (\ref{NFMF})
is satisfied. Now we consider the case when all $\zeta ,\zeta ^{\prime
},\ldots ,\zeta ^{\left( m\right) }$ are finite. We denote for brevity $%
\omega _{n_{0}}=\omega $, $\mathbf{\mathbf{k}}_{\ast l}=\mathbf{\mathbf{k}}%
_{\ast }$ and $\phi _{n,\zeta ,\vec{n},\vec{\zeta}}=\phi $. Since (\ref{kmik}%
) holds we get \ from (\ref{phim}) that%
\begin{gather*}
\left\vert \phi \left( \mathbf{\mathbf{k}},\vec{k}\right) -\phi \left( 
\mathbf{\mathbf{k}},\vec{k}_{\ast }\right) \right\vert \leq \left\vert
\omega \left( \mathbf{k}^{\prime }\right) -\omega \left( \zeta ^{\prime }%
\mathbf{k}_{\ast }\right) \right\vert +\ldots +\left\vert \omega \left( 
\mathbf{k}^{\left( m\right) }\right) -\omega \left( \zeta ^{\left( m\right) }%
\mathbf{k}_{\ast }\right) \right\vert \\
\leq \max_{\left\vert \mathbf{k}_{\ast }-\mathbf{k}\right\vert \leq \delta
_{NFM}}\left\vert \nabla \omega \left( \mathbf{\mathbf{k}}\right)
\right\vert \sum_{j=1}^{m}\left\vert \mathbf{k}^{\left( j\right) }-\zeta
^{\left( j\right) }\mathbf{k}_{\ast }\right\vert \leq \delta
_{NFM}\max_{\left\vert \mathbf{k}_{\ast }-\mathbf{k}\right\vert \leq \delta
_{NFM}}\left\vert \nabla \omega \left( \mathbf{\mathbf{k}}\right)
\right\vert .
\end{gather*}%
Using (\ref{desmall}) we conclude that 
\begin{equation}
\left\vert \phi \left( \mathbf{\mathbf{k}},\vec{k}\right) \right\vert \geq
\left\vert \phi \left( \mathbf{\mathbf{k}},\vec{k}_{\ast }\right)
\right\vert -\frac{1}{4}\left\vert \omega _{\ast }\right\vert .  \label{fi1}
\end{equation}%
By (\ref{NFMz}) the condition (\ref{FMeq}) is not satisfied, therefore (\ref%
{zomomz}) \ holds and implies that 
\begin{equation}
\left\vert \phi \left( \mathbf{k}_{\vec{\zeta}},\vec{k}_{\ast }\right)
\right\vert \geq \omega _{\ast }.  \label{zom}
\end{equation}%
Using (\ref{zom}), (\ref{kmik}) and (\ref{desmall}) we conclude that 
\begin{equation}
\left\vert \phi \left( \mathbf{\mathbf{k}},\vec{k}_{\ast }\right)
\right\vert \geq \omega _{\ast }-\left\vert \omega \left( \mathbf{k}\right)
-\omega \left( \mathbf{k}_{\vec{\zeta}}\right) \right\vert \geq \omega
_{\ast }-\delta _{NFM}\max_{\left\vert \mathbf{k}_{\ast }-\mathbf{k}%
\right\vert \leq \delta _{NFM}}\left\vert \nabla \omega \left( \mathbf{%
\mathbf{k}}\right) \right\vert \geq \frac{3}{4}\omega _{\ast }.  \label{fige}
\end{equation}%
Together with (\ref{fi1}) this inequality implies that when (\ref{kmik})
holds we have the estimate 
\begin{equation}
\left\vert \phi \left( \mathbf{\mathbf{k}},\vec{k}\right) \right\vert \geq 
\frac{1}{2}\omega _{\ast }.  \label{fi2}
\end{equation}%
Now we note that the oscillatory factor in (\ref{Fm}) 
\begin{equation*}
\exp \left\{ \mathrm{i}\phi \left( \mathbf{\mathbf{k}},\vec{k}\right) \frac{%
\tau _{1}}{\varrho }\right\} =\frac{\varrho }{\mathrm{i}\phi \left( \mathbf{%
\mathbf{k}},\vec{k}\right) }\partial _{\tau _{1}}\exp \left\{ \mathrm{i}\phi
\left( \mathbf{\mathbf{k}},\vec{k}\right) \frac{\tau _{1}}{\varrho }\right\}
.
\end{equation*}%
Integrating (\ref{Fm})\ by parts with respect to $\tau _{1}$ we obtain%
\begin{gather}
\mathcal{F}_{\zeta ,\vec{\zeta}}^{\left( m\right) }\left( \mathbf{\tilde{u}}%
_{1}\ldots \mathbf{\tilde{u}}_{m}\right) \left( \mathbf{k},\tau \right) =
\label{NFMintbp} \\
\int_{\mathbb{D}_{m}}\frac{\varrho \exp \left\{ \mathrm{i}\phi \left( 
\mathbf{\mathbf{k}},\vec{k}\right) \frac{\tau }{\varrho }\right\} }{\mathrm{i%
}\phi \left( \mathbf{\mathbf{k}},\vec{k}\right) }\chi _{\zeta ,\vec{\zeta}%
}^{\left( m\right) }\left( \mathbf{\mathbf{k}},\vec{k}\right) \tilde{u}%
_{1,\zeta ^{\prime }}\left( \mathbf{k}^{\prime },\tau \right) \ldots \tilde{u%
}_{m,\zeta ^{\prime }}\left( \mathbf{k}^{\left( m\right) }\left( \mathbf{k},%
\vec{k}\right) ,\tau \right) \mathrm{\tilde{d}}^{\left( m-1\right) d}\vec{k}
\notag \\
-\int_{\mathbb{D}_{m}}\frac{\varrho }{\mathrm{i}\phi \left( \mathbf{\mathbf{k%
}},\vec{k}\right) }\chi _{\zeta ,\vec{\zeta}}^{\left( m\right) }\left( 
\mathbf{\mathbf{k}},\vec{k}\right) \tilde{u}_{1,\zeta ^{\prime }}\left( 
\mathbf{k}^{\prime },0\right) \ldots \tilde{u}_{m,\zeta ^{\prime }}\left( 
\mathbf{k}^{\left( m\right) }\left( \mathbf{k},\vec{k}\right) ,0\right) 
\mathrm{\tilde{d}}^{\left( m-1\right) d}\vec{k}  \notag \\
-\int_{0}^{\tau }\int_{\mathbb{D}_{m}}\frac{\varrho }{\mathrm{i}\phi \left( 
\mathbf{\mathbf{k}},\vec{k}\right) }\exp \left\{ \mathrm{i}\phi \left( 
\mathbf{\mathbf{k}},\vec{k}\right) \frac{\tau _{1}}{\varrho }\right\}  \notag
\\
\chi _{\zeta ,\vec{\zeta}_{\left( m\right) }}^{\left( m\right) }\left( 
\mathbf{\mathbf{k}},\vec{k}\right) \partial _{\tau _{1}}\left[ \tilde{u}%
_{1,\zeta ^{\prime }}\left( \mathbf{k}^{\prime }\right) \ldots \tilde{u}%
_{m,\zeta ^{\prime }}\left( \mathbf{k}^{\left( m\right) }\left( \mathbf{k},%
\vec{k}\right) \right) \right] \mathrm{\tilde{d}}^{\left( m-1\right) d}\vec{k%
}d\tau _{1}.  \notag
\end{gather}%
Estimating the denominator by (\ref{fi2}) \ and using (\ref{Yconv}) we
obtain (\ref{NFMF}). Finally, we consider the case when $\zeta =\infty $ and
all remaining indices $\zeta ^{\left( j\right) }$ equal $\pm $. We expand $%
\Pi _{\infty }$ into sum of $\Pi _{n,\zeta }$ as in (\ref{Pinfk}). In this
case $\chi _{\zeta ,\vec{\zeta}_{\left( m\right) }}^{\left( m\right) }\left( 
\mathbf{\mathbf{k}},\vec{k}\right) $ involves a projection $\Pi _{n,\zeta }$
with $n\neq n_{0}$ (the oscillatory integral may involve $N_{h}-1$ terms
with such $n$). For a fixed $n$ the corresponding phase function $\phi
\left( \mathbf{\mathbf{k}},\vec{k}\right) $ takes the form%
\begin{equation*}
\phi \left( \mathbf{\mathbf{k}},\vec{k}\right) =\phi _{n,\zeta ,\vec{n},\vec{%
\zeta}}\left( \mathbf{\mathbf{k}},\vec{k}\right) =\zeta \omega _{n}\left( 
\mathbf{k}\right) -\zeta ^{\prime }\omega _{n_{0}}\left( \mathbf{k}^{\prime
}\right) -\ldots -\zeta ^{\left( m\right) }\omega _{n_{0}}\left( \mathbf{k}%
^{\left( m\right) }\right) .
\end{equation*}%
Using again (\ref{zomomz}) (now with $n\neq n_{0}$) we obtain \ that (\ref%
{zom}) holds. This implies (\ref{fi2}) as above provided $\delta _{NFM}$ is
small enough. Hence, the relation (\ref{NFMintbp}) holds, implying readily
the desired bound (\ref{NFMF}).
\end{proof}

\subsubsection{FM and NFM monomials for SI oscillatory integrals}

The following below theorem shows that NFM monomials are of the order $%
O\left( \varrho \right) $ as $\varrho \rightarrow 0$. We begin first with
the following statement.

\begin{lemma}
\label{Lemma NFMmin} Assume that Condition \ref{Condition generic1} holds.
Let a monomial $S=\mathcal{F}_{\zeta }^{\left( s\right) }\left( M_{1,\zeta
^{\left( 1\right) }}\ldots M_{s,\zeta ^{\left( s\right) }}\right) $ have all
submonomials $M_{1,\zeta ^{\left( 1\right) }}\ldots M_{s,\zeta ^{\left(
s\right) }}$ which satisfy FM condition (\ref{FMz}), but $S$ itself is not
FM. Assume that $S$ is applied to wavepackets $\mathbf{h}_{l}$ which satisfy
Definition \ref{dwavepack}. and%
\begin{equation}
\mathbf{\tilde{h}}_{l,\zeta }\left( \mathbf{k},\beta \right) =0\text{ if }%
\left\vert \mathbf{k}-\zeta \mathbf{k}_{\ast l}\right\vert \geq \pi
_{0}\beta ^{1-\epsilon },\zeta =\pm .  \label{hl0}
\end{equation}%
Then 
\begin{gather}
\left\Vert S\right\Vert _{E}\leq \frac{4\varrho \left\Vert \chi ^{\left(
s\right) }\right\Vert }{\left\vert \omega \left( \mathbf{\mathbf{k}}_{\ast
}\right) \right\vert }C_{\Xi }^{2s+1}\dprod\limits_{j}\left\Vert M_{j,\zeta
^{\left( j\right) }}\right\Vert _{E}+  \label{SNFM} \\
+\frac{4\varrho \tau _{\ast }\left\Vert \chi ^{\left( s\right) }\right\Vert 
}{\left\vert \omega \left( \mathbf{\mathbf{k}}_{\ast }\right) \right\vert }%
C_{\Xi }^{2s+1}\sum_{i=1}^{s}\left\Vert \partial _{\tau }M_{i,\zeta ^{\left(
i\right) }}\right\Vert _{E}\dprod\limits_{j\neq i}\left\Vert M_{j,\zeta
^{\left( j\right) }}\right\Vert _{E},\;E=C\left( \left[ 0,\tau _{\ast }%
\right] ,L_{1}\right) .  \notag
\end{gather}
\end{lemma}

\begin{proof}
Since $M_{1,\zeta ^{\left( 1\right) }}\ldots M_{s,\zeta ^{\left( s\right) }}$
are decorated FM \ submonomials we can use Lemma \ref{Lemma Fsupport} \ and
Corollary \ref{Corollary 100}. Applying Corollary \ref{Corollary 100} \ and
using (\ref{j0}) we obtain that 
\begin{equation}
M_{l,\zeta ^{\left( l\right) }}\left( \mathbf{k}^{\left( l\right) },\tau
_{1}\right) =0\text{ \ when }\left\vert \mathbf{k}^{\left( l\right) }-\zeta
^{\left( l\right) }\mathbf{k}_{\ast }\right\vert >\nu \left( M_{l,\zeta
^{\left( l\right) }}\right) \beta ^{1-\epsilon }\pi _{0},l=1,\ldots ,s.
\label{Meq00}
\end{equation}%
where $\nu \left( M\right) $ is homogeneity index of $M$. Consider now the
oscillatory integral (\ref{Fm}) which determines $S$, namely 
\begin{gather}
\mathcal{F}_{\zeta ,\vec{\zeta}}^{\left( s\right) }\left( M_{1,\zeta
^{\left( 1\right) }}\ldots M_{s,\zeta ^{\left( s\right) }}\right) \left( 
\mathbf{k},\tau \right) =\int_{0}^{\tau }\int_{\mathbb{D}_{s}}\exp \left\{ 
\mathrm{i}\phi _{\ \zeta ,\vec{\zeta}}\left( \mathbf{\mathbf{k}},\vec{k}%
\right) \frac{\tau _{1}}{\varrho }\right\}  \label{FS} \\
\chi _{\ \zeta ,\vec{\zeta}}^{\left( s\right) }\left( \mathbf{\mathbf{k}},%
\vec{k}\right) M_{1,\zeta ^{\left( 1\right) }}\left( \mathbf{k}^{\prime
},\tau _{1}\right) \ldots M_{s,\zeta ^{\left( s\right) }}\left( \mathbf{k}%
^{\left( s\right) }\left( \mathbf{\mathbf{k}},\vec{k}\right) ,\tau
_{1}\right) \mathrm{\tilde{d}}^{\left( s-1\right) d}\vec{k}\mathrm{d}\tau
_{1}.  \notag
\end{gather}%
We apply Lemma \ref{Lemma NFM} where, according to (\ref{Meq00}) and (\ref%
{sumnu}) $\delta _{NFM}=m\beta ^{1-\epsilon }\pi _{0}$. According to (\ref%
{NFMF}) \ 
\begin{gather}
\left\Vert S\right\Vert _{E}=\left\Vert \mathcal{F}_{\zeta ,\vec{\zeta}%
_{\left( s\right) }}^{\left( s\right) }\left( M_{1,\zeta ^{\left( 1\right)
}}\ldots M_{s,\zeta ^{\left( s\right) }}\right) \left( \mathbf{k},\tau
\right) \right\Vert _{E}\leq \frac{4\varrho \left\Vert \chi ^{\left(
s\right) }\right\Vert }{\left\vert \omega \left( \mathbf{\mathbf{k}}_{\ast
}\right) \right\vert }C_{\Xi }^{2s+1}\dprod\limits_{j}\left\Vert M_{j,\zeta
^{\left( j\right) }}\right\Vert _{E}  \label{normS} \\
+\frac{4\varrho \tau _{\ast }\left\Vert \chi ^{\left( s\right) }\right\Vert 
}{\left\vert \omega \left( \mathbf{\mathbf{k}}_{\ast }\right) \right\vert }%
C_{\Xi }^{2s+1}\sum_{i=1}^{s}\left\Vert \partial _{\tau }M_{i,\zeta ^{\left(
i\right) }}\right\Vert _{E}\dprod\limits_{j\neq i}\left\Vert M_{j,\zeta
^{\left( j\right) }}\right\Vert _{E},\;E=C\left( \left[ 0,\tau _{\ast }%
\right] ,L_{1}\right) ,  \notag
\end{gather}%
that implies (\ref{SNFM}).
\end{proof}

\begin{theorem}
\label{Theorem NFM} Suppose that (i) the inequalities (\ref{omne0}) hold;
(ii) $\mathbf{\tilde{h}}_{l}$ are wavepackets in the sense of Definition \ref%
{dwavepack}; (iii) the relations (\ref{hl0}) hold; (iv) the projections are
defined by (\ref{PPosc}); (v) Condition \ref{Condition generic1} holds.\
Then a NFM decorated monomial based on oscillatory integrals $\mathcal{F}$
defined by (\ref{Fu}) satisfies the estimate \ 
\begin{gather}
\left\Vert M\left( \mathcal{F},T,\vec{\lambda},\vec{\zeta}\right) \left( 
\mathbf{\tilde{h}}_{1}\ldots \mathbf{\tilde{h}}_{m}\right) \right\Vert
_{C\left( \left[ 0,\tau _{\ast }\right] ,L_{1}\right) }\leq  \label{GmRm} \\
\frac{4\varrho \tau _{\ast }^{i-1}\left[ 1+m\right] }{\left\vert \omega
\left( \mathbf{\mathbf{k}}_{\ast }\right) \right\vert }C_{\Xi
}^{2i+e}\dprod\limits_{N\in T,r\left( N\right) >0}\left\Vert \chi ^{\left(
\mu \left( N\right) \right) }\right\Vert \dprod\limits_{l=1}^{m}\left\Vert 
\mathbf{\tilde{h}}_{l,\zeta ^{\left( l\right) }}\right\Vert _{C\left( \left[
0,\tau _{\ast }\right] ,L_{1}\right) },  \notag
\end{gather}%
where $i,m$ and $e$ are respectively the incidence number, the homogeneity
index and the number of edges of $T$.
\end{theorem}

\begin{proof}
Let $M\left( \mathcal{F},T,\vec{\lambda}_{\left( q\right) },\vec{\zeta}%
_{\left( m\right) }\right) \left( \mathbf{\tilde{h}}_{1}\ldots \mathbf{%
\tilde{h}}_{m}\right) $ be a NFM decorated $m$-homogenious monomial. We find
a decorated submonomial $S=M\left( \mathcal{F},T\left( N_{0}\right) ,\vec{%
\lambda}_{\left( q\right) },\vec{\zeta}_{\left( m\right) }\right) $ of $%
M\left( \mathcal{F},T,\vec{\lambda}_{\left( q\right) },\vec{\zeta}_{\left(
m\right) }\right) $ with such $N_{0}$ that $S$ \ is NFM \ \ and has minimal
rank of all NFM submonomials. We denote by $r_{0}$ the rank of $S$, by $%
i^{\prime }$ its incidence number \ and by $s=\nu \left( S\right) =\nu
\left( T\left( N_{0}\right) \right) $ its homogeneity index. This monomial
has the form $S=\mathcal{F}_{\zeta }^{\left( s\right) }\left( M_{1,\zeta
^{\left( 1\right) }}\ldots M_{s,\zeta ^{\left( s\right) }}\right) $. Since
the rank is minimal, all decorated submonomials $M_{l,\zeta ^{\left(
l\right) }}$ \emph{are FM} and their ranks do not exceed $r_{0}-1$. Then
according to (\ref{sumhom}) their homogeneity indices satisfy 
\begin{equation}
\nu \left( M_{1,\zeta ^{\left( 1\right) }}\right) +\ldots +\nu \left(
M_{s,\zeta ^{\left( s\right) }}\right) =s\leq m.  \label{sumnu}
\end{equation}%
Applying Lemma \ref{Lemma NFMmin} we obtain (\ref{SNFM}). Now we use Lemma %
\ref{Lemma 8} \ and \ref{Corollary normM}. Applying Lemma \ref{Norm
submonomial} we obtain \ 
\begin{gather*}
\left\Vert M\left( \left\{ \mathcal{F}\right\} ,T,\Gamma \right) \left( 
\mathbf{\tilde{h}}_{1}\ldots \mathbf{\tilde{h}}_{m}\right) \right\Vert
_{E}\leq \\
\left\Vert S\right\Vert _{E}\dprod\limits_{N\in T\setminus T^{\prime }\left(
N_{0}\right) ,r\left( N\right) >0}\left\Vert \mathcal{F}_{\Gamma \left(
N\right) }^{\left( \mu \left( N\right) \right) }\right\Vert
\dprod\limits_{l<\varkappa \ }\left\Vert \mathbf{\tilde{h}}_{l,\zeta
^{\left( l\right) }}\right\Vert _{E}\dprod\limits_{l\geq \varkappa +\nu
\left( T^{\prime }\left( N_{0}\right) \right) \ }\left\Vert \mathbf{\tilde{h}%
}_{l,\zeta ^{\left( l\right) }}\right\Vert _{E}.
\end{gather*}%
Note that the norm of $\left\Vert \mathcal{F}_{\Gamma \left( N\right)
}^{\left( \mu \left( N\right) \right) }\right\Vert $ is estimated by (\ref%
{dtf}) and norm of $S$ by (\ref{SNFM}). In turn, we estimate right-hand side
of (\ref{SNFM}) using (\ref{dtf}) and (\ref{dtf1}) Taking into account that $%
s\leq m$ in\ the sum in (\ref{normS}) we get the estimate (\ref{GmRm}).
\end{proof}

We also consider the case when Condition \ref{Condition generic1} does not
hold and Condition \ref{Condition generic2} holds. In this case we give an
alternative definition of FM\ and NFM\ decorated monomials.

\begin{definition}[Alternative Frequency Matching]
\label{Definition AFM} We call a decorated composition monomial $M\left( 
\mathcal{F},T,\vec{\lambda},\vec{\zeta}\right) $ \emph{alternatively} \emph{%
frequency matched} (AFM) if (i) every node of $T$ has an odd number of child
nodes (at least three); (ii) for every non-end node $N\in T$ the
corresponding decorated submonomial $M^{\prime }\left( \mathcal{F},T\left(
N\right) ,\vec{\lambda},\vec{\zeta}\right) =\mathcal{F}_{\lambda }^{\left(
m^{\prime }\right) }\left( M_{1,\zeta ^{\prime }}\ldots M_{m^{\prime },\zeta
^{\left( m^{\prime }\right) }}\right) $ satisfies (\ref{neinf}) and%
\begin{equation}
\text{sign}\left( \sum_{j=1}^{m^{\prime }}\zeta ^{\left( j\right) }\right)
=\lambda ,  \label{signlam}
\end{equation}%
where $\lambda ,\zeta ^{\left( j\right) }\in \Lambda $ defined by (\ref{PPPI}%
), we identify $\pm $ with $\pm 1$. A decorated composition monomial which
is not AFM is called \emph{alternatively\ not frequency matched} (ANFM)
monomial.
\end{definition}

Now we prove a statement analogous to Theorem \ref{Theorem NFM} \ when
Condition \ref{Condition generic2} holds.

\begin{theorem}
\label{Theorem ANFM} \ Assume that assumptions of Theorem NFM hold with
Condition \ref{Condition generic1} \ replaced by Condition \ref{Condition
generic2}. Then (\ref{GmRm}) holds.
\end{theorem}

\begin{proof}
According to Corollary \ref{Corollary 100},\ if $\mathbf{\tilde{h}}%
_{l_{1}}=\ldots =\mathbf{\tilde{h}}_{l_{m}}=\mathbf{\tilde{h}}_{l}$ satisfy
Definition \ref{dwavepack} and (\ref{hl0}), then $M\left( \mathcal{F},T,\vec{%
\lambda},\vec{\zeta}\right) \left( \mathbf{\tilde{h}}_{l_{1}}\ldots \mathbf{%
\tilde{h}}_{l_{m}}\right) =M\left( \mathcal{F},T,\vec{\lambda},\vec{\zeta}%
\right) \left( \mathbf{\tilde{h}}_{l_{1},\zeta ^{\prime }}\ldots \mathbf{%
\tilde{h}}_{l_{m},\zeta ^{\left( m\right) }}\right) $ has support in a $%
m\beta ^{1-\epsilon }\ $ vicinity of \ \ $\mathbf{k}_{\vec{\zeta}}=\nu 
\mathbf{k}_{\ast }$ defined by (\ref{kz}), $\nu $ and $m$ are odd integers,$%
m\geq 3$, $\nu =\zeta ^{\prime }+\ldots +\zeta ^{\left( m\right) }$. Let $%
S=M\left( \mathcal{F},T^{\prime },\vec{\lambda},\vec{\zeta}\right) $ be
minimal ANFM submonomial of $M\left( \mathcal{F},T,\vec{\lambda},\vec{\zeta}%
\right) $, that is if $T^{\prime \prime }\subset T^{\prime }$ then \ $%
M\left( \mathcal{F},T^{\prime \prime },\vec{\lambda},\vec{\zeta}\right) $ is
AFM\ submonomial of $S$. The monomial $S$ has the form of (\ref{FS}) with
the interaction phase function%
\begin{equation}
\phi _{\zeta ,\vec{\zeta}}\left( \mathbf{\mathbf{k}},\vec{k}\right) =\zeta
\omega _{n}\left( \mathbf{k}\right) -\zeta ^{\prime }\omega _{n_{0}}\left( 
\mathbf{k}^{\prime }\right) -\ldots -\zeta ^{\left( s\right) }\omega
_{n_{0}}\left( \mathbf{k}^{\left( s\right) }\right) .  \label{fiz}
\end{equation}%
The integrand is non-zero near $\ \mathbf{k}^{\left( j\right) }=\nu _{j}%
\mathbf{k}_{\ast },$and applying (\ref{signlam}) to every AFM submonomial we
get 
\begin{equation}
\zeta ^{\left( l\right) }=\text{sign}\left( \nu _{l}\right) .  \label{zsign}
\end{equation}%
Using (\ref{nd23}) and (\ref{omeven}) we obtain 
\begin{gather}
\phi _{\zeta ,\vec{\zeta}}\left( \nu \mathbf{k}_{\ast },\vec{k}_{\ast
}\right) =\zeta \omega _{n}\left( \nu \mathbf{k}_{\ast }\right) -\zeta
^{\prime }\omega _{n_{0}}\left( \nu _{1}\mathbf{k}_{\ast }\right) -\ldots
-\zeta ^{\left( s\right) }\omega _{n_{0}}\left( \nu _{s}\mathbf{k}_{\ast
}\right)  \label{fiAFM} \\
=\zeta \omega _{n}\left( \nu \mathbf{k}_{\ast }\right) -\text{sign}\left(
\nu _{1}\right) \left\vert \nu _{1}\right\vert \omega _{n_{0}}\left( \mathbf{%
k}_{\ast }\right) -\ldots -\text{sign}\left( \nu _{s}\right) \left\vert \nu
_{s}\right\vert \omega _{n_{0}}\left( \mathbf{k}_{\ast }\right)  \notag \\
=\zeta \left\vert \nu \right\vert \omega _{n}\left( \mathbf{k}_{\ast
}\right) -\left( \nu _{1}+\ldots +\nu _{s}\right) \omega _{n_{0}}\left( 
\mathbf{k}_{\ast }\right) ,\;\nu =\nu _{1}+\ldots +\nu _{s}.  \notag
\end{gather}%
Therefore, since $S$ is ANFM, $\zeta \neq $sign$\left( \nu \right) $ and
since $\ \nu $ is odd, 
\begin{equation}
\phi _{\zeta ,\vec{\zeta}}\left( \nu \mathbf{k}_{\ast },\vec{k}_{\ast
}\right) =-2\nu \omega _{n_{0}}\left( \mathbf{k}_{\ast }\right) \neq 0,
\label{fineq0}
\end{equation}%
therefore (\ref{zom}) holds. We can repeat the proofs of Lemma \ref{Lemma
NFM} and Lemma \ref{Lemma NFMmin} and obtain (\ref{SNFM}). From (\ref{SNFM})
\ we obtain (\ref{GmRm}) \ as in the proof of Theorem \ref{Theorem NFM}.
\end{proof}

Below we give estimations for the derivatives with respect to $\mathbf{k}$
of a composition monomial applied to a wavepacket. Note that (\ref{hbder}) \
admits a singular dependence on $\beta $ of wavepackets $\tilde{h}_{\zeta
}\left( \beta ,\mathbf{k}\right) $. This type of dependence also naturally
comes from explicit formulas as (\ref{h0}) which yield that the first
derivative with respect to $\mathbf{k}$ has a factor $\beta ^{-1}$. Below we
estimate dependence on $\beta $ of monomials applied to wavepackets and will
show that they have the same type of singularity.

Observe that by (\ref{kpi0}) all the points $\mathbf{k}_{\ast l}$ are at the
distance at least $\ 2\pi _{0}$ from $\sigma $. Hence, according to
Definition \ref{Definition band-crossing point}, and (\ref{gradchi}) 
\begin{equation}
\max_{\left\vert \mathbf{k\pm k}_{\ast l}\right\vert \leq \pi _{0},\
l=1,\ldots ,N_{h},}\left( \left\vert \nabla _{\mathbf{k}}^{2}\omega
\right\vert +\left\vert \nabla _{\mathbf{k}}\omega \right\vert \right) \leq
C_{\omega ,2},\   \label{Com2}
\end{equation}%
\begin{equation}
\max_{\left\vert \mathbf{k\pm k}_{\ast l}\right\vert \leq \pi _{0},\
l=1,\ldots ,N_{h}}\left\vert \nabla \chi _{\zeta ,\vec{\zeta}}^{\left(
m\right) }\left( \mathbf{\mathbf{k}},\mathbf{k}^{\prime },\ldots ,\mathbf{k}%
^{\left( m\right) }\right) \right\vert \leq C_{\chi }C_{\Xi }^{m+1}.
\label{grchi}
\end{equation}%
The following seemingly technical Lemma describes a very important property
of solutions. It shows that the $\mathbf{k}$\textbf{-}gradient of solutions
behaves, roughly speaking, as the gradient of initial data. Corresponding
estimates play a crucial role in the control of smallness of interaction of
different wavepackets.

\begin{lemma}
\label{Lemma gradF}\ Let $M\left( \mathcal{F},T,\vec{\lambda},\vec{\zeta}%
\right) \left( \mathbf{\tilde{h}}_{l_{1}}\ldots \mathbf{\tilde{h}}%
_{l_{m}}\right) $ be a decorated monomial which is SI. Assume that $\mathbf{%
\tilde{h}}_{l_{j}}=\mathbf{\tilde{h}}_{l}$ are wavepackets satisfying
Definition \ref{dwavepack}, (\ref{hl0}) and (\ref{Nhh}), that (\ref{scale1})
holds \ and 
\begin{equation}
\beta ^{1-\epsilon }m\leq \pi _{0}.  \label{delm}
\end{equation}%
Assume that either Condition \ref{Condition generic1} holds \ and the
monomial is FM or Condition \ref{Condition generic2} holds \ and the
monomial is AFM. Then 
\begin{equation}
\left\Vert \nabla _{\mathbf{k}}M\left( \mathcal{F},T,\vec{\lambda},\vec{\zeta%
}\right) \left( \mathbf{\tilde{h}}_{l_{1}}\ldots \mathbf{\tilde{h}}%
_{l_{m}}\right) \right\Vert _{E}\leq CC_{\chi }\tau _{\ast }^{i}C_{\Xi
}^{2i+e}C_{\chi }^{i-1}R^{m-1}\beta ^{-1-\epsilon }m^{2},  \label{gradk}
\end{equation}%
where $E=C\left( \left[ 0,\tau _{\ast }\right] ,L_{1}\right) $, $\tau _{\ast
}\leq 1,$ with $i=i\left( T\right) $ and $e=e\left( T\right) $ being
respectively the incidence number and the number of edges of $T$.
\end{lemma}

\begin{proof}
We use the induction with respect to the incidence number $i$ of a tree $T$.
First, we consider the case when Condition \ref{Condition generic1} holds
and $M\left( \mathcal{F},T,\vec{\lambda},\vec{\zeta}\right) $ is FM. For $i=0
$ (\ref{gradk}) \ follows from (\ref{hbder}). Now we assume that (\ref{gradk}%
) \ holds for the incidence number less than $i$ and prove it when the
incidence number equals $i$. Since arguments of $M\left( \mathcal{F},T,\vec{%
\lambda},\vec{\zeta}\right) \ $ are SI, according Definition \ref{Definition
GVM} $l_{1}=\ldots =l_{m}=l$. It is sufficient to prove the boundedness of 
\begin{equation*}
\left\Vert \nabla _{\mathbf{k}}M\left( \mathcal{F},T,\vec{\lambda},\vec{\zeta%
}\right) \left( \mathbf{\tilde{h}}_{l}^{m}\right) \right\Vert
_{E}=\left\Vert \nabla _{\mathbf{k}}\mathcal{F}_{\lambda }^{\left( s\right)
}\left( M_{1,\zeta ^{\prime }}\ldots M_{s,\zeta ^{\left( s\right) }}\right)
\right\Vert _{E},
\end{equation*}%
where $M_{1}\ldots M_{s}$ are decorated submonomials, $M_{j,\zeta }=\Pi
_{\zeta }M_{j}$. Let the submonomials have incidence numbers $i_{1},\ldots
,i_{s}$ and homogeneities $m_{1},\ldots ,m_{s}$ respectively satisfying 
\begin{equation}
i_{1}+\ldots +i_{s}=i-1,\ m_{1}+\ldots +m_{s}=m.  \label{rm}
\end{equation}%
We have by (\ref{Fm}) 
\begin{gather}
\nabla _{\mathbf{k}}\mathcal{F}_{\lambda }^{\left( s\right) }\left(
M_{1,\zeta ^{\prime }}\ldots M_{s,\zeta ^{\left( s\right) }}\right) \left( 
\mathbf{k},\tau \right) =\nabla _{\mathbf{k}}\int_{0}^{\tau }\int_{\left[
-\pi ,\pi \right] ^{\left( s-1\right) d}}\exp \left\{ \mathrm{i}\phi _{\
\lambda ,\vec{\zeta}}\left( \mathbf{\mathbf{k}},\vec{k}\right) \frac{\tau
_{1}}{\varrho }\right\}   \label{gradF} \\
\chi _{\lambda ,\vec{\zeta}}^{\left( s\right) }\left( \mathbf{\mathbf{k}},%
\vec{k}\right) M_{1,\zeta ^{\prime }}\left( \mathbf{k}^{\prime }\right)
\ldots M_{s,\zeta ^{\left( s\right) }}\left( \mathbf{k}^{\left( s\right)
}\left( \mathbf{k},\vec{k}\right) \right) \,\mathrm{\tilde{d}}^{\left(
s-1\right) d}\vec{k}d\tau _{1}.  \notag
\end{gather}%
By Leibnitz formula%
\begin{equation}
\nabla _{\mathbf{k}}\mathcal{F}_{\lambda }^{\left( s\right) }\left(
M_{1,\zeta ^{\prime }}\ldots M_{s,\zeta ^{\left( s\right) }}\right) \left( 
\mathbf{k},\tau \right) =I_{1}+I_{2}+I_{3},  \label{Leib1}
\end{equation}%
where%
\begin{gather*}
I_{1}=\int_{0}^{\tau }\int_{\left[ -\pi ,\pi \right] ^{\left( s-1\right) d}}%
\left[ \nabla _{\mathbf{k}}\exp \left\{ \mathrm{i}\phi _{\ \lambda ,\vec{%
\zeta}}\left( \mathbf{\mathbf{k}},\vec{k}\right) \frac{\tau _{1}}{\varrho }%
\right\} \right]  \\
\chi _{\lambda ,\vec{\zeta}}^{\left( s\right) }\left( \mathbf{\mathbf{k}},%
\vec{k}\right) M_{1,\zeta ^{\prime }}\left( \mathbf{k}^{\prime }\right)
\ldots M_{s,\zeta ^{\left( s\right) }}\left( \mathbf{k}^{\left( s\right)
}\left( \mathbf{k},\vec{k}\right) \right) \,\mathrm{\tilde{d}}^{\left(
s-1\right) d}\vec{k}d\tau _{1},
\end{gather*}%
\begin{gather*}
I_{2}=\int_{0}^{\tau }\int_{\left[ -\pi ,\pi \right] ^{\left( s-1\right)
d}}\exp \left\{ \mathrm{i}\phi _{\ \lambda ,\vec{\zeta}}\left( \mathbf{%
\mathbf{k}},\vec{k}\right) \frac{\tau _{1}}{\varrho }\right\}  \\
\left[ \nabla _{\mathbf{k}}\chi _{\lambda ,\vec{\zeta}}^{\left( s\right)
}\left( \mathbf{\mathbf{k}},\vec{k}\right) \right] M_{1,\zeta ^{\prime
}}\left( \mathbf{k}^{\prime }\right) \ldots M_{s,\zeta ^{\left( s\right)
}}\left( \mathbf{k}^{\left( s\right) }\left( \mathbf{k},\vec{k}\right)
\right) \,\mathrm{\tilde{d}}^{\left( s-1\right) d}\vec{k}d\tau _{1},
\end{gather*}%
\begin{gather*}
I_{3}=\int_{0}^{\tau }\int_{\left[ -\pi ,\pi \right] ^{\left( s-1\right)
d}}\exp \left\{ \mathrm{i}\phi _{\ \lambda ,\vec{\zeta}}\left( \mathbf{%
\mathbf{k}},\vec{k}\right) \frac{\tau _{1}}{\varrho }\right\}  \\
\chi _{\lambda ,\vec{\zeta}}^{\left( s\right) }\left( \mathbf{\mathbf{k}},%
\vec{k}\right) M_{1,\zeta ^{\prime }}\left( \mathbf{k}^{\prime }\right)
\ldots \nabla _{\mathbf{k}}M_{s,\zeta ^{\left( s\right) }}\left( \mathbf{k}%
^{\left( s\right) }\left( \mathbf{k},\vec{k}\right) \right) \,\mathrm{\tilde{%
d}}^{\left( s-1\right) d}\vec{k}d\tau _{1}.
\end{gather*}%
By (\ref{gq}) 
\begin{equation}
\left\Vert M_{j,\zeta ^{\left( j\right) }}\left( \mathbf{k}^{\left( j\right)
}\right) \right\Vert _{L_{1}}\leq C\tau ^{i^{\left( j\right) }}C_{\Xi
}^{2i^{\left( j\right) }+e^{\left( j\right) }}C_{\chi }^{i^{\left( j\right)
}}R^{m_{j}},j=1,\ldots ,s.  \label{ucr}
\end{equation}%
Using (\ref{Yconv}), (\ref{ucr}), (\ref{rm}) and the induction assumption we
get%
\begin{equation}
\left\vert I_{3}\right\vert \leq \left\Vert \chi ^{\left( s\right)
}\right\Vert \dprod\limits_{j=1}^{s-1}\left\Vert M_{j,\zeta ^{\left(
j\right) }}\left( \mathbf{k}^{\left( j\right) }\right) \right\Vert
_{E}\int_{0}^{\tau }\left\Vert \nabla _{\mathbf{k}}M_{s,\zeta ^{\left(
s\right) }}\right\Vert _{E}d\tau _{1}\leq CC_{1}^{m}R^{m-1}\tau ^{i}C_{\Xi
}^{2i+e}C_{\chi }^{i}\beta ^{-1-\epsilon }.  \label{I3}
\end{equation}%
From (\ref{ucr}) and the smoothness of $\chi ^{\left( s\right) }\left( 
\mathbf{\mathbf{k}},\vec{k}\right) $ we get 
\begin{equation}
\left\vert I_{2}\right\vert \leq C\beta ^{-1-\epsilon }\tau
^{i}C_{1}^{m}C_{\Xi }^{2i+e}C_{\chi }^{i}R^{m}.  \label{I2}
\end{equation}%
\bigskip Now we estimate $I_{1}$. Using (\ref{phim}) we obtain%
\begin{gather}
I_{1}=\int_{0}^{\tau }\int_{\left[ -\pi ,\pi \right] ^{\left( s-1\right) d}}%
\left[ \exp \left\{ \mathrm{i}\phi _{\ \lambda ,\vec{\zeta}}\left( \mathbf{%
\mathbf{k}},\vec{k}\right) \frac{\tau _{1}}{\varrho }\right\} \right] 
\label{I1int} \\
\frac{\tau _{1}}{\varrho }\left[ -\lambda \nabla _{\mathbf{k}}\omega \left( 
\mathbf{k}\right) +\zeta ^{\left( s\right) }\nabla _{\mathbf{k}}\omega
\left( \mathbf{k}^{\left( s\right) }\left( \mathbf{k},\vec{k}\right) \right) %
\right]   \notag \\
\chi _{\lambda ,\vec{\zeta}}^{\left( s\right) }\left( \mathbf{\mathbf{k}},%
\vec{k}\right) M_{1,\zeta ^{\prime }}\left( \mathbf{k}^{\prime }\right)
\ldots M_{s,\zeta ^{\left( s\right) }}\left( \mathbf{k}^{\left( s\right)
}\left( \mathbf{k},\vec{k}\right) \right) \mathrm{\tilde{d}}^{\left(
s-1\right) d}\vec{k}d\tau _{1}.  \notag
\end{gather}%
The difficulty in the estimation of the integral $I_{1}$ comes from the
factor $\frac{\tau _{1}}{\varrho }$ since $\varrho $ is small. \ Note that
according to (\ref{scale1}) $\beta ^{2}/\varrho \leq C$. Since $M\left( 
\mathcal{F},T,\vec{\lambda},\vec{\zeta}\right) $ is FM, its every
submonomial is FM too\ and we can apply to them Corollary \ref{Corollary 100}%
, which yields 
\begin{equation*}
M_{j,\zeta ^{\left( j\right) }}\left( \mathbf{k}^{\left( j\right) }\right) =0%
\text{ for }\left\vert \mathbf{k}^{\left( j\right) }-\zeta ^{\left( j\right)
}\mathbf{k}_{\ast }\right\vert >m_{j}\pi _{0}\beta ^{1-\epsilon
},\;j=1,\ldots ,s.
\end{equation*}%
Hence, it is sufficient to estimate $I_{1}$for$.$%
\begin{equation}
\left\vert \mathbf{k}^{\left( j\right) }-\zeta ^{\left( j\right) }\mathbf{k}%
_{\ast }\right\vert \leq \delta _{1}=m\pi _{0}\beta ^{1-\epsilon }\text{ for
all }j.  \label{allj}
\end{equation}%
According to Lemma \ref{Lemma FMgr}, since $\vec{\lambda},\vec{\zeta}$ are
FM 
\begin{equation}
\nabla _{\mathbf{\mathbf{k}}}\phi _{\lambda ,\vec{\zeta}}\left( \lambda 
\mathbf{\mathbf{k}}_{\ast },\vec{k}_{\ast }\right) =\left[ -\lambda \nabla _{%
\mathbf{k}}\omega \left( \mathbf{\mathbf{k}}_{\ast }\right) +\zeta ^{\left(
s\right) }\nabla _{\mathbf{k}}\omega \left( \left( \mathbf{k}^{\left(
s\right) }\left( \mathbf{\mathbf{k}}_{\ast },\vec{k}_{\ast }\right) \right)
\right) \right] =0.  \label{grk0}
\end{equation}%
Using (\ref{Com2}) we conclude that in a vicinity of $\mathbf{k}_{\ast }$
defined by (\ref{allj}) we have 
\begin{equation*}
\left\vert \left[ -\lambda \nabla _{\mathbf{k}}\omega \left( \mathbf{k}%
\right) +\zeta ^{\left( s\right) }\nabla _{\mathbf{k}}\omega \left( \mathbf{k%
}^{\left( s\right) }\left( \mathbf{k},\vec{k}\right) \right) \right]
\right\vert \leq 2\left( s+1\right) C_{\omega ,2}\delta _{1}.
\end{equation*}%
This yields the estimate%
\begin{equation}
\left\vert I_{1}\right\vert \leq CC_{\Xi }^{2i+e}\tau ^{i}C_{\chi
}^{i}C_{1}^{m}\beta ^{-1-\epsilon }m^{2}R^{m}.  \label{I1}
\end{equation}%
Combining (\ref{I1}), (\ref{I2}) and (\ref{I3}) we obtain (\ref{gradk}) and
the induction is completed. \ Now we consider the case when Condition \ref%
{Condition generic2} holds \ and the monomial is AFM. Note that according to
Corollary \ref{Corollary 100} \ the submonomials $M_{j,\zeta ^{\left(
j\right) }}$ have supports near $\nu _{j}\mathbf{k}_{\ast }$, with an odd $%
\nu _{j}$. By Lemma \ref{Lemma Fsupport} the monomial itself is non-zero
near $\nu \mathbf{k}_{\ast },$ $\mathbf{\ }\nu =\nu _{1}+\ldots +\nu _{s}$;
since $s$ is odd $\nu $ is odd too. Obviously, one of $\nu _{j}$ has the
same sign as $\nu $, we assume that $j=s$, that is 
\begin{equation}
\text{sign}\left( \nu _{s}\right) =\text{sign}\left( \nu _{1}+\ldots +\nu
_{s}\right) =\text{sign}\left( \nu \right) ,  \label{sigsig}
\end{equation}%
the general case can be reduced to this by a relabeling of variables. \ The
interaction phase function is given by (\ref{fiz}) and since the
submonomials are AFM (\ref{zsign}) holds. According to (\ref{omeven}) $%
\nabla _{\mathbf{k}}\left( \omega \left( -\mathbf{k}\right) \right) =-\left(
\nabla _{\mathbf{k}}\omega \right) \left( \mathbf{k}\right) $. Therefore,
using (\ref{nd22}) we obtain 
\begin{gather*}
\nabla _{\mathbf{\mathbf{k}}}\phi _{\lambda ,\vec{\zeta}}\left( \nu \mathbf{%
\mathbf{k}}_{\ast },\vec{k}_{\ast }\right) =\lambda \nabla _{\mathbf{k}%
}\omega \left( \nu \mathbf{k}_{\ast }\right) -\zeta ^{\left( s\right)
}\nabla _{\mathbf{k}}\omega \left( \nu _{s}\mathbf{k}_{\ast }\right) = \\
\lambda \left( \nabla _{\mathbf{k}}\omega \right) \left( \text{sign}\left(
\nu \right) \left\vert \nu \right\vert \mathbf{k}_{\ast }\right) -\zeta
^{\left( s\right) }\nabla _{\mathbf{k}}\omega \left( \text{sign}\left( \nu
_{s}\right) \left\vert \nu _{s}\right\vert \mathbf{k}_{\ast }\right)  \\
=\lambda \left( \nabla _{\mathbf{k}}\omega \right) \left( \text{sign}\left(
\nu \right) \mathbf{k}_{\ast }\right) -\zeta ^{\left( s\right) }\nabla _{%
\mathbf{k}}\omega \left( \text{sign}\left( \nu _{s}\right) \mathbf{k}_{\ast
}\right) =\left( \lambda \text{sign}\left( \nu \right) -\zeta ^{\left(
s\right) }\text{sign}\left( \nu _{s}\right) \right) \left( \nabla _{\mathbf{k%
}}\omega \right) \left( \mathbf{k}_{\ast }\right) .
\end{gather*}%
Using (\ref{zsign}) \ we conclude that 
\begin{equation}
\nabla _{\mathbf{\mathbf{k}}}\phi _{\lambda ,\vec{\zeta}}\left( \nu \mathbf{k%
}_{\ast },\vec{k}_{\ast }\right) =0,\;\vec{k}_{\ast }=\left( \nu _{1}\mathbf{%
k}_{\ast },\ldots ,\nu _{s}\mathbf{k}_{\ast }\right) .  \label{fi0}
\end{equation}%
Using (\ref{fi0}) \ instead of (\ref{grk0}) we conclude as in the first half
of the proof that (\ref{gradk}) holds in the AFM case too.
\end{proof}

\subsection{Properties of minimal CI monomials}

Here we consider CI evaluated monomials with arguments involving different
wavepackets $\mathbf{\tilde{h}}_{l}$. Since the group velocities of
wavepackets are different, namely (\ref{NGVM}) is satisfied, there exists $%
p_{0}>0$ such that 
\begin{equation}
\left\vert \nabla \omega \left( \mathbf{k}_{\ast l_{1}}\right) -\nabla
\omega \left( \mathbf{k}_{\ast l_{2}}\right) \right\vert \geq p_{0}>0\text{
if \ }l_{1}\neq l_{2}.  \label{NGVMp}
\end{equation}%
The next lemma is a standard implication of the Stationary Phase Method in
the case when the phase function has no critical points in the domain of
integration, namely when (\ref{NGVM}) holds.

\begin{lemma}
\label{Lemma GVM0}\textbf{\ }Let $\mathbf{k}_{\ast l}$ \ and $\omega _{n}$
be generic in the sense of Definition \ref{Condition generic22}. Let $%
\mathcal{F}^{\left( m\right) }$ be defined by (\ref{Fu}), $m\left( \beta
\right) $ be as in (\ref{mofrho}). We assume that (\ref{gradchi}) and (\ref%
{NGVM}) hold. We also assume that (\ref{Nhh}), (\ref{hl0}), (\ref{sourloc}),
(\ref{hbder}) and (\ref{scale1}) hold. We assume that $M\left( \mathcal{F}%
,T\right) \left( \mathbf{\tilde{h}}_{l_{1}}\ldots \mathbf{\tilde{h}}%
_{l_{m}}\right) $ is a monomial with homogeneity index $m$ evaluated \ at
arguments with CI multiindex $l_{1},\ldots ,l_{m}$, but every evaluated
submonomial of $M\left( \mathcal{F},T\right) \left( \mathbf{\tilde{h}}%
_{l_{1}}\ldots \mathbf{\tilde{h}}_{l_{m}}\right) $ is SI. Then for $m\leq
m\left( \beta \right) $ and small $\beta $ 
\begin{equation}
\left\Vert M\left( \mathcal{F},T\right) \left( \mathbf{\tilde{h}}%
_{l_{1}}\ldots \mathbf{\tilde{h}}_{l_{m}}\right) \right\Vert _{E}\leq \frac{C%
}{p_{0}}\tau _{\ast }^{i-1}C_{\Xi }^{2i+e}3^{2m}C_{\chi }^{i}\left[ \frac{%
\varrho \left\vert \ln \beta \right\vert }{\beta ^{1+\epsilon }}+\beta %
\right] m^{2}R^{m-1},  \label{NGVMest0}
\end{equation}%
where $i$ and $e$ are respectively the incidence number and number of edges
of $T$, $R$ is as in (\ref{Nhh}).
\end{lemma}

\begin{proof}
Since $\mathbf{k}_{\ast l}$ are not band-crossing points, the relations (\ref%
{grchi}) and\ (\ref{Com2}) hold. \ We expand $M\left( \mathcal{F},T\right) $
into a sum of decorated monomials $M\left( \mathcal{F},T,\vec{\lambda},\vec{%
\zeta}\right) $ as in (\ref{MFsplit}), which contains no more than $%
3^{i\left( T\right) +m}$ terms, and $i\left( T\right) +m\leq 2m$. The
submonomials of every decorated monomial are SI by the assumption of the
theorem. If Condition \ref{Condition generic1} holds, the submonomials are
either FM or NFM; if Condition \ref{Condition generic2} holds, the
submonomials are either AFM or ANFM. If a decorated submonomial $M\left( 
\mathcal{F},T^{\prime },\vec{\lambda}^{\prime },\vec{\zeta}^{\prime }\right) 
$ is NFM we use Theorem \ref{Theorem NFM} \ and obtain from (\ref{GmRm}) the
inequality 
\begin{equation}
\left\Vert M\left( \mathcal{F},T^{\prime },\vec{\lambda}^{\prime },\vec{\zeta%
}^{\prime }\right) \left( \mathbf{\tilde{h}}_{l_{j^{\prime }+1}}\ldots 
\mathbf{\tilde{h}}_{l_{j^{\prime }+m^{\prime }}}\right) \right\Vert _{E}\leq
C\varrho \tau _{\ast }^{i^{\prime }-1}\left[ 1+m\right] C_{\Xi }^{2i^{\prime
}+e^{\prime }}C_{\chi }^{i^{\prime }}R^{m^{\prime }},  \label{MNFM}
\end{equation}%
where $i^{\prime }$ and $e^{\prime }$ are the incidence number and number of
edges of the subtree $T^{\prime }$. Alternatively, if Condition \ref%
{Condition generic2} holds, and .a decorated monomial $M\left( \mathcal{F}%
,T^{\prime },\vec{\lambda}^{\prime },\vec{\zeta}^{\prime }\right) \ $is
ANFM, we use Theorem \ref{Theorem ANFM} \ and obtain from (\ref{GmRm}) the
inequality (\ref{MNFM}). Using (\ref{MNFM}) in both cases \ we obtain 
\begin{equation}
\left\Vert M\left( \mathcal{F},T,\vec{\lambda},\vec{\zeta}\right) \left( 
\mathbf{\tilde{h}}_{l_{1}}\ldots \mathbf{\tilde{h}}_{l_{m}}\right)
\right\Vert _{E}\leq C\varrho \tau _{\ast }^{i-1}C_{\Xi }^{2i+e}C_{\chi
}^{i}mR^{m}.  \label{NGVM1}
\end{equation}%
Now we consider the case when Condition \ref{Condition generic1} holds and
every submonomial of $M\left( \mathcal{F},T,\vec{\lambda},\vec{\zeta}\right) 
$ is FM. We write the integral with respect to $\tau _{1}$ in (\ref{FS}) as
a sum of two integrals from $0$ to $\beta $ and from $\beta $ to $\tau $,
namely 
\begin{gather}
\mathcal{F}_{\zeta ,\vec{\zeta}}^{\left( s\right) }\left( M_{1}\ldots
M_{s}\right) \left( \mathbf{k},\tau \right) =F_{1}+F_{2},  \label{Fbet} \\
F_{1}=\int_{\beta }^{\tau }\int_{\mathbb{D}_{m}}\exp \left\{ \mathrm{i}\phi
_{\zeta ,\vec{\zeta}}\left( \mathbf{\mathbf{k}},\vec{k}\right) \frac{\tau
_{1}}{\varrho }\right\} A_{\zeta ,\vec{\zeta}}^{\left( s\right) }\left( 
\mathbf{k},\vec{k}\right) \mathrm{\tilde{d}}^{\left( s-1\right) d}\vec{k}%
d\tau _{1},F_{2}=\int_{0}^{\beta }\ldots d\tau _{1}  \notag
\end{gather}%
where 
\begin{equation}
A_{\zeta ,\vec{\zeta}_{\left( m\right) }}^{\left( s\right) }\left( \mathbf{k}%
,\vec{k}\right) =\chi _{\zeta ,\vec{\zeta}}^{\left( s\right) }\,\left( 
\mathbf{k},\vec{k}\right) M_{1}\left( \mathbf{k}^{\prime }\right) \ldots
M_{s}\left( \mathbf{k}^{\left( s\right) }\right) ,  \label{Akk}
\end{equation}%
$M_{j}$ are submonomials of $M$. According to Corollary \ref{Corollary normM}
with $\tau _{\ast }=\beta $ 
\begin{equation}
\left\Vert F_{2}\right\Vert _{L_{1}}\leq 2C_{\Xi }^{1+2s}C_{\chi }\beta
\dprod\limits_{j=1}^{s}\left\Vert M_{j}\right\Vert _{E}\leq \beta C_{\Xi
}^{e+2i}\tau _{\ast }^{i-1}C_{\chi }\dprod\limits_{j=1}^{m}\left\Vert 
\mathbf{\tilde{h}}_{l_{j}}\right\Vert _{E}\leq \beta C_{\chi }C_{\Xi
}^{e+2i}\tau _{\ast E}^{i-1}R^{m}.  \label{0beta2}
\end{equation}%
Now we estimate $F_{1}$. Since $M\left( \mathcal{F},T\right) $ is CI, there
are two SI submonomials $M_{j_{1}}$ and $M_{j_{2}}$ applied to $\left( 
\mathbf{\tilde{h}}_{l_{j_{1}}}\right) ^{m_{1}}$ and $\left( \mathbf{\tilde{h}%
}_{l_{j_{2}}}\right) ^{m_{2}}$ \ with $l_{j_{1}}\neq l_{j_{2}}$. Let us
assume that $l_{j_{1}}=l_{1}$, $l_{j_{2}}=l_{s}$ \ (the general case $\ $
can be easily reduced to it by a relabeling of variables). We denote 
\begin{equation}
\mathbf{\phi }^{\prime }=\nabla _{\mathbf{k}^{\prime }}\phi _{\zeta ,\vec{%
\zeta}}\left( \mathbf{\mathbf{k}},\vec{k}_{\ast }\right) =\nabla _{\mathbf{k}%
^{\prime }}\omega \left( \mathbf{k}_{\ast l_{1}}\right) -\nabla _{\mathbf{k}%
^{\left( s\right) }}\omega \left( \mathbf{k}_{\ast l_{s}}\right) \neq 0,%
\mathbf{p}=\mathbf{\phi }^{\prime }/\left\vert \mathbf{\phi }^{\prime
}\right\vert .  \label{fipp}
\end{equation}%
By (\ref{NGVMp}) and (\ref{pmom})\ we obtain%
\begin{equation}
\left\vert \mathbf{p\cdot }\nabla _{\mathbf{k}^{\prime }}\phi _{\zeta ,\vec{%
\zeta}}\left( \mathbf{\mathbf{k}},\vec{k}_{\ast }\right) \right\vert \geq
p_{0}>0\text{ for }\vec{k}=\vec{k}_{\ast }=\left( \mathbf{k}_{\ast
l_{1}},\ldots ,\mathbf{k}_{\ast l_{s}}\right) .  \label{pstar}
\end{equation}%
Note that 
\begin{equation*}
\exp \left\{ \mathrm{i}\phi _{\zeta ,\vec{\zeta}}\left( \mathbf{\mathbf{k}},%
\vec{k}\right) \frac{\tau _{1}}{\varrho }\right\} =\frac{\varrho }{\mathrm{i}%
\mathbf{p\cdot }\nabla _{\mathbf{k}^{\prime }}\phi _{\zeta ,\vec{\zeta}%
}\left( \mathbf{\mathbf{k}},\vec{k}\right) \tau _{1}}\mathbf{p\cdot }\nabla
_{\mathbf{k}^{\prime }}\exp \left\{ \mathrm{i}\phi _{\zeta ,\vec{\zeta}%
}\left( \mathbf{\mathbf{k}},\vec{k}\right) \frac{\tau _{1}}{\varrho }%
\right\} .
\end{equation*}%
Using this identity, (\ref{kkar}) and integrating by parts the integral \
which defines $F_{1}$ \ in (\ref{Fbet}) we obtain 
\begin{gather}
F_{1}=\int_{\beta }^{\tau }I\left( \mathbf{\mathbf{k}},\tau _{1}\right)
d\tau _{1},\;I\left( \mathbf{\mathbf{k}},\tau _{1}\right) =\int_{\mathbb{D}%
_{m}}\exp \left\{ \mathrm{i}\phi _{\zeta ,\vec{\zeta}}\left( \mathbf{\mathbf{%
k}},\vec{k}\right) \frac{\tau _{1}}{\varrho }\right\} A_{\zeta ,\vec{\zeta}%
}^{\left( s\right) }\left( \mathbf{k},\vec{k}\right) \mathrm{\tilde{d}}%
^{\left( s-1\right) d}\vec{k}=  \label{Ik} \\
-\int_{\mathbb{D}_{s}}\frac{\varrho \exp \left\{ \mathrm{i}\phi _{\zeta ,%
\vec{\zeta}}\left( \mathbf{\mathbf{k}},\vec{k}\right) \frac{\tau _{1}}{%
\varrho }\right\} }{\mathrm{i}\tau _{1}}\mathbf{p\cdot }\nabla _{\mathbf{k}%
^{\prime }}\frac{A_{\zeta ,\vec{\zeta}}^{\left( s\right) }\,\left( \mathbf{k}%
,\vec{k}\right) }{\nabla _{\mathbf{k}^{\prime }}\phi _{\zeta ,\vec{\zeta}%
}\left( \mathbf{\mathbf{k}},\vec{k}\right) \cdot \mathbf{p}}\mathrm{\tilde{d}%
}^{\left( s-1\right) d}\vec{k}.  \notag
\end{gather}%
From (\ref{hl0}), \ Lemma \ref{Lemma Fsupport} and Corollary \ref{Corollary
100} we see that in the integral$\;I\left( \mathbf{\mathbf{k}},\tau
_{1}\right) $ the integrands are nonzero only if 
\begin{equation}
\left\vert \mathbf{\mathbf{k}}^{\left( j\right) }-\zeta ^{\left( j\right) }%
\mathbf{\mathbf{k}}_{\ast }^{\left( j\right) }\right\vert \leq m_{j}\pi
_{0}\beta ^{1-\epsilon },\left\vert \mathbf{\mathbf{k}}-\zeta \mathbf{%
\mathbf{k}}_{\ast }\right\vert \leq m\pi _{0}\beta ^{1-\epsilon
},m_{1}+\ldots +m_{s}\leq m,  \label{nzero}
\end{equation}%
where $\pi _{0}\leq 1$. Using the Taylor remainder estimate \ for $\phi
_{\zeta ,\vec{\zeta}}$ at $\vec{k}_{\ast }$ we obtain the inequality 
\begin{equation}
\left\vert \nabla _{\mathbf{k}^{\prime }}\phi _{\zeta ,\vec{\zeta}}\left( 
\mathbf{\mathbf{k}},\vec{k}\right) -\mathbf{\phi }^{\prime }\right\vert \leq
3m\beta ^{1-\epsilon }C_{\omega ,2}\text{ if (\ref{nzero}) holds.\ }
\label{Tayest}
\end{equation}%
Suppose that $\beta $ is small and satisfies%
\begin{equation}
3m\beta ^{1-\epsilon }C_{\omega ,2}\leq \frac{p_{0}}{2}.  \label{ses}
\end{equation}%
Condition \ (\ref{ses}) is satisfied for small $\beta $ if $m\leq m\left(
\beta \right) $ as in (\ref{mofrho}).\ Using (\ref{Tayest}) we derive from (%
\ref{pstar}), (\ref{ses}) and (\ref{hl0}) that \ 
\begin{equation}
\left\vert \mathbf{p\cdot }\nabla _{\mathbf{k}^{\prime }}\phi _{\zeta ,\vec{%
\zeta}}\left( \mathbf{\mathbf{k}},\vec{k}\right) \right\vert \geq \frac{p_{0}%
}{2}>0\text{ if (\ref{nzero}) holds}.  \label{gradom}
\end{equation}%
Now we use (\ref{gradom}) to estimate denominators, \ (\ref{Com2}) \ to
estimate second $\mathbf{k}^{\prime }$-derivatives of $\omega $ \ and (\ref%
{grchi}) to estimate $\nabla _{\mathbf{k}^{\prime }}\chi $. We conclude that 
\begin{gather}
\left\vert I\left( \mathbf{\mathbf{k}},\tau _{1}\right) \right\vert \leq
C_{\Xi }^{2s+1}\int_{\mathbb{D}_{s}}\left[ \frac{\varrho }{\tau _{1}p_{0}}%
\left\vert \nabla _{\mathbf{k}^{\prime }}A_{\zeta ,\vec{\zeta}}^{\left(
s\right) }\left( \mathbf{k},\vec{k}\right) \right\vert +\frac{8\varrho
C_{\omega ,2}}{\tau _{1}p_{0}^{2}}\left\vert A_{\zeta ,\vec{\zeta}_{\left(
m\right) }}^{\left( s\right) }\left( \mathbf{k},\vec{k}\right) \right\vert %
\right] \mathrm{\tilde{d}}^{\left( s-1\right) d}\vec{k}  \label{Ikt1} \\
\leq \frac{\varrho }{\tau _{1}p_{0}}\left[ \left\Vert \left( \nabla
_{k^{\prime }}-\nabla _{k^{\left( s\right) }}\right) \chi ^{\left( s\right)
}\left( \mathbf{\mathbf{k}},\cdot \right) \right\Vert +\frac{8C_{\omega ,2}}{%
p_{0}}\left\Vert \chi ^{\left( m\right) }\left( \mathbf{\mathbf{k}},\cdot
\right) \right\Vert \right] C_{\Xi }^{2s+1}\dprod\limits_{j=1}^{s}\left\Vert
M_{j}\right\Vert _{L_{1}}+  \notag \\
\frac{\varrho C_{\Xi }^{2s+1}\left\Vert \chi ^{\left( s\right) }\left( 
\mathbf{\mathbf{k}},\cdot \right) \right\Vert }{\tau _{1}p_{0}}\left[
\dprod\limits_{j=2}^{s}\left\Vert M_{j}\right\Vert _{L_{1}}\left\Vert \nabla
_{\mathbf{k}^{\prime }}M_{1}\right\Vert
_{L_{1}}+\dprod\limits_{j=1}^{s-1}\left\Vert M_{j}\right\Vert
_{L_{1}}\left\Vert \nabla _{\mathbf{k}^{\left( s\right) }}M_{s}\right\Vert
_{L_{1}}\right] .  \notag
\end{gather}%
To estimate $\nabla M_{i}$ we use Lemma \ref{Lemma gradF}. We also use \ (%
\ref{dtf}) and (\ref{gq}) \ to estimate $\left\Vert M_{j}\right\Vert _{L_{1}}
$. Therefore, using (\ref{rm}), we obtain 
\begin{equation}
\left\vert I\left( \mathbf{\mathbf{k}},\tau _{1}\right) \right\vert \leq 
\frac{C}{\tau _{1}}\tau _{\ast }^{i-1}C_{\Xi }^{2i+e}C_{\chi }^{i}\frac{%
\varrho }{\beta ^{1+\epsilon }p_{0}}m^{2}R^{m-1}.  \label{II}
\end{equation}%
Finally, we consider the case when \ the alternative Condition \ref%
{Condition generic2} holds. In this case $M_{1}$ and $M_{s}$ \ according to
Lemma \ref{Lemma Fsupport} are localized near $\nu _{1}\mathbf{k}_{\ast
l_{1}}$ and $\nu _{2}\mathbf{k}_{\ast l_{s}}$ \ with some $\nu _{1}$ \ and $%
\nu _{2}$; we use (\ref{csg1}) to obtain (\ref{pstar}) both for AFM and ANFM
submonomials. Therefore (\ref{gradom}) holds \ and we again get (\ref{Ikt1})
and (\ref{II}). So, we proved (\ref{II}) in all cases. Integrating (\ref{II}%
) in $\tau _{1}$ we obtain 
\begin{equation}
\left\Vert F_{1}\right\Vert _{E}\leq C\tau _{\ast }^{i-1}C_{\Xi
}^{2i+e}C_{\chi }^{i}\frac{\varrho }{\beta ^{1+\epsilon }p_{0}}%
m^{2}\left\vert \ln \beta \right\vert R^{m-1}.  \label{NGVM2}
\end{equation}%
Using summation over all $\vec{\lambda},\vec{\zeta}$ (the sum involves no
more than $3^{2m}$ terms) we obtain (\ref{NGVMest0}) from (\ref{NGVM1}) and (%
\ref{NGVM2}).
\end{proof}

\section{Proof of the superposition theorems}

In this section we prove Theorems \ref{Theorem Superposition}, \ref{Theorem
Superposition1} on the approximate modal superposition principle.

\subsection{Proof of the Superposition principle for lattice equations}

\ Here we prove Theorem \ref{Theorem Superposition}. First we note that
according to Lemma \ref{Lemma 1minpsi} we can replace $\mathbf{\tilde{h}}%
_{l} $ by $\mathbf{\tilde{h}}_{l}^{\Psi }\mathbf{\ }$in the statement of
Theorem \ref{Theorem Superposition}, in particular in (\ref{Gsum}), (\ref%
{rem}). Hence we can assume that (\ref{hl0}) holds.

Based on Theorem \ref{Imfth1} \ we expand the solution of (\ref{eqFur}) into
series (\ref{uTj}) \ and then into the sum of composition monomials $M\left( 
\mathcal{F},T\right) $ as in (\ref{Gtreeexp}): 
\begin{equation}
\mathcal{G}\left( \mathcal{F},\mathbf{\tilde{h}}\right) =\mathbf{\tilde{h}}%
+\sum_{m=2}^{\infty }\sum_{T\in T_{m}}c_{T}M\left( \mathcal{F},T\right)
\left( \mathbf{\tilde{h}}^{m}\right) ,  \label{Goscmon}
\end{equation}%
where%
\begin{equation}
\mathbf{\tilde{h}}=\sum_{l=1}^{N_{h}}\mathbf{\tilde{h}}_{l},\ \left\Vert 
\mathbf{\tilde{h}}_{l}\right\Vert _{E}\leq R,\ l=1,\ldots ,N_{h},
\label{hRm}
\end{equation}%
and the relation (\ref{Nhh}) \ (that is $N_{h}R<R_{\mathcal{G}}$) holds,
where $R_{\mathcal{G}}$ is the radius of convergence\ from Theorem \ref%
{Imfth1}, \ $R$ will be specified below. Using Lemma \ref{Lemma Truncation}
we conclude that 
\begin{equation}
\mathcal{G}\left( \mathcal{F},\mathbf{\tilde{h}}\right) =\mathbf{\tilde{h}}%
+\sum_{m=2}^{m\left( \beta \right) }\sum_{T\in T_{m}}c_{T}M\left( \mathcal{F}%
,T\right) \left( \mathbf{\tilde{h}}^{m}\right) +g,\ \left\Vert g\right\Vert
_{E}\leq \beta ,  \label{Gtr}
\end{equation}%
where $m\left( \beta \right) $ is defined by (\ref{mofrho}). Then we expand
every monomial $M\left( \mathcal{F},T\right) \left( \mathbf{\tilde{h}}%
^{m}\right) $ according to (\ref{Msum}) into the sum of the terms $M\left( 
\mathcal{F},T\right) \left( \mathbf{\tilde{h}}_{l_{1}}\ldots \mathbf{\tilde{h%
}}_{l_{m}}\right) $. Note that since $m\left( \beta \right) \leq C\left\vert
\ln \beta \right\vert $, conditions (\ref{ses}), (\ref{delm}), (\ref{desmall}%
) are satisfied\ if $\beta $ is small enough \ for every $m\leq m\left(
\beta \right) $. The monomials $M\left( \mathcal{F},T\right) \left( \mathbf{%
\tilde{h}}_{l_{1}}\ldots \mathbf{\tilde{h}}_{l_{m}}\right) $ belong to two
classes, SI and CI (according to Definition \ref{Definition GVM}) and the
class is determined by the multiindex $\left( l_{1},\ldots ,l_{m}\right) =%
\bar{l}$. Using (\ref{Gtr}) we conclude that 
\begin{gather}
\mathcal{G}\left( \mathcal{F},\sum_{l=1}^{N_{h}}\mathbf{\tilde{h}}%
_{l}\right) =\sum_{l=1}^{N_{h}}\mathcal{G}\left( \mathcal{F},\mathbf{\tilde{h%
}}_{l}\right) +\mathbf{\tilde{D}},  \label{GG} \\
\mathbf{\tilde{D}}=\sum_{m=2}^{m\left( \beta \right) }\sum_{T\in T_{m}}\sum_{%
\text{CI\ }l_{1},\ldots ,l_{m}}c_{T}M\left( \mathcal{F},T\right) \left( 
\mathbf{\tilde{h}}_{l_{1}}\ldots \mathbf{\tilde{h}}_{l_{m}}\right) +g_{1},\
\left\Vert g_{1}\right\Vert _{E}\leq C\beta .  \notag
\end{gather}%
To obtain (\ref{rem}) we have to estimate the sum in $\mathbf{\tilde{D}}$
and show that it is small. It follows from (\ref{ctes}) that 
\begin{gather*}
\left\Vert \sum_{m=2}^{m\left( \beta \right) }\sum_{T\in T_{m}}\sum_{\text{CI%
}l_{1},\ldots ,l_{m}}c_{T}M\left( \mathcal{F},T\right) \left( \mathbf{\tilde{%
h}}_{l_{1}}\ldots \mathbf{\tilde{h}}_{l_{m}}\right) \right\Vert _{E}\leq \\
\leq \sum_{m=2}^{m\left( \beta \right) }N_{h}^{m}\sum_{T\in
T_{m}}c_{T}\sup_{T\in T_{m},\text{CI}\bar{l}}\left\Vert M\left( \mathcal{F}%
,T\right) \left( \mathbf{\tilde{h}}_{l_{1}}\ldots \mathbf{\tilde{h}}%
_{l_{m}}\right) \right\Vert _{E}\leq \\
\leq \sum_{m=2}^{m\left( \beta \right) }N_{h}^{m}c_{0}c_{1}^{m}\sup_{T\in
T_{m},\text{CI}\bar{l}}\left\Vert M\left( \mathcal{F},T\right) \left( 
\mathbf{\tilde{h}}_{l_{1}}\ldots \mathbf{\tilde{h}}_{l_{m}}\right)
\right\Vert _{E}.
\end{gather*}%
Now we consider an evaluated monomial $M\left( \mathcal{F},T\right) \left( 
\mathbf{\tilde{h}}_{l_{1}}\ldots \mathbf{\tilde{h}}_{l_{m}}\right) $ with
arguments given by CI \ multiindex $\bar{l}=\left( l_{1},\ldots
,l_{m}\right) $. To prove that this monomial has a small norm, according to
Lemma \ref{Norm submonomial} it is sufficient to show that one of its
submonomials is small and the relevant operators are bounded. According to
Proposition \ref{minsub} the monomial $M\left( \mathcal{F},T\right) \left( 
\mathbf{\tilde{h}}_{l_{1}}\ldots \mathbf{\tilde{h}}_{l_{m}}\right) $
contains a submonomial $M\left( \mathcal{F},T^{\prime }\right) \left( 
\mathbf{\tilde{h}}_{l_{s^{\prime }}}\ldots \mathbf{\tilde{h}}_{l_{s^{\prime
\prime }}}\right) $ with the homogeneity index $s=s^{\prime \prime
}-s^{\prime }+1$, the incidence number $i^{\prime }$ and the rank $r^{\prime
}$ which is minimal in the following sense. The monomial $M\left( \mathcal{F}%
,T^{\prime }\right) \left( \mathbf{\tilde{h}}_{l_{s^{\prime }}}\ldots 
\mathbf{\tilde{h}}_{l_{s^{\prime \prime }}}\right) $ is CI, but every its
submonomial $M\left( \mathcal{F},T^{\prime \prime }\right) \left( \mathbf{%
\tilde{h}}_{l_{s^{\prime \prime }}}\ldots \mathbf{\tilde{h}}_{l_{s^{\prime
\prime \prime }}}\right) $ is SI. Now we use the space decomposition (\ref%
{PPosc}) and expand $M\left( \mathcal{F},T^{\prime }\right) $ as in (\ref%
{Fsumlamz}) into a sum of no more than $3^{2m}$ decorated monomials $M\left( 
\mathcal{F},T^{\prime },\vec{\lambda},\vec{\zeta}\right) \left( \mathbf{%
\tilde{h}}_{l_{s^{\prime }}}\ldots \mathbf{\tilde{h}}_{l_{s^{\prime \prime
}}}\right) $.\ The decorated submonomials of every decorated monomial are
SI. We apply Lemma \ref{Lemma GVM0} \ and conclude that 
\begin{equation}
\left\Vert M\left( \mathcal{F},T^{\prime },\vec{\lambda},\vec{\zeta}\right)
\left( \mathbf{\tilde{h}}_{l_{s^{\prime }}}\ldots \mathbf{\tilde{h}}%
_{l_{s^{\prime \prime }}}\right) \right\Vert _{E}\leq C\left[ \frac{\varrho 
}{\beta ^{1+\epsilon }}\left\vert \ln \beta \right\vert +\beta \right] \frac{%
s^{2}}{p_{0}}\tau _{\ast }^{i^{\prime }-1}C_{\Xi }^{e^{\prime }+2i^{\prime
}}C_{\chi }^{i^{\prime }}R^{s^{\prime \prime }-s^{\prime }}.  \label{FMc}
\end{equation}%
Hence, there is a submonomial of $M\left( \mathcal{F},T\right) \left( 
\mathbf{\tilde{h}}_{l_{1}}\ldots \mathbf{\tilde{h}}_{l_{m}}\right) $ with a
small norm. Namely, since (\ref{scale1}) and (\ref{tau1}) are assumed, this
small submonomial provides the smallness of the norm of the whole monomial $%
M\left( \mathcal{F},T\right) \left( \mathbf{\tilde{h}}_{l_{1}}\ldots \mathbf{%
\tilde{h}}_{l_{m}}\right) $ according to Lemma \ref{Norm submonomial}. We
also use Corollary \ref{Corollary normM} and (\ref{chiCR}) to estimate norms
of remaining submonomials of rank $r$ and apply (\ref{MFr}) and (\ref{rm})
to obtain 
\begin{equation}
\left\Vert M\left( \mathcal{F},T\right) \left( \mathbf{\tilde{h}}%
_{l_{1}}\ldots \mathbf{\tilde{h}}_{l_{m}}\right) \right\Vert \leq 3^{2m}%
\left[ \frac{\varrho }{\beta ^{1+\epsilon }}\left\vert \ln \beta \right\vert
+\beta \right] C_{1}m^{2}\tau _{\ast }^{i-1}C_{\Xi }^{e+2i}C_{\chi
}^{i}R^{m-1}.  \label{NGVMest}
\end{equation}%
Since $e=i+m-1$, using (\ref{mofs}) and the inequalities $i\left( T\right)
=i\geq m/m_{F}$, $i\leq m-1$ we get 
\begin{eqnarray}
&&\sum_{m=2}^{m\left( \beta \right) }\sum_{T\in T_{m}}\sum_{\text{CI}%
l_{1},\ldots ,l_{m}}c_{T}M\left( \mathcal{F},T\right) \left( \mathbf{\tilde{h%
}}_{l_{1}}\ldots \mathbf{\tilde{h}}_{l_{m}}\right)  \label{in1} \\
&\leq &C_{2}\left[ \frac{\varrho }{\beta ^{1+\epsilon }}\left\vert \ln \beta
\right\vert +\beta \right] \sum_{m=2}^{\infty }\tau _{\ast
}^{m/m_{F}-1}m^{2}N_{h}^{m}c_{1}^{m}R^{m-1}  \notag
\end{eqnarray}%
with $c_{1}=9C_{\Xi }^{5}C_{\chi }$. The series converges if, in addition to
(\ref{Nhh}), $R$ satisfies the inequality%
\begin{equation*}
RN_{h}c_{1}\tau _{\ast }^{1/m_{F}}<1.
\end{equation*}%
For such $R$ and $\tau _{\ast }$, combining (\ref{in1}) with (\ref{Gtr}) and
using (\ref{scale1}) we obtain (\ref{rem}), and the Theorem \ref{Theorem
Superposition} is proved.

\subsection{Proof of the Superposition principle for PDE}

Here we prove \ Theorem \ref{Theorem Superposition11} (and its particular
case Theorem \ref{Theorem Superposition1}). The proof is completely similar
to the above proof of Theorem \ref{Theorem Superposition} up to every
detail. One only have to replace $\mathbb{D}_{m}$ given by (\ref{Dm}) \ by $%
\mathbb{D}_{m}$ given by (\ref{DmR}) \ and the space $L_{1}$ is now defined
by (\ref{L1R}) instead of (\ref{L1}).

\begin{remark}
Note that smallness of CI terms is essential and is based on different group
velocities of single band wavepackets. Note that separation of different
wavepackets based only on FM and NFM\ arguments as in Lemma \ref{Lemma
NFMmin} is impossible since there are always FM monomials with different $l$
because of the symmetry conditions (\ref{invsym}), (\ref{omeven}), for
example FM condition%
\begin{equation*}
\zeta \omega _{n,\zeta }\left( \zeta \mathbf{k}_{\ast }\right) -\zeta
^{\prime }\omega _{n^{\prime }}\left( \zeta ^{\prime }\mathbf{k}_{\ast
1}\right) -\zeta ^{\prime \prime }\omega _{n^{\prime \prime }}\left( \zeta
^{\prime \prime }\mathbf{k}_{\ast 2}\right) -\zeta ^{\prime \prime \prime
}\omega _{n^{\prime \prime \prime }}\left( \zeta ^{\prime \prime \prime }%
\mathbf{k}_{\ast 3}\right) =0
\end{equation*}%
is fulfilled if 
\begin{equation*}
n=n^{\prime },\ \zeta =\zeta ^{\prime },\ \mathbf{k}_{\ast }=\mathbf{k}%
_{\ast 1},\ n^{\prime \prime }=n^{\prime \prime \prime },\ \zeta ^{\prime
\prime }=-\zeta ^{\prime \prime \prime },\ \mathbf{k}_{\ast 2}=\mathbf{k}%
_{\ast 3}
\end{equation*}%
independently of the values of $\mathbf{k}_{\ast }$, $\mathbf{k}_{\ast 3}$
and independently of a particular form of functions $\omega _{n}\left( 
\mathbf{k}\right) $.
\end{remark}

\section{Examples and possible generalizations}

\subsection{ Fermi-Pasta-Ulam equation}

FPU equation on the infinite lattice has the form%
\begin{gather}
\partial _{t}^{2}x_{n}=\left( x_{n-1}-2x_{n}+x_{n+1}\right) +\alpha
_{3}\left( \left( x_{n+1}-x_{n}\right) ^{3}-\left( x_{n}-x_{n-1}\right)
^{3}\right)  \label{FPM2} \\
+\alpha _{2}\left( \left( x_{n+1}-x_{n}\right) ^{2}-\left(
x_{n}-x_{n-1}\right) ^{2}\right) .  \notag
\end{gather}%
It can be reduced to the following first-order equation 
\begin{equation}
\partial _{t}x_{n}=y_{n}-y_{n-1},\ \partial _{t}y_{n}=x_{n+1}-x_{n}+\alpha
_{3}\left( x_{n+1}-x_{n}\right) ^{3}+\alpha _{2}\left( x_{n+1}-x_{n}\right)
^{2}.  \label{FPM11}
\end{equation}%
We introduce lattice Fourier transforms $\tilde{x}\left( k\right) $ and $%
\tilde{y}\left( k\right) $ by (\ref{Fourintr}), namely%
\begin{equation*}
\tilde{x}\left( k\right) =\sum_{n}x_{n}\mathrm{e}^{-\mathrm{i}nk},\ k\in %
\left[ -\pi ,\pi \right] .
\end{equation*}%
First we write Fourier transform of the linear part of (\ref{FPM11}) (that
is with $\alpha _{3}=\alpha _{2}=0$). Multiplying by $\mathrm{e}^{-\mathrm{i}%
nk}$ \ and doing summation we obtain%
\begin{equation*}
\partial _{t}\tilde{x}\left( k\right) =\tilde{y}\left( k\right) -\mathrm{e}%
^{-\mathrm{i}k}\tilde{y}\left( k\right) ,\ \partial _{t}\tilde{y}\left(
k\right) =\mathrm{e}^{\mathrm{i}k}\tilde{x}\left( k\right) -\tilde{x}\left(
k\right) .
\end{equation*}%
that can be recast in the matrix form as follows 
\begin{equation*}
\partial _{t}\left[ 
\begin{array}{c}
\tilde{x} \\ 
\tilde{y}%
\end{array}%
\right] =\left[ 
\begin{array}{cc}
0 & -\left( \mathrm{e}^{\mathrm{i}k}-1\right) ^{\ast } \\ 
\mathrm{e}^{\mathrm{i}k}-1 & 0%
\end{array}%
\right] \left[ 
\begin{array}{c}
\tilde{x} \\ 
\tilde{y}%
\end{array}%
\right] .
\end{equation*}%
The eigenvalues of the matrix are purely imaginary and equal $\mathrm{i}%
\omega _{\zeta }\left( k\right) $ with 
\begin{equation*}
\omega _{\zeta }\left( k\right) =\zeta \left\vert \mathrm{e}^{\mathrm{i}%
k}-1\right\vert =2\zeta \left\vert \sin \frac{k}{2}\right\vert ,\ \zeta =\pm
,\ -\pi \leq k\leq \pi .
\end{equation*}%
The eigenvectors are orthogonal and are given explicitly by 
\begin{equation}
\mathbf{g}_{\zeta }\left( k\right) =\frac{1}{\sqrt{2}\left\vert \mathrm{e}^{%
\mathrm{i}k}-1\right\vert }\left[ 
\begin{array}{c}
i\zeta \left\vert \mathrm{e}^{\mathrm{i}k}-1\right\vert \\ 
\mathrm{e}^{\mathrm{i}k}-1%
\end{array}%
\right] =\frac{1}{\sqrt{2}}\left[ 
\begin{array}{c}
i\zeta \\ 
\frac{\mathrm{e}^{\mathrm{i}k}-1}{\left\vert \mathrm{e}^{\mathrm{i}%
k}-1\right\vert }%
\end{array}%
\right] ,\ \zeta =\pm ,k\neq 0.  \label{Gzet}
\end{equation}%
Now let us consider nonlinear terms. Note that the lattice Fourier transform
of the product $x\left( \mathbf{\mathbf{n}}\right) z\left( \mathbf{\mathbf{n}%
}\right) $, $\mathbf{\mathbf{n}}\in \mathbb{Z}^{d}$ \ is given by the
following convolution formula 
\begin{equation}
\widetilde{xz}\left( \mathbf{\mathbf{k}}\right) =\frac{1}{\left( 2\pi
\right) ^{d}}\int_{\left[ -\pi ,\pi \right] ^{d}}\tilde{x}\left( \mathbf{%
\mathbf{s}}\right) \tilde{z}\left( \mathbf{\mathbf{k}}-\mathbf{\mathbf{s}}%
\right) \,\mathrm{d}\mathbf{s}  \label{convtor}
\end{equation}%
as in the case of the continuous Fourier transform. Note that 
\begin{equation*}
\widetilde{x_{n+1}-x_{n}}\left( k\right) =\left( \mathrm{e}^{\mathrm{i}%
k}-1\right) \tilde{x}\left( k\right) ,
\end{equation*}%
and, hence, the Fourier transform of the cubic term of the nonlinearity in (%
\ref{FPM11}) is 
\begin{gather}
\widetilde{\left( x_{n+1}-x_{n}\right) ^{3}}=  \label{chi3} \\
=\frac{1}{\left( 2\pi \right) ^{2}}\int_{\substack{ k^{\prime }+k^{\prime
\prime }+k^{\prime \prime \prime }=\mathbf{k}  \\ \left[ -\pi ,\pi \right]
^{2}}}\left( \mathrm{e}^{\mathrm{i}k^{\prime }}-1\right) \left( \mathrm{e}^{%
\mathrm{i}k^{\prime \prime }}-1\right) \left( \mathrm{e}^{\mathrm{i}%
k^{\prime \prime \prime }}-1\right) \tilde{x}\left( k^{\prime }\right) 
\tilde{x}\left( k^{\prime \prime }\right) \tilde{x}\left( k^{\prime \prime
\prime }\right) \mathrm{d}k^{\prime }\mathrm{d}k^{\prime \prime },  \notag
\end{gather}%
and similar convolution for the quadratic term.

\subsection{Examples of wavepacket data}

Here we give examples of initial data for PDE in $\mathbb{R}^{d}$ and on the
lattice $\mathbb{Z}^{d}$ which are wavepackets in the sense of Definition %
\ref{dwavepack}. We define a wavepacket by (\ref{hbold}) where $h_{\zeta }$
is chosen to satisfy (\ref{hbder}) and (\ref{sourloc}).

Recall that a \emph{Schwartz function} is an infinitely smooth function $%
\Phi \left( \mathbf{r}\right) $, $\mathbf{r}\in \mathbb{R}^{d}$ which
rapidly decays and satisfies for every $s\geq 0$ the inequality\ 
\begin{equation}
\sup_{\mathbf{r}}\sum_{\left\vert \alpha \right\vert +p\leq s}\left\vert 
\mathbf{r}\right\vert ^{p}\left\vert \partial _{\mathbf{r}}^{\alpha }\Phi
\left( \mathbf{r}\right) \right\vert d\mathbf{r}\leq C_{1}\left( s\right) ,
\label{fismooth}
\end{equation}%
where 
\begin{equation*}
\partial _{\mathbf{r}}^{\alpha }\Phi \left( \mathbf{r}\right) =\partial
_{r_{1}}^{\alpha _{1}}\ldots \partial _{r_{d}}^{\alpha _{d}}\Phi \left( 
\mathbf{r}\right) ,\ \alpha =\left( \alpha _{1},\ldots ,\alpha _{d}\right)
,\left\vert \ \alpha \right\vert =\alpha _{1}+\ldots +\alpha _{d}.
\end{equation*}%
It is well known that Fourier transform of a Schwartz function remains to be
a Schwartz function and that its derivatives satisfy the inequality 
\begin{equation}
\sup_{\mathbf{k}}\sum_{\left\vert \alpha \right\vert +p\leq s}\left\vert
\left\vert \mathbf{k}\right\vert ^{p}\partial _{\mathbf{k}}^{\alpha }\hat{%
\Phi}\left( \mathbf{k}\right) \right\vert \leq C_{2}\left( s\right) .\ \text{
}  \label{fihat}
\end{equation}

\paragraph{Example 1.}

We consider equation in $\mathbb{R}^{d}$ as in Subsection 1.2. The simplest
example of a wavepacket in the sense of Definition \ref{dwavepack} is a
function of the form \ (\ref{h0}) where 
\begin{equation}
\int_{\mathbb{R}^{d}}\left\vert \hat{h}_{\zeta }\left( \mathbf{k}\right)
\right\vert +\left\vert \nabla _{\mathbf{k}}\hat{h}_{\zeta }\left( \mathbf{k}%
\right) \right\vert +\left\vert \mathbf{k}\right\vert ^{1/\epsilon
}\left\vert \hat{h}_{\zeta }\left( \mathbf{k}\right) \right\vert \mathrm{d}%
\mathbf{k}<\infty ,\ .  \label{condh}
\end{equation}%
and $\mathbf{g}_{n,\zeta }\left( \mathbf{k}\right) $ is an eigenvector from (%
\ref{OmomL}). Note that $\beta ^{-d}\hat{h}_{\zeta }\left( \mathbf{k}/\beta
\right) $ is the Fourier transform of a function $h_{\zeta }\left( \beta 
\mathbf{r}\right) $.

\begin{lemma}
\label{Lemma NLSwave} Let $\mathbf{\hat{h}}\left( \beta ,\mathbf{k}\right) ,$
$\mathbf{k}\in \mathbb{R}^{d}$ be defined by (\ref{h0}), (\ref{condh}). Then 
$\mathbf{\hat{h}}_{l,\zeta }\left( \beta ,\mathbf{k}\right) $ is a
wavepacket with wavepacket center $\mathbf{k}_{\ast }$ in the sense of
Definition \ref{dwavepack} \ with $L_{1}=L_{1}\left( \mathbb{R}^{d}\right) $.
\end{lemma}

\begin{proof}
First, condition (\ref{L1b}) holds since 
\begin{equation*}
\left\Vert \mathbf{\hat{h}}_{\zeta }\left( \beta ,\mathbf{\cdot }\right)
\right\Vert _{L_{1}}=\int_{\mathbb{R}^{d}}\beta ^{-d}\left\vert \hat{h}%
_{\zeta }\left( \frac{\mathbf{k}-\zeta \mathbf{k}_{\ast }}{\beta }\right) 
\mathbf{g}_{n,\zeta }\left( \mathbf{k}_{\ast }\right) \right\vert \mathrm{d}%
\mathbf{k}=\left\vert \mathbf{g}_{n,\zeta }\left( \mathbf{k}_{\ast }\right)
\right\vert \int_{\mathbb{R}^{d}}\left\vert \hat{h}_{\zeta }\left( \mathbf{k}%
\ \right) \right\vert \mathrm{d}\mathbf{k}.
\end{equation*}%
Condition (\ref{hbold}) is obviously fulfilled since 
\begin{equation*}
\mathbf{\hat{h}}_{\zeta }\left( \beta ,\mathbf{k}\right) =\Pi _{n,\zeta
}\left( \mathbf{\mathbf{k}}\right) \mathbf{\tilde{h}}_{\zeta }\left( \beta ,%
\mathbf{k}\right) .
\end{equation*}%
Inequality (\ref{sourloc}) follows from \ the estimate 
\begin{equation}
\beta ^{-d}\int_{\left\vert \mathbf{k-}\zeta \mathbf{k}_{\ast }\right\vert
\geq \beta ^{1-\epsilon }}\left\vert \hat{h}_{\zeta }\left( \frac{\mathbf{k}%
-\zeta \mathbf{k}_{\ast }}{\beta }\right) \right\vert \mathrm{d}\mathbf{k}%
\leq \beta \int_{\left\vert \mathbf{k}\right\vert \geq \beta ^{-\epsilon }}\
\left\vert \mathbf{k}\right\vert ^{1/\epsilon }\left\vert \hat{h}_{\zeta
}\left( \mathbf{k}\right) \right\vert \mathrm{d}\mathbf{k}\leq C\beta .
\label{p3}
\end{equation}%
To verify (\ref{hbder}) we note that since $\Pi _{n,\zeta }\left( \mathbf{%
\mathbf{k}}\right) $ smoothly depend on $k$ near $\zeta \mathbf{k}_{\ast }$
we have 
\begin{gather*}
\int_{\left\vert \mathbf{k-}\zeta \mathbf{k}_{\ast }\right\vert \leq \beta
^{1-\epsilon }}\left\vert \nabla _{\mathbf{k}}\mathbf{\hat{h}}_{\zeta
}\left( \beta ,\mathbf{k}\right) \right\vert \mathrm{d}\mathbf{k} \\
\leq C\int_{\left\vert \mathbf{k-}\zeta \mathbf{k}_{\ast }\right\vert \leq
\beta ^{1-\epsilon }}\beta ^{-d-1}\left\vert \nabla _{\mathbf{k}}\hat{h}%
_{l}\left( \frac{\mathbf{k}-\zeta \mathbf{k}_{\ast }}{\beta }\right)
\right\vert +\beta ^{-d}\left\vert \hat{h}_{l}\left( \frac{\mathbf{k}-\zeta 
\mathbf{k}_{\ast }}{\beta }\right) \right\vert \mathrm{d}\mathbf{k} \\
\leq C\beta ^{-1}\int_{\mathbb{R}^{d}}\left\vert \nabla _{\mathbf{k}}\hat{h}%
_{\zeta }\left( \mathbf{k}\right) \right\vert \mathrm{d}\mathbf{k}+C
\end{gather*}%
and (\ref{condh}) implies (\ref{hbder}).
\end{proof}

\paragraph{Example 2.}

Let us consider a lattice equation in $\mathbb{Z}^{d}$ as in Section 1.1. We
\ would like to give a sufficient condition for functions defined on the
lattice which ensures that their Fourier transforms satisfy all requirements
of Definition \ref{dwavepack}. We pick a Schwartz function $\Phi \left( 
\mathbf{r}\right) $ (see (\ref{fismooth})), a vector $\mathbf{k}_{\ast }\in %
\left[ -\pi ,\pi \right] ^{d}$ and introduce 
\begin{equation}
h\left( \beta ,\mathbf{r}\right) =\mathrm{e}^{-\mathrm{i}\mathbf{\mathbf{r}}%
\cdot \mathbf{\mathbf{\mathbf{\mathbf{k}}}}_{\ast }}\Phi \left( \beta 
\mathbf{r}\right) ,\ \mathbf{r}\in \mathbb{R}^{d}.  \label{he}
\end{equation}%
Then we restrict the above function to the lattice $\mathbb{Z}^{d}$ by
setting $\mathbf{r}=\mathbf{m}$. The following lemma is similar to Lemma \ref%
{Lemma NLSwave}.

\begin{lemma}
\label{Lemma 201}\textbf{\ }Let $\Phi \left( \mathbf{r}\right) \ $be a
Schwartz function, $h_{\zeta }\left( \beta ,\mathbf{r}\right) $ be \ defined
by (\ref{he}), $\tilde{h}_{\zeta }\left( \beta ,\mathbf{k}\right) $ $\ $ be
its lattice Fourier transform. Then the function\ $\tilde{h}_{\zeta }\left(
\beta ,\mathbf{k}\right) $ extended to $\mathbb{R}^{d}$ as a periodic
function with period $2\pi $ \ satisfies all requirements of Definition \ref%
{dwavepack} with $L_{1}=L_{1}\left( \left[ -\pi ,\pi \right] ^{d}\right) $.
\end{lemma}

\begin{proof}
The lattice Fourier transform of $h\left( \beta ,\mathbf{r}\right) $ equals
\ \ 
\begin{equation}
\tilde{h}\left( \beta ,\mathbf{k}\right) =\sum_{\mathbf{\mathbf{m\in }}%
\mathbb{Z}^{d}}\mathrm{e}^{-\mathrm{i}\mathbf{\mathbf{m}}\cdot \mathbf{%
\mathbf{\mathbf{\mathbf{k}}}}_{\ast }}\Phi \left( \beta \mathbf{\mathbf{m}}%
\right) e^{-\mathrm{i}\mathbf{\mathbf{m}}\cdot \mathbf{\mathbf{\mathbf{%
\mathbf{k}}}}}=\sum_{\mathbf{\mathbf{m\in }}\mathbb{Z}^{d}}\Phi \left( \beta 
\mathbf{\mathbf{m}}\right) \mathrm{e}^{-\mathrm{i}\mathbf{\mathbf{m}}\cdot
\left( \mathbf{\mathbf{\mathbf{\mathbf{k-k}}}}_{\ast }\right) }.  \label{hfi}
\end{equation}%
Since the above expression naturally defines $\tilde{h}\left( \beta ,\mathbf{%
k}\right) $ as a function of $\mathbf{\mathbf{\mathbf{\mathbf{k}}}}-\mathbf{%
\mathbf{\mathbf{\mathbf{k}}}}_{\ast }$, it is sufficient to take $\mathbf{%
\mathbf{\mathbf{\mathbf{k}}}}_{\ast }=0$. To get (\ref{sourloc}), we use the
representation of $\Phi \left( \mathbf{\mathbf{r}}\right) $ in terms of
inverse Fourier transform (\ref{Finv}) 
\begin{equation}
\Phi \left( \mathbf{\mathbf{r}}\right) =\frac{1}{\left( 2\pi \right) ^{d}}%
\int_{\mathbb{R}^{d}}\hat{\Phi}\left( \mathbf{\mathbf{k}}\right) \mathrm{e}^{%
\mathrm{i}\mathbf{r}\cdot \mathbf{\mathbf{k}}}\mathrm{d}\mathbf{k},\mathbf{\ 
}\Phi \left( \beta \mathbf{\mathbf{m}}\right) =\frac{1}{\left( 2\pi \beta
\right) ^{d}}\int_{\mathbb{R}^{d}}\hat{\Phi}\left( \frac{1}{\beta }\mathbf{%
\mathbf{k}}\right) \mathrm{e}^{\mathrm{i}\mathbf{m}\cdot \mathbf{\mathbf{k}}}%
\mathrm{d}\mathbf{k}.  \label{fif}
\end{equation}%
We split $\Phi \left( \beta \mathbf{\mathbf{m}}\right) $ into two terms: 
\begin{eqnarray}
\Phi \left( \beta \mathbf{\mathbf{m}}\right) &=&\frac{1}{\left( 2\pi \beta
\right) ^{d}}\int_{\mathbb{R}^{d}}\Psi \left( \frac{1}{\beta ^{1-\epsilon }}%
\xi \right) \hat{\Phi}\left( \frac{1}{\beta }\xi \right) \mathrm{e}^{\mathrm{%
i}\mathbf{m}\cdot \mathbf{\xi }}\mathrm{d}\mathbf{\xi }+\Phi _{1}\left( 
\mathbf{m}\right) ,  \label{fifi1} \\
\Phi _{1}\left( \mathbf{m}\right) &=&\frac{1}{\left( 2\pi \beta \right) ^{d}}%
\int_{\mathbb{R}^{d}}\left( 1-\Psi \left( \frac{1}{\beta ^{1-\epsilon }}\xi
\right) \right) \hat{\Phi}\left( \frac{1}{\beta }\xi \right) \mathrm{e}^{%
\mathrm{i}\mathbf{m}\cdot \mathbf{\xi }}\mathrm{d}\mathbf{\xi }  \notag
\end{eqnarray}%
with $\Psi \left( \xi \right) $ defined by (\ref{j0}). The first term in (%
\ref{fifi1}) coincides with the inverse lattice Fourier transform, its
lattice Fourier transform is explicitly given and can be treated as in Lemma %
\ref{Lemma NLSwave}. The second term gives $O\left( \beta ^{N}\right) $ with
large $N$ for Schwartz functions $\hat{\Phi}$. Using these observations we
check all points of Definition \ref{dwavepack} \ as in Lemma \ref{Lemma
NLSwave}.
\end{proof}

\subsection{The Nonlinear Maxwell equation}

We expect that the approximate superposition principle can be generalized to
the Nonlinear Maxwell equations (NLM)\ in periodic media studied in \cite%
{BF4} . A concise operator form of the NLM is%
\begin{equation}
\partial _{\tau }\mathbf{U}=-\frac{\mathrm{i}}{\varrho }\mathbf{MU}+\mathcal{%
F}_{\text{NL}}\left( \mathbf{U}\right) -\mathbf{J}_{0},\ \mathbf{U}\left(
\tau \right) =0\ \text{for }\tau \leq 0,  \notag
\end{equation}%
where the excitation current 
\begin{equation*}
\mathbf{J}\left( \tau \right) =0\ \text{for }\tau \leq 0.
\end{equation*}%
We were studying the properties of nonlinear wave interactions as described
by the Nonlinear Maxwell equations in series of papers \cite{BF1}-\cite{BF6}%
. Our analysis of the solutions to the NLM uses an expansion in terms of
orthonormal Floquet-Bloch basis $\mathbf{\tilde{G}}_{n,\zeta }\left( \mathbf{%
r},\mathbf{k}\right) $, $n=1,\ldots .$, namely 
\begin{equation}
\mathbf{\tilde{U}}\left( \mathbf{k},\mathbf{r},\tau \right) =\sum_{\zeta
=\pm 1}\sum_{n=1}^{\infty }\tilde{U}_{n,\zeta }\left( \mathbf{k},\tau
\right) \mathbf{\tilde{G}}_{n,\zeta }\left( \mathbf{r},\mathbf{k}\right) ,\ 
\mathbf{k}\in \left[ -\pi ,\pi \right] ^{d}.  \label{Utild}
\end{equation}%
This expansion is similar to (\ref{Uboldj}) with $J$ replaced by $\infty $,
since the linear Maxwell operator in a periodic medium has infinitely many
bands. The excitation currents take the form similar to forcing term in (\ref%
{varc}), namely 
\begin{gather*}
\mathbf{\tilde{J}}\left( \mathbf{r},\mathbf{k},\tau \right) =\tilde{j}%
_{n,+}\left( \mathbf{k},\tau \right) \mathbf{\tilde{G}}_{n,+}\left( \mathbf{r%
},\mathbf{k}\right) \mathrm{e}^{-\frac{\mathrm{i}}{\varrho }\omega
_{n}\left( \mathbf{k}\right) \tau }+\tilde{j}_{n,-}\left( \mathbf{k},\tau
\right) \mathbf{\tilde{G}}_{n,-}\left( \mathbf{r},\mathbf{k}\right) \mathrm{e%
}^{\frac{\mathrm{i}}{\varrho }\omega _{n}\left( \mathbf{k}\right) \tau }, \\
\ \mathbf{\tilde{J}}_{n}\left( \mathbf{r},\mathbf{k},\tau \right) =0,\ n\neq
n_{0},
\end{gather*}%
with a fixed $n=n_{0}$. The difference with (\ref{varc})\ is that
time-independent $\mathbf{h}_{n,\zeta }\left( \mathbf{k}\right) $ is
replaced by $\tilde{j}_{n,\zeta \ }\left( \mathbf{k},\tau \right) $. The
functions $\tilde{j}_{n,\zeta \ }\left( \mathbf{k},\tau \right) $ for every $%
\tau $ have the form of wavepackets in the sense of Definition \ref%
{dwavepack}, or in particular the form similar to (\ref{h0}) with fixed $%
\mathbf{k}_{\ast }$.

The Existence and uniqueness Theorem for the NLM is proven in \cite{BF4}, in
particular function-analytic representation of the solution as a function of
the excitation current. The results of this paper can be extended to the
NLM\ equations provided that certain technical difficulties are addressed.
Particularly, the classical NLM\ equation allows for the time dispersion
with consequent time-convolution integration in the nonlinear term. This
complication can be addressed by approximating it with a nonlinearity of the
form (\ref{Fmintr}) with an error $O\left( \varrho \right) =O\left( \beta
^{2}\right) $, see \cite{BF6}. Then the derivation of the approximate linear
superposition principle for wavepackets can be done as in this paper.
Another complication with the NLM is that it has infinite number of bands.

\subsection{Dissipative terms in the linear part}

Equations (\ref{eqFur}) and (\ref{difF}) \ involve linear operators $\mathrm{%
i}\mathbf{L}\left( \mathbf{k}\right) $ with \ purely imaginary spectrum.
Quite similarly we can consider equations of the form 
\begin{equation}
\partial _{\tau }\mathbf{\hat{U}}\left( \mathbf{\mathbf{k}},\tau \right) =%
\left[ -\mathbf{G}\left( \mathbf{k}\right) -\frac{\mathrm{i}}{\varrho }%
\mathbf{L}\left( \mathbf{\mathbf{k}}\right) \right] \mathbf{\hat{U}}\left( 
\mathbf{\mathbf{k}},\tau \right) +\hat{F}\left( \mathbf{\hat{U}}\right)
\left( \mathbf{\mathbf{k}},\tau \right) ,  \label{eqdiss}
\end{equation}%
where a Hermitian \ matrix $\mathbf{G}\left( \mathbf{k}\right) $ commutes
with the Hermitian matrix $\mathbf{L}\left( \mathbf{k}\right) $ and $\mathbf{%
G}\left( \mathbf{k}\right) $ is non-negative. In this case the approximate
superposition principle also holds. The proofs are quite similar. In the
case (\ref{difF}), which corresponds to of PDE, $\mathbf{G}\left( \mathbf{k}%
\right) $ determines a dissipative term, for example $\mathbf{G}\left( 
\mathbf{k}\right) =\left\vert \mathbf{k}\right\vert ^{2}I,\mathbf{k}\in 
\mathbb{R}^{d}$, where $I$ is the identity matrix, corresponds to Laplace
operator $\Delta $. When such a dissipative term is introduced, we can
consider nonlinearities $\hat{F}$ which involve derivatives, see \cite{BMN1}%
, \cite{BMN2} in a similar situation. For such nonlinearities our framework
remains the same, but some statements and proofs have to be modified. We
will consider this case in a separate paper.

\section{Appendix A: Structure of a composition monomial based on
oscillatory integrals}

Every composition monomial $M\left( \mathcal{F},T,\vec{\lambda}_{\left( \hat{%
s}\right) },\vec{\zeta}_{\left( m\right) }\right) \left( \mathbf{\tilde{h}}%
_{1}\ldots \mathbf{\tilde{h}}_{m}\right) $ based on oscillatory integral
operators $\mathcal{F}^{\left( m\right) }$ as defined by (\ref{Fmdel}) and
the space decomposition as defined by (\ref{PPosc}) has the following
structure. Let $T$ be the tree corresponding to the monomial $M$. The
monomial involves integration with respect to time variables $\tau _{\left(
N\right) }$ where $N\in T$ are the nodes of the tree $T$. The monomial also
involves integration with respect to variables $\mathbf{k}_{N}$, $N\in T$.
The argument of the integral operator $M\left( \mathcal{F},T,\vec{\lambda}%
_{\left( \hat{s}\right) },\vec{\zeta}_{\left( m\right) }\right) $ involves
only end nodes (of zero rank) and has the form 
\begin{equation*}
\dprod\limits_{\limfunc{rank}\left( N\right) =0}\mathbf{\tilde{h}}_{N}\left( 
\mathbf{k}_{N}\right) .
\end{equation*}%
The kernel of the integral operator involves the composition monomial $%
M\left( \chi ,T,\vec{\lambda}_{\left( \hat{s}\right) },\vec{\zeta}_{\left(
m\right) }\right) $ based on the susceptibilities tensors $\chi _{\zeta ,%
\vec{\zeta}_{\left( m\right) }}^{\left( m\right) }\left( \mathbf{\mathbf{k}},%
\vec{k}_{\left( m\right) }\right) $ with the same tree $T$. Note that the
phase matching condition (\ref{Phmc}) takes the form 
\begin{equation*}
\mathbf{k}_{N}=\mathbf{k}_{N}^{\prime }+\ldots +\mathbf{k}_{N}^{\left( \mu
\left( N\right) \right) }=\sum_{i=1}^{\mu \left( N\right) }\mathbf{k}%
_{c_{i}\left( N\right) }.
\end{equation*}%
Recall that if $c_{i}\left( N\right) $, $i=1,\ldots ,\mu \left( N\right) $
is the $i$-th child node of $N$, then the arguments in (\ref{Fmdel}) are
determined by the formula 
\begin{equation*}
\mathbf{k}_{c_{i}\left( N\right) }=\mathbf{k}_{N}^{\left( c_{i}\right) }.
\end{equation*}%
Hence, the kernel of the integral operator $M\left( \mathcal{F},T,\vec{%
\lambda}_{\left( \hat{s}\right) },\vec{\zeta}_{\left( m\right) }\right)
\left( \mathbf{\tilde{h}}_{1}\ldots \mathbf{\tilde{h}}_{m}\right) $ involves
the product of normalized delta functions 
\begin{equation*}
\dprod\limits_{\limfunc{rank}\left( N\right) >0}\delta \left( \mathbf{k}_{N}-%
\mathbf{k}_{c_{1}\left( N\right) }-\ldots -\mathbf{k}_{c_{\mu \left(
N\right) }\left( N\right) }\right) ,
\end{equation*}%
and the integration with respect to $\mathbf{k}_{N}$ is over the torus 
\begin{equation*}
\left( \dprod\limits_{N\neq N_{\ast }}\int_{\left[ -\pi ,\pi \right] ^{\mu
\left( N\right) d}}\right) \left[ \ldots \right] \dprod\limits_{N\neq
N_{\ast }}\,\mathrm{d}\mathbf{k}_{N},
\end{equation*}%
and, obviously, the variable $\mathbf{k}_{N_{\ast }}$ corresponding to the
root node $N_{\ast }$ is not involved into the integration.

Since every operator $\mathcal{F}^{\left( m\right) }$ at a node $N$ of the
monomial $M\left( \mathcal{F},T,\vec{\lambda}_{\left( \hat{s}\right) },\vec{%
\zeta}_{\left( m\right) }\right) $ contains the oscillatory factor 
\begin{gather*}
\exp \left\{ \mathrm{i}\phi _{\ \zeta ,\vec{\zeta}_{\left( m\right)
,N}}\left( \mathbf{\mathbf{k}},\vec{k}_{\left( m\right) }\right) \frac{\tau
_{\left( N\right) }}{\varrho }\right\} = \\
\exp \left\{ \mathrm{i}\left[ \zeta _{N}\omega \left( \mathbf{k}_{N}\right)
-\zeta _{N}^{\prime }\omega \left( \mathbf{k}_{N}^{\prime }\right) -\ldots
-\zeta _{N}^{\left( m\right) }\omega \left( \mathbf{k}_{N}^{\left( m\right)
}\right) \right] \frac{\tau _{\left( N\right) }}{\varrho }\right\} ,
\end{gather*}%
we obtain the following total oscillatory factor%
\begin{equation}
\exp \left\{ \mathrm{i}\frac{1}{\varrho }\Phi _{\ \zeta ,\vec{\zeta}_{\left(
m\right) ,T}}\left( \mathbf{\mathbf{k}},\vec{k}_{\left( m\right) }\right)
\right\} ,  \label{expFi}
\end{equation}%
where the phase function $\Phi _{\ T,\vec{\zeta}}\left( \vec{k}\right) $ of
the monomial is defined by the formula 
\begin{equation}
\Phi _{T,\vec{\zeta}}\left( \vec{k},\vec{\tau}\right) =\sum_{N\in T}\left[
\zeta _{N}\omega \left( \mathbf{k}\right) -\sum_{i=1}^{\mu \left( N\right)
}\zeta _{N}^{\left( c_{i}\left( N\right) \right) }\omega \left( \mathbf{k}%
_{c_{i}\left( N\right) }\right) \right] \tau _{\left( N\right) }.
\label{FiT}
\end{equation}%
The vectors $\vec{k}$, $\vec{\tau}$ and $\vec{\zeta}$ are composed of $%
\mathbf{k}_{N},$ $\tau _{N}$ and $\zeta _{N}$ using the standard labeling of
the nodes.

\emph{Notice then that the oscillatory exponent (\ref{expFi}) is the only
expression in the composition monomial which involves parameter} $\varrho $.
Observe also that the FM condition takes here the form 
\begin{equation*}
\zeta _{N}=\sum_{i=1}^{\mu \left( N\right) }\zeta _{N}^{\left( c_{i}\left(
N\right) \right) }.
\end{equation*}%
The domain of integration with respect to time variables is given in terms
of the tree $T$ by the following inequalities 
\begin{equation}
D_{T}=\left\{ \tau _{\left( N\right) }:0\leq \tau _{\left( N\right) }\leq
\tau _{\left( p\left( N\right) \right) },N\in T\setminus N_{\ast }\right\}
\label{DT}
\end{equation}%
where $p\left( N\right) $ is the parent node of the node $N$. Using
introduced notations we can write the action of the monomial $M\left( 
\mathcal{F},T,\vec{\lambda}_{\left( \hat{s}\right) },\vec{\zeta}_{\left(
m\right) }\right) $ in the form%
\begin{gather}
M\left( \mathcal{F},T,\vec{\lambda},\vec{\zeta}\right) \left( \dprod\limits_{%
\text{rank}\left( N\right) =0}\mathbf{\tilde{h}}_{N}\right) \left( \mathbf{k}%
_{N_{\ast }},\tau _{N_{\ast }}\right) =\int_{D_{T}}\left(
\dprod\limits_{N\neq N_{\ast }}\int_{\left[ -\pi ,\pi \right] ^{\mu \left(
N\right) d}}\right)  \label{Monint} \\
\exp \left\{ \mathrm{i}\frac{1}{\varrho }\Phi _{\ T,\vec{\zeta}}\left( \vec{k%
},\vec{\tau}\right) \right\} M\left( \chi ,T,\vec{\lambda},\vec{\zeta},\vec{k%
}\right) \dprod\limits_{\text{rank}\left( N\right) =0}\mathbf{\tilde{h}}%
_{N}\left( \mathbf{k}_{N}\right)  \notag \\
\dprod\limits_{\text{rank}\left( N\right) >0}\delta \left( \mathbf{k}_{N}-%
\mathbf{k}_{c_{1}\left( N\right) }-\cdots -\mathbf{k}_{c_{\mu \left(
N\right) }\left( N\right) }\right) \dprod\limits_{N\neq N_{\ast }}\,\mathrm{d%
}\mathbf{k}_{N}\dprod\limits_{N\neq N_{\ast }}\mathrm{d}\tau _{\left(
N\right) }.  \notag
\end{gather}%
Note that $m$ equals the number of end nodes, that is nodes with zero rank
and they are numerated using the standard labeling of the nodes, that is 
\begin{equation*}
\mathbf{\tilde{h}}_{1}\left( \mathbf{k}_{1}\right) \cdots \mathbf{\tilde{h}}%
_{m}\left( \mathbf{k}_{m}\right) =\dprod\limits_{\text{rank}\left( N\right)
=0}\mathbf{\tilde{h}}_{N}\left( \mathbf{k}_{N}\right) .
\end{equation*}%
The formula (\ref{Monint}) gives a closed form of a composition monomial
based on oscillatory integral operators $\mathcal{F}^{\left( m\right) }$
with an arbitrary large rank.

\section{Appendix B: Proof of the refined implicit function theorem}

Here we give the proof of Theorem \ref{Theorem Monomial Convergence}.

First, we consider the following elementary problem which provides majorants
for the problem of interest. Let a function of one complex variable be
defined by the formula 
\begin{equation}
\mathcal{\check{F}}\left( u\right) =C_{\mathcal{F}}\sum_{m=2}^{\infty
}u^{m}R_{\mathcal{F}}^{-m}=C_{\mathcal{F}}\left[ \frac{u^{2}/R_{\mathcal{F}%
}^{2}}{1-u/R_{\mathcal{F}}}\right] ,C_{\mathcal{F}}>0,R_{\mathcal{F}}>0.
\label{fu}
\end{equation}%
In this case $\mathcal{\check{F}}^{\left( m\right) }\left( x_{1}\ldots
x_{m}\right) =C_{\mathcal{F}}R_{\mathcal{F}}^{-m}x_{1}\ldots x_{m}$. Let us
introduce the equation 
\begin{equation}
u=\mathcal{\check{F}}\left( u\right) +x,\ u,x\in \mathbb{C}  \label{fux}
\end{equation}%
which is a particular case of (\ref{eqF}). A small solution $u\left(
x\right) $ of this equation such that $u\left( 0\right) =0$ is given by the
series 
\begin{equation*}
u=\check{G}\left( x\right) =\sum_{m=1}^{\infty }\check{G}^{\left( m\right)
}x^{m},
\end{equation*}%
which is a particular case of formula (\ref{Gser}). The terms $\check{G}%
^{\left( m\right) }x^{m}$ of this problem are determined from (\ref{recG0})
and can be written in the form (\ref{treesum})%
\begin{equation}
\check{G}^{\left( m\right) }x^{m}=\sum_{T\in T_{m}}c_{T}M\left( \mathcal{%
\check{F}},T\right) x^{m}.  \label{Ghat}
\end{equation}%
Obviously, 
\begin{equation}
M\left( \mathcal{\check{F}},T\right) x^{m}=C_{\mathcal{F}}^{i\left( T\right)
}R_{\mathcal{F}}^{-e\left( T\right) }x^{m}  \label{mfhat}
\end{equation}%
where $i\left( T\right) $ is the incidence number of the tree $T$, $e\left(
T\right) $ is the number of edges of $T$. Now we compare solution of the
general equation (\ref{eqF}). It is given by the formula (\ref{Gser}) with
operators $\ \mathcal{G}^{\left( m\right) }\left( \mathbf{\mathbf{u}}%
^{m}\right) $ admitting expansion (\ref{treesum}). Since 
\begin{equation*}
\left\Vert \mathcal{F}^{\left( m\right) }\right\Vert \leq C_{\mathcal{F}}R_{%
\mathcal{F}}^{-m},
\end{equation*}%
where the constants are the same as in (\ref{fu}) we have 
\begin{equation*}
\left\Vert M\left( \mathcal{F},T\right) \left( \mathbf{\mathbf{x}}_{1}\ldots 
\mathbf{\mathbf{x}}_{\nu }\right) \right\Vert \leq M\left( \mathcal{\check{F}%
},T\right) \left\Vert \mathbf{\mathbf{x}}_{1}\right\Vert \ldots \left\Vert 
\mathbf{\mathbf{x}}_{\nu }\right\Vert ,
\end{equation*}%
implying 
\begin{equation}
\sum_{T\in T_{m}}c_{T}\left\Vert M\left( \mathcal{F},T\right) \left( \mathbf{%
\mathbf{x}}_{1}\ldots \mathbf{\mathbf{x}}_{m}\right) \right\Vert \leq
\sum_{T\in T_{m}}c_{T}M\left( \mathcal{\check{F}},T\right) \left\Vert 
\mathbf{\mathbf{x}}_{1}\right\Vert \ldots \left\Vert \mathbf{\mathbf{x}}%
_{m}\right\Vert =\check{G}^{\left( m\right) }\left\Vert \mathbf{\mathbf{x}}%
_{1}\right\Vert \ldots \left\Vert \mathbf{\mathbf{x}}_{m}\right\Vert .
\label{cGMG}
\end{equation}%
Solving (\ref{fux}) we get explicitly 
\begin{equation*}
u=\frac{R_{\mathcal{F}}}{2c}\left( 1-\sqrt{1-4c\frac{x}{R_{\mathcal{F}}}}%
\right) =\check{G}\left( x\right) ,\,c=\frac{C_{\mathcal{F}}}{R_{\mathcal{F}}%
}+1.
\end{equation*}%
We have the following estimate of the coefficients$\ $ 
\begin{equation}
\check{G}^{\left( m\right) }\leq \frac{R_{\mathcal{F}}^{2}}{2\left( C_{%
\mathcal{F}}+R_{\mathcal{F}}\right) }\left( 4\frac{C_{\mathcal{F}}+R_{%
\mathcal{F}}}{R_{\mathcal{F}}^{2}}\right) ^{m},\,m=1,2,\ldots ,
\label{Gmhat}
\end{equation}%
(see \cite{BF4} for details in a similar situation). From (\ref{mfhat}) and (%
\ref{Gmhat}) we infer the following inequality%
\begin{equation*}
\sum_{T\in T_{m}}c_{T}C_{\mathcal{F}}^{i\left( T\right) }R_{\mathcal{F}%
}^{-e\left( T\right) }\leq \frac{R_{\mathcal{F}}^{2}}{2\left( C_{\mathcal{F}%
}+R_{\mathcal{F}}\right) }\left( 4\frac{C_{\mathcal{F}}+R_{\mathcal{F}}}{R_{%
\mathcal{F}}^{2}}\right) ^{m}
\end{equation*}%
which hods for all $C_{\mathcal{F}}$, $R_{\ \mathcal{F}}>0$. We set $C_{%
\mathcal{F}}=R_{\mathcal{F}}=1$ and obtain the desired bound (\ref{ctes}).

\section{Notations and abbreviations}

For reader's convenience we provide below a list of notations and
abbreviations used in this paper.

\textbf{AFM}-- alternatively frequency matched, see Definition \ref%
{Definition AFM}

\textbf{ANFM} -- alternatively non-frequency-matched see Definition \ref%
{Definition AFM}

\textbf{band-crossing points} -- see Definition \ref{Definition
band-crossing point}

\textbf{cc} -- complex conjugate to the preceding terms in the formula

\textbf{composition monomial -- }see Definition \ref{Definition Rank}

\textbf{decoration projections} -- see (\ref{Ps}), (\ref{PP})

\textbf{decorated monomial }-- see Definition \ref{Generalized monomial of a
tree}

\textbf{CI} \textbf{monomials} -- \ cross-interacting monomials, see
Definition \ref{Definition GVM}

\textbf{FPU, Fermi-Pasta-Ulam equation -- }(\ref{FPMint}), (\ref{FPMin2}), (%
\ref{FPM2})

\textbf{Floquet-Bloch modal decomposition -- }see (\ref{Utild})

\textbf{Fourier transform} -- see (\ref{FoR})

\textbf{FM} -- frequency matched, see Definition \ref{Definition FM} see
also (\ref{sumz})

\textbf{homogeneity index of a monomial } Definition \ref{Definition Rank}

\textbf{homogeneity index of a tree} -- Definition \ref{Tree}

\textbf{incidence number of a monomial -- }number of occurrences of
operators $\ \mathcal{F}^{\left( l\right) }$ in the composition monomial

\textbf{incidence number of a monomial -- }see Definition \ref{defsubmon}

\textbf{incidence number of a tree -- }Definition \ref{Definition incidence}

\textbf{lattice Fourier transform} -- see (\ref{Fourintr})

\textbf{monomial} -- Definition \ref{Definition Rank}

\textbf{NFM} -- non-frequency-matched see Definition \ref{Definition FM},
see also (\ref{NFMz})

\textbf{oscillatory integral operator} -- see (\ref{Fm}), (\ref{varcu})

\textbf{rank of monomial} -- see Definition \ref{Definition Rank}

\textbf{root operator} (\ref{submon})

\textbf{SI monomials} -- self-interacting monomials, see Definition \ref%
{Definition GVM}

\textbf{Schwartz functions -- }infinitely smooth functions on $\mathbb{R}%
^{d} $ \ which decay faster than any power, see (\ref{fismooth})

\textbf{single-mode wavepacket -- }see Definition \ref{dwavepack}

\textbf{submonomial} (\ref{defsubmon})

\textbf{wavepacket} see Definition \ref{dwavepack}

$\mathrm{\tilde{d}}^{\left( m-1\right) }\vec{k}=\frac{1}{\left( 2\pi \right)
^{\left( m-1\right) d}}\,\mathrm{d}\mathbf{k}^{\prime }\ldots \,\mathrm{d}%
\mathbf{k}^{\left( m-1\right) }$ -- see (\ref{dtild})

$\mathbb{D}_{m}=\left[ -\pi ,\pi \right] ^{\left( m-1\right) d}$ -- see (\ref%
{Dm}) or $\mathbb{D}_{m}=\mathbb{R}^{\left( m-1\right) d}$ see (\ref{DmR})

$E=C\left( \left[ 0,\tau _{\ast }\right] ,L_{1}\right) $ see (\ref{Elat})

$\hat{F}^{\left( m\right) }$ -- $m$-linear operator in $L_{1}$, see (\ref%
{Fmintr}), (\ref{FY})

$\mathcal{F}_{n,\zeta ,\vec{n},\vec{\zeta}}^{\left( m\right) }$ -- basis
element of the $\ m$-linear operator $\mathcal{F}_{\ }^{\left( m\right) }$
in $E$ see (\ref{Fm})

$\mathcal{F}_{\lambda ,\vec{\zeta}}^{\left( n\right) }$ -- see (\ref{Flamz})

$\hat{h}_{\zeta }\left( \beta ,\mathbf{k}\right) ,\ \zeta =\pm $ -- Fourier
transform of the wavepacket initial data $h_{\zeta }\left( \beta ,\mathbf{r}%
\right) $, see Definition \ref{dwavepack}

$\hat{h}_{\zeta }\left( \frac{1}{\beta }\mathbf{\xi }\right) ,\ \zeta =\pm $
-- Fourier transform of the wavepacket $h_{\zeta }\left( \beta \mathbf{r}%
\right) $initial data, see Definition \ref{dwavepack}

$\mathbf{\tilde{h}}_{l}^{\Psi }\left( \mathbf{k},\beta \right) $ -- a
function nullified outside $\beta ^{1-\epsilon }$ vicinity of $\pm \mathbf{k}%
_{\ast }$, see (\ref{hpsi})

$\mathbf{k}=\left( k_{1},\ldots ,k_{d}\right) \in \left[ -\pi ,\pi \right]
^{d}$ - quasimomentum (wave vector) variable, (\ref{Fourintr}), (\ref{kkar}).

$\mathbf{k}=\left( k_{1},\ldots ,k_{d}\right) \in \mathbb{R}^{d}\ $--
Fourier wave vector variable, (\ref{FoR}), (\ref{kkar}).

$\mathbf{k}_{\ast }=\left( k_{\ast 1},\ldots ,k_{\ast d}\right) $ -- center
of the wavepacket see Definition \ref{dwavepack}

$\mathbf{k}_{\ast l}\ $-- center of $l$-th wavepacket

$\vec{k}=\left( \mathbf{k}^{\prime },\ldots ,\mathbf{k}^{\left( m\right)
}\right) $, -- interaction multi-wave vector, (\ref{kkar}), (\ref{zetaar}) .

$\mathbf{k}^{\left( s\right) }\left( \mathbf{k},\vec{k}\right) =\mathbf{k}-%
\mathbf{k}^{\prime }-\ldots -\mathbf{k}^{\left( s-1\right) }$ -- see (\ref%
{kkar})

$L_{1}$ -- Lebesgue space $L_{1}\left( \left[ -\pi ,\pi \right] ^{d}\right) $
or $L_{1}\left( \mathbb{R}^{d}\right) $ - see (\ref{L1}) and (\ref{L1R})

$n$ -- band number

$\vec{n}=\left( n^{\prime },\ldots ,n^{\left( m\right) }\right) $ -- band
interaction index, (\ref{zetaar})

$\nabla _{\mathbf{r}}=\left( \frac{\partial }{\partial r_{1}},\frac{\partial 
}{\partial r_{2}},\cdots ,\frac{\partial }{\partial r_{d}}\right) $ --
spatial gradient

$O\left( \mu \right) $ -- any quantity having the property that $\frac{%
O\left( \mu \right) }{\mu }$ is bounded as $\mu \rightarrow 0$.

$\omega _{\bar{n}}\left( \mathbf{k}\right) =\zeta \omega _{n}\left( \mathbf{k%
}\right) $ -- dispersion relation of the band $\left( \zeta ,n\right) $, see
(\ref{OmomL})

$\omega _{n_{0}}^{\prime }\left( \mathbf{k}\right) =\nabla _{\mathbf{k}%
}\omega _{n_{0}}\left( \mathbf{k}\right) $ -- group velocity vector

$\omega _{n}\left( \mathbf{k}\right) $-- $n$-th eigenvalue of $\mathbf{L}%
\left( \mathbf{\mathbf{k}}\right) $, see (\ref{OmomL}); dispersion relation
of $n$-th band

$\Psi $ -- cutoff function in quasimomentum domain, see (\ref{j0})

$\phi _{\vec{n}}\left( \mathbf{k},\vec{k}\right) =\zeta \omega _{n}\left( 
\mathbf{k}\right) -\zeta ^{\prime }\omega _{n^{\prime }}\left( \mathbf{k}%
^{\prime }\right) -\ldots -\zeta ^{\left( m\right) }\omega _{n^{\left(
m\right) }}\left( \mathbf{k}^{\left( m\right) }\right) $ -- interaction
phase function, (\ref{phim})

$\pi _{0}$ -- see (\ref{kpi0})

$\Pi _{n,\zeta }\left( \mathbf{\mathbf{k}}\right) \ $-- projection in $%
\mathbb{C}^{2J}$ onto direction of $\mathbf{g}_{n,\zeta }\left( \mathbf{k}%
\right) $; see (\ref{Pin})

$\mathbf{r}=\left( r_{1},\ldots ,r_{d}\right) $ -- spatial variable

$\varrho =\beta ^{2}$ -- (\ref{scale1})

$\sigma $ -- the set of band-crossing points, see Definition \ref{Definition
band-crossing point}

$\mathbf{\hat{U}}\left( \mathbf{k}\right) \ $-- Fourier transform of $%
U\left( \mathbf{r}\right) $, see (\ref{FoR})

$\mathbf{\tilde{U}}_{n,\zeta }\left( \mathbf{k},\tau \right) =\mathbf{\tilde{%
u}}_{n,\zeta }\left( \mathbf{k},\tau \right) \mathrm{e}^{-\frac{\mathrm{i}%
\tau }{\varrho }\zeta \omega _{n}\left( \mathbf{k}\right) }$ -- amplitudes,
see (\ref{Uu})

$\zeta =\pm $ or $\zeta =\pm 1$ -- band binary index.

$\vec{\zeta}=\left( \zeta ^{\prime },\ldots ,\zeta ^{\left( m\right)
}\right) $ -- binary band index vector, see (\ref{zetaar})

$Z^{\ast }$ -- complex conjugate to $Z$\newline

\textbf{Acknowledgment:} Effort of A. Babin and A. Figotin is sponsored by
the Air Force Office of Scientific Research, Air Force Materials Command,
USAF, under grant number FA9550-04-1-0359.\newline


\begin{thebibliography}{99}
\bibitem{BF1} Babin A. and Figotin A., \textsl{Nonlinear Photonic Crystals:
I. Quadratic nonlinearity}, Waves in Random Media, 11, R31-R102, (2001).

\bibitem{BF2} Babin A. and Figotin A., \textsl{Nonlinear Photonic Crystals:
II. Interaction classification for quadratic nonlinearities}, Waves in
Random Media, 12, R25-R52, (2002).

\bibitem{BF3} Babin A. and Figotin A., \textsl{Nonlinear Photonic Crystals:
III. Cubic Nonlinearity}, Waves in Random Media, \textbf{13}, pp. R41-R69
(2003).

\bibitem{BF4} Babin A. and Figotin A., \textsl{Nonlinear Maxwell Equations
in Inhomogenious Media}, Commun. Math. Phys. 241, 519-581 (2003).

\bibitem{BF5} Babin A. and Figotin A., \textsl{Polylinear spectral
decomposition for nonlinear Maxwell equations}, in Partial Differential
Equations, M.S. Agranovich and M.A. Shubin eds, Advances in Mathematical
Sciences, American Mathematical Society Translations -Series 2, Vol. 206,
2002, p. 1-28.

\bibitem{BF6} Babin A. and Figotin A., \textsl{Nonlinear Photonic Crystals:
IV Nonlinear Schrodinger Equation Regime}, Waves in Random and Complex
Media, Vol. 15, No. 2 (2005), pp. 145-228.

\bibitem{BF7} Babin A. and Figotin A., \ \textsl{Wavepacket preservation
under nonlinear evolution}, submitted; e-print available online at arxiv.org
arXiv:math.AP/0607723

\bibitem{BMN1} Babin A., Mahalov A. and Nicolaenko B., \textsl{Global
regularity of 3D rotating Navier-Stokes equations for resonant domains},
Indiana University Mathematics Journal vol. 48 no. 3 (1999), p.1133-1176.

\bibitem{BMN2} Babin A., Mahalov A. and Nicolaenko B, \textsl{Fast Singular
Oscillating Limits and Global Regularity for the 3D Primitive Equations of
Geophysics}, M2AN,v.34,no.2, 2000, p.201-222.

\bibitem{Bambusi03} Bambusi, D., \textsl{Birkhoff normal form for some
nonlinear PDEs}, Comm. Math. Phys. 234 (2003), no. 2, 253--285.

\bibitem{BenYoussefLannes02} Ben Youssef, W.; Lannes, D., \textsl{The long
wave limit for a general class of 2D quasilinear hyperbolic problems}, Comm.
Partial Differential Equations 27 (2002), no. 5-6, 979--1020.

\bibitem{BermanI} Berman G.P, Izrailev F.M., \textsl{The Fermi-Pasta-Ulam
problem: 50 years of progress}, arXiv:nlin.CD

\bibitem{BM} Bogoliubov N. N. and Mitropolsky Y. A., \textsl{Asymptotic
Methods In The Theory Of Non-Linear Oscillations, }Delhi, Hindustan Pub.
Corp., 1961.

\bibitem{BonaCL05} Bona, J. L.; Colin, T.; Lannes, D., \textsl{Long wave
approximations for water waves}, Arch. Ration. Mech. Anal. 178 (2005), no.
3, 373--410.

\bibitem{Bourgain} Bourgain, J., \textsl{Global solutions of nonlinear Schr%
\"{o}dinger equations.} American Mathematical Society Colloquium
Publications, 46. American Mathematical Society, Providence, RI, 1999.

\bibitem{Caz} Cazenave T., Semilinear Schr\"{o}dinger equations. Courant
Lecture Notes in Mathematics, 10., New York University, Courant Institute of
Mathematical Sciences, New York; American Mathematical Society, Providence,
RI, 2003.

\bibitem{Colin} Colin, T., \textsl{Rigorous derivation of the nonlinear Schr%
\"{o}dinger equation and Davey-Stewartson systems from quadratic hyperbolic
systems}, Asymptot. Anal. 31 (2002), no. 1, 69--91.

\bibitem{ColinLannes} Colin, T.; Lannes, D., \textsl{Justification of and
long-wave correction to Davey-Stewartson systems from quadratic hyperbolic
systems.}, Discrete Contin. Dyn. Syst. 11 (2004), no. 1, 83--100.

\bibitem{Craig} Craig W.; Groves M. D., \textsl{Normal forms for wave motion
in fluid interfaces}, Wave Motion 31 (2000), no. 1, 21--41.

\bibitem{CraigSulemS92} Craig, W.; Sulem, C.; Sulem, P.-L., \textsl{%
Nonlinear modulation of gravity waves: a rigorous approach}, Nonlinearity 5
(1992), no. 2, 497--522.

\bibitem{Din} Dineen S., \textsl{Complex Analysis on Infinite Dimensional
Spaces}, Springer, 1999.

\bibitem{GW} Gallay T.; Wayne C. E., \textsl{Invariant manifolds and the
long-time asymptotics of the Navier-Stokes and vorticity equations on} $%
\mathbf{R}^{2}$., Arch. Ration. Mech. Anal. 163 (2002), no. 3, 209--258.

\bibitem{GiaMielke} Giannoulis, J.; Mielke, A.,\textsl{The nonlinear Schr%
\"{o}dinger equation as a macroscopic limit for an oscillator chain with
cubic nonlinearities.}, Nonlinearity 17 (2004), no. 2, 551--565.

\bibitem{Hayashi03} N. Hayashi and P. Naumkin, \textsl{Asymptotics of small
solutions to nonlinear Schr\"{o}dinger equations with cubic nonlinearities.}
Int. J. Pure Appl. Math. 3 (2002), no. 3, 255--273.

\bibitem{HPh} Hille E. and Phillips R. S., \textsl{Functional Analysis and
Semigroups}, AMS, 1991.

\bibitem{InfeldR} Infeld, E. and Rowlands, G. \textsl{Nonlinear Waves,
Solitons, and Chaos}, 2nd ed. Cambridge, England: Cambridge University
Press, 2000.

\bibitem{Iooss} Iooss, G.; Lombardi, E., \textsl{Polynomial normal forms
with exponentially small remainder for analytic vector fields.} J.
Differential Equations 212 (2005), no. 1, 1--61.

\bibitem{JolyMR98} Joly, J.-L.; Metivier, G.; Rauch, J., \textsl{Diffractive
nonlinear geometric optics with rectification}, Indiana Univ. Math. J. 47
(1998), no. 4, 1167--1241.

\bibitem{KalyakinUMN} Kalyakin, L. A.,\textsl{\ Long-wave asymptotics.
Integrable equations as the asymptotic limit of nonlinear systems.}, Uspekhi
Mat. Nauk 44 (1989), no. 1(265), 5--34, 247; translation in Russian Math.
Surveys 44 (1989), no. 1, 3--42.

\bibitem{Kalyakin2} Kalyakin L.A., Asymptotic decay of a one-dimensional
wave packet in a nonlinear dispersive medium, Math. USSR Sb. Surveys 60 (2)
(1988) 457--483.

\bibitem{Kuksin} Kuksin S. B., \textsl{Fifteen years of KAM for PDE.}
Geometry, topology, and mathematical physics, 237--258, Amer. Math. Soc.
Transl. Ser. 2, 212, Amer. Math. Soc., Providence, RI, 2004.

\bibitem{KSM} Kirrmann P.; Schneider G.; Mielke A., The validity of
modulation equations for extended systems with cubic nonlinearities, \ Proc.
Roy. Soc. Edinburgh Sect. A 122 (1992), no. 1-2, 85--91.

\bibitem{Lax} Lax P.D., \ Integrals of nonlinear equations of evolution and
solitary waves, Comm. Pure Appl. Math. 21 (1968), 467-490.

\bibitem{Maslov83} Maslov V.P., \textsl{Non-standard characteristics in
asymptotic problems}, Uspekhi Mat. Nauk \textbf{38}:6 (1983), 3-36,
translation in Russian Math. Surveys \textbf{38}:6 (1983),1-42.

\bibitem{MSZ} Mielke A., Schneider G., Ziegra A., \textsl{Comparison of
inertial manifolds and application to modulated systems}, Math. Nachr. 214
(2000), 53--69.

\bibitem{Nayfeh} Nayfeh, A. H., \textsl{Perturbation Methods}, New York:
Wiley, 1973.

\bibitem{Pankov} A. Pankov, \textsl{Travelling Waves And Periodic
Oscillations In Fermi-Pasta-Ulam Lattices}, Imperial College Press, 2005.

\bibitem{PW} Pierce R. D.; Wayne C. E., \textsl{On the validity of
mean-field amplitude equations for counterpropagating wavetrains},
Nonlinearity 8 (1995), no. 5, 769--779.

\bibitem{Schneider98a} Schneider, G., \textsl{Justification of modulation
equations for hyperbolic systems via normal forms}, NoDEA Nonlinear
Differential Equations Appl. 5 (1998), no. 1, 69--82.

\bibitem{Schneider05} Schneider, G., Justification and failure of the
nonlinear Schr\"{o}dinger equation in case of non-trivial quadratic
resonances. J. Differential Equations 216 (2005), no. 2, 354--386.

\bibitem{SU} Schneider G., Uecker H. \textsl{Existence and stability of
modulating pulse solutions in Maxwell's equations describing nonlinear optics%
}, Z. Angew. Math. Phys. 54 (2003), no. 4, 677--712.

\bibitem{Sulem} Sulem C. and Sulem P.-L. , \textsl{The Nonlinear Schrodinger
Equation}, Springer , 1999.

\bibitem{Weinstein} Soffer A., Weinstein M. I., \textsl{Resonances,
radiation damping and instability in Hamiltonian nonlinear wave equations},
Invent. Math. 136 (1999), no. 1, 9--74.

\bibitem{Weissert97} Weissert T.P.,\textsl{The Genesis of Simulation in
Dynamics: pursuing the Fermi-Pasta-Ulam problem, Springer-Verlag}, New York,
1997.

\bibitem{W} Whitham G., \textsl{Linear and Nonlinear Waves}, John Wiley \&
Sons, 1974.
\end{thebibliography}
\end{document}